\documentclass{article}
\usepackage[utf8]{inputenc}
\usepackage[margin=1in,left=1in]{geometry}
\usepackage{setspace}
\usepackage{slashed}
\usepackage{graphicx}
\usepackage{amssymb, amsmath, amsthm, graphics}
\usepackage{mathabx,epsfig}
\usepackage[mathscr]{euscript}
\usepackage{cancel}
\usepackage{tikz}
\usetikzlibrary{calc}
\usetikzlibrary{matrix}
\usepackage{mathtools}

\DeclarePairedDelimiter\floor{\lfloor}{\rfloor}

\usepackage[hyphens]{url}
\usepackage{hyperref}

\usepackage{enumerate} 

\usepackage{diagbox, eqparbox, hhline}    

\usepackage{stmaryrd}    

\usepackage{latexsym}
\usepackage[all]{xy}
\usepackage{color}

\usepackage{bbold}

\usepackage{tikz}
\input xy
\xyoption{all}
\usepackage{euscript}
\usepackage{multirow}
\usepackage{etex, pictexwd,dcpic}
\usetikzlibrary{positioning}
\usetikzlibrary{shapes.geometric}
\usetikzlibrary{shapes.misc}
\usetikzlibrary{calc}
\usetikzlibrary{positioning}

\usepackage{pgflibraryarrows}
\usepackage{pgflibrarysnakes}

\usetikzlibrary{trees} 
\usetikzlibrary[trees] 

\newtheorem{theorem}{Theorem}
\newtheorem{definition}[theorem]{Definition}
\newtheorem{lemma}[theorem]{Lemma}
\newtheorem{corollary}[theorem]{Corollary}
\newtheorem{proposition}[theorem]{Proposition}

\newtheorem{remark}{Remark}

\newtheorem{example}{Example}

\numberwithin{equation}{section}

\usepackage{cleveref}

\crefformat{section}{\S#2#1#3} 
\crefformat{subsection}{\S#2#1#3}
\crefformat{subsubsection}{\S#2#1#3}

\renewcommand{\(}{\begin{equation*}}
\renewcommand{\)}{\end{equation*}}
\newcommand{\bea}{\begin{eqnarray*}}
\newcommand{\eea}{\end{eqnarray*}}

\renewcommand{\a}{\alpha}
\renewcommand{\b}{\beta}

\def\endofproof {\hfill{$\Box$}\\}

\newcommand{\Exterior}{\mathchoice{{\textstyle\bigwedge}}%
{{\bigwedge}}%
{{\textstyle\wedge}}%
{{\scriptstyle\wedge}}}

\def\S{\ensuremath{\ES{S}}}

\def\S{\ensuremath{\ES{S}}}

\newcommand{\N}{\ensuremath{\mathcal N}}
\newcommand{\beq}{\begin{equation}}
\newcommand{\eeq}{\end{equation}}

\newcommand{\into}{\hookrightarrow}

\newcommand{\op}[1]{\ensuremath{\operatorname{#1}}}

\newcommand{\ES}[1]{\ensuremath{\EuScript{#1}}}

\newcommand{\theproof}{\noindent {\bf Proof.\ }}

\numberwithin{equation}{section}

\renewcommand{\(}{\begin{equation}}
\renewcommand{\)}{\end{equation}}

\def\ch{{\rm  ch}}

\def\1{{\bf 1}}

\def\<{\langle}
\def\>{\rangle}

\def\a{\alpha}
\def\b{\beta}

\numberwithin{equation}{section}


\newcommand{\R}{\ensuremath{\mathbb R}}
\newcommand{\RR}{\ensuremath{\mathbb R}}
\newcommand{\NN}{\ensuremath{\mathbb N}}

\newcommand{\ZZ}{\ensuremath{\mathbb Z}}
\newcommand{\Z}{\ensuremath{\mathbb Z}} 
\newcommand{\QQ}{\ensuremath{\mathbb Q}}
\newcommand{\HH}{\ensuremath{\mathbb H}}
\newcommand{\BB}{\ensuremath{\mathbf B}}

\newcommand{\CC}{\ensuremath{\mathbb C}}

\newcommand{\chp}{\ensuremath{\mathscr{C}\mathrm{h}^{+}}}

\newcommand{\sset}{\ensuremath{s\mathscr{S}\mathrm{et}}}

\newcommand{\sab}{\ensuremath{s\mathscr{A}\mathrm{b}}}

\newcommand{\sh}{\ensuremath{\mathscr{S}\mathrm{h}}}

\newcommand{\E}{\ensuremath{\mathscr{E}}}

\newcommand{\map}{\mathrm{Map}}

\begin{document}

\title{
Differential KO-theory: \\ Constructions, computations, and applications}

\author{
Daniel Grady and Hisham Sati\\
}

\maketitle

\begin{abstract} 

We provide several constructions in differential KO-theory. First, we construct a differential refinement of the $\hat{A}$-genus and a pushforward leading to a Riemann-Roch theorem. We set up a differential refinement of the Atiyah-Hirzebruch spectral sequence (AHSS) for differential KO-theory and 
explicitly identify the differentials, including ones which mix geometric and topological data. We highlight the power of these explicit identifications by providing a characterization of forms in the image of the Pontrjagin character. Along the way, we fill gaps in the literature where K-theory is usually worked out leaving KO-theory essentially untouched.
We also illustrate with examples and applications, including higher tangential structures, 
Adams operations, and a differential Wu formula.

\end{abstract}

\tableofcontents

\section{Introduction}

The initial goal of this project was to explicitly develop twisted differential KO-theory.
However, as we embarked on the project, we realized that doing all of this at once would be
a massive task, in a large part due to the fact that the differential (untwisted) real theory
itself has not been fully developed, at least not in a systematic way that we could find in 
the literature. Moreover, even in the untwisted differential case, there are many constructions which have not been fully worked out even for the differential refinement of complex K-theory and we would like to have these constructions available for both KO and its twists. Therefore, we have decided to first develop the untwisted
differential theory in the current paper and leave the twisted case for a follow-up treatment. 
We hope that the current treatment will serve both as a stepping stone for that original 
goal, as well as being useful and interesting in its own right. For that purpose we have
aimed to make the discussions as complete as possible, highlighting similarities 
and differences with the complex case, both for K-theory and its differential incarnation,
as well as with the underlying bare topological real case, i.e. that of KO-theory. We also note that the two main results of this paper -- the Riemann-Roch theorem and the identification of the differentials in the refined AHSS -- have not been developed even for the differential refinement of complex K-theory.

\medskip
There are two main definitions of KO-theory of a space $X$. 
The first is the one by Atiyah and Hirzebruch (AH) using real vector bundles over $X$ \cite{AH}. 
The second is via the KO-spectrum, which is an $\Omega$-spectrum with 
every term $E_{8r}$ being the space $\ZZ \times \op{BO}$, 
and every space $E_{8r+4}$ being the space $\ZZ \times \op{BSp}$, and
other terms arising from Bott's periodicity theorem for the orthogonal 
group. 
For $X$ a finite-dimensional CW-complex, the definition as homotopy classes of maps
$[X, \ZZ \times \op{BO}]$ agrees with the AH definition. See 
\cite{Ad} for remarks  on this comparison.

\medskip
Differential KO-theory has not attracted as much explicit attention as its 
complex counterpart, i.e. differential K-theory.
Freed in \cite{Fr} proposes considering smooth extensions of 
real versions of K-theory, with an indication of how cycle
models  could look like.  Further brief assertions are indicated and 
used by him in applications in \cite{Fr-lec}. Brief remarks by Bunke and Schick 
appear in \cite[Section 4.9]{BS} with some expectations and 
the need  is raised for fully working out the various models. 
On the applications side, in addition to Freed \cite{Fr-lec},
Belov and Moore \cite{BeM} use differential KO-theory
for the construction of 
Chern-Simons action functionals in type II string theory. 
\footnote{As we mentioned above, we will consider the twisted case separately,
where we will give extensive references for motivations and existing constructions in that case.}

\medskip
The two main goals of this paper are to provide explicit constructions of differential 
KO-theory and to provide computational tools
which allow for explicit computations of this cohomology theory for various smooth manifolds. 
The latter will be in the form of an Atiyah-Hirzebruch spectral sequence (AHSS), in 
a treatment that is parallel the constructions in the series \cite{GS3}\cite{GS4}\cite{GS5}\cite{GS6}.
One point to highlight in contrast with the complex case is that the coefficients of 
KO-theory have torsion, which makes the discussion here more subtle. Another major point is that the identification of the `geometric' differentials in the spectral sequence has not been achieved even in the (refinement of the) complex case. The reason for this has less to do with the differences between KO and K and more to do with fact that we did not know what these differentials were during the writing of \cite{GS3} and \cite{GS4}.

\medskip
Note that one can consider KR-theory \cite{Real}, which can be viewed as a generalization 
of both KO-theory and K-theory. However,  for our purposes (as we will see) we prefer to 
deal with  KO-theory directly without passing  through KR-theory. 
While applications of the latter are quite interesting
(see e.g. \cite{DFM1}\cite{DFM2}\cite{Fo}\cite{DMDR}\cite{HMSV}), 
the former suffices for the ones we have in mind. Note that in the topological case the AHSS for KR-theory is constructed in \cite{Dug},
with ways of how to reduce to KO-theory when the involution is trivial.

\medskip
What does it mean to provide a differential refinement of a cohomology theory? 
A priori one would like to add geometric data to the otherwise topological 
structure. In implementing that, one might mean three different notions:
\begin{enumerate}
\item A cohomology theory on manifolds, i.e., in the presence of smooth structure.
\item A cohomology theory with geometric data included: connections.
\item A cohomology theory with geometric data included: Chern character forms.
\end{enumerate}
So one might consequently ask the following questions: 

\begin{enumerate}[(i)]
\item Which one of the above does one really mean when talking 
about differential cohomology?
\item  Are the above equivalent, or at least related? 
\end{enumerate} 
In \cite{BNV}, it was shown that every sheaf of spectra fits into a hexagon diagram, similar 
to the differential cohomology hexagon  of Simons and Sullivan \cite{SSu2}. Hence, every 
sheaf of spectra can be thought of as representing a differential cohomology theory. 
Such generality leads to many differential refinements which are quite different from the 
Hopkins-Singer type refinements which were originally considered in \cite{HS}. 
Here we will clarify what we mean by a differential KO-theory, presenting several approaches with 
enough generality to thereby also clarifying a bit what is meant by a differential cohomology theory
in general. After producing the various possibilities, we explore interrelations which we
make explicit.  

\medskip
One conceptual approach is to view a differential cohomology theory as an amalgam of an 
underlying (topological) cohomology theory and the data of differential forms, so that one 
has a schematic picture 
$$
\xymatrix{
\text{Differential cohomology} \ar[d] \ar[rr] && \text{Closed forms} \ar[d]
\\
\text{Cohomology} \ar[rr] && \text{de Rham cohomology} \;.
}
$$
That is, we have a fiber product or twisted product 
$$
`` \text{Differential cohomology}= \text{Cohomology}\times_{\text{de Rham}} \text{Closed forms}"\;.
$$
This has various explicit presentations; for instance as a triple in the original model
\cite{HS}, as well as using sheaves of spectra \cite{BNV} and toposes and smooth stacks \cite{Urs}. 

\medskip
\noindent {\bf Notation.} In order to highlight the various notions of refinements of 
cohomology that we introduce, we adopt the following notation.

\medskip
\hspace{-10mm}
\begin{tabular}{|lll|}
\hline
{\bf Topological:} & ${\rm BO}$ & The topological space classifying stable orthogonal bundles.
\\ 
& KO & The KO-theory spectrum. 
\\
& $\op{KO}^\ast(-)$ & The generalized cohomology theory associated to the KO-theory spectrum.
\\
\hline
{\bf Smooth:} & $\mathbf{B}{\rm O}$ & The smooth stack (on the site of smooth manifolds)
\\
&& that classifies smooth stable orthogonal bundles.
\\
& $\mathbf{K}{\rm O}$ & The sheaf of spectra (on the site of smooth manifolds) 
\\
&& obtained by $\infty$-group completion of the commutative monoid $\coprod_{n\in \NN}\mathbf{B}{\rm O(n)}$.
\\
& $\mathbf{K}{\rm O}^*(-)$ & The smooth generalized cohomology theory corresponding to 
\\
&& the sheaf of spectra for 
$\mathbf{K}{\rm O}$. This is a functor from the homotopy category 
\\
&& of smooth stacks to abelian groups (or rings). 
\\
\hline
{\bf Smooth}  & $\mathbf{B}{\rm O}_\nabla$ & The smooth stack (on the site of smooth manifolds) classifying 
\\
{\bf with connection:} && smooth stable orthogonal bundles with
compatible connection $\nabla$.
\\
& $\mathbf{K}{\rm O}_\nabla$ & The sheaf of spectra (on the site of smooth manifolds) 
\\
&& obtained by $\infty$-group completion of the commutative monoid $\coprod_{n\in \NN}\mathbf{B}{\rm O(n)}_{\nabla}$.
\\
& $\mathbf{K}{\rm O}^*_\nabla(-)$ & The smooth generalized cohomology theory corresponding to 
\\
&& the sheaf of spectra for 
$\mathbf{K}{\rm O}_{\nabla}$. This is a functor from the homotopy category 
\\
&& of smooth stacks to abelian groups (or rings). 
\\
\hline
{\bf Differential:} & $\widehat{\rm BO}$ &  The smooth stack (on the site of smooth manifolds)
\\
&& that classifies stable orthogonal bundles, equipped with differential form data.
\\
& $\widehat{\rm KO}$ & The sheaf of spectra (on the site of smooth manifolds) 
\\
&& obtained by taking a fiber product of the KO-theory spectrum with 
\\
&& (closed) differential forms, over 4-periodic real cohomology. 
\\
& $\widehat{\rm KO}^*(-)$ & The differential cohomology theory corresponding to the sheaf of spectra for 
$\widehat{\rm KO}$.
\\
\hline
\end{tabular}

\medskip
The three versions of differential KO-theory, namely  $\mathbf{K}{\rm O}$,
$\mathbf{K}{\rm O}_\nabla(-)$, and $\widehat{\rm KO}$ 
are related via diagram \eqref{shspdfko} in section \ref{Ch-topKO}:
\( \label{3-related}
\xymatrix@R=1.3em{
\mathbf{K}\op{O}_{\nabla}\ar[rr]^-{\widehat{\rm cyc}}\ar[d]_-{\mathcal{F}} && 
\widehat{\op{KO}} \ar[d]^-{\mathcal{I}}
\\
\mathbf{K}\op{O}\ar[rr]_-{\rm cyc} && \op{KO}\;.
}
\)
The map $\mathcal{F}:\mathbf{K}\op{O}_{\nabla}\to \mathbf{K}{\rm O}$ forgets the data of 
connections of vector bundles, while the map on the right $\mathcal{I}:\widehat{\op{KO}}\to \op{KO}$ 
is the canonical map  which topologically realizes the geometric data. The cycle map cyc and its 
differential refinement ${\widehat{\rm cyc}}$ are described later in that section and in 
Proposition \ref{csconvb}.

\medskip
We also encounter other variants related to the spectrum KO, namely the ones arising from 
taking the underlying discrete form ${\rm O}^\delta$ of the stable orthogonal group ${\rm O}$, leading to 
the (orthogonal) algebraic K-theory of $\RR$ which we denote $\op{K}({\rm O}^\delta)$. This is related to $\mathbf{K}{\rm O}$ 
and to the  geometrically discrete spectrum $\delta({\rm KO})$ via diagram \eqref{agtpkocn1}, i.e., 
\(
\xymatrix@R=4pt@C=4em{
& \mathbf{K}{\op{O}} \ar[rd]^-{\mathcal{I}} &
\\
\op{K(O^{\delta})}\ar[ru]^{j}\ar[rr] & & \delta(\op{KO})\;.  &
}
\)
Similarly, the spectrum $\mathbf{K}{\rm O}_{\nabla}$ satisfies diagram \eqref{agtpkocn2}, i.e., 
\(
\xymatrix@R=4pt@C=4em{
& \mathbf{K}{\op{O}}_{\nabla} \ar[rd]^-{\mathcal{I}} &
\\
\op{K(O^{\delta})}\ar[ru]^{j}\ar[rr] & & \delta(\op{KO})\;.  &
}
\)

The paper is organized as follows. In Section \ref{Ch-topKO} we provide the setting of the underlying topological 
KO-theory, starting with the coefficients, the description of the spectrum, and the product structure.
In Subsection \ref{Sec-PontCh} we recall the Pontrjagin classes $p_i$ and the Pontrjagin character ${\rm Ph}$, together 
with their properties and corresponding differential form representatives. Then, in Subsection \ref{Sec-orgen},
we recall orientations and genera associated to KO-theory, mainly the $\hat{A}$-genus and to a lesser extent 
the somewhat related signature. Then we describe the Thom class in relation to the index. We end this section in 
Subsection \ref{Sec-cplx} by 
highlighting the basics of the relation to complex K-theory, focusing mainly on the Bott sequence which 
mixes the two theories, and which we refine later in Subsection \ref{Sec-cplx-diff}.

\medskip
The main section of the paper is Section \ref{Ch-diffKO}, where we construct the three  variants of 
differential KO-theory (summarized above). This uses smooth stacks, starting with classifying stacks of 
real orthogonal bundles with connections in Subsection \ref{Sec-Pont-class}. We also consider differential refinements 
of Pontrjagin classes as natural transformations between the category of smooth vector bundles with 
connections and differential integral cohomology in degree a multiple of 4. We also consider the 
differential Pontrjagin characters, which end up being polynomials in differential Pontrjagin classes, 
except that multiplication of these is given by the Deligne-Beilinson cup product in differential cohomology.
(A more fundamental construction via sheaves of spectra is given later in Subsection \ref{Sec-Diff-Pont-Cha}
in Section \ref{Ch-diffKO} -- see Proposition \ref{phchrfuni}).
Similarly, we describe the genera, particularly, the differential A-genus, which we denote $\widehat{\mathbb{A}}$. 
In Subsection \ref{Sec-construct} we use the smooth stacks to describe the three differential variant theories in detail, 
including the properties they enjoy. In Subsection \ref{Sec-cplx-diff} we provide a differential refinement of the 
Bott sequence from Subsection \ref{Sec-cplx} above. 
This ends up having the very interesting and rich feature of relating and mixing 
$\widehat{\rm KO}$-theory to the underlying KO-theory, to the complex K-theory and even to the 
differential refinement $\widehat{\rm K}$ of the latter. One main ingredient in the differential 
refinement is the flat theory $\op{KO}_{\rm flat}(-)$ which is related to the theory with $U(1)$ coefficients,
$\widetilde{\op{KO}}^{*}(M; {\rm U}(1))$, both of which we describe in Subsection \ref{Sec-flat}, 
where we also describe the corresponding properties and provide a vanishing theorem.
The cycle map cyc and its differential refinement $\widehat{\rm cyc}$ appearing in diagram 
\eqref{3-related} are described explicitly in Subsection \ref{Sec-cyclemap}. In Subsection \ref{Sec-Diff-Pont-Cha}
we describe the differential Pontrjagin character via sheaves of spectra, refining the topological description 
from Subsection \ref{Sec-PontCh}, and explain how the character behaves
under various maps in the differential cohomology diamond or hexagon diagram \eqref{kodfdiam} 
for $\widehat{\rm KO}$.

\medskip
Further constructions are provided in Section \ref{Ch-further}.  In order to describe the orientation 
and push-forward, in Subsection \ref{Sec-compact} we introduce differential KO-theory with compact  
vertical supports, $\widehat{\rm KO}_c(V)$, for a vector bundle $V$. There we also define 
differential forms with vertical 
compact support on the Thom space, 
integration over the fiber $\int_D$, and a corresponding  differential desuspension map $\hat{\sigma}_n^{-1}$. 
These are used in Subsection \ref{Sec-or-gen} to connect the differential A-genus $\widehat{\mathbb{A}}$
with a Thom isomorphism theorem for $\widehat{\rm KO}$-theory. This requires differential Thom classes
$\hat{\nu}\in \widehat{\op{KO}}(V)$, the set of which we study and a distinguished one of which is identified via naturality in the Spin case. Then we prove a differential KO Thom injection formula, which we use to 
provide an expression for $\widehat{\mathbb{A}}$ in terms of the Thom injection $\Phi_{\widehat{H}}$. 
The constructions in this section then culminate in Subsection \ref{Sec-push} 
with establishing a pushforward and a Riemann-Roch theorem for $\widehat{\rm KO}$-theory.

\medskip
Having constructed differential KO-theory, we would like to be able to compute. So in Section \ref{Ch-AHSS}
we describe our main computational tool, which is the Atiyah-Hirzebruch spectral sequence (AHSS).
We start in Subsection \ref{Sec-AHSS-KObare} by elucidating the AHSS for topological KO-theory, summarizing the main
results scattered in the literature, which we hope will be useful in its own right. 
Here torsion is very subtle and is what distinguishes the real theory from its complex counterpart. 
Improving on the existing results 
in the case of KO-theory, we identify a large portion of the differentials in analogy to the complex case. 
This serves as a basis for the differential version in Subsection \ref{Sec-AHSS-KOhat}. Here we use 
our earlier general constructions in \cite{GS3} to construct the differential AHSS and identify its main pages.
We show in Subsection \ref{Sec-Ident} that here there are two main differentials, one on the 1st quadrant 
and another on the 4th quadrant,
both of which we identify using the refinements of Steenrod squares $\widehat{Sq}$ we established earlier in  \cite{GS2}. 
Being a differential theory, $\widehat{\rm KO}$ will have an AHSS the differential within which should
mix geometric and topological data. Indeed, we explicitly identify such mixing differentials, which are not
differential refinements of those of the underlying topological theory, as images of certain differentials forms,
depending on dimension ($8k$ vs. $8k+4$). This then allows us in Subsection \ref{Sec-IntPont} to provide
integrality and realizability results for the Pontrjagin character Ph, presenting  answers for when a differential 
form can arise as the image of Ph.  Then we end this section in Subsection \ref{Sec-comp} by providing explicit 
computations for $\widehat{\rm KO}$-groups of spaces, particularly spheres. 
These serve as illustrations of the utility of the machinery 
developed earlier in this section. 

\medskip
Equipped with the properties and computational tools for differential KO-theory, we provide in 
Section \ref{Ch-apps} some applications of these. First, in Subsection \ref{Sec-higher} we study integrality 
conditions arising from higher structures, such as String, Fivebrane, and Ninebrane structures (see \cite{SSS2}\cite{SSS3}\cite{9}). These
can be obtained efficiently from the results in Subsection \ref{Sec-AHSS-KOhat}.
In Subsection \ref{Sec-Adams} we construct Adams operations for $\widehat{\rm KO}$-theory, extending those
of KO-theory in such a way that is compatible with ones constructed by Bunke \cite{Bu-Ad} 
for $\widehat{\rm K}$-theory. We also do so for (the degree zero part of) 
${\bf K}{\rm O}_\nabla$, whose complex counterpart 
does not seem to have been considered before. We further describe main properties, uniqueness, 
and the relation between the two versions of the operations via the differential cycle maps. 
The third application in Subsection \ref{Sec-Classicalw} 
is describing differential refinements of Stiefel-Whitney classes $w_j$, which are 
useful for those refined differentials given by $\widehat{Sq}^i$'s, leading to a Wu formula at the level
of differential cohomology. 

\medskip
Finally, in  the Appendix,
we provide basic higher categorical framework
in just enough abstraction to allow our constructions in the bulk of the paper to go through. This we hope 
could also be of independent interest as a summary of the homotopy background for differential 
cohomology in general. 

\section{Recollections on topological ${\rm KO}$-theory}
\label{Ch-topKO}

In this section we provide a recollection of basic material on KO-theory (see 
\cite{AH}\cite{Bo}\cite{Ad}\cite{Ka} \cite{Hu}\cite{Sw}). 
For each space $X$, the ${\rm KO}$-theory of $X$ can be defined as the 
Grothendiek group completion of the monoid of real vector bundles on 
$X$ (with monoidal operation given by Whitney sum). When $X$ is 
pointed, the reduced theory can be defined as 
$$
\widetilde{{\rm KO}}(X):=\op{ker}({\rm KO}(X) \longrightarrow {\rm KO}(x_0))\;,
$$
so that ${\rm KO}(X) \cong \widetilde{{\rm KO}}(X) \oplus \ZZ$.

\medskip
\noindent {\bf The spectrum definition of KO-theory.}
The reduced functor $\widetilde{{\rm KO}}(-)$ is represented by
the infinite loop space $\ZZ \times \op{BO}$, with
$\ZZ \times \op{BO} \simeq \Omega^8(\ZZ \times \op{BO})$,
where $\op{BO}$ is the classifying space of the 
orthogonal group $\op{O}$ and $\Omega^p \op{BO}$ is the $p$-th based loop space on $\op{BO}$.
The cohomology theory $\op{KO}^*(-)$ is represented by 
the $\Omega$-spectrum $\op{KO}$ with the $n$th space 
$\underline{\op{KO}}_n=\Omega^m(\ZZ \times \op{BO})$, where
$m$ is chosen so that $n - m \equiv 0$ mod $8$
and $0 \leq m <8$. The tensor product pairing is represented by 
pairings of spectra that makes $\op{KO}$ into a 
an $E_\infty$ ring specturm.  
The unit $\mathbb{S} \to \op{KO}$ is 
generated by the map $S^0 \to \ZZ \times \op{BO}$ that takes the 
non-basepoint to trivial rank 1 bundle on the point.

\medskip
Explicitly, Bott periodicity implies 
\bea
\op{O}\simeq \Omega (\ZZ \times \op{BO}), \qquad \quad
\op{O}/\op{U} \simeq \Omega \op{O}, \qquad \quad
\op{U}/\op{Sp} \simeq \Omega (\op{O}/\op{U}), \qquad \quad
\ZZ \times \op{BSp} \simeq \Omega(\op{U}/\op{O})\;, 
\\
\op{Sp} \simeq \Omega (\ZZ \times \op{BSp}), \qquad \quad
\op{Sp}/\op{U} \simeq \Omega \op{Sp}, \qquad \quad
\op{U}/\op{O} \simeq \Omega(\op{Sp}/\op{U}), \qquad \quad
\ZZ \times \op{BO} \simeq \Omega (\op{U}/\op{O})\;. 
\eea
Thus the corresponding $\op{KO}$-spectrum 
$\op{KO}=(\varepsilon_n: \Sigma \op{KO}_n \to \op{KO}_{n+1})_{n \in \ZZ}$
is given as follows  (see e.g. \cite{Yam})
\bea
\op{KO}_{8n}= \Z \times \op{BO}, \qquad \quad
\op{KO}_{8n+1}=\op{U}/\op{O}, \qquad \quad
\op{KO}_{8n+2}=\op{Sp}/\op{U}, \qquad \quad
\op{KO}_{8n+3}=\op{Sp}
\\
\op{KO}_{8n+4}= \Z \times \op{BSp}, \qquad \quad
\op{KO}_{8n+5}=\op{U}/\op{Sp}, \qquad \quad
\op{KO}_{8n+6}=\op{O}/\op{U}, \qquad \quad
\op{KO}_{8n+7}=\op{O}\;.
\eea
A topological symmetric $E_\infty$ ring 
spectrum that represents (periodic) topological KO-theory
is constructed in \cite{Jo} using an appropriate subspace of
the space of Fredholm  operators. This is developed further in
\cite[Appendix]{AGG}.

\medskip
There are elements $\eta$, $\alpha$, and $\beta$, with $|\eta|=-1$, $|\alpha|=-4$, 
$|\beta|=-8$, as generators of $\op{KO}^{-1}=\pi_1(\op{KO})\cong \ZZ/2$,
$\op{KO}^{-4}=\pi_4(\op{KO})\cong \ZZ$, and $\op{KO}^{-8}=\pi_8(\op{KO})\cong \ZZ$,
respectively. The elements $\eta$, $\alpha$, $\beta$ are represented by the real Hopf line bundle 
over $S^1$, the symplectic Hopf line bundle over $S^4$, and 
the canonical bundle over $S^8$, respectively \cite{Ka}. 
Then the homotopy groups of the spectrum are given by Bott periodicity as
$$
\pi_i(\op{KO})=
\left\{
\begin{array}{ll}
\ZZ\{\beta^k\} & \text{for}\; i=8k,
\\
\ZZ/2\{\eta \beta^k\} & \text{for}\; i=8k+1,
\\
\ZZ/2\{\eta^2\beta^k\} & \text{for} \; i=8k+2,
\\
\ZZ\{\alpha \beta^k\} & \text{for}\;  i=8k+4,
\\
0 & \text{otherwise}.
\end{array}
\right.
$$
As a graded ring, the coefficient ring is given by (\cite{ABS}; see also \cite[p. 73]{Bott})
\(
\label{eq-coeff}
\op{KO}^*({\rm pt})=\pi_*(\op{KO})=\ZZ[\eta, \alpha, \beta^{\pm 1}]/
(2\eta, \eta^3, \eta \a, \a^2 -4\b)\;.
\)
For example, $\op{KO}^{-10}({\rm pt})=\left<\beta \eta^2 \right>\cong \Z/2$
and $\op{KO}^6({\rm pt})=\left<\beta^{-1} \eta^2 \right>\cong \Z/2$.


\subsection{The Pontrjagin character}
\label{Sec-PontCh} 

Consider the element ${\rm Ph}\in H^*(\op{BO}(n);\QQ)$ determined by the even characteristic series 
$$
{\rm Ph}=\frac{1+(-1)^{n+1}}{2}+\sum_{i=1}^{[n/2]}(e^{x_i}+e^{-x_i})\;.
$$
Writing this characteristic series in terms of elementary symmetric functions, using Newton's identities, one can expand  in terms of Pontrjagin classes. This gives a characteristic class whose first few terms are 
\(
\label{Eq-Php}
\op{Ph}=\op{rank} + p_1 + \tfrac{1}{12}(p_1^2 - 2 p_2) + \cdots \;,
\)
called the Pontrjagin character. Since each $p_n$ represents an integral cohomology class, each homogeneous component of the formal power series is rational. The forms $p_n$ are natural with respect to pullback of vector bundles with connection. From the identity $p_n(V)=(-1)^nc_{2n}(V\otimes \CC)$, it follows at once that ${\rm Ph}$ gives a natural transformation factoring through $K$-theory
\(\label{pchrtrko}
\op{Ph}: \op{KO}^0(-) \longrightarrow \op{K}^0(-) \longrightarrow 
\bigoplus_{i \geq 0} H^{4i}(-; \QQ)\;.
\)
Moreover, one can show that this transformation preserves the ring structure. Thus, the Pontrjagin character is the analogue for KO-theory of the Chern character transformation and, in fact, it factors through this transformation.

\begin{remark}
[Properties of the Pontrjagin character]\label{RemPropPont}
The following are satisfied:

\item {\bf (i)} $\op{Ph}$ is a natural transformation of cohomology theories. 

\item {\bf (ii)} It is graded by $\op{Ph}_k=\op{ch}_k \circ (\otimes_\RR \CC)$.
It is also the composition of the Chern character and the complexification
$$
\op{Ph}=\op{ch} \circ c: 
\op{KO} \overset{c}{\longrightarrow} \op{K}
\overset{\op{ch}}{\longrightarrow} \mathscr{H}\QQ[u,u^{-1}]\;,
$$
i.e. $\op{Ph}(E)=\op{Ch}(E \otimes_\RR \CC)$. In components,
$\op{Ph}_k(E)=(-1)^k\op{ch}_{2k}(E \otimes_\RR \CC)$.

\item {\bf (iii)} It is an isomorphism when tensored with $\QQ$.  The rational
Pontrjagin character defines an isomorphism 
$$
\op{Ph}: \op{KO}(X) \otimes \QQ\cong [X, \ZZ \times \op{BO}] \otimes \QQ
\longrightarrow  H^{4*}(X; \QQ)\;.
$$
In components, $\op{KO}^i(-) \otimes \QQ \cong \oplus_{k=0}^{\infty}H^{4k+i}(-; \QQ)$,
$i=0, 1, 2, 3$. 

\item {\bf (iv)}  It satisfies 
$$
\op{Ph}(E \otimes E')=\op{Ph}(E) \cdot \op{Ph}(E')
\qquad
\text{and}
\qquad 
\op{Ph}(E \oplus E')=\op{Ph}(E) + \op{Ph}(E')\;.
$$

\item {\bf (v)} In degree 4, $\pi_4(\op{KO})=\op{KO}(S^4) \to H^*(S^4; \QQ)$,
the Pontrjagin character of the generator $\alpha$ can be evaluated
on the fundamental class $[S^4]\in H_4(S^4; \ZZ)$ giving 
$\langle \op{ph}(\alpha), [S^4] \rangle =2$. 
For $X=\op{pt}$, one has $\op{KO}^*(X)\cong \op{KO}_*$ and $H^{**}(X; \QQ)\cong \QQ$
and $\op{ph}$ is entirely determined by 
$$
\op{Ph}(\eta)=0\;, \qquad \op{Ph}(\alpha)=2\;.
$$
In particular, $\op{Ph}$ is \emph{integral} (see \cite{Och})
$
\op{Ph}: \op{KO}_* \to \ZZ
$.
\end{remark}

\begin{example}[Pontrjagin character of the universal bundle]
Let $\xi$ be the universal oriented 2-plane bundle over 
$\CC P^\infty=BSO(2)$, and let $u=\xi-2 \in KO(\CC P^\infty)=\ZZ[u]$,
the ring of formal power series in $u$ (see \cite{AW}).  
For $x\in H^2(\CC P^\infty)$
the generator, one has (see \cite{CJ}) 
$\op{Ph}(u)=\sum_{i \geq 1} 
\frac{1}{(2i)!}{x^{2i}}$.
\end{example}

\medskip
\noindent {\bf A more general take on the Pontrjagin character.} The rationalized coefficients of $\op{KO}$ are given by (as utilized, e.g., in \cite{BZ})
\(
\label{Rat-coeff}
\op{KO}_*({\rm pt})\otimes \QQ \cong \QQ[\a, \b, \b^{-1}]/\langle \a^2 - 4 \b \rangle
\cong \QQ[\alpha,\alpha^{-1}]
\)
with $\a \in \op{KO}_4({\rm pt})$ and $\beta \in \op{KO}_8({\rm pt})$. 
This implies that rational KO-theory is 4-periodic. It is a module over the above
coefficients with the module multiplication with $\alpha/2$ implementing 
the 4-periodicity. 
The complexification $c: \op{KO}^*({\rm pt}) \to \op{K}^*({\rm pt})\cong \ZZ[u, u^{-1}]$
satisfies $c(\a)=2 u^2$ and $c(\b)=u^4$, where $u \in \op{K}^2({\rm pt})$. 
This implies that, rationally, the coefficients for complex K-theory split as two 
copies of the coefficients of KO-theory (see \cite{BZ})
$$
{\bf c}:=c + u^{-1} c: (\op{KO}_* \oplus \op{KO}_*) \otimes \QQ 
\overset{\cong}{\longrightarrow} \op{K}_* \otimes \QQ\;.
$$
More generally, there is an isomorphism 
$\op{KO}^*_\QQ(X) \cong H^*(X; \pi_{-*}(\op{KO}) \otimes \QQ)$ and the 
rationalization gives a map
$$
{\rm Ph}: \op{KO}^*(X) \longrightarrow H^*(X; \pi_{-*}(KO)\otimes \QQ)
$$
which coincides with the above Pontrjagin character, and the diagram 
$$
\xymatrix{
\op{K}^*(X)\ar[rr] \ar@/^2pc/[rrrr]^-{\rm Ch} && 
\op{K}^*(X)\otimes \QQ\ar[rr]^-{\cong} && H^{*}(X;\QQ[u,u^{-1}])
\\
\op{KO}^*(X)\ar[rr]\ar[u]^-{c} \ar@/_2pc/[rrrr]^-{\rm Ph} &&
\op{KO}^*(X)\otimes \QQ\ar[rr]^-{\cong}\ar[u]^-{c\otimes \QQ} 
&&  H^{*}(X;\QQ[\alpha,\alpha^{-1}])\ar[u]^-{\alpha\mapsto 2u^2}
}
$$
commutes.
The Pontrjagin character defines an injection 
\cite[Theorem 4.29]{MM}
$$
[\op{BSO}, \op{BO}] \longrightarrow H^*(\op{BSO}; \QQ)\;.
$$
Note that the real representation ring of the special orthogonal group is given by the polynomial ring
$$
RO(\op{SO}(2n+1))=\ZZ[\gamma_1, \gamma_2, \cdots, \gamma_n]
$$
where $\gamma_i=\gamma^i(V - \dim V)$, with $V$ the standard representation of 
$\op{SO}(2n+1)$ on $\RR^{2n+1}$.
The augmentation idea $IO(\op{SO}(2n+1)$ is the usual 
augmentation ideal so the completion 
$
RO(\op{SO}(2n+1))^\wedge =\ZZ[\![\gamma_1, \gamma_2, \cdots, \gamma_n]\!]
$
is the usual power series ring. Now each finite skeleton of $\op{BSO}$ is contained
in some $BSO(2n+1)$ so 
$$
[\op{BSO}, \op{BO}]=\ZZ[\![\gamma_1, \gamma_2, \cdots]\!]\;.
$$
Then $[\op{BSO}, \op{BO}]$ maps indicatively into 
$[\op{BSO}, \op{BO}]\otimes \QQ \cong [\op{BSO}, \op{BO}[\QQ]]$ where 
$\op{BO}[\QQ]$ denotes the rational type of $\op{BO}$. However, 
$\op{Ph}: \op{BO}[\QQ] \to \prod_{n \geq 1} K(\QQ, 4n)$ is 
a homotopy equivalence.

\medskip
\noindent {\bf Geometric representatives.} 
The Pontrjagin classes admit differential form representatives via Chern-Weil theory. Taking the above approach one then immediately has a \emph{Pontrjagin character form} which depends on a choice of connection, as the right hand side 
of \eqref{Eq-Php} a formal power series of certain closed differential forms. As such, we review the construction of the Pontjagin forms via Chern-Weil theory
in \cite[Section 9.6]{GHV2}. Let $V\to M$ be an 
$n$-dimensional real vector bundle on $M$ and let $\mathfrak{gl}_n$ denote the Lie 
algebra of endomorphisms of the fiber. Consider the homogeneous functions
$$
\sigma_p(\varphi):={\rm tr}(\Exterior^p\varphi)\;,\ \ \varphi\in \mathfrak{gl}_n\;.
$$
Then we have the formula
$$
{\rm det}(\varphi+\lambda\mathbb{1})=\sum_{p=0}^n\sigma_p(\varphi)\lambda^{n-p}\;.
$$
Fix a connection $\nabla$ on $V\to M$ and let $\mathcal{F}_{\nabla}$ be the associated local 
curvature forms on $M$. The homogeneous functions $\sigma_p$ are 
${\rm Ad}_{\mathfrak{gl}_n}$-invariant and the following are well-defined
via 
Chern-Weil theory.

\begin{definition} [Pontrjagin forms]
{\bf (i)} The 
\emph{Pontrjagin 
forms} are defined by 
$$
p_k(\mathcal{F}_{\nabla}):=\sigma_{2k}(\mathcal{F}_{\nabla})\;.
$$
\item {\bf (ii)} The \emph{total Pontrjagin  form} is the sum
$$
p(\mathcal{F}_{\nabla}):=p_1(\mathcal{F}_{\nabla})+p_2(\mathcal{F}_{\nabla})+
p_3(\mathcal{F}_{\nabla})+\hdots\;.
$$
\end{definition} 
The forms $p_n(\mathcal{F}_{\nabla})$ can also be obtained as the $2n$-th Chern form 
of the corresponding  complexified bundle. From the splitting principle, one immediately sees the integrality of  the integrals
$$
\frac{1}{2\pi }\oint_cp_n(\mathcal{F}_{\nabla})\in \ZZ
$$ 
for every cycle $c:\Delta^n\to M$. Hence, each $p_n$ is in the image of the canonical map $H^{4n}(M;\ZZ)\to H^{4n}(M;\RR)$.

\begin{remark}[Realizability of forms by Pontrjagin/Chern characters]
One can ask which differential form can arise as the Pontrjagin 
character form of some real vector bundle. Some partial work along these lines has been done for complex K-theory.
\begin{enumerate}[{\bf (i)}]
%

\vspace{-2mm}
\item Every even exact form (and odd exact form for $\op{K}^1$) arises as the component of a Chern character form \cite{TWZ1} (Corollary 2.7 and Corollary 5.8). See also \cite{PT} for the case of Hermitian holomorphic vector bundles. 
Earlier computational results on the image of the Chern character appear in \cite[Theorem 5.4]{Da}

\vspace{-2mm} 
\item 
For any bundle over a closed Riemannian manifold after stabilizing, there is a unitary connection
on the bundle whose Chern-Weil form is the harmonic representative of the Chern character
of the bundle \cite[Corollary 5.1]{SSu1}.

\vspace{-2mm} 
\item 
From \cite{SSu3} one has the following results involving the Todd genus Td. 
A cohomology class $c$ in $H^{\rm even}(X, \QQ)$ 
is the Chern character of a complex bundle over $X$ if and only 
for every closed even-dimensional stable almost complex mapping to $X$, 
$V \xrightarrow{f}X$, the integral 
$\int_X f^* c {\rm Td} V$
is an integer. A similar statement holds for the transgressed $\ch$ in $U$, 
odd-dimensional closed stable almost complex structures in
$X$, elements in $H^{\rm odd}(X, \QQ)$ and maps $X\to U$.
\end{enumerate}
We will address the image of the Pontrjagin character form using the Atiyah-Hirzebruch spectral sequence
(AHSS) in Section \ref{Sec-AHSS-KOhat}. Our techniques work equally well in the complex case and can be used to give a characterization of Chern character forms. In the holomorphic setting and for holomorphic K-theory, this question is central to Hodge theory. Indeed, one version of the Hodge conjecture asserts that all harmonic forms of type $(p,p)$, taking rational values on cycles, can be realized as a component of the Chern character of some holomorphic vector bundle. We wonder if our spectral sequence detects all such obstructions on forms, or whether there are additional obstructions in moving to the holomorphic category.
\end{remark}

\medskip
As 2-torsion is central to KO-theory, we include the following discussion
which seems to have been somewhat ignored in the literature. 

\begin{remark}[Whitney sum formula for Pontrjagin classes]
The total Pontrjagin class $p(E \oplus E')$ of a Whitney sum is congruent to 
$p(E)p(E')$ modulo elements of order 2, i.e.,
$2\big(p(E \oplus E')- p(E) p( E')\big)=0$. The proof 
(e.g., from \cite[Theorem 15.3]{MiS}) simply ignores the 
2-torsion arising from the odd Chern classes altogether. In fact, 
keeping the 2-torsion is not so bad to identify.  This has been worked out in \cite{Th62}\cite{Br}. 
Let $E$ be a real $n$-dimensional vector bundle over a paracompact 
space $X$. Associated with $E$ are the Pontrjagin classes
$$
p(E)= 1 + p_1(E) + \cdots + p_q(E)    
$$
where $p_j(E)\in H^{4j}(X; \ZZ)$ and $q=[n/2]$ the largest integer in $n/2$. 
For $E_\CC=E \otimes \CC$ the complexification of $E$, the homogeneous 
Pontrjagin classes $p_i(E)$ are defined by 
$$
p_i(E)=(-1)^ic_{2i}(E_\CC) \qquad 1 \leq i \leq [n/2]\;.
$$
In order to keep track of all degrees, one sets
$$
P_{2k}(E)=(-1)^q c_k(E_\CC) \qquad q=[k/2], \; 1 \leq k \leq 2n\;,
$$
so that $p_j(E)=P_{4j}(E)$. The classes 
$P_{4k+2}(E)$, for $k \geq 0$, are called the torsion Pontrjagin classes of 
$E$ as they satisfy $2P_{4k+2}(E)=0$. They are given by \cite{Th62}
$$
P_{4k+2}(E)= \beta_2 (w_{2k}(E) w_{2k+1}(E))=
p_k(E) \beta_2 w_1(E) + (\beta_2 w_{2k}(E))^2=\beta_2 Sq^{2k}w_{2k+1}(E)\;,
$$
where $\beta_2$ is the Bockstein corresponding to the mod reduction. 
The 2-torsion correction to the multiplicative property of the Pontrjagin 
classes
$$
D_k(E, E')=p_k(E \oplus E') - \sum_{i+j=k}p_i(E) p_j(E')
$$

\vspace{-2mm} 
\noindent is given by 
$$
D_k(E, E')=\sum_{i+j=k-1} \beta_2 [w_{2i}(E) w_{2i+1}(E')]
\beta_2 [w_{2i}(E') w_{2i+1}(E)]\;.
$$
One observes that this correction is, however, trivial in many relevant instances, especially in applications:
\begin{itemize}
\item If $E$ and $E'$ are both orientable bundles, meaning 
$w_1(E)=0=w_1(E')$, then $D_1(E, E')=0=D_2(E, E')$.  
\item If $E$ and $E'$ are Spin bundles then 
$D_k(E, E')=0$ for $1 \leq k \leq 6$. 
\end{itemize}
\end{remark}

\subsection{Orientation and genus}
\label{Sec-orgen}

In this section we discuss orientation and the genus associated with KO-theory and which will be useful for us in later constructions. We begin with the A-genus,
which takes values in KO-theory for Spin manifolds. 
%



\medskip
\noindent {\bf The $\widehat{A}$-genus.} 
By Borel-Hirzebruch \cite{BH}, if $M$ is a Spin manifold, then 
$\widehat{A}(M)$ is an integer. The reason for this is given by 
Atiyah and Singer via the interpretation as the index of a Dirac operator
(see \cite{AH-diff}\cite{ASi} and \cite{LM} for exposition).
By using families of operators, the index of a family operators is 
now a difference of vector bundles, which is a KO-element. Then 
the index provides a ring homomorphism 
$$
\widehat{A}: \op{M}{\rm Spin}_* \longrightarrow \op{KO}_*
$$
This can be lifted to a morphism of spectra, as follows; see the nice summary in \cite{Kn}. 
As orientation is a stable operation, it is enough to construct a Thom class for Spin bundles $E \to X$ of dimension $8k$
and recover the other case by suspension and shifting degrees. 
Let $P_{\rm Spin}(E)$ be the associated principal bundle. 
Using the periodicity element $\beta$ in the Grothendieck group of 
$C\ell_8$-modules, one forms the bundle 
$$
E=P_{\rm Spin}(E) \times_{{\rm Spin}(8k)} \beta^k\;.
$$
Over the disk bundle $D(E)$ of $E$, there is a map 
$$
\xymatrix@R=1.5em{
\pi^*E^0  \ar[rd] \ar[rr]^\sigma && \pi^*E^1 \ar[ld]
\\
& D(E) \ar[d]^\pi& 
\\
& X\;, &
}
$$
defined by $\sigma_{(p, v)}(x)=v\cdot x$. Since $\sigma$ restricted to the 
sphere bundle is an isomorphism, it defines a difference class 
$[\pi^*E^1, \pi^*E^0, \sigma]\in \op{KO}(D(E), S(E))\cong \widetilde{\op{KO}}({\rm Th}(E))$, 
which is $\op{KO}$-Thom class for $E$. The construction is functorial and multiplicative,
and so can be refined to a map of ring spectra $\op{M}{\rm Spin} \to \op{KO}$ lifting the $A$-genus. This construction is essentially the same as the one discussed in Example \ref{Ex-LMbasic}.

\medskip
Let $E$ be a real vector bundle over $X$ with a $\op{Spin}(n)$
structure, $n=\dim E=8k$. For such bundles there is a Thom
isomorphism in KO-theory
$$
\Phi: \op{KO}^*(X) \longrightarrow \widetilde{\op{KO}}({\rm Th}(E))\;.
$$
Let $\Phi_H: H^*(X; \QQ) \to \widetilde{H}^*({\rm Th}(E); \QQ)$
be the Thom isomorphism in cohomology, which is uniquely determined by 
the orientation of the bundle $E$. Then we have the compatibility
\(\label{arfcmpc}
\widehat{A}(E)^{-1}=\Phi_H^{-1} \op{Ph}(\Phi(1))\;.
\)
This compatibility immediately implies the Grothendieck-Riemann-Roch theorem for $\op{KO}$-theory (see 
\cite{AH-diff}). Indeed, let $f:W\to M$ be a proper map between Spin oriented manifolds and for simplicity, 
assume $M$ to be connected and that all connected components of $W$ have the same dimension. 
Set $n={\rm dim}(M)-{\rm dim}(W)$ for the codimension. Choose a fiberwise embedding 
$i:W\into M\times \RR^k$ for sufficiently large $k$. We can extend $i$ to a tubular neighborhood 
$\mathcal{N}$ of the embedding and this gives rise to a map of Thom spaces
$$
c(i):M_+\wedge S^k\longrightarrow {\rm Th}({\mathcal{N}})\;.
$$
Using the Thom isomorphism and the rationalized Pontrjagin character, we 
form the diagram
$$
\xymatrix{
\op{KO}^{*}(W)\ar[d]^-{\rm Ph}\ar[r]_-{\cong}^-{\Phi} & \widetilde{\op{KO}}^{*-n+k}({\rm Th}({\mathcal{N}}))\ar[d]^-{\rm Ph}\ar[r]^-{c(i)^*} & \op{KO}^{*-n}(M)\ar[d]^-{{\rm Ph}}
\\
H^{*}(W;\QQ[\alpha,\alpha^{-1}])\ar[r]^-{\Phi_{H}}_-{\cong} & \widetilde{H}^{*-n+k}({\rm Th}({\mathcal{N}});\QQ[\alpha,\alpha^{-1}])\ar[r]^-{c(i)^*}  & H^{*-n}(M;\QQ[\alpha,\alpha^{-1}])\;.
}
$$
Now for $V\to W$ a real vector bundle on $W$, consider the two horizontal compositions
$f_!:=c(i)^*\circ \Phi$ and $f_*:=c(i)^*\circ \Phi_{H}$. Then, from the compatibility \eqref{arfcmpc} and using the projection formula for $\Phi_{H}$, we immediately have the formula
\begin{align*}
f_*({\rm Ph}(V)\cup \hat{A}(\mathcal{N})) &=
f_*({\rm Ph}(V)\cup \Phi^{-1}_{H}{\rm Ph}(\Phi(1))) 
\\&= 
f_*\Phi_{H}^{-1}(p^*{\rm Ph}(V)\cup {\rm Ph}(\Phi(1))) \nonumber
\\
&= {\rm Ph}(c(i)^*(p^*V\cdot \Phi(1)))
\\&={\rm Ph}(f_!(V))\;.
\end{align*}
By definition, the A-genus satisfies $\hat{A}(E\oplus F)=\hat{A}(E)\cup \hat{A}(F)$. Since $[\mathcal{N}]=f^*[TM]-[TW]$ as classes  in $\op{KO}$, we have $\hat{A}(\mathcal{N})^{-1}=\hat{A}(TW)\cup f^*\hat{A}(TM)^{-1}$. Hence,
\begin{align*} 
f_*({\rm Ph}(V)\cup \hat{A}(TW))\cup \hat{A}(TM)^{-1}&=f_*({\rm Ph}(V)\cup \hat{A}(TW)\cup f^*\hat{A}(TM)^{-1})
\\& ={\rm Ph}(f_!(V))
\end{align*} 
which is the Grothendieck-Riemann-Roch theorem for $\op{KO}$ (see \cite{AH-diff}).


%

\medskip
\noindent {\bf The ABS orientation:} Let $\pi: E \to X$ denote a real vector bundle over a compact space $X$. 
Via the projection $\pi: E \to X$,  the compactly supported version $\op{KO}^{-*}_{\rm cpt}(E)$ is canonically 
a $\op{KO}^{-*}(X)$-module. A class $u \in {\rm KO}_{\rm cpt}(E)={\rm KO}^0_{\rm cpt}(E)$
is a \emph{KO-theory orientation} for the bundle $E$ if 
${\rm KO}^{-*}_{\rm cpt}(E)$ is  a free $\op{KO}^{-*}(X)$-module with generator $u$. 
An extensive modern discussion of KO-orientation can be found in \cite{AHR}. 

\medskip
One way to approaching orientability of bundles with respect to KO-theory is via the one introduced in \cite{ABS}. Another approach 
is to consider the Atiyah-Hirzebruch spectral sequence (AHSS) \cite{AH}; see \cite{Kn}
for a nice description, which we follow. Starting with
a bundle $E \to X$ which is orientable with respect to integral cohomology, 
we have a Thom class $\mu_E^\ZZ$. 
The Thom isomorphism gives an isomorphism of $E_2$-pages
$$
\xymatrix{
H^p(X; \pi_{-q}\op{KO}) 
\ar@{=>}[r] \ar[d] & 
\op{KO}^{p+q}(X)
\ar@{-->}[d]
\\
\widetilde{H}^p(X^E; \pi_{-q}\op{KO}) 
\ar@{=>}[r] & 
\widetilde{\op{KO}}^{p+q}(X^E) 
}
$$ 
induced by cup product with the Thom classes
$\mu^\ZZ_V$ and $\mu_E^{\ZZ/2}$. When $n=\dim E$, the
Leibniz rule implies that the first obstruction to the existence of 
the dashed arrow is 
$$
d_2^{n, -q} \mu_E^{\pi_q KO}=0\;.
$$
One way to construct the AHSS is via the Postnikov tower of  KO, 
in which case the differentials are exactly the $k$-invariants. In 
particular, $d_2^{n, -1}$ is a nontrivial  stable 
cohomology operation $H \ZZ/2 \to \Sigma^2 H \ZZ/2$,
i.e., given by $Sq^2$. Hence a necessary condition for 
orientability is 
$$
w_2(E):= Sq^2 \mu_E^{\ZZ/2}=0\;.
$$
That is, there is a nullhomotopy of the composite in the 
following diagram, hence a lift to the homotopy fiber
$$
\xymatrix{
&& B{\rm Spin}(n) \ar[d] && 
\\
X \ar@{-->}[rru]^E \ar[rr] && B{\rm SO}(n)  \ar[r]^{w_2} & K(\ZZ/2, 2)\;.
}
$$
One can show that this is also sufficient, so that a vector bundle is 
KO-orientable if and only if it is a Spin bundle -- see \cite{ABS} for the hard work. 
The following example 
provides a concrete geometric description of $KO$-orientation for 
a bundle $\pi:E\to X$ which we will use later.

\begin{example}[Spinorial description of KO-orientation]\label{Ex-LMbasic}
See \cite[Example C.4]{LM}.
Let $E=X \times \RR^{8m} \xrightarrow{\pi} X$  be a trivialized Riemannian vector
bundle and $X$ compact. Using the Clifford algebra, define 
$$
u=[\S^+, \S^- ; \mu]\in \op{KO}_{\rm c}(E)\;,
$$
where $\S=\S^+ \oplus \S^-=\pi^*\S (E)$ is the irreducible real graded $C\ell_{8m}$-module
(extended trivially over $E$) and where 
$
\mu_{x, v}(\varphi)=v \cdot \varphi
$
is given by Clifford multiplication. By Bott periodicity, $u$ is a KO-orientation 
for $E$. This can be generalized to nontrivial bundles \cite[Theorem C.9]{LM}.
Let $\pi: E \to X$ be a real $8m$-dimensional bundle with a Spin structure, over
a compact space $X$. Consider the class
$$
\mathbf{S}(E)=[\pi^*\S^+(E), \pi^*\S^-(E); \mu] \in \op{KO}_{\rm c}(E)\;,
$$
where $\S(E)=\S^+(E) \oplus \S^-(E)$ is the irreducible graded real spinor bundle of $E$ and for an element $e\in E$, the map $\mu$ is defined fiberwise by $\mu_e(\varphi)=e \cdot \varphi$ via Clifford multiplication.
Then $\mathbf{S}(E)$ is a KO-theory orientation 
for $E$. In particular, the map $i_!: \op{KO}(X) \to \op{KO}_{\rm c}(E)$ given by 
$$
i_!(a)=(\pi^* a) \cdot \mathbf{S}(E)
$$
is an isomorphism. 
\end{example}


\medskip
\noindent {\bf The index and the Thom class}. The Dirac operator on $X$ has an index in $\op{KO}^m({\rm pt})$, given
by $f_!(1)$, where $f: X \to {\rm pt}$ is the terminal map. Depending on 
$m$, one has  \cite[Sec. 4.2]{Hit}
\begin{align*} 
\op{KO}^{-(8m+1)}({\rm pt}) &\cong \ZZ/2    \hspace{12mm} f_!(1)= {\rm dim}_{\CC} ~h_D \; \text{mod}\; 2,
\\
\op{KO}^{-(8m+2)}({\rm pt})& \cong \ZZ/2    \hspace{12mm} f_!(1)= {\rm dim}_{\HH} ~h_{D} \; \text{mod}\; 2,
\\
\op{KO}^{-(8m+4)}({\rm pt}) &\cong \ZZ    \hspace{16mm} f_!(1)=\tfrac{1}{2}\widehat{A}(X)[X],
\\
\op{KO}^{-8m}({\rm pt}) &\cong \ZZ    \hspace{16mm} f_!(1)=\widehat{A}(X)[X],
\end{align*} 
where $D$ is the dirac operator on the spinor bundle and $h_D={\rm ker}(D)$ is the space of harmonic spinors. The $\alpha$-invariant is defined as $\alpha(X)=f_!(1) \in \op{KO}^{-n}(\ast )$ for a Spin 
manifold of dimension $n$. This is related to bilinear forms -- 
applications are given in \cite[Sec. 5.2.2]{pf}. 

\medskip
To connect to the index for twisted Dirac operators, let $X$ be an $n$-dimensional compact Riemannian Spin manifold 
and let $E$ be a real vector bundle over $X$ with an orthogonal 
connection. Then the spinor bundle $\S$ can be twisted by $E$ and the index
of the corresponding Dirac operator $D_E$ on $\S \otimes E$  is
$$
{\rm ind}_E(X)={\rm ind}_n(D_E) \in \op{KO}^{-n}({\rm pt})\;.
$$
When $n\equiv 1$ or 2 (mod 8), let 
$h_E\equiv \ker D_E$ denote the space of harmonic $E$-valued spinors. Then 
(see \cite[Theorem 7.13]{LM})
$$
{\rm ind}_E(X)=
\left\{
\begin{array}{lll}
\dim_\CC h_E  \; ({\rm mod}\;2) && {\rm if}\; n \equiv 1 \; ({\rm mod}\;8),
\\
\dim_{\mathbb{H}} h_E \;  ({\rm mod}\;2) && {\rm if}\; n \equiv 2 \; ({\rm mod}\;8),
\\
\tfrac{1}{2}\{\ch(E) \cdot \widehat{A}(X) \}[X] && {\rm if}\; n \equiv 4 \; ({\rm mod}\;8),
\\
\{\ch(E) \cdot \widehat{A}(X) \}[X] && {\rm if}\; n \equiv 0 \; ({\rm mod}\;8).
\end{array}
\right.
$$
The topological interpretation of the index is as follows: 
$X^n$ can be embedded in $S^{n + 8k}$ for sufficiently large $k$, and the normal 
bundle $\N$ of $X$ can be identified with a tubular neighborhood of the embedding. 
The Spin structure on $X$ determines a unique Spin structure on $\N$. Let $\S^\pm(\N)$
denote the canonical real spinor bundles of $\N$ of dimension $8k$. Via the 
projection $\pi: \N \to X$ one can lift $\S^\pm(\N)$ to the total space of $\N$.
At each nonzero vector ${\bf n} \in \N$ consider the isomorphism 
$\mu_n: \pi^*\S^+(\N) \xrightarrow{\cong} \pi^* \S^-(\N)$ given by the 
Clifford multiplication by ${\bf n}$. Then the difference element 
$$
\tau_{\N}\equiv [\pi^*\S^+(\N) , \pi^*\S^-(\N); \mu]
$$
represents a class in the relative KO-group $\op{KO}(\N, \N - X)$, where
$X \subset \N$ is the zero-section, and is the KO-theory 
\emph{Thom class} of $\N$.
Given a real bundle $E$ over $X$, one can consider the class
$$
\tau_{\N}(E)=\tau_\N \cdot [\pi^*E] \in \op{KO}(\N, \N - X)\;.
$$
Since $\N$ is embedded as a domain $\N \subset S^{n +8k}$, there is an excision 
isomorphism $j: \op{KO}(\N, \N - X) \cong \op{KO}(S^{n + 8k}, S^{n + 8k} -X)$. Composing this with 
the natural map $i: \op{KO}(S^{n + 8k}, S^{n + 8k} -X)\to \widetilde{\rm KO}(S^{n + 8k})$ and applying
Bott periodicity $\beta: \widetilde{\op{KO}}(S^{n + 8k}) \cong \widetilde{\rm KO}(S^{n })\equiv \op{KO}^{-n}({\rm pt})$,
one obtains a class 
$$
\widehat{A}_E(X)= \beta \circ i \circ j (\tau_\N (E)) \in \op{KO}^{-n}({\rm pt})\;.
$$
The Atiyah-Singer index theorem gives ${\rm ind}_E(X)=\widehat{A}_E(X)$. 

\medskip We will consider the corresponding orientations and genera in the differential case in Section \ref{Ch-diffKO}.

\subsection{Relation to the complex theory}
\label{Sec-cplx}

Here we consider some basic relations between KO-theory and K-theory. We will refine this discussion in 
Section \ref{Sec-cplx-diff} and towards the end of Subsection \ref{Sec-cyclemap}. 

\medskip
\noindent {\bf Operations on vector bundles.}
The following operations on vector bundles are inherited from corresponding ones on the vector space fibers:
\begin{enumerate}
\item  {\bf Conjugation:}  Given a complex vector bundle $E$ we can take the conjugate bundle 
$\overline{E}=\tau^{-1}(E)$.
\item {\bf Complexification: } Given a real vector bundle $E$ we can define a complex vector
bundle $c(E)=E_\CC=E \otimes_\RR \CC$.
\item {\bf Realification: } Given a complex vector bundle $E$ we can forget the complex structure, 
obtaining a real vector bundle $r(E)$. 
\end{enumerate} 
These operations define respective natural transformations at the level of the corresponding K-theories
$$
\tau^{-1}: \widetilde{\op{K}}^* \longrightarrow \widetilde{\op{K}}^*\;, \qquad
c: \widetilde{\op{KO}}^* \longrightarrow \widetilde{\op{K}}^*\;, \qquad 
r: \widetilde{\op{K}}^* \longrightarrow \widetilde{\op{KO}}^*\;, \qquad
$$
and even maps of spectra
$$
\tau^{-1}:  \op{K}  \longrightarrow  \op{K}\;, \qquad
c:  \op{KO} \longrightarrow  \op{K} \;, \qquad
r:  \op{K} \longrightarrow \op{KO}\;.
$$
They further satisfy the relations (see \cite[13.93]{Sw})
$$
rc=2\;, \qquad cr=1+ \tau^{-1}\;, \qquad \tau^{-1}c=c\;, \qquad 
r \tau^{-1}=r\;, \qquad \tau^{-1} \tau^{-1}=1\;, \quad r(x \cdot c(z))=r(x) z\;,
$$
for all $x \in \op{K}(X)$ and $z \in \op{KO}(X)$. 
Note that $c$ and $\tau$ are ring homomorphisms, hence multiplicative,
while  $r$ is only a homomorphism of the underlying abelian groups. 
This gives that $\op{KO}^*(X) \cong \{ \op{K}^*(X)\}^{\ZZ/2}$, the fixed elements 
under the conjugation automorphism of $K^*(X)$ \cite{Bott}. Note that this would
also produce the coefficients via the homotopy fixed point spectral sequence
corresponding to this action. 

\medskip
The map $c$ is realized by a fibration
$U/O \to BO \to BU$ and by Bott periodicity $U/O \simeq \Omega^{-1}BO$. 
This then leads to the following long exact sequence \cite{Bott}
\cite[Theorem III.5.18]{Ka}.
\footnote{
The existence of such a sequence is noted by Bott \cite{Bott}
but with a typo. Atiyah \cite{At-K} gave further clarifications writing it as
$\Sigma KO \xrightarrow{\eta} KO \xrightarrow{c} KI \xrightarrow{R} \Sigma^2KO$, but
without determining $R$ (see \cite[Remark 4.1.1]{BG}). In \cite[Lemma 5.11]{MT} the following 
essentially equivalent form
is given: With $\eta$ a generator of $\widetilde{KO}(S^1)\cong \Z/2$ and $g$
a generator of $\widetilde{K}(S^2)\cong \Z$, the following sequence is exact for $n\geq 1$:
$\widetilde{KO}(S^{n+1})\xrightarrow{\wedge \eta}
\widetilde{KO}(S^{n+2})\xrightarrow{(\wedge g)^{-1}\circ c}
\widetilde{K}(S^{n})\xrightarrow{r}
\widetilde{KO}(S^{n})\xrightarrow{\wedge \eta}
\widetilde{KO}(S^{n+1})$.
}

\medskip
\noindent 
{\bf The Bott sequence:} 
Real and complex K-theory are related by the Bott exact sequence
\(
\label{Bott-seq}
\xymatrix{
\cdots \ar[r] & \op{KO}^{* +1}(X) \ar[r]^{\cdot \eta} & \op{KO}^*(X) \ar[r]^{\chi} &
\op{K}^{*+2}(X) \ar[r]^{r} & \op{KO}^{*+2}(X) \ar[r] & \cdots 
}
\)
where the homomorphism $\cdot \eta$ is multiplication by $\eta$ and $\chi$ is complexification followed
by multiplication  by  the Bott generator $u^{-1}$.
Note that for any element $g \in \op{K}^*(X)$, the difference $g - \overline{g}$ lies in the kernel
of $r$, hence in the image of $\chi$, and 
$\chi(r(ug))=g - \overline{g}$. Several consequences follow: 

\begin{itemize}
\item This immediately leads to $\op{KO}(S^{8n}) \overset{c}{\cong} \op{K}^*(S^{8n})$. 
\item The product of any element of $\op{KO}^*(X)$ with $\eta$ is zero if and only if the 
element is in the image of the realification homomorphism. 
\item Whenever $\op{KO}^*(X)$ is a free abelian group, 
the sequence shows that $r$ is an epimorphism and $c$ is a monomorphism. 
\item Taking $X=S^0$ and using the fact that the first Chern class $c_i(-): \op{K}^0(S^2) \to H^2(S^2)$ 
is an isomorphism that commutes with conjugation, so that $c_i(\overline{u})=\overline{c_i(u)}$, 
leads to the following characterization of the image of $u^i \in \op{K}^{-2i}(S^0)$ under 
the realification homomorphism
$$
r(u^i)=\left\{
\begin{array}{lll}
2 \beta^{\tfrac{i}{4}}  &&   i \equiv 0  \; ({\rm mod}\; 4), \\
\eta^2 \beta^{\tfrac{i-1}{4}} & &   i \equiv 1 \; ({\rm mod}\; 4), \\
\alpha \beta^{\tfrac{i-2}{4}}  &&   i \equiv 2 \; ({\rm mod}\; 4), \\
0  &&   i \equiv 3 \; ({\rm mod}\; 4),
\end{array}
\right.
$$
as well as the images of $\eta, \alpha, \beta$ under complexification
$$
c(\eta)=0\;, \qquad c(\alpha)=2u^2\;, \qquad c(\beta)= u^4\;. 
$$
\end{itemize}
We will make explicit use of the Bott sequence in identifying differentials in the AHSS 
for KO-theory (see Proposition \ref{diff-Ko-ahss}).

\medskip
There is one further relationship between K-theory and KO-theory after inverting 2:
$$
\op{KO}^0(X) \otimes_\ZZ \ZZ[\tfrac{1}{2}] \cong \op{K}^0(X) \otimes_\ZZ \ZZ[\tfrac{1}{2}]\;,
$$
which highlights the importance of 2-torsion. 
We will discuss KO-theory with coefficients further, 
first in relation to orientations and genera towards the end of Section \ref{Sec-orgen}
and then in Section \ref{Sec-Adams} in the context of Adams operations.

\begin{remark}[Involutions at the level of spectra]
The involutions can be lifted to the level of spectra (see \cite[Chapter VI, Lemma 3.3]{Ru}
and can be described explicitly as follows.  
A map ${\bf c}=\{ {\bf c}_k\}_{k\in \ZZ}$ of $\Omega$-spectra that represents the 
complexification $c: \op{KO}^*(X) \to \op{K}^*(X)$ is given in the first few levels
by the following (see \cite{Wat1}\cite{Wat}):
\begin{enumerate}
\item ${\bf c}_0=BC \times 1: \op{BO} \times \ZZ \to \op{BU} \times \ZZ$. 

\item ${\bf c}_1=\xi_c: \op{U}/\op{O}  \to \op{U}$.   Let $(G, \sigma)$ be a symmetric pair, 
where $G$ is a Lie group and $\sigma$ is an 
involutive automorphism of $G$. Let $G^\sigma$ be the fixed point 
subgroup of $\sigma$. Then we have a map 
$\xi: G/G^\sigma \to G$ defined by 
$\xi(x G^\sigma)=x \sigma(x)^{-1}$ for $xG^{\sigma}\in G/G^\sigma$. 
For $G=\op{U}$ with the appropriate involution, one has 
$\xi: \op{U}/\op{O} \to \op{U}$. 

\item ${\bf c}_2=(j_q, 0): \op{Sp}/\op{U} \to \op{BU} \times \ZZ$.
The map $q: \op{U} \to \op{Sp}$ leads to the following fiber sequence 
$\op{U} \xrightarrow{q} \op{Sp} \xrightarrow{\pi_q} \op{Sp}/\op{U} \xrightarrow{j_q} \op{BU}$. 
\end{enumerate} 
\end{remark}


We will work out the differential refinement of the above sequence later
in 
Subsection \ref{Sec-cplx-diff}.

\section{Differential ${\rm KO}$-theory}
\label{Ch-diffKO}

\subsection{Differential refinements of Pontrjagin classes/characters}
\label{Sec-Pont-class}

In this section, we introduce differential refinements of the Pontrjagin classes. General methods for 
constructing differential refinements of characteristic classes have been studied extensively from various 
points of view (e.g. \cite{CS}\cite{Urs}\cite{Bun}\cite{Bun10}\cite{Cech}). Our starting point 
will follow closely both \cite{Bun} and \cite{Cech} (see also \cite{SSS3}\cite{FSS1}\cite{FSS2}), combining elements of both perspectives. We will use the quasi-categorical treatment of $\infty$-categories due to Lurie \cite{Lur}.

\begin{remark}
To circumvent size issues, we will always identify the underlying set of a smooth manifold with some subset of $\RR$. We emphasize that this identification is not compatible with the smooth structure. 
\end{remark}

\begin{definition}[Category of vector bundles with connections on smooth manifolds] 
Let $\mathscr{M}{\rm an}$ denote the category of smooth manifolds and let $\mathscr{C}{\rm at}$ denote the $(\infty,1)$-category of (small) categories, obtained by localizing at categorical equivalences. Define the $\infty$-functors 
$$
{\rm Vect}:\mathscr{M}{\rm an}^{\rm op}\longrightarrow \mathscr{C}{\rm at}\;,
\qquad 
{\rm Vect}_{\nabla}:\mathscr{M}{\rm an}^{\rm op}\longrightarrow \mathscr{C}{\rm at}
$$
which associate to each smooth manifold $M$, the category whose objects are vector bundles (respectively, with connection) on $M$ and morphisms are the usual morphisms of vector bundles. The operation of pullback of vector bundles (respectively, with connection) makes both associations functorial (in the $(\infty,1)$-sense).
\end{definition}

The category of smooth manifolds can be topologized by taking good open covers as covering families. Since bundles 
(respectively, with connection) are glued from local data, both functors ${\rm Vect}$ and ${\rm Vect}_{\nabla}$ define smooth stacks. For each smooth manifold $M\in \mathscr{M}{\rm an}$, we can take the core (maximal sub-groupoid)
of the categories ${\rm Vect}(M)$ and ${\rm Vect}_{\nabla}(M)$, yielding groupoids. This operation preserves sheaves and we get corresponding sheaves of groupoids ${\rm Iso}({\rm Vect})$ and ${\rm Iso}({\rm Vect}_{\nabla})$. Passing to $\infty$-groupoids via the nerve operation (which also preserves sheaves) we get corresponding sheaves of $\infty$-groupoids. 

\begin{definition}[Classifying stack for real bundles]
We define the smooth stack $\BB{\rm O}(k)$ as the orbit stack (i.e., the stackification of the corresponding quotient in prestacks)
$$
\BB{\rm O}(k):=\ast/\!/C^{\infty}(-;{\rm O}(k))\;,
$$ 
and the smooth stack $\BB{\rm O}=\coprod_{k\in \NN}\BB{\rm O}(k)$,
where the maps $\BB{\rm O}(k)\to \BB{\rm O}(k+1)$ are induced by the obvious inclusions as block matrices.
\end{definition} 

\begin{definition}[Classifying stack for real bundles with connection]
We define the smooth stack $\BB{\rm O}(k)_{\nabla}$ as the orbit stack 
$$
\BB{\rm O}(k)_{\nabla}:=\Omega^1(-;\mathfrak{o}(k))/\!/C^{\infty}(-;{\rm O}(k))\;,
$$ 
where locally $g\in C^{\infty}(U;{\rm O}(k))$ acts on ${\mathcal A}\in \Omega^1(U;\mathfrak{o}(k))$ by gauge transformations ${\mathcal A}\mapsto g^{-1}{\mathcal A}g+g^{-1}dg\;.$ We define the smooth stack 
$\BB{\rm O}_{\nabla}:=\coprod_{k\in \NN}\BB{\rm O}(k)_{\nabla}$.
\end{definition}

We now summarize some well known relationships between these smooth stacks 
(see \cite{Cech}\cite{FSS1} \cite{FSS2}\cite{Urs}). For each $n\in \mathbb{N}$, there is an equivalence of smooth stacks ${\rm O}(n)$ bundles ${\rm O}(n)\text{-}\mathscr{B}{\rm un}
\simeq \BB {\rm O}(n)$, where ${\rm O}(n)\text{-}\mathscr{B}{\rm un}$ is the smooth stack of principal orthogonal bundles. The operation of taking associated vector bundles gives 
an equivalence $
{\rm Iso}({\rm Vect}^g)\simeq \BB{\rm O}\;,
$
where the stack ${\rm Vect}^g$ is the smooth stack of vector bundles equipped with a smooth fiberwise inner product \footnote{Note that the associated bundle comes equipped with a canonical inner product inherited from the standard inner product on $\RR^n$. A choice of metric is required to define the inverse construction.}. Similarly, we have an equivalence of smooth stacks ${\rm O}(n)\text{-}\mathscr{B}{\rm un}_{\nabla}\simeq \BB{\rm O}(n)_{\nabla}$ where the stack on the right is the moduli stack of principal bundles ${\rm O}(n)$-bundles equipped with connection. Again, taking associated bundles yields an equivalence
$
{\rm Iso}({\rm Vect}_{\nabla}^g)\simeq \BB{\rm O}_{\nabla}\;,
$ 
where the stack ${\rm Vect}_{\nabla}^g$ is that of vector bundles equipped with metric connections. 
At the level of the \emph{connected components} of these stacks, the data provided by the metrics become
redundant -- every vector bundle can be equipped with a (noncanonical) metric. 
Indeed, we have the following. 

\begin{proposition}[Stacky realization of vector bundles (with connection)]
We have isomorphisms of presheaves of commutative monoids
$$\pi_0{\rm Vect}^g(-)\cong \pi_0\map(-,\BB {\rm O})\;.$$
Adding connections, the following are also isomorphic:  
$$
\pi_0{\rm Vect}_{\nabla}^g(-)\cong \pi_0\map(-,\BB{\rm O}_{\nabla})\;.
$$
Moreover, the canonical map $\pi_0{\rm Vect}^g(-)\to \pi_0{\rm Vect}(-)$ which forgets the fiberwise metric is an objectwise isomorphism. 
\end{proposition}
\theproof
From the preceding discussion, we only need to prove the last assertion. Since every vector bundle admits a fiberwise metric and 
any isomorphism between two such bundles can be improved to an isometric isomorphism, the claim is immediate.
\endofproof

We now discuss the differential refinement of the Pontrjagin classes. Via pullback of vector bundles, the Pontrjagin classes $p_n\in H^{4n}(BO;\ZZ)$ give a natural transformation
$$
p_n:\pi_0{\rm Vect}(-)\longrightarrow H^{4n}(-;\ZZ)\;.
$$ 
Similarly, by Chern-Weil theory, the Pontrjagin forms give rise to natural transformations 
$$
p_n(\mathcal{F}_{(-)}):\pi_0{\rm Vect}_{\nabla}(-)\longrightarrow \Omega_{\rm cl}^{\rm 4n}(-)\;,
$$
where the sheaf on the right is the sheaf of closed differential forms of degree $4n$. Now we can ask whether there is a unique differential refinement $\hat{p}_n:\pi_0{\rm Vect}_{\nabla}(-)\to \widehat{H}^{4n}(-;\ZZ)$ which makes the diagrams 
\(\label{dfrfpcls}
\xymatrix{
\pi_0{\rm Vect}_{\nabla}(-)\ar[rr]^-{\hat{p}_n} \ar[d]^-{U}&&  \widehat{H}^{4n}(-;\ZZ)\ar[d]^-{\mathcal{I}}
\\
\pi_0{\rm Vect}(-)\ar[rr]^-{p_n} && H^{4n}(-;\ZZ)
}
\qquad \text{and} \qquad 
\xymatrix{
\pi_0{\rm Vect}_{\nabla}(-)\ar[d]^-{\mathcal{F}_{(-)}}\ar[rr]^-{\hat{p}_n} && \widehat{H}^{4n}(-;\ZZ)\ar[d]^-{\mathcal{R}}
\\
\Omega^2(-;\mathfrak{o})\ar[rr]^-{p_n(\mathcal{F}_{(-)})} && \Omega^{4n}_{\rm cl}(-)
}
\)
commute. Indeed, this is the case.
\begin{proposition}
[Differential $n$-th Pontrjagin class]
\label{pclprop}
There is a unique natural transformation $\hat{p}_n:\pi_0{\rm Vect}_{\nabla}(-)$ $\to \widehat{H}^{4n}(-;\ZZ)$ making the diagrams \eqref{dfrfpcls} commute. We call this transformation the $n$-th \emph{differential Pontrjagin class}.
\end{proposition}
\theproof
Fix a vector bundle $V\to M$ and a connection $\nabla$ on $V$. The Pontrjagin classes $p_n(V)\in H^{4n}(M;\ZZ)$ satisfy the compatibility $r(p_n(V))=[p_n(\nabla)]$, where $r:H^{4n}(M;\ZZ)\to H^{4n}(M;\RR)$ is induced by tensoring with $\RR$ and $p_n(\nabla)$ is the Pontrjagin form defined by Chern-Weil theory. The ambiguity in defining the differential refinement $\hat{p}_n(V,\nabla)$ is measured by the image if the map $a:H^{4n-1}(M;\RR)\to \widehat{H}^{4n}(M;\ZZ)$ (see Appendix, Proposition \ref{prop-setref}). Let $\op{BO}:=\coprod_{n\in \NN}\op{BO}(n)$ and $V:M\to \op{BO}$ be the classifying map. Since $\op{BO}$ is a countable CW complex, it admits an approximation by finite dimensional smooth manifolds. More precisely, for each $n\in \NN$ there is a smooth manifold $N$, a classifying map $W:N\to \op{BO}$ which is $r=\max\{n,{\rm dim}(M)\}+1$ connected and a smooth map $f:M\to N$ such that $f^*W\cong V$. Fix a connection $\nabla^{\prime}$ on $W$. Since the rational cohomology of $\op{BO}$ vanishes in degrees which are not a multiple of 4, it follows that the ambiguity in defining $\hat{p}_n(W,\nabla^{\prime})$ vanishes. Then by naturality, we are forced to define $\hat{p}_n(V,f^*\nabla^{\prime})=f^*\hat{p}_n(W,\nabla^{\prime})$, whence
$$\hat{p}_n(V,\nabla)=\hat{p}_n(V,f^*\nabla^{\prime})+a({\rm cs}(\nabla,f^*\nabla^{\prime}))\;.$$
This proves uniqueness, but it remains to show that this map is well defined and independent of the choices.

Let $W^{\prime}:N^{\prime}\to \op{BO}$ be another $r$-connected map and $g:M\to N^{\prime}$ be such that $g^*W^{\prime}\cong V$. Fix a connection $\nabla^{\prime\prime}$ on $W^{\prime}$. We must show that 
$$
\hat{p}_n(V,f^*\nabla^{\prime})+a({\rm cs}(\nabla,f^*\nabla^{\prime}))=
\hat{p}_n(V,g^*\nabla^{\prime\prime})+a({\rm cs}(\nabla,g^*\nabla^{\prime\prime}))\;.
$$
Choose a third manifold $N^{\prime\prime}$, a map $W^{\prime\prime}:N^{\prime\prime}\to \op{BO}$ and maps $h:N^{\prime}\to N^{\prime\prime}$, $k:N^{\prime}\to N^{\prime\prime}$ giving rise to the commutative diagram 
$$
\xymatrix@C=4em{
& N\ar[dr]^-{h} \ar@/^1pc/[drr]^-{W} &
\\
M\ar[ru]^-{f}\ar[dr]_-{g} \ar[dr] && N^{\prime\prime} \ar[r]^-{W^{\prime\prime}} & \op{BO}\;.
\\
& N^{\prime}\ar[ru]_-{k}\ar@/_1pc/[urr]_-{W^{\prime}} &
}
$$ 
Then we have
\begin{align*}
\hat{p}_n(V,f^*\nabla^{\prime})+a({\rm cs}(\nabla,f^*\nabla^{\prime}))  - &
\hat{p}_n(V,g^*\nabla^{\prime\prime})+a({\rm cs}(\nabla,g^*\nabla^{\prime\prime})) 
\\
&= f^*(\hat{p}_n(W,h^*\nabla^{\prime\prime\prime}))+f^*a({\rm cs}(\nabla^{\prime},h^*\nabla^{\prime\prime\prime}))+a({\rm cs}(\nabla,f^*\nabla^{\prime}))
\\
& \;\;  \;\;\;  \;\;\;-  g^*(\hat{p}_n(W^{\prime},k^*\nabla^{\prime\prime\prime}))+g^*a({\rm cs}(\nabla^{\prime\prime},k^*\nabla^{\prime\prime\prime}))+a({\rm cs}(\nabla,g^*\nabla^{\prime\prime}))
\\
&=  f^*(\hat{p}_n(W,h^*\nabla^{\prime\prime\prime}))-g^*(\hat{p}_n(W^{\prime},k^*\nabla^{\prime\prime\prime})) -a({\rm cs}(f^*h^*\nabla^{\prime\prime\prime},g^*k^*\nabla^{\prime\prime\prime}))\;.
\end{align*}
The last expression is zero by the homotopy formula (see Appendix). This shows that the map is well defined and independent of the choices made and naturality is easily verified. 
\endofproof

\begin{remark}[Real vs. complex]
The method of proof of the previous proposition closely follows the corresponding proof for the Chern classes in the complex case (see \cite[Theorem 3.42]{Bun}), adapted to the real setting. Using the formula for the Pontrjagin classes in terms of Chern classes (and their corresponding characteristic forms), one could deduce the previous proposition from the complex case. However, following our general methodology, we have chosen to work completely in the real setting. 
\end{remark}

Notice that the unique morphisms of presheaves $\hat{p}_n:\pi_0{\rm Vect}_{\nabla}\to \widehat{H}^{4n}(-;\ZZ)$ lift to morphisms of prestacks $\hat{p}_n:{\rm Vect}_{\nabla}\to \widehat{H}^{4n}(-;\ZZ)$, since such a lift is completely determined by its value on an isomorphism class. Using the identifications ${\rm Iso}({\rm Vect}^g)\simeq \BB{\rm O}$ and ${\rm Iso}({\rm Vect}^g_{\nabla})\simeq \BB{\rm O}_{\nabla}$, the Yoneda lemma implies that there is a unique homotopy class 
$$
\hat{p}_n\in \pi_0\map\big(\BB {\rm O}_{\nabla},\BB^{4n-1}U(1)_{\nabla}\big)
\cong 
\prod_{k\in \mathbb{N}}\pi_0\map\big(\BB{\rm O}(k)_{\nabla},\BB^{4n-1}U(1)_{\nabla}\big)\;,
$$ 
given by precomposing $\hat{p}_n$ with the canonical forgetful map ${\rm Vect}_{\nabla}^g(-)\to {\rm Vect}_{\nabla}(-)$. For a rank-$k$ bundle, post-composition with the canonical map $\mathcal{R}:\BB^{4n-1}U(1)_{\nabla}\to \Omega^{4n}$ gives a corresponding element of 
$\mathcal{R}(\hat{p}_n)\in \pi_0\map(\BB{\rm O}(k)_{\nabla},\Omega^{\rm 4n}_{\rm cl})$
which is identified with the Chern-Weil form associated to the $n$-th Pontrjagin-class. Post-composition with the canonical map $\mathcal{I}:\BB^{4n-1}U(1)_{\nabla}\to \BB^{4n}\ZZ$ produces the underlying topological class in $\pi_0\map(\BB{\rm O}(k)_{\nabla},\BB^{4n}\ZZ)\simeq \pi_0\map(\op{BO}(k),B^{4n}\ZZ)$. Summarizing:
\begin{lemma}[Differential Pontrjagin classes for metric connections]
Precomposition of $\hat{p}_n$ with the canonical map $\pi_0{\rm Vect}^g_{\nabla}(-)\to \pi_0{\rm Vect}_{\nabla}(-)$ determines an equivalence class of natural transformations
\(\label{pclcech}
\hat{p}_n:\BB{\rm O}_{\nabla} \longrightarrow \BB^{4n}U(1)_{\nabla}\;.
\)
\end{lemma}

\begin{remark}[Various approaches]
We could have constructed the class \eqref{pclcech} directly via Chern-Weil theory (see \cite{Cech}) and defined this class to be the Pontrjagin class. However, this definition is slightly more restrictive, since it requires the vector bundles to be equipped with a metric connection. Therefore, we have chosen the more general definition in Proposition \ref{pclprop}.
\end{remark}

Recall that differential cohomology admits a cup product operation  (see \cite{FSS1}\cite{FSS2})
$$
\cup_{\rm DB}:\BB^{n-1}U(1)_{\nabla}\times \BB^{m-1}U(1)_{\nabla}
\; \longrightarrow \;
\BB^{n+m-1}U(1)_{\nabla}\;,
$$
induced by the Deligne-Beilinson cup product. At the level of cohomology this gives the graded group 
$\widehat{H}^*(M;\ZZ)$, for each smooth manifold $M$, a graded commutative ring structure which
is natural with respect to pullback. We can consider rational differential cohomology in an analogous way, 
via the hypercohomology of the positively graded complex
$$
\xymatrix{
\hdots \ar[r] & 0\ar[r] &  \QQ \; \ar@{^{(}->}[r] & \Omega^0\ar[r]^-{d} & \Omega^1\ar[r]^-{d}& \hdots \ar[r]^-{d}& \Omega^n}\;.
$$
The Deligne-Beilinson cup product is defined on local sections via 
$$
\alpha\cup_{\rm DB} \beta:=\left\{\begin{array}{cc}
\alpha\beta & \vert \alpha\vert=n,
\\
\alpha\wedge d\beta & \vert \beta\vert=0,
\\
0 & \text{otherwise}.
\end{array}\right.
$$ 
This product gives $\widehat{H}^*(M;\QQ)$ the structure of a graded commutative ring. 

\medskip
Recall that the Pontrjagin character arises as an element of the ring 
$H^{*}(\op{BO}(n); \QQ)$, determined by the even
series $\sum_{i=1}^{[n/2]} 2 \op{cosh}(x_i)$. 
Given any formal power series in the Pontrjagin classes, we can define differential cohomology classes which refine these underlying classes in an obvious way -- replacing the Pontrjagin classes by their differential refinements and using the multiplicative structure of the Deligne-Beilinson cup product.  With expansion \eqref{Eq-Php} in mind, 
we thus have the following natural definition for the \emph{differential Pontrjagin character}.

\begin{definition}[Differential Pontrjagin character]\label{Def-DiffPontChar} Let $V\to M$ be a real vector bundle with connection $\nabla$. Consider the series
$$
g(x_1,\hdots,x_{[n/2]})=\sum_{i=1}^{[n/2]}2\cosh(x_i)\;.
$$
Let $g=\sum_{j} f_j$ be a presentation of $g$ in terms of the elementary symmetric functions, with $f_j$ the homogeneous component of degree $j$. We define the \emph{differential Pontrjagin character} via 
$$
\widehat{{\rm Ph}}(V,\nabla)=\sum_{j}f_{j}
\big(\hat{p}_1(V,\nabla),\hat{p}_2(V,\nabla),\hdots,\hat{p}_j(V,\nabla)\big)\;,
$$
where products are taken with respect to the Deligne-Beilinson cup product, i.e., 
$\hat{p}_1^2:=\hat{p}_1 \cup_{\rm DB}\hat{p}_1$ and so on. The first few terms are
$$
\widehat{{\rm Ph}}(V,\nabla):={\rm rank}+
\hat{p}_1(V,\nabla) +\tfrac{1}{12}(\hat{p}_1^2\big(V,\nabla)-2\hat{p}_2(V,\nabla)\big)+\hdots\;.
$$
\end{definition}
We will provide characterizations in Section 
\ref{Sec-Diff-Pont-Cha} below (see Proposition \ref{phchrfuni}).
Similarly, for the $\hat{A}$-genus, we provide the following.
\footnote{Another refinement $\hat{\hat{A}}$ in the context of differential character is provided in \cite{CS}
as the lift in $\hat{H}^{\rm odd}(M; \RR/\QQ)$ of the Chern-Weil representative $\hat{A}(TM, \nabla_{TM})$.}

\begin{definition}[Differential $A$-genus]\label{Def-diffA}
We define the refinement of the $\hat{A}$-genus as 
$$
\hat{\mathbb{A}}(V,\nabla)=
\sum_{j}K_j\big(\hat{p}_1(V,\nabla),\hdots ,\hat{p}_j(V,\nabla)\big)\;,
$$ 
where $K_j$ is the multiplicative sequence associated to
the formal power series 
$$
Q(z)=\frac{\sqrt{z}/2}{{\rm sinh}(\sqrt{z}/2)}\;.
$$
The first few terms are
$$
\hat{\mathbb{A}}(V,\nabla)=1-\tfrac{1}{24}\hat{p}_1(V,\nabla)+
\tfrac{1}{5760}\big(-4\hat{p}_2(V,\nabla)+7\hat{p}_1^2(V,\nabla)\big)+\hdots\;,
$$
and the products are Deligne-Beilinson cup products,  i.e., 
$\hat{p}_1^2:=\hat{p}_1 \cup_{\rm DB}\hat{p}_1$ and so on. 
\end{definition}
The following observation will be useful in later sections.
\begin{lemma}
[The differential A-genus as a unit]
\label{Lemm-a-roof}
The class $\hat{\mathbb{A}}(V,\nabla)$ is a unit in the ring $\widehat{H}^*(M;\mathbb{Q})$. Thus, the multiplicative inverse $\hat{\mathbb{A}}(V,\nabla)^{-1}$ is well-defined. 
\end{lemma}
\theproof
Since $M$ is finite-dimensional, $\hat{\mathbb{A}}(V,\nabla)=1+x$, with $x$ a nilpotent element. Hence
$$\frac{1}{1+x}=1-x+x^2-x^3+\hdots +(-1)^{k-1}x^k$$
is a multiplicative inverse. 
\endofproof

\begin{remark}[Refinements of classical genera]
Formally, it is clear that these classes refine their corresponding topological counterparts. 
Indeed, this follows immediately from the fact that $\mathcal{I}(\hat{p}_j)=p_j$ and $\mathcal{I}$ 
is a ring homomorphism. Differential form representatives of the $A$-genus have been used in 
\cite{Bu1} in discussing geometric index theorems. We will connect these `naive' definitions 
with more refined versions appearing at the level of spectra in Section \ref{Sec-or-gen}.
\end{remark} 
\subsection{Constructing differential KO-theory} 
\label{Sec-construct}

In this section we study two sheaves of spectra which deserve to be called differential refinement of ${\rm KO}$. One spectrum is the analogue of the algebraic $K$-theory spectrum in the smooth context, capturing the data of locally trivial vector bundles and connections. The other is the Hopkins-Singer construction \cite{HS} applied to $\op{KO}$, which captures the key ingredients in Chern-Weil theory. As a prerequisite to the algebraic approach we also give a sheaf of spectra which represents \emph{smooth} $\op{KO}$-theory. These theories are related by maps in the following ways
\(\label{shspdfko}
\xymatrix{
\mathbf{K}{\rm O}_{\nabla}\ar[rr]^-{\widehat{\rm cyc}}\ar[d]_-{\mathcal{F}} && 
\widehat{{\rm KO}} \ar[d]^-{\mathcal{I}}
\\
\mathbf{K}{\rm O}\ar[rr]_-{\rm cyc} && \op{KO}\;.
}
\)
The map $\mathcal{F}:\mathbf{K}{\rm O}_{\nabla}\to \mathbf{K}{\rm O}$ forgets the data of connections of vector bundles, while the map on the right $\mathcal{I}:\widehat{\op{KO}}\to \op{KO}$ is canonical map in the 
Hopkins-Singer model which topologically realizes the geometric data. We study the cycle map cyc and
its differential refinement ${\widehat{\rm cyc}}$
explicitly  later in Subsection \ref{Sec-cyclemap} and in Proposition \ref{csconvb}.

\medskip
In order to present these various theories in a systematic way, we need to briefly review some of the theory of sheaves of spectra. More details are found in the appendix, which closely follows \cite{BNV} and \cite{Urs}. The $\infty$-category of sheaves of spectra $\mathscr{S}{\rm h}_{\infty}(\mathscr{M}{\rm an};\mathscr{S}{\rm p})$ admits several adjoints relating it to the $\infty$-category of spectra. More precisely, we have a quadruple adjunction $\Pi\dashv \delta \dashv \Gamma \dashv \delta^{\dagger}$,
\(
\label{quad-adj}
\xymatrix{
\mathscr{S}{\rm h}_{\infty}(\mathscr{M}{\rm an};\mathscr{S}{\rm p})\ar@<-.5em>[rr]|-{\Gamma} 
\ar@<1em>[rr]|-{\Pi} && \mathscr{S}{\rm p} \ar@<1.2em>[ll]|-{ \delta^{\dagger}} \ar@<-.3em>[ll]|-{\delta} 
}.
\)
The functor $\delta:\mathscr{S}{\rm p}\into \mathscr{S}{\rm h}_{\infty}(\mathscr{M}{\rm an};\mathscr{S}{\rm p})$ takes a topological spectrum and forms the associated constant sheaf of spectra. The functor $\Gamma$ evaluates a sheaf of spectra on the point manifold and the functor $\Pi$ geometrically realizes, recovering a topological spectrum from a sheaf of spectra.
It will turn out that the only differences in the sheaves of spectra depicted in diagram \eqref{shspdfko}
are in the geometry. More precisely, the topological  
\footnote{\label{footgeom}also often (confusing) called `geometric' realization. This name arose from the need to 
distinguish between simplicial constructions and topological space constructions, and the latter was 
broadly and vaguely referred to as geometric. We believe, especially in light of recent interest in actual 
geometric constructions via differential cohomology, that the second classical construction should be referred 
to as \emph{topological} realization, because that is really what it is, and reserve geometric realization 
for differential contexts as the one we have (elsewhere) in this paper.} 
realization $\mathcal{I}$
of all maps there are equivalences. For any sheaf of spectra $\E$, the quadruple adjunction allows us to form the differential cohomology hexagon diagram (see the Appendix)
\(\label{spdggcoh}
\xymatrix@C=1pt @!C{
&{\rm fib}(\mathcal{I}) \ar[rd]\ar[rr] & & {\rm cofib}(j) \ar[rd] &
\\
\Sigma^{-1}{\delta \Pi}{\rm cofib}(j)\ar[ru]\ar[rd] & & {\E}\ar[rd]^-{\mathcal{I}} \ar[ru]^-{\mathcal{R}}& &  {\delta \Pi}{\rm cofib}(j)\;,
\\
& \delta\Gamma{\E}\ar[ru]^{j}\ar[rr] & & {\delta \Pi}{ \E}\ar[ru] &
}
\)
where both the right and left square are Cartesian, the two diagonal sequences are exact and the bottom and top sequences are exact. We will spell out what the various ingredients that go into this diagram look like for the sheaves of spectra depicted in diagram \eqref{shspdfko}.

\begin{remark}[Topological spectra as sheaves of spectra]
The functor $\delta:\mathscr{S}{\rm p}\into \mathscr{S}{\rm h}_{\infty}(\mathscr{M}{\rm an};\mathscr{S}{\rm p})$ is fully faithful. As should be clear from the context, whenever a topological spectrum appears as the source or target of a map of sheaves a spectra, the reader should assume that it is embedded via $\delta$.  
\end{remark}

The usual functorial constructions one can perform in the category of spectra carry through without difficulty to the category of sheaves of spectra. In particular, given a sheaf of spectra $\E$, one can form the infinite loop stack \footnote{Note that sheaves of spectra are defined via pointed smooth stacks, just as ordinary spectra are defined via pointed topological spaces.} $\Omega^{\infty}\E$, with the "loop stack" operation given by using the intrinsic 
homotopy theory of smooth stacks. Similarly, one defines the infinite suspension spectrum of a sheaf of spectrum in the usual way (see Appendix).

\medskip
We now turn our attention to constructing the two sheaves of spectra $\mathbf{K}{\op{O}}$ and $\mathbf{K}{\op{O}}_{\nabla}$. Recall the smooth stacks ${\rm Vect}$, ${\rm Vect}_{\nabla}$, $\BB {\rm O}$ and $\BB{\rm O}_{\nabla}$ from Section \ref{Sec-Pont-class}. It is well known that we have a natural bijective correspondence
$$
\pi_0{\rm Iso}({\rm Vect})(M)\cong \pi_0\map(M,\op{BO})
$$
between isomorphism classes of locally trivial vector bundles on $M$ and homotopy classes of maps to $B{\rm O}$. 
Note, however, that it is not true that we have an equivalence at the level of the corresponding $\infty$-groupoids (i.e.,
${\rm Vect}(M) \not\simeq \map(M,\op{BO})$). This fails for obvious 
reasons -- for example, ${\rm Vect}(M)$ is the nerve of a category (hence 1-truncated), while 
$\map(M,\op{BO})$ may have complicated homotopy type. In contrast, we saw in Section
\ref{Sec-Pont-class} that we have an equivalence of groupids 
\(
{\rm Iso}({\rm Vect}^g)(M)\simeq \map(M,\BB{\rm O})\;.
\)
Thus, we see that working in the smooth setting has the advantage of capturing the entire categorical structure of vector bundles (with isomorphisms as morphisms) and not just their underlying set of homotopy classes. 

\medskip
Let $C\mathscr{M}{\rm on}(\infty\mathscr{G}{\rm pd})$ denote the sub $\infty$-category of commutative monoids in $\infty\mathscr{G}{\rm pd}$ and let $C\mathscr{G}{\rm rp}(\infty\mathscr{G}{\rm pd})$ be the subcategory of $\infty$-abelian groups (i.e. connected spectra). The inclusion $i:C\mathscr{G}{\rm rp}(\infty\mathscr{G}{\rm pd})\into C\mathscr{M}{\rm on}(\infty\mathscr{G}{\rm pd})$ admits a left adjoint $\mathscr{K}$ (see the Appendix) which can be thought of as taking the group completion. The functor $\mathscr{K}$ prolongs to a functor between presheaves of $\infty$-monoids and $\infty$-abelian groups. We have the following natural definition in harmony with Quillen's definition for the algebraic $K$-theory spectrum \cite{Qu}. 

\begin{definition}[Smooth KO-spectrum]
Let $L:\mathscr{P}\mathscr{S}{\rm h}_{\infty}(\mathscr{M}{\rm an};\mathscr{S}{\rm p})\to \mathscr{S}{\rm h}_{\infty}(\mathscr{M}{\rm an};\mathscr{S}{\rm p})$ denote the stackification functor (left adjoint to the inclusion $i$). We define the smooth $\op{KO}$-spectrum as the connected sheaf of spectra defined by 
$$
\Omega^{\infty}{\bf K}{\rm O}:=L\circ \mathscr{K}\Big(\coprod_{n\in \NN}\BB{\rm O}(n)\Big)\;.
$$
\end{definition}

As pointed out in \cite{BNV}, application of the sheafification functor $L$ is needed, since $\mathcal{K}$ 
is a left adjoint and hence does not preserve sheaves in general. Given the discussion in Section
\ref{Sec-Pont-class}, it is immediate from the definition
that we have the following. 
\begin{lemma}[Equivalence of smooth and topological KO-theory]
There is  a natural isomorphism
$$
{\bf K}{\rm O}(M)\cong {\rm Gr}\big({\rm Iso}({\rm Vect}^g)(M)\big)\cong \op{KO}(M)\;.
$$ 
\end{lemma}
Thus, as a cohomology theory, the smooth variant of $\mathbf{K}{\rm O}$-theory contains 
the same information as topological $\op{KO}$-theory. At the level of spectra, we have the following. 

\begin{proposition}
[Topological realization of  $\mathbf{K}{\rm O}$]
\label{smkotpko}
The unit $\eta:\mathbb{1}\to \delta\Pi$ of the adjunction $\Pi\dashv \delta$ gives rise to a morphism of sheaves of spectra
$$
{\rm cyc}:=\eta_{\mathbf{K}{\rm O}}:\mathbf{K}{\rm O}\longrightarrow \op{KO}\simeq \delta \Pi(\mathbf{K}{\rm O}).
$$
\end{proposition}
\theproof
The monoidal structure on $\coprod_{n\in \mathbb{N}}\mathbf{B}{\rm O}(n)$ corresponds to 
the direct sum of vector bundles, and pullback by a smooth map preserves this operation. Hence the diagram 
\(\label{dbomngrp}
\coprod_{n\in \mathbb{N}}\mathbf{B}{\rm O}(n):\mathscr{M}{\rm an}^{\rm op}
\longrightarrow \infty\mathscr{G}{\rm pd}
\)
factors through $C\mathscr{M}{\rm on}(\infty\mathscr{G}{\rm pd})$. At the level of presheaves, $\Pi$ is presented by the colimit operation. We claim that the colimit over diagram \eqref{dbomngrp} in $C\mathscr{M}{\rm on}(\infty\mathscr{G}{\rm pd})$ can be computed by $\Pi$. To see this, note that by \cite[Proposition 4.3.36]{Pav}, $\Pi$ can also be presented by the colimit over the simplicial diagram indexed by the smooth simplicies $\Delta^n$ -- a diagram of shape $N(\Delta^{op})$. By \cite[Lemma 5.5.8.4]{Lur},  $N(\Delta^{op})$ is a sifted simplicial set. By \cite[Corollary 3.2.3.2]{LurHA}, the category $C\mathscr{M}{\rm on}(\infty\mathscr{G}{\rm pd})$ admits such colimits and these colimits are detected by the inclusion $C\mathscr{M}{\rm on}(\infty\mathscr{G}{\rm pd})\into \infty\mathscr{G}{\rm pd}$. Hence the colimit is 
$$
\Pi\Big(\coprod_{n\in \mathbb{N}}\mathbf{B}{\rm O}(n)\Big)\simeq \coprod_{n\in \mathbb{N}}\op{BO}(n)\;.
$$
Finally, since $\mathscr{K}$ is a left adjoint and $\Pi$ sends local equivalences in $\mathscr{S}{\rm h}_{\infty}(\mathscr{M}{\rm an})$ to equivalences, we have 
$$
\Pi\Big(L\circ \mathscr{K}\Big(\coprod_{n\in \mathbb{N}}\mathbf{B}{\rm O}(n)\Big)\Big)\simeq \Pi\Big(\mathscr{K}\Big(\coprod_{n\in \mathbb{N}}\mathbf{B}{\rm O}(n)\Big)\Big)\simeq \mathscr{K}\Big(\coprod_{n\in \mathbb{N}}\op{BO}(n)\Big)\;.
$$
The claim then follows by passing to associated spectra and the well-known fact that (see the beginning of Section 
\ref{Ch-topKO})
$$
\Omega^{\infty}\op{KO}\simeq \mathscr{K}\Big(\coprod_{n\in \mathbb{N}}\op{BO}(n)\Big)\simeq \op{BO}\times \mathbb{Z}\;.
$$

\vspace{-6mm}
\endofproof

We can add connections into the picture in a similar way, giving the following analogous definition.
\begin{definition}
[Smooth $\op{KO}$-spectrum with connections]
We define the smooth $\op{KO}$-spectrum with connections as the connected sheaf of spectra defined by 
$$
\Omega^{\infty}{\bf K}{\rm O}_{\nabla}:=L\circ \mathscr{K}\Big(\coprod_{n\in \NN}\BB{\rm O}(n)_{\nabla}\Big)\;.
$$
\end{definition}
Again we immediately see that we have an equivalence $L\circ \mathscr{K}\circ {\rm Iso}\big({\rm Vect}_{\nabla}^g\big)\simeq \Omega^{\infty}{\bf K}{\rm O}_{\nabla}\;,$
which gives rise to a natural isomorphism $\mathbf{K}{\rm O}_{\nabla}(M)\cong {\rm Gr}\big({\rm Iso}({\rm Vect}^g_{\nabla})(M)\big)$. In this case, the theory $\mathbf{K}{\op{O}}_{\nabla}$ differs from $\op{KO}$, as it detects the data of connection. However, the two theories agree upon topological realization.

\begin{proposition}
[Topological realization of smooth $\op{KO}$-spectrum (with connections)]
\label{grekocn}
We have equivalences in $\mathscr{S}{\rm p}$
$$
\Pi({\bf K}{\rm O}_{\nabla})\simeq \Pi({\bf K}{\rm O})\simeq \op{KO}.
$$
\end{proposition}
\theproof
The same argument as in the proof of Proposition \ref{smkotpko} applies to ${\bf K}{\rm O}_{\nabla}$, since 
$$
\Pi\big(\mathbf{B}{\op{O}}(n)_{\nabla}\big)\simeq \op{BO}(n)
$$
and pullback of bundles with connections commutes with direct sums. We deduce that we have an equivalence of infinite loop spaces 
$$
\Pi\Big(L\circ \mathscr{K}\Big(\coprod_{n\in \mathbb{N}}\mathbf{B}{\rm O}(n)_{\nabla}\Big)\Big)\simeq \Pi\Big(\mathscr{K}\Big(\coprod_{n\in \mathbb{N}}\mathbf{B}{\rm O}(n)_{\nabla}\Big)\Big)\simeq \mathscr{K}\Big(\coprod_{n\in \mathbb{N}}\op{BO}(n)\Big)\;.
$$ 
We, therefore, have an equivalence of infinite loop spaces
$$
\Pi(\Omega^{\infty}{\bf K}{\rm O}_{\nabla})\simeq \Pi\circ L\circ \mathscr{K}\Big(\coprod_{n\in \NN}\BB{\rm O}(n)_{\nabla})\Big)\simeq \mathscr{K}\Big(\coprod_{n\in \NN}\Pi(\BB{\rm O}(n)_{\nabla})\Big)\;.
$$
Combining this with Proposition \ref{smkotpko} gives the result.
\endofproof

Interestingly, both theories $\op{{\bf K}O}$ and $\op{\mathbf{K}O}_{\nabla}$ are related to algebraic orthogonal $\op{K}$-theory in a nontrivial way. Here, by algebraic orthogonal $\op{K}$-theory, we mean the cohomology theory obtained by applying the Quillen plus construction to the classifying space of orthogonal group $\op{O}^{\delta}$ equipped with the  \emph{discrete} topology. This theory is related to the algebraic $\op{K}$-theory of $\RR$ via the inclusion of discrete groups $\op{O}(n)\into \op{GL}_n(\RR)$. 

\begin{proposition}
[Relations among theories] We have the following:
\vspace{-2mm}
\begin{enumerate}[{\bf (i)}]
\item The bottom triangle in the differential cohomology diamond diagram for $\mathbf{K}{\rm O}$-theory takes the following form
\(\label{agtpkocn1}
\xymatrix@R=4pt@C=4em{
& \mathbf{K}{\op{O}} \ar[rd]^-{\mathcal{I}} &
\\
\op{K(O^{\delta})}\ar[ru]^{j}\ar[rr] & & \delta(\op{KO})\;,  &
}
\)
where $\op{K(O^{\delta})}$ is the orthogonal algebraic $K$-theory spectrum of $\mathbb{R}$,
and $\delta (\op{KO})$ is the geometrically discrete spectrum corresponding to the 
topological KO-spectrum.

\item For the spectrum $\mathbf{K}{\rm O}_{\nabla}$, it takes the similar form
\(\label{agtpkocn2}
\xymatrix@R=4pt@C=4em{
& \mathbf{K}{\op{O}}_{\nabla} \ar[rd]^-{\mathcal{I}} &
\\
\op{K(O^{\delta})}\ar[ru]^{j}\ar[rr] & & \delta(\op{KO})\;.  &
}
\)
Hence the differences between the two spectra are purely geometric.
\end{enumerate}
\end{proposition}
\theproof
Proposition \ref{grekocn} already identifies the right diagonal map $\mathcal{I}$ in diagrams \eqref{agtpkocn1} and \eqref{agtpkocn2}. 
We need only compute $\delta\Gamma(\mathbf{K}{\op{O}}_{\nabla})$ and $\delta\Gamma(\mathbf{K}{\op{O}})$. 
The global sections functor $\Gamma$ 
commutes with $\infty$-group completion $\mathscr{K}$. Indeed, the latter is computed by taking the loop space of the $\infty$-colimit over the simplicial object associated to the commutative monoid (i.e. the underlying simplicial object associated to a $\Gamma$-space). Since $\Gamma$ is both a left and right adjoint, it commutes with these operations. Thus, we have
\(\label{dscgrbo}
\Gamma\circ L\circ \mathscr{K}\Big(\coprod_{n\in \mathbb{N}}\mathbf{B}{\op{O}}(n)\Big)\simeq \mathscr{K}\Big(\Gamma\Big(\coprod_{n\in \mathbb{N}}\mathbf{B}{\op{O}}(n)\Big)\Big)\;.
\)
Now, since 
$$
\Gamma\big(\mathbf{B}{\op{O}}(n)\big)\simeq\ast/\!/C^{\infty}(\ast,{\op{O}(n)})=\ast/\!/{\rm O}(n)^{\delta}\simeq {\rm BO}(n)^{\delta}
$$
and since $\Gamma$ commutes with coproducts, we identify the right hand side of \eqref{dscgrbo} with $\mathscr{K}\Big(\coprod_{n\in \mathbb{N}}{\rm B}{\op{O}}(n)^{\delta}\Big)$.
Since we also have 
$\Gamma\big(\mathbf{B}{\op{O}}(n)_{\nabla}\big)\simeq \ast/\!/{\rm O}(n)^{\delta}
\simeq {\rm BO}(n)^{\delta}$, a similar argument shows that 
$$
\Gamma\circ L\circ \mathscr{K}\Big(\coprod_{n\in \mathbb{N}}\mathbf{B}{\op{O}}(n)_{\nabla}\Big)\simeq \mathscr{K}\Big(\coprod_{n\in \mathbb{N}}{\rm B}{\op{O}}(n)^{\delta}\Big)\;.
$$

\vspace{-5mm}
\endofproof

We now would like to construct the Hopkins-Singer differential $\widehat{\op{KO}}$-spectrum. 
We will need the following ingredients:

\begin{enumerate}[{\bf (i)}]
\item Recall 
that we have the Eilenberg-MacLane functor 
$$
\mathscr{H}:\mathscr{C}{\rm h}\longrightarrow \mathscr{S}{\rm p}\;,
$$
which sends an unbounded chain complex to a corresponding spectrum (see the Appendix). 
\item Let $\Omega^*(-;\pi_*(\op{KO}))$ be the complex of forms with coefficients in $\pi_*(\op{KO})$. Explicitly, by rationalizing the 
coefficients of KO,  this complex is 4-periodic and looks as follows
$$
\Omega^*(-;\pi_*(\op{KO}))=\big(\vcenter{\xymatrix@=1.5em{
\hdots \ar[r] &\bigoplus_n \Omega^{4n}\ar[r] & \bigoplus_n \Omega^{4n+1}\ar[r] & \bigoplus_n \Omega^{4n+2}\ar[r] & \bigoplus_n \Omega^{4n+3}\ar[r] & \bigoplus_n \Omega^{4n} \ar[r] & \hdots 
}}\big).
$$
\item We can truncate the complex $\Omega^*(-;\pi_*(\op{KO}))$ at degree zero, removing all forms in negative degrees. We denote this truncated complex by 
$$
\tau_{\leq 0}\Omega^*(-;\pi_*(\op{KO}))=\big(\vcenter{\xymatrix{
\hdots \ar[r] & 0 \ar[r] &  0 \ar[r] & \bigoplus_n \Omega^{4n}\ar[r] & \bigoplus_n \Omega^{4n+2}\ar[r] & \bigoplus_n \Omega^{4n+3} \ar[r] & \hdots 
}}\big)\;,
$$
where the first nonzero component appears in degree zero. 
\item The Pontrjagin character form gives a morphism of smooth stacks 
$$
{\rm Ph}:\coprod_{n\in \mathbb{N}}\mathbf{B}{\op{O}}(n)_{\nabla}
\xymatrix{\ar[r]&} \Omega^0(-;\pi_*(\op{KO}))\;.
$$
\end{enumerate} 
Since $\Omega^0(-;\pi_*(\op{KO}))$ is already a sheaf of abelian groups, we have $\mathscr{K}(\Omega^0(-;\pi_*(\op{KO})))=\mathscr{H}(\Omega^0(-;\pi_*(\op{KO})))$. Postcomposing with the canonical map 
$$
i^*: \mathscr{H}\big(\Omega^0(-;\pi_*(\op{KO}))\big)\xymatrix{\ar[r]&} 
\mathscr{H}\big(\tau_{\leq 0}\Omega^*(-;\pi_*(\op{KO}))\big)
$$
induced by the inclusion $i:\Omega^0(-;\pi_*(\op{KO}))\into \tau_{\leq 0}\Omega^*(-;\pi_*(\op{KO}))$, 
we get an induced map on completions:
\(\label{pchrefsp}
{\rm Ph}:\mathbf{K}{\op{O}}(n)_{\nabla}:=
\mathscr{K}\Big(\coprod_{n\in \mathbb{N}}\mathbf{B}{\op{O}}(n)_{\nabla}\Big)\xymatrix{\ar[r]&} \mathscr{H}(\tau_{\leq 0}\Omega^*(-;\pi_*(\op{KO}))\;.
\)
By Proposition \ref{grekocn}, geometrically realizing gives a map 
$$
\Pi({\rm Ph}):\op{KO}\xymatrix{\ar[r]&} \mathscr{H}(\mathbb{R}[\alpha,\alpha^{-1}])\;,
$$
where $\vert \alpha\vert=4$. 
Setting $\widetilde{\rm Ph}:=\Pi({\rm Ph})$, we now have the following definition. 

\begin{definition} [Hopkins-Singer differential KO-theory]
\item{\bf (i)} The differential KO-theory spectrum is defined via the pullback in sheaves of spectra
$$
\xymatrix@R=1.7em{
{\rm diff}\big({\rm KO}, \widetilde{\rm Ph}, \pi_*(\op{KO})\big)
\ar[r]\ar[d] & \mathscr{H}\big(\tau_{\leq 0}\Omega^*(-;\pi_*(\op{KO}))\big)\ar[d]
\\
\op{KO}\ar[r]^-{\widetilde{\rm Ph}} &  \mathscr{H}(\pi_*(\op{KO})\otimes \RR)
}.
$$
This pullback depends on the map $\widetilde{\rm Ph}$ and the graded ring $\pi_*({\rm KO})$. 
We fix this data once and for all and denote the sheaf of spectra simply as
$$
\widehat{\rm KO}:= {\rm diff}\big({\rm KO}, \widetilde{\rm Ph} ,\pi_*(\op{KO})\big)\;.
$$
{\bf (ii)} The differential KO-spectrum refining higher degree KO-groups is given by the pullback 
$$
\xymatrix@R=1.7em{
{\rm diff}\big(\Sigma^n{\rm KO},\Sigma^n(\widetilde{\rm Ph}), \pi_*(\op{KO})[n]\big)
\ar[r]\ar[d] & \mathscr{H}\big(\tau_{\leq 0}\Omega^*(-;\pi_*(\op{KO})[n])\big)\ar[d]
\\
\Sigma^n\op{KO}\ar[r]^-{\Sigma^n(p)} &  \mathscr{H}(\pi_*(\op{KO})[n]\otimes \RR)
}.
$$
where $\Sigma^n$ denotes the $n$-fold suspension and $\pi_*(\op{KO})[n]$ denotes the shift of the complex $\pi_*(\op{KO})$ up $n$-units. Again we fix this data once and for all and define 
\footnote{These sheaves of spectra are not to be confused with the notation for homology, which we do not consider in 
this paper.} 
$$
\widehat{\rm KO}_n:= {\rm diff}\big(\Sigma^n{\rm KO},\Sigma^n(\widetilde{\rm Ph}), \pi_*(\op{KO})[n]\big)\;.
$$
\item {\bf (iii)} Differential KO-cohomology of a manifold $M$ is defined as 
$$
\widehat{\rm KO}^n(M):=\pi_0\map(M, \widehat{\rm KO}_n)\;.
$$
\end{definition}

One has the following properties, as for any differential cohomology theory.   

\begin{remark}
[Basic properties of $\widehat{\rm KO}$]
\item [{\bf (i)}] {\rm (Diamond)}
From \cite[Lemma 6.8]{BNV}, we see that the differential cohomology hexagon diagram takes the following form
\(
\label{kodfdiam}
\xymatrix @C=-34pt @!C{
&\Omega^{*-1}(M;\pi_*(\op{KO}))/{\rm im}(d) \ar[rd]^{a}\ar[rr]^{d} & &
\Omega_{\rm cl}^*(M;\pi_*(\op{KO}))\ar[rd] &  
\\
H^{*-1}(M;\pi_*(\op{KO})\otimes \RR)\ar[ru]\ar[rd] & &\widehat{\rm KO}(M)\ar[rd]^{\mathcal I} \ar[ru]^{\mathcal R}& &  H^*(M;\pi_*(\op{KO})\otimes \RR) \;,
\\
&\op{KO}^{*-1}(M;U(1))\ar[ru]^-j\ar[rr]^{\beta} & & \op{KO}^*(M) \ar[ru] &
}
\) 
which related differential KO-theory to the underlying topological theory and differential form representatives for the rationalization. 
\item [{\bf (ii)}] {\rm (Coefficients)}
Both diagonals in the diagram are exact and the bottom sequence is exact -- induced from the cofiber/fiber sequence
$$
\op{KO}\simeq \op{KO}\wedge \mathbb{S}\longrightarrow \op{KO}\wedge \mathbb{S}\RR\longrightarrow 
\op{KO}\wedge \mathbb{S}U(1)\;,
$$
where $\mathbb{S}\RR$ and $\mathbb{S}U(1)$ are Moore spectra for $\RR$ and $U(1)$, respectively. 
These correspond to the cohomology theories with coefficients, namely 
${\rm KO}^*(-)$, ${\rm KO}^*(-; \R)$, and ${\rm KO}^*(-; U(1))$, respectively. 

\item [{\bf (iii)}] {\rm (Mayer-Vietoris)}
Again applying the general construction of \cite{BNV} to our case, if $M$ a smooth manifold  and $\{U,V\}$ an open cover, we also have a Mayer-Vietoris sequence
$$
\xymatrix@R=1.6em{
\cdots \ar[r] & \op{KO}^{n-2}(U\cap V;U(1)) \ar[r] & \widehat{\rm KO}^n(M) \ar[r]&\widehat{\rm KO}^n(U)\oplus \widehat{\rm KO}^n(V)
\ar@{->} `r/8pt[d] `/12pt[l] `^d[ll]+<-10ex> `^r/8pt[dll] [dll] 
\\
&  \widehat{\rm KO}^n(U\cap V)\ar[r] & \op{KO}^{n+1}(M)\ar[r] & \cdots .
}
$$
\end{remark}

We now discuss the ring structure on $\widehat{\op{KO}}$. We recall that the $\infty$-category $\mathscr{S}{\rm p}$ admits the structure of a symmetric monoidal $\infty$-category \cite{HSS}\cite{MMSS}\cite[4.8.2.19]{LurHA}.
The $\infty$-category of sheaves with values in $\mathscr{S}{\rm p}$,
being a localization of a functor category into a symmetric monoidal category, 
itself admits a symmetric monoidal    structure from $\mathscr{S}{\rm p}$.

\begin{proposition}
[Ring structure on $\widehat{\op{KO}}$]\label{Prop-ringstr}
The sheaves of spectra $\widehat{\rm KO}_n$ fit together to form an $E_{\infty}$-algebra $\widehat{\rm KO}_\ast:=\bigvee_{n\in \ZZ}\widehat{\rm KO}_n$ in the symmetric monoidal $\infty$-category $\sh_{\infty}(\mathscr{M}{\rm an},\mathscr{S}{\rm p})$. At the level of cohomology, this gives a bigraded algebra 
$$
\pi_0\map\Big(M,\bigvee_{n\in \ZZ}\widehat{\rm KO}_n\Big)\cong 
\bigoplus_{n\in \ZZ}\pi_m\map(M,\widehat{\rm KO}_n)\;.
$$
In particular, if we restrict to the diagonal, we have natural multiplication maps
$$
\hat{\mu}:\widehat{\rm KO}^n(M)\otimes \widehat{\rm KO}^m(M)
\xymatrix{\ar[r]&}
\widehat{\rm KO}^{n+m}(M)
$$
which make the diagrams 
\(\label{rngcmdko1}
\xymatrix{
\widehat{\rm KO}^n(M)\otimes \widehat{\rm KO}^m(M)\ar[rr]^-{\hat{\mu}}\ar[d]^-{\mathcal{I}\otimes \mathcal{I}} && 
\widehat{\rm KO}^{n+m}(M)\ar[d]^-{\mathcal{I}}
\\
{\rm KO}^n(M)\otimes {\rm KO}^m(M)\ar[rr]^-{\mu} && {\rm KO}^{n+m}(M)
}
\)
and
\(\label{rngcmdko2}
\xymatrix{
\widehat{\rm KO}^n(M)\otimes \widehat{\rm KO}^m(M)\ar[rr]^-{\hat{\mu}}\ar[d]^-{\mathcal{R}\otimes \mathcal{R}} &&
\widehat{\rm KO}^{n+m}(M)\ar[d]^-{\mathcal{R}}
\\
\Omega^n(M;\pi_*(\op{KO}))\otimes \Omega^m(M;\pi_*(\op{KO}))\ar[rr]^-{\wedge} && 
\Omega^{n+m}(M;\pi_*(\op{KO}))
}
\)
commute.
\end{proposition}
\theproof
The subcategory of $E_{\infty}$-algebras in $\sh_{\infty}(\mathscr{M}{\rm an};\mathscr{S}{\rm p})$ admits pullbacks. Moreover, this pullback can be computed as the pullback in $\sh_{\infty}(\mathscr{M}{\rm an};\mathscr{S}{\rm p})$ and the algebra structure is inherited from its associated cospan diagram \cite[Section 3.2]{Lur2}. Since ${\rm KO}$ admits the structure of an $E_{\infty}$-algebra, and the map $\widetilde{\rm Ph} :\op{KO}\to \op{KO}\wedge \mathscr{H}\RR\simeq \mathscr{H}(\pi_*(\op{KO})\otimes \RR)$ is a map of $E_{\infty}$-ring spectra, we need only 
focus on the differential form part of the cospan. The wedge product of forms gives maps
$$
\wedge :\;\tau_{\leq 0}\Omega^*\big(-;\pi_*(\op{KO}))[n]\otimes \tau_{\leq 0}\Omega^*(-;\pi_*(\op{KO})\big)[m]
\xymatrix{\ar[r]&}\tau_{\leq 0}\Omega^*(-;\pi_*(\op{KO}))[n+m]
$$
which are strictly graded commutative and associative. Thus, with these multiplication maps, 
the wedge product 
$\bigvee_{n\in \ZZ}\big(\tau_{\leq 0}\Omega^*(-;\pi_*(\op{KO}))[n]\big)$ admits the structure of an $E_{\infty}$-algebra. Since coproducts commute with pullbacks in $\sh_{\infty}(\mathscr{M}{\rm an};\mathscr{S}{\rm p})$, we immediately compute 
$$
\xymatrix{
\bigvee_{n\in \ZZ}\widehat{\rm KO}_n\ar[r]\ar[d] & \bigvee_{n\in \ZZ}\big(\tau_{\leq 0}\Omega^*(-;\pi_*(\op{KO}))[n]\big)\ar[d]
\\
\bigvee_{n\in \ZZ}\Sigma^n\op{KO}\ar[r] & \bigvee_{n\in \ZZ}\Sigma^n \mathscr{H}(\pi_*(\op{KO})\otimes \RR)
\;.}
$$
Since all the maps are morphisms of $E_{\infty}$-algebras, diagram \eqref{rngcmdko1} and 
diagram \eqref{rngcmdko2} manifestly commute.
\endofproof

Using the above multiplicative properties for $\mathcal{I}$ and $\mathcal{R}$, 
we also deduce the following multiplicative property for the third 
natural transformation in the diamond. 

\begin{proposition}
[Multiplicative property of the  secondary curvature map $a$ for $\widehat{\op{KO}}(M)$]
The map $a$ in the diamond (or hexagon) diagram \eqref{kodfdiam}, satisfies 
$$
a(\lambda\wedge \mathcal{R}(x))=a(\lambda)\cdot x\;,
$$
for all $\lambda\in \Omega^{-1}(M;\pi_*(\op{KO}))/{\rm im}(d)$.
\end{proposition}
\theproof
By the commutativity of diagram \eqref{kodfdiam} and the fact that the ring structure on
$\widehat{\op{KO}}^*(M)$ refines the ring structure on $\Omega^*(M;\pi_*(\op{KO}))$ 
(by Proposition \ref{Prop-ringstr} above), we have 
$$
\mathcal{R}(a(\lambda\wedge \mathcal{R}(x)))=d(\lambda\wedge \mathcal{R}(x))=
d\lambda\wedge \mathcal{R}(x)=\mathcal{R}(a(\lambda)\cdot x)\;.
$$
By exactness of the diagonal in the diamond, we have $a(\lambda\wedge \mathcal{R}(x))-a(\lambda)\cdot x=j(y)$ for some $y\in \op{KO}(M;U(1))$. By commutativity of the bottom triangle and the left square of \eqref{kodfdiam}, 
$j(y)=a(\gamma)$ for some class $\gamma\in H^{-1}(M;\pi_*(\op{KO})\otimes \mathbb{R})$ and
\(\label{crvprdop}
a(\lambda\wedge \mathcal{R}(x))-a(\lambda)\cdot x=a(\gamma)
\quad \Longrightarrow \quad
a(\lambda)\cdot x= a(\lambda\wedge \mathcal{R}(x)-\gamma)\;.
\)
To prove the claim, it suffices to show that  $\gamma$ is 0. 
Since the product on $\widehat{\op{KO}}$, the wedge product of forms and maps $a$ and $\mathcal{R}$ are natural with respect to smooth maps $f:M\to N$, it follows that $\gamma$ must define a natural operation in $\lambda$ and $x$. By \eqref{crvprdop}, one easily checks that $\gamma$ is additive in both variables. In particular, restricting to classes with $\mathcal{R}(x)=\mathcal{I}(x)=0$, it defines a natural operation 
$$
\gamma:H^{4m-1}(-;\mathbb{R})\times H^{4n-1}(-;\mathbb{R})\longrightarrow H^{4n+4m-1}(-;\mathbb{R})\;,
$$
for some $m,n\in \mathbb{Z}$. By the K\"unneth Theorem and the Representability Theorem, it follows that $\gamma=0$. 
\endofproof

Just as $\op{KO}$ has as infinite loop space $\Omega^{\infty}\op{KO}\simeq \op{BO}\times \ZZ$, the sheaf of spectra $\widehat{\op{KO}}$ has a related infinite loop stack. Define the smooth stack $\widehat{\op{BO}}$  as follows. Let $\RR[\alpha]$ be the graded algebra with $\vert\alpha \vert=4$ and consider the canonical monomorphism of sheaves of differentially graded algebras 
\(\label{inc4pder}
\xymatrix{
i:\RR[\alpha] \; \ar@{^{(}->}[r] & \tau_{\geq 0}\Omega^*(-;\pi_*(\op{KO}))
}
,
\)
where $\tau_{\geq 0}\Omega^*(-;\RR[\alpha,\alpha^{-1}])$ is the truncated, positively graded complex of differential forms with values in $ \RR[\alpha,\alpha^{-1}]$. We also consider the sheaves of chain complexes $\tau_{\geq 0}\Omega_{\rm cl}^*(-;\RR[\alpha,\alpha^{-1}])$ which is truncated, but has closed forms in degree zero, and the complex $\tau_{\leq 0}\Omega^*(-;\pi_*(\op{KO}))$ appearing in the Hopkins-Singer model. 

\medskip
Note also that from the definition of the Eilenberg-MacLane functor $\mathscr{H}:\mathscr{C}{\rm h}\to \mathscr{S}{\rm p}$, we have a commutative diagram 
$$
\xymatrix@R=1.5em{
\mathscr{C}{\rm h}\ar[rr]^-{\mathscr{H}} && \mathscr{S}{\rm p}
\\
\mathscr{C}{\rm h}_+\ar[rr]^-{\rm DK}\ar@{^{(}->}[u] && \mathscr{S}{\rm p}_+\ar@{^{(}->}[u]\;,
}
$$
where $\mathscr{S}{\rm p}_+$ is the category of connected spectra, or infinite loop spaces, and ${\rm DK}:\mathscr{C}{\rm h}_+\to \sab\to \mathscr{S}{\rm p}_+$ is the Dold-Kan functor. This commutative diagram prolonges to sheaves. Using these ingredients, we make the following definition. 

\begin{definition}[Classifying stack for real bundles with form data]
Define $\widehat{\op{BO}}$ by the $\infty$-pullback diagram 
\(\label{bohatpul}
\xymatrix{
\widehat{\op{BO}}\ar[rr] \ar[d]&& \Omega^0_{\rm cl}(-;\RR[\alpha,\alpha^{-1}])\ar[d]
\\
\delta(\op{BO})\ar[rr]^-{\rm Ph^{\prime}} &&
\op{DK}(\tau_{\geq 0}\Omega_{\rm cl}^*(-;\RR[\alpha,\alpha^{-1}]))\;,
}
\)
where ${\rm DK}:\sh_{\infty}(\mathscr{M}{\rm an};\mathscr{C}{\rm h}_+)\to \sh_{\infty}(\mathscr{M}{\rm an})$ is the Dold-Kan functor, prolonged to sheaves, and the map ${\rm Ph}^{\prime}:\op{BO}\to  \op{DK}(\tau_{\geq 0}\Omega^*(-;\pi_*(\op{KO}))$ is defined by the composite of the Pontrjagin character ${\rm Ph}:\op{BO}\to \bigvee_{n>0}K(\RR,4n)\simeq {\rm DK}(\RR[\alpha])$ with $\op{DK}(i)$ (see \eqref{inc4pder}). The right vertical map is induced by the canonical inclusion.
\end{definition}

\begin{remark}[Corresponding differential cocycles]
A map $M\to \widehat{\op{BO}}$ is equivalently a triple, $(\omega, V\to M, h)$, where
\item {\bf (i)}  $V\to M$ is a real vector bundle with orthogonal structure, classified by a map $V:M\to BO$, 
\item {\bf (ii)} $\omega=\omega_0+\omega_4+\omega_6+\hdots$ is a differential form, and
\item  {\bf (iii)} the map $h:M\times \Delta[1]\to \op{DK}(\tau_{\geq 0}\Omega_{\rm cl}^*(-;\RR[\alpha,\alpha^{-1}])$ whose restriction to $i_0:\Delta[0]\into \Delta[1]$ the Pontrjagin character ${\rm DK}(i)\circ {\rm Ph}(V)$ and the restriction to the endpoint $i_1:\Delta[0]\into \Delta[1]$ is $\omega$. 

\noindent Note that, by resolving $M$ by the {\v C}ech nerve of a good open cover $\{U_{\alpha}\}$, $h$ can be identified with a {\v C}ech cocycle with values in the sheaf $\tau_{\geq 0}\Omega_{\rm cl}^*(-;\RR[\alpha,\alpha^{-1}])$, such that $D(h)=(d-\delta)(h)=\omega-i{\rm Ph}(V)$, where ${\rm Ph}(V)$ is a {\v C}ech cocycle representative of the Pontrjagin character. 
\end{remark}

\begin{proposition}[Infinite loop stack for differential KO]
The zero stack (or infinite loop stack) of the sheaf of spectra $\widehat{\op{KO}}$ can be identified with $\widehat{\op{BO}}\times \ZZ$
\end{proposition}
\theproof
From the definition of the Hopkins-Singer $\widehat{\op{KO}}$-theory, and using the fact that passing to infinite loop stacks preserves $\infty$-pullbacks, we have the $\infty$-pullback diagram 
$$
\xymatrix{
\Omega^{\infty}\widehat{\op{KO}}\ar[rr]\ar[d] &&
{\rm DK}(\tau_{\leq 0}\Omega^*(-;\pi_*(\op{KO}))\ar[d]
\\
\Omega^{\infty}\op{KO}\ar[rr]^-{\widetilde{\rm Ph}} && 
{\rm DK}(\pi_*(\op{KO}))\;.
}
$$
Since the map $\RR[\alpha]\into \tau_{\geq 0}\Omega_{\rm cl}^*(-;\RR[\alpha,\alpha^{-1}])=\tau_{\geq 0}\Omega_{\rm cl}^*(-;\pi_*(\op{KO}))$ and the inclusion $i:\Omega_{\rm cl}^0(-;\pi_*(\op{KO}))\into \tau_{\leq 0}\Omega^*(-;\pi_*(\op{KO}))$ are quasi-isomorphisms and ${\rm DK}$ sends quasi-isomorphisms
to equivalences, it follows that $\Omega^{\infty}\widehat{\op{KO}}$ can be presented by 
the $\infty$-pullback \eqref{bohatpul}.
\endofproof

\subsection{Relation to the complex theory: Differential refinement of the Bott sequence}
\label{Sec-cplx-diff}

We would like to have differential refinements of some of the relationships between $\op{KO}$-theory 
and $\op{K}$-theory. 
The main emphasis will be on the following refinement of the Bott sequence from Section \ref{Sec-cplx}. 
Using the notation from that section, recall the realification 
$r:\op{K}^*\to \op{KO}^*$ and complexification $c:\op{KO}^*\to \op{K}^*$ operations as described there. 
Recall also the inclusion of flat classes $j:\op{KO}(-;U(1))\into \widehat{\op{KO}}(-)$, the topological 
realization map $\mathcal{I}:\widehat{\op{KO}}(-)\to \op{KO}(-)$, and the Bockstein map 
$\beta_{U(1)}:\op{KO}(-;U(1))\to \op{KO}(-)$ (see the diamond diagram \eqref{kodfdiam}). 
The maps $r$ and $c$ can be refined to the differential theories $\widehat{\op{K}}$ and 
$\widehat{\op{KO}}$ in a fairly obvious way. For example, the realification map $r:\op{K}\to \op{KO}$ 
induces an epimorphism on coefficients $r:\pi_*(\op{K})\to \pi_*(\op{KO})$. This induces a morphism 
on differential forms
$$
\overline{r}:\Omega^*(-;\pi_*(\op{K}))\longrightarrow \Omega^*(-;\pi_*(\op{KO}))
$$
which carries through to a homotopy commutative diagram of sheaves of spectra
\(
\xymatrix@R=1.5em{
\mathscr{H}(\tau_{\leq i}\Omega^{*}(-;\pi_*(\op{K})))\ar[r]^{\overline{r}_i}\ar[d] & \mathscr{H}(\tau_{\leq i}\Omega^{*}(-;\pi_*(\op{KO}))) \ar[d]
\\
\mathscr{H}(\pi_{*-i}(\op{K})) \ar[r]^{r_i} & \mathscr{H}(\pi_{*-i}(\op{KO}))
\\
\Sigma^{i}\op{K}\ar[r]^-{r_i}\ar[u] & \Sigma^{i}\op{KO} \ar[u]
}
\)
which induces a morphism of sheaves of spectra 
\(\label{rea-morsp}
\hat{r}:{\rm diff}(\Sigma^{i}\op{K},\widetilde{{\rm Ch}}, \pi_{*+i}(\op{K}))\longrightarrow 
{\rm diff}(\Sigma^{i}\op{KO},\widetilde{{\rm Ph}}, \pi_{*+i}(\op{KO}))\;.
\)
Similar considerations gives rise to a complexification morphism of sheaves of spectra
\(\label{cpx-morsp}
\hat{c}:{\rm diff}(\Sigma^{i}\op{KO},\widetilde{{\rm Ph}}, \pi_{*+i}(\op{KO}))\longrightarrow
{\rm diff}(\Sigma^{i}\op{K},\widetilde{{\rm Ch}}, \pi_{*+i}(\op{K}))\;.
\)
\begin{definition}[Differential realification and complexification maps]\label{Def-diff-r-c} 
{\bf (i)} We define the \emph{differential} realification map 
$$
\hat{r}_i:\widehat{\op{K}}^{i}(M)\longrightarrow \widehat{\op{KO}}^{i}(M)
$$
as the natural map induced by \eqref{rea-morsp}.
\\
\noindent {\bf (ii)} We define the \emph{differential} complexification map 
$$
\hat{c}_i:\widehat{\op{KO}}^{i}(M)\longrightarrow \widehat{\op{K}}^{i}(M)
$$
as the natural map induced by \eqref{cpx-morsp}.
\end{definition}

\begin{proposition}[Differential refinement of the Bott sequence]\label{prop-diffBott}
There is a long exact sequence
$$
\xymatrix@R=1.6em{
\cdots \ar[r] & \op{K}^{n-2}(M;U(1)) \ar[r]^-{r\beta_{U(1)}} & {\rm KO}^{n+1}(M) \ar[r]^-{j\tilde{\eta} }&\widehat{\rm KO}^n(M)
\ar@{->} `r/8pt[d] `/12pt[l]`^d[ll]+<-10ex>_-{\hat{c}} `^r/8pt[dll] [dll] 
\\
&  \widehat{\rm K}^n(M)\ar[r]^-{\hat{r}\mathcal{I}} & \op{KO}^{n+2}(M)\ar[r] & \cdots 
}
$$
refining the Bott sequence \eqref{Bott-seq}, where the map $\tilde{\eta}:\op{KO}^*(M)\to \op{KO}^*(M;U(1))$ satisfies $\beta_{U(1)}(\tilde{\eta})=\eta:\op{KO}^*(M)\to \op{KO}^{*-1}(M)$. 
\end{proposition}
\theproof
The fiber sequence 
$$\Sigma\op{KO}\overset{\eta}{\longrightarrow} \op{KO}\overset{c}{\longrightarrow} \op{K}$$
induces a long exact sequence on homotopy groups, which (rationally) gives an exact sequence of graded abelian groups
\(\label{cefkoksq}
0  \longrightarrow \pi_*(\op{KO})\otimes \RR \xrightarrow{\;\;c\;\;} \pi_*(\op{K})\otimes \RR \longrightarrow (\pi_{*}(\op{KO})\otimes \RR)[2]\longrightarrow 0\;.
\)
Applying the Eilenberg-MacLane functor $\mathscr{H}:\sh_{\infty}(\mathscr{M}{\rm an};\mathscr{C}{\rm h})\to \sh_{\infty}(\mathscr{M}{\rm an};\mathscr{S}{\rm p})$, we have homotopy commutative diagrams 
\(\label{thsqdgdif}
\xymatrix{
\ast \ar[d]\ar[rr]^-{0}_{\ }="s2" &&
\mathscr{H}(\tau_{\leq -n}\Omega^*(-;\pi_*(\op{KO})))\ar[rr]_{\ }="s3" \ar[d]&&
\mathscr{H}(\tau_{\leq -n}\Omega^*(-;\pi_*(\op{K})))\ar[d]
\\
\ast \ar[rr]^-{0}="t2"_{\ }="s1" &&
\mathscr{H}(\pi_{*+n}(\op{KO}))\ar[rr]^{\ }="t3"_{\ }="s4" && 
\mathscr{H}( \pi_{*+n}(\op{K}))
\\
\Sigma^{-n-1} \op{KO}\ar[u] \ar[rr]_-{\eta}^{\ }="t1" && 
\Sigma^{-n}\op{KO}\ar[rr]^{\ }="t4"_-{c}\ar[u]_-{\widetilde{\rm Ph}} && \Sigma^{-n}\op{K}\ar[u]_-{\widetilde{\rm Ch}}
\ar@{=>}_{\bf 2}"s1"+<4ex,-1ex>;"t1"+<-6ex,1ex>
\ar@{=>}_{\rm id}"s2"+<2ex,-1ex>;"t2"+<-4ex,1ex>
\ar@{=>}_{\rm id}"s3"+<4ex,-1ex>;"t3"+<-5ex,1ex>
\ar@{=>}_{\bf 1}"s4"+<4ex,-1ex>;"t4"+<-5ex,1ex>
}.
\)
The homotopy {\bf 1} exists from the compatibility of complexification followed by the Chern character and the Pontrjagin character, the vertical maps on the left are the induced maps on fibers and the homotopy {\bf 2} is induced from the universal property of the pullback. Since taking pullbacks commutes with 
taking fibers, we have an induced homotopy commutative diagram
\(\label{infskok}\xymatrix{
& \mathscr{H}(\tau_{\leq -n}\Omega^*(-;\pi_*(\op{KO})) &
\\
\Sigma^{-n-1}\op{KO} \ar[r]^-{\eta^{\prime}}\ar@{-->}[dr]&  
{\rm diff}\big(\Sigma^{-n}\op{KO},\widetilde{{\rm Ph}}, \pi_{*+n}(\op{KO})\big)
\ar[r]^-{\hat{c}} \ar[u]^-{\mathcal{R}} &  
{\rm diff}\big(\Sigma^{-n}\op{K},\widetilde{{\rm Ch}}, \pi_{*+n}(\op{K})\big)
\\
& \Sigma^{-n+1}\op{KO}\wedge \mathbb{S}U(1) \ar[u]^-{j}& 
}
\)
where the horizontal and vertical sequences are fiber sequence and $\mathbb{S}{\rm U}(1)$ is the Moore spectrum for $U(1)$. The map $\eta^{\prime}$ refines $\eta$ and factors (up to homotopy) through the fiber of the map $\mathcal{R}$, hence factorizes through $j$,
filling the dashed arrow. Let $\tilde{\eta}$ denote such a filler. Then at the level of homotopy classes, we have 
$$\beta_{U(1)}(\tilde{\eta})=\mathcal{I}(j(\tilde{\eta}))=\mathcal{I}(\eta^{\prime})=\eta\;.$$ 
Recall that for any sheaf of spectra $\mathscr{E}$, with refinement $\widehat{\mathscr{E}}:={\rm diff}(\mathscr{E},c,A)$, we have $\Sigma^{n}{\rm diff}(\Sigma^n\mathscr{E},c,A[n])\simeq \Sigma^n\mathscr{E}\wedge \mathbb{S}U(1) $ for $n>0$. Then using this, applying $\mathcal{I}$ to the horizontal fiber sequence in \eqref{infskok} and using the homotopy commutative diagram 
\vspace{-2mm} 
$$
\xymatrix@R=1.5em{
\Sigma^{-n+1}\op{K}\wedge \mathbb{S}U(1) \ar[rr]^-j\ar[drr]_-{\beta} && {\rm diff}(\Sigma^{-2i}\op{K},\widetilde{{\rm Ph}}, \pi_{*+n}(\op{K}))
\ar[d]^-{\mathcal{I}}
\\
&&\Sigma^{n}\op{K}
}
$$ 

\vspace{-2mm} 
\noindent gives the identification $\hat{r}\beta_{U(1)}$ for the connecting map. 
\endofproof

%

We highlight that the above differential Bott sequence involves an interesting interplay of various 
theories from topological to flat to differential, as well as real and complex flavors. We, therefore,
believe that it has a very rich content and hence deserves further study.
For now we consider direct consequences of the differential Bott sequence. 
In general we know that $j$ is always injective while $\mathcal{I}$ is always 
surjective. For the other maps we have the following statements which 
generalize classical statements (see e.g. \cite{Do}). Note that the nature of the maps
is different from the ones in the classical case. 

\begin{corollary}[Maps in the differential Bott sequence]
In the differential Bott sequence
\begin{enumerate}[{\bf (i)}]
\vspace{-1mm}
\item The product of any element of $\op{KO}^*(X)$ with $j \widetilde{\eta}$ is zero if and only if
that element is in the image of $r\beta_{U(1)}$.
\vspace{-2mm}
\item Whenever $\op{KO}^*(X)$ is free abelian, we have that $r\mathcal{I}$ is 
an epimorphism and $c$ is a monomorphism. 
\end{enumerate}
\end{corollary}

In \cite{Zi}, a systematic computation of
the KO-groups of all full flag varieties is carried out using Witt rings 
and the relation with K-theory. The Witt group of a space $X$
is defined to be $\mathcal{W}^i(X):={\rm KO}^{2i}(X)/r_i$, where 
$r_i: \op{K}^{2i}(X) \to {\rm KO}^{2i}(X)$ is the realification map. These
capture the 2-torsion of the KO-groups of $X$.
The Witt ring is $\mathcal{W}^*(X):= \bigoplus_{i \in \Z/4}\op{KO}^{2i}(X)/r_i$.
We consider the corresponding ring in differential KO-theory:

\begin{definition}[Differential Witt group and ring]
With $\hat{r}_i$ the realification map in Definition \ref{Def-diff-r-c},
\vspace{-1mm}
\item {\bf (i)} We define the \emph{differential Witt group} $\widehat{\mathcal{W}}_i(M)$ as the group 
$\widehat{\op{KO}}^{2i}(M)/\hat{r}_i$.
\vspace{-1mm}
\item {\bf (ii)} Similarly, we define the \emph{differential Witt ring} as 
$\widehat{\mathcal{W}}^*(X):= \bigoplus_{i \in \Z/4}\widehat{\op{KO}}^{2i}(X)/\hat{r}_i$.
\end{definition}

The following shows explicitly that differential refinement does not capture any new
geometric information. This is in contrast to torsion cohomology operations retaining 
some geometry under differential refinement \cite{GS1}\cite{GS2}.

\begin{proposition}[Correspondence of differential Witt groups and rings]
We have a natural isomorphism $\widehat{\mathcal{W}}_i(M)\cong \mathcal{W}_i(M)$, where  $\mathcal{W}_i(M)$ is the Witt group corresponding to the topological theory $\op{KO}$. Moreover, this isomorphism lifts to an isomorphism on corresponding Witt rings. 
\end{proposition}
\theproof
From the differential cohomology diamond and the definition of $\hat{r}_i$, we have a commutative diagram 
\vspace{-1mm}
\(\label{wittdiag}
\xymatrix{
\Omega^{2i-1}(M;\pi_*(\op{K}))/{\rm im}(d)\ar[r]^{\overline{r}_i}\ar[d]^-{a} & \Omega^{2i-1}(M;\pi_*(\op{KO}))/{\rm im}(d)\ar@{-->}[r]\ar[d]^-{a} & 0 \ar[d]
\\
\widehat{\op{K}}^{2i}(M)\ar[r]^{\hat{r_i}}\ar[d]\ar[d]^-{\mathcal{I}} & \widehat{\op{KO}}^{2i}(M)\ar[d]^-{\mathcal{I}}\ar@{-->}[r] & \widehat{\mathcal{W}}_i(M)\ar[d]^-{\overline{\mathcal{I}}}
\\
\op{K}^{2i}(M)\ar[r]^-{r_i} & \op{KO}^{2i}(M)\ar@{-->}[r] & \mathcal{W}_i(M)
}
\)
where the two columns on the left are exact. The map $\overline{r}_i$ is the map on forms induced by the 
realification on coefficients $r_i:\pi_*(\op{K})\to \pi_*(\op{KO})$. Since this map is surjective, it follows that 
$\overline{r}_i$ is surjecitve. Passing to the corresponding quotients gives the dashed arrows in diagram
\eqref{wittdiag}. A quick diagram chase reveals that the induced map $\overline{\mathcal{I}}$ is a both 
an epimorphism and a monomorphism, hence an isomorphism. Since the map $\mathcal{I}$ is a ring homomorphism, it is immediate that the isomorphism holds even at the level of rings.
\endofproof

We will encounter further interactions between the real and the complex differential theories 
in the flat case below.

\subsection{The flat theory} 
\label{Sec-flat}
KO-theory with arbitrary coefficients has been studied in \cite{Mei}, the results of which, based on \cite{A-thesis}, 
we will apply for the case when the coefficients take values in $U(1)$. The spectrum $\op{KO}$ with 
coefficients in the 1-dimensional unitary group $U(1)$ is the smash product with the Moore spectrum 
with $U(1)$-coefficients, $\mathbb{S}{\rm U}(1)$, i.e., 
$\op{KO}_{U(1)}=\op{KO} \wedge \mathbb{S}{\rm U}(1)$, and similarly for the complex case. 
The flat theory will be important in our identifications of the differentials in the Atiyah-Hirzebruch 
spectral sequence in Section \ref{Ch-AHSS}. 

\medskip
As in the complex case \cite[(2.2)]{FL}, the flat real theory $\op{KO}_{\rm flat}(-)$ participates in the following 
short exact sequence
$$
\xymatrix{
0 \ar[r] & \op{KO}_{\rm flat}(X)\ar[r] &
\widehat{\op{KO}}(X) \ar[r]^{\rm Ph} &
\Omega^{4n}_{\op{BO}}(X) \ar[r] & 0\;,
}
$$
where 
$$
\Omega^{4n}_{\op{BO}}(X):=\left\{ \omega \in \Omega^{4n}_{d=0}(X)\;|\; [\omega]\in {\rm Im}\big(r \circ {\rm Ph}: \op{KO}^0 \longrightarrow H^{4n}(X; \RR)\big) \right\}\;.
$$
In \cite[Section 2.3]{BS}, it was shown that for any Hopkins-Singer type differential cohomology theory $\widehat{{E}}$, one always has a natural isomorphism \footnote{In \cite{BS}, there is a condition on the torsion subgroup of $E^*({\rm pt})$ which is not satisfied by $\op{KO}$. However, in \cite[Proposition 4.57]{Bun} there is a more general proof which does not rely on this condition.}
$$\widehat{{E}}_{\rm flat}(M)\cong {E}(M;U(1))\;,$$
where the cohomology theory on the right is the theory represented by the spectrum $E_{U(1)}:=E\wedge \mathbb{S}{\rm U}(1)$. One always has a long exact sequence (Bockstein sequence) associated to the fiber/cofiber sequence
$$
E\longrightarrow E\wedge \mathbb{S}\mathbb{R}\longrightarrow E\wedge \mathbb{S}{\rm U}(1)\;.
$$
This long exact sequence appears as the bottom long fiber sequence in diagram \eqref{spdggcoh}.

\medskip
Let $\eta: S^3 \to S^2$ be the
second Hopf map and let $H$ denote the complex Hopf bundle over the 
complex projective plane $\CC P^2$. In \cite{Mei} it was shown that 
there is an equivalence of spectra
$\op{K} \simeq \op{KO} \wedge h$, where $h=S^0 \cup_\eta e^2$ is the suspension
spectrum whose second term is $\CC P^2$ (see \cite[p. 206]{Ad}).
For $X$ a pointed CW-complex there is a cofibration of spectra
$$
\op{KO} \wedge h \longrightarrow \op{KO} \wedge \Sigma^2 S^0 \longrightarrow 
\op{KO} \wedge \Sigma S^0\;.
$$
After smashing with the Moore spectrum $\mathbb{S}{\rm U}(1)$, this induces a long exact sequence 
$$
\xymatrix{
\cdots \ar[r] & \widetilde{\op{K}}^n(X; {\rm U}(1)) 
\ar[r] &
\widetilde{\op{KO}}^{n+2}(X; {\rm U}(1))
\ar[r] &
\widetilde{\op{KO}}^{n+1}(X; {\rm U}(1))
\ar[r] &
\widetilde{\op{K}}^{n+1}(X; {\rm U}(1))
\ar[r] &
\cdots\;.
}
$$
The following is then an immediate application of \cite[Theorem 1.7]{Mei}
adapted from pointed CW-complexes to smooth manifolds (hence we omit the proof). 

\begin{proposition} 
[Vanishing of flat KO-theory]
Let $M$ be a smooth manifold. Then 
$\widetilde{\op{KO}}^{*}(M_+; {\rm U}(1))\cong \widehat{\op{KO}}_{\rm flat}^{*+1}(M)=0$ if and only if $\widetilde{\op{K}}^{*}(M_+; {\rm U}(1))\cong \widehat{\op{K}}_{\rm flat}^{*+1}(M)=0$. 
\end{proposition} 

We now use the known results about flat K-theory to readily deduce 
results about flat KO-theory \footnote{In the journal version of this article, the next proposition contains a calculation error. This error was already corrected when the proposition was used for spectral sequence calculations, however we failed to correct the original error. It has been corrected in the present version.}

\begin{proposition}
[Coefficients of KO-theory with U(1)-coefficients]\label{codfrng}
We have the identification
$$
{\widetilde{\rm KO}}_{U(1)}^{-i}({\rm pt})=
\left\{
\begin{array}{ll}
U(1)\{\beta^k\} & \text{for}\; \; i=8k,
\\
\ZZ/2\{\eta \beta^k\} & \text{for}\;\; i=8k+2,
\\
\ZZ/2\{\eta^2\beta^k\} & \text{for} \; \; i=8k+3,
\\
U(1)\{\alpha \beta^k\} & \text{for}\; \;  i=8k+4,
\\
0 & \text{otherwise}.
\end{array}
\right.
$$  
where for a generator $\beta$, we make the convention $U(1)\{\beta\}:=\RR\{\beta\}/\ZZ\{\beta\}$.
\end{proposition}
\theproof
We have a long exact sequence on coefficients
$$
\cdots \longrightarrow  \widetilde{\rm KO}_{U(1)}^{*-1}({\rm pt})\xrightarrow{\beta_{U(1)}} \widetilde{\rm KO}^*({\rm pt}) \longrightarrow \widetilde{\rm KO}^*({\rm pt})\otimes \mathbb{R} \longrightarrow \widetilde{\rm KO}_{U(1)}^{*}({\rm pt})\longrightarrow \cdots.
$$
Using the identification of $\pi_*(\widetilde{\rm KO})$, we immediately see that we have isomorphisms
$$
\widetilde{\rm KO}_{U(1)}^{-2}({\rm pt})\overset{\beta_{U(1)}}{\cong} \widetilde{\rm KO}^{-1}({\rm pt})\cong \mathbb{Z}/2\{\eta\}
\qquad \text{and} \qquad
\widetilde{\rm KO}_{U(1)}^{-3}({\rm pt})\overset{\beta_{U(1)}}{\cong} \widetilde{\rm KO}^{-2}({\rm pt})\cong \mathbb{Z}/2\{\eta^2\}\;.
$$
We have an exact sequence $\widetilde{\rm KO}_{U(1)}^{-1}({\rm pt})\into \ZZ\into \RR,$ so $\widetilde{\rm KO}_{U(1)}^{-1}({\rm pt})\cong 0$. From the above coefficient sequence, we have an exact sequence
$$
\mathbb{Z}\{\alpha\}\longrightarrow \mathbb{R}\{\alpha\}\longrightarrow \widetilde{\rm KO}_{U(1)}^{-4}({\rm pt})\longrightarrow 0\;,
$$
where $\ZZ\{\alpha\}$ denotes the free abelian group generated by $\alpha$ and $\RR\{\alpha\}$ is the free 1-dimensional vector space, generated by $\alpha$. Hence, 
$$\widetilde{\rm KO}^{-4}_{U(1)}({\rm pt})=U(1)\{\alpha\}=\RR\{\alpha\}/\ZZ\{\alpha\}.$$
Finally, we have $\widetilde{\rm KO}_{U(1)}^{-5}({\rm pt})\cong \widetilde{\rm KO}_{U(1)}^{-6}({\rm pt})\cong \widetilde{\rm KO}^{-7}({\rm pt})\cong 0$, 
as well as exact sequences
$$
0\longrightarrow \widetilde{\rm KO}^{-7}_{U(1)}({\rm pt})\longrightarrow \ZZ\{\beta\}\longrightarrow \RR\{\beta\} \quad\text{and} \quad \mathbb{Z}\{\beta\}\longrightarrow \mathbb{R}\{\beta\}\longrightarrow \widetilde{\rm KO}_{U(1)}^{-8}({\rm pt}) \longrightarrow 0.
$$
Hence, $\widetilde{\rm KO}_{U(1)}^{-7}({\rm pt})=0$ and $\widetilde{\rm KO}_{U(1)}^{-8}({\rm pt})\cong U(1)\{\beta\}$. The claim follows by 8-periodicity.
\endofproof

Using the identification $\widehat{\op{KO}}_{\rm flat}^*(M)\cong \op{KO}^{*-1}_{U(1)}(M)$, 
\footnote{Note that in the differential theory we are still using the underlying reduced 
topological theory. However, in order to not overburden the notation we avoid writing 
$\widehat{\widetilde{\rm KO}}$ and assume it to be understood.} 
we immediately have the following identification for coefficients of the flat theory. This will be useful for the AHSS in Section \ref{Ch-AHSS}.

\begin{corollary}[Coefficients of flat KO] We have the following coefficients
$$
\widehat{{\op{KO}}}_{\rm flat}^{-i}({\rm pt})=
\left\{
\begin{array}{ll}
U(1)\{\beta^k\} & \text{for}\; \; i=8k-1,
\\
\ZZ/2\{\eta \beta^k\} & \text{for}\;\; i=8k+1,
\\
\ZZ/2\{\eta^2\beta^k\} & \text{for} \; \; i=8k+2,
\\
U(1)\{\alpha \beta^k\}& \text{for}\; \;  i=8k+3,
\\
0 & \text{otherwise}.
\end{array}
\right.
$$  
\end{corollary}

\begin{remark}[Pairing with KO-homology]
There is a pairing between KO-theory and KO-homology defined as 
a map $\op{KO}^p \times \op{KO}_p \to \ZZ$. This factors as
$\op{KO}^p \times \op{KO}_p \to \op{KO}_0({\rm pt})\cong \ZZ$, where the isomorphism
between $\op{KO}_0({\rm pt})$ and $\ZZ$  is given by a form of an index map
(in geometric or analytic settings).  For coefficients in $U(1)$, 
there should be an index pairing $\op{KO}^p(X; U(1)) \times \op{KO}_p(X) \to U(1)$, which 
again should involve an isomorphism $\op{KO}({\rm pt}; U(1))\cong U(1)$
that requires an appropriate index. The complex case $\op{K}(-; U(1))$ 
was realized analytically 
via the reduced $\eta$-invariant \cite{Lott} and geometrically via 
KK-theory involving $II_1$-factors \cite{Dee}. It would be interesting to develop
this aspect, which we believe holds, but we will not consider this here. 
\end{remark}

%

\subsection{The cycle map}
\label{Sec-cyclemap}

We have already seen that theory $\mathbf{K}{\rm O}_{\nabla}$ represents real orthogonal vector bundles 
equipped with metric connections (see Section \ref{Sec-construct} in Section \ref{Ch-diffKO}). 
In this section, we describe a cycle map which relates this theory to the 
Hopkins-Singer differential $\widehat{\op{KO}}$ model, completing diagram \eqref{shspdfko}. 

\begin{definition}[Cycle maps for KO and $\widehat{\op{KO}}$]\label{def-cyc-dif}
{\bf (i)} Fix an equivalence of infinite loop spaces $\Omega^{\infty}\op{KO}\simeq \op{BO}\times \ZZ$. 
This equivalence induces a natural transformation of functors 
\(\label{cyckth}
{\rm cyc}:\pi_0{\rm Iso}({\rm Vect})\longrightarrow \op{KO}(-)\;.
\)
We call such a map a \emph{cycle map} for $\op{KO}$-theory. 
\item {\bf (ii)} Let $\mathscr{F}:\pi_0{\rm Iso}({\rm Vect}_{\nabla})\to \pi_0{\rm Iso}({\rm Vect})$ 
be the forgetful map which sends an isomorphism class of vector bundles with connection to the 
underlying class of bundles. We call any map ${\widehat{{\rm cyc}}}$ filling the top arrow in the
commutative diagram 
$$
\xymatrix{
\pi_0{\rm Iso}({\rm Vect}_{\nabla})\ar[rr]^-{\widehat{{\rm cyc}}}\ar[d]_-{\mathscr{F}}&& 
\widehat{\op{KO}}(-)\ar[d]^-{\mathcal{I}}
\\
\pi_0{\rm Iso}({\rm Vect})\ar[rr]^-{\rm cyc} && \op{KO}(-)
}
$$
a \emph{refinement of the cycle map}.
\end{definition}

The next proposition shows that there is a unique refinement of a cycle maps which is compatible 
with the underlying topological and geometric data. We use the notation from diagram \eqref{kodfdiam}.

\begin{proposition}
[Unique refinement of cycle maps]
\label{dfcycmp}
The set of natural transformations 
\(\label{hclvcnko}
\widehat{{\rm cyc}}:\pi_0{\rm Iso}({\rm Vect}_{\nabla})\longrightarrow \widehat{\op{KO}}(-)
\)
which satisfy $\mathcal{R}(\widehat{{\rm cyc}}(V,\nabla))={\rm Ph}(\mathcal{F}_{\nabla})$ and 
$\mathcal{I}(\widehat{{\rm cyc}}(V,\nabla))={\rm cyc}(V)$ is either empty or contains a unique element.
\end{proposition}
\theproof
Let $V$ be a $k$-dimensional real vector bundle on a smooth manifold $M$. First, observe that if $\nabla$ and
$\nabla^{\prime}$ are connections on $V$ then, from the affine structure on the space of connections, 
the homotopy formula for $\widehat{\op{KO}}$ (see \eqref{hofmdfco} in the Appendix) implies that 
$$
\phi(V,\nabla)-\phi(V,\nabla^{\prime})=a\big({\rm cs}_{\phi}(\nabla,\nabla^{\prime})\big)
$$
for any natural map $\phi:\pi_0{\rm Vect}(-)\to \widehat{\op{KO}}$.  It is always possible 
(see \cite{Bun}\cite{Bun10} in the spirit of \cite{NR}) to find an $(n+1)$-connected map 
$W:N\to \op{BO}(k)$ (with $N$ a smooth manifold and $n$ sufficiently large), and a smooth 
map $f:M\to N$ such that $f^*W\cong V$. Now, since the rational cohomology of $\op{BO}$ 
is generated by Pontrjagin classes (see, e.g., \cite{MT}), it follows from connectivity that 
\footnote{The notation $k:4k-1>n$ here indicates that we are summing over all $k$'s up to 
the greatest integer $\floor*{\tfrac{1}{4}(n+1)}$.}
$$
\bigoplus_{k: 4k-1\leq n}H^{4k-1}(N;\RR)=0\;.
$$
Fix an element $x\in \widehat{\op{KO}}(N)$ with $\mathcal{R}(x)={\rm Ph}(\nabla^{\prime})$ 
and $\mathcal{I}(x)={\rm cyc}(V)$. Such choices are parametrized by (see Proposition
\ref{prop-setref} in the Appendix)
\vspace{-2mm}
$$
a\big(H^{-1}(N;\pi_*(\op{KO})\otimes \RR)\big)=
a\bigg(\bigoplus_{k: 4k-1 >n}H^{2k+1}(N;\RR)\bigg)\;.
$$

\vspace{-1mm}
\noindent Let $\phi$ and $\psi$ be elements of the set described in the proposition (if nonempty). 
Since, $\phi$ and $\psi$ both refine the Pontrjagin character form, it follows via transgression that 
${\rm cs}_{\phi}(\nabla,\nabla^{\prime})={\rm cs}_{\psi}(\nabla,\nabla^{\prime})$
for any two connections $\nabla$, $\nabla'$ on the same bundle. By the above, it follows that there is a 
form $\rho_N:=\sum_{k:4k-1>n}\rho_{4k-1}$ on $N$ such that, by naturality, 
\vspace{-2mm}
\bea
a(\rho_M)=\phi(V,\nabla)-\psi(V,\nabla) &=& \phi(V,f^*\nabla^{\prime})+a\big({\rm cs}_{\phi}(\nabla,f^*\nabla^{\prime})\big)-\psi(V,f^*\nabla^{\prime})-a\big({\rm cs}_{\psi}(\nabla,f^*\nabla^{\prime})\big)
\\
&=& f^*(\phi(W,\nabla^{\prime})-\psi(V,\nabla^{\prime}))
\\
&=& a(f^*\rho_N)\;.
\eea 

\vspace{-1mm}
\noindent
Hence, up to forms in the image of ${\rm Ph}$, the components of $\rho_M$ in degrees $\leq n$ must vanish. By finite dimensionality of $M$, we can take $n$ sufficiently large so that $a(\rho_M)$ vanishes. 
\endofproof

\newpage 

We now turn our attention to showing the existence of a differential refinement of the cycle map. We first 
construct this map as a full morphism of sheaves of spectra. We will make use of the adjunction between 
the Grothendieck group completion Gr and the functor $i$ that forgets the group structure but retains the 
monoid structure (see the Appendix). 

\begin{proposition}
[Refinement of cycle map via sheaves of spectra]\label{csconvb}
Given an equivalence of infinite loop spaces ${\rm cyc}:\Omega^{\infty}\op{KO}\simeq \op{BO}\times \ZZ$ 
giving rise to a cycle map, there is a morphism of sheaves of spectra
$$
\widehat{{\rm cyc}}:{\bf K}{\op{O}}_{\nabla}\longrightarrow \widehat{\op{KO}}\;,
$$
such that $\mathcal{I}(\widehat{{\rm cyc}})\simeq {\rm cyc}\circ \mathscr{F}$ and 
$\mathcal{R}(\widehat{{\rm cyc}})\simeq {\rm Ph}$, where ${\rm Ph}$ is the map 
defined by \eqref{pchrefsp}. Hence, ${\widehat{\rm cyc}}$ gives rise to a natural transformation 
$$
\widehat{\rm cyc}:{\rm Gr}\big(\pi_0({\rm Iso}({\rm Vect}^g_{\nabla}))\big)\longrightarrow
\widehat{\op{KO}}(-)\;,
$$
and precomposition with the unit of the adjunction $({\rm Gr}\dashv i)$ 
gives the unique transformation satisfying the conditions in Proposition \ref{dfcycmp}.
\end{proposition}
\theproof
The proof is essentially contained in the construction of the Hopkins-Singer $\widehat{\op{KO}}$-theory. 
The map $\widetilde{\rm Ph}$, given as part of the data in the differential function spectrum ${\rm diff}(\op{KO}, \widetilde{\rm Ph},\pi_*(\op{KO}))=:\widehat{\op{KO}}$, was defined by applying the topological realization 
$\delta\Pi$ to the morphism ${\rm Ph}$ in Proposition \ref{dfcycmp}. Thus, by definition, we have a homotopy commutative diagram 
\(\label{excycmpko}
\xymatrix{
{\bf K}{\op{O}}_{\nabla} \ar@{-->}[dr]^-{\widehat{\rm cyc}}\ar@/^2pc/[rrd]^-{\rm {Ph}}\ar@/_2pc/[ddr]_-{{\rm cyc}\circ \mathscr{F}}& & 
\\
&{\rm diff}\big(\op{KO}, \widetilde{\rm Ph},\pi_*(\op{KO})\big)\ar[d]^-{\mathcal{I}}\ar[r]^-{\mathcal{R}} &
\mathscr{H}\big(\tau_{\geq 0}\Omega^*(-;\pi_*(\op{KO})\big)\ar[d]
\\
& \op{KO}\ar[r]^-{\widetilde{\rm Ph}} & \mathscr{H}\big(\Omega^*(-;\pi_*(\op{KO}))\big)
\;,}
\)
where $\mathscr{F}$ is the canonical map induced by the unit $\eta:\mathbb{1}\to \delta \Pi$ (i.e., the map that forgets the geometric data). \footnote{Actually this map uses the geometric data to build a topological spectrum -- it is the constant sheaf of spectra on the topological realization. We have called it 'forgetful' in order to emphasize that it reduces to an underlying theory (in this case $\op{KO}$).} The dashed arrow is the map induced by the universal property and depends on the choices of homotopy filling the diagram. From the homotopy commutativity, we immediately read off the desired compatibility.
\endofproof

\begin{remark}[Independence of choice of homotopy]
Note that the universal map \eqref{excycmpko} is already unique up to a contractible choice (as part of the universal property). However, this uniqueness depends on a choice of homotopy filling the outer square in the diagram. Proposition \ref{dfcycmp} shows that after mapping out of a smooth manifold and passing to the induced map on connected components, this map is unique, independent of the choice of homotopy.
\end{remark}

From the above result, namely Proposition \ref{csconvb} and its proof, we see that we have completed the homotopy commutative diagram \eqref{shspdfko}. We now show that, upon topological realization, each of those sheaves 
of spectra gives a model for topological $\op{KO}$.
\begin{proposition}
[Equivalences via topological realization]
The topological realization of the maps in diagram \eqref{shspdfko} are equivalences.
\end{proposition}
\theproof
We have already shown in Proposition \ref{grekocn} that we have equivalences $\Pi(\mathbf{K}\op{O}_{\nabla})\overset{\mathscr{F}}{\simeq} \Pi(\mathbf{K}\op{O})\simeq \op{KO}$. By the  two-out-of-three property, 
it remains to show that the map $\mathcal{I}$ induces an equivalence upon topological realization. But this holds by construction. Indeed, the map $\mathcal{I}$ is defined as the unit $\mathbb{1}\to \delta\Pi$ applied to the sheaf of spectra $\widehat{\op{KO}}$. Since $\Pi\delta \Pi=\Pi$, it immediately follows that  $\Pi(\mathcal{I})=\mathbb{1}$.
\endofproof

We now consider the relationship of the cycle map with complex $K$-theory (see \cite{Bun10}\cite{Bun}). In \cite{Bun10}, it was shown by similar means that there is exists a cycle map
$$
\widehat{{\rm cyc}}_{\CC}:{\rm Vect}_{\nabla}^{\CC}\longrightarrow \op{K}(-)\;,
$$
which satisfies $\mathcal{R}(V,\nabla)={\rm Ch}(\nabla)$ and $\mathcal{I}(V,\nabla)={\rm cyc}^{\CC}$. A similar argument as the one in Proposition \ref{dfcycmp} shows that this map is also unique. At the topological level, we have commutative diagrams 
$$
\xymatrix{
{\rm Iso}({\rm Vect})\ar[d]_-{\otimes \CC}\ar[rr]^-{{\rm cyc}} && \op{KO}(-)\ar[d]^-{\otimes \CC}
\\
{\rm Iso}({\rm Vect}^{\CC})\ar[rr]^-{{\rm cyc}^{\CC}} && \op{K}(-)
}
\qquad  \quad\text{and} \qquad \quad
\xymatrix{
{\rm Iso}({\rm Vect})\ar[rr]^-{{\rm cyc}} && \op{KO}(-)
\\
{\rm Iso}({\rm Vect}^{\CC})\ar[rr]^-{{\rm cyc}^{\CC}}\ar[u]^-{F} && \op{K}(-)\ar[u]^-{F}\;,
}
$$
where the right vertical maps are induced by identifications of infinite loop spaces $\Omega^{\infty} \op{ KO}\simeq {\rm BO}\times \ZZ$ and $ \Omega^{\infty}\op{K}\simeq {\rm BU}\times \ZZ$, along with the complexification map $\otimes \CC: {\rm BO}(k)\to {\rm BU}(k)$ and the forgetful representation $F: {\rm BU}(k)\to {\rm BO}(2k)$, both 
stabilized. 

\begin{proposition}
[Refinement of cycle maps -- real vs. complex]\label{Prop-RvsC}
Let $c:\widehat{\op{KO}}(-)\to \widehat{\op{K}}(-)$ be the complexification map and $F:\widehat{\op{K}}\to \widehat{\op{KO}}$ be the forgetful map. We have commutative diagrams 
\(\label{forcpxdf}
\xymatrix{
{\rm Iso}({\rm Vect}_{\nabla})\ar[d]_-{\otimes \CC}\ar[rr]^-{\widehat{{\rm cyc}}} && 
\widehat{\op{KO}}(-)\ar[d]^-{c}
\\
{\rm Iso}({\rm Vect}^{\CC}_{\nabla})\ar[rr]^-{\widehat{{\rm cyc}}^{\CC}} &&
\widehat{\op{K}}(-)
}
\qquad \text{and} \qquad
\xymatrix{
{\rm Iso}({\rm Vect}_{\nabla})\ar[rr]^-{\widehat{{\rm cyc}}} && \widehat{\op{KO}}(-)
\\
{\rm Iso}({\rm Vect}_{\nabla}^{\CC})\ar[rr]^-{\widehat{{\rm cyc}}^{\CC}}\ar[u]^-{\mathscr F} && 
\widehat{\op{K}}(-)\ar[u]^-{F}\;.
}
\)
\end{proposition}
\theproof
With minor modification, the proof of Proposition \ref{dfcycmp} implies that there are unique maps
$$
\phi: \pi_0 {\rm Iso}({\rm Vect}_{\nabla})\longrightarrow  \widehat{\op{K}}(-)
\qquad \text{and} \qquad
\psi: \pi_0 {\rm Iso}({\rm Vect}_{\nabla}^{\CC})\longrightarrow \widehat{\op{KO}}(-)\;,
$$
satisfying $\mathcal{R}(\phi)={\rm Ch}(\nabla\otimes \CC)$, $\mathcal{I}(\phi)={\rm cyc}^{\CC}\otimes \CC$, $\mathcal{R}(\psi)={\rm Ph}(\mathscr{F}(\nabla))$ and $\mathcal{I}(\psi)={\rm cyc}\circ \mathscr{F}$. Since the compositions in diagrams \eqref{forcpxdf} satisfy the given properties, uniqueness implies commutativity.
\endofproof

\subsection{The differential Pontrjagin character}
\label{Sec-Diff-Pont-Cha}

In this section, we discuss the Pontrjagin character as a morphism of sheaves of spectra. Recall the discussion 
on the Deligne-Beilinson cup product in the latter part of Section \ref{Sec-Pont-class}. Given this product operation,
one has an induced product operation on the tensor product of graded rings 
$\widehat{H}^*(M;\mathbb{Q})\otimes \pi_*(\op{KO})$. Alternatively, one could consider taking the 
Hopkins-Singer differential refinement of the corresponding spectrum $\mathscr{H}(\mathbb{Q}\otimes \pi_*(\op{KO}))$. 
In particular, we have the differential function spectrum 
$$
\widehat{H}\big(\mathbb{Q}\otimes \pi_*(\op{KO})\big):=
{\rm diff}\big(\mathbb{Q}\otimes \pi_*(\op{KO}),i,\mathbb{R}\otimes \pi_*(\op{KO})\big)\;,
$$
where $i:\mathbb{Q}\otimes \pi_*(\op{KO})\into \mathbb{R}\otimes \pi_*(\op{KO})$ is the canonical 
inclusion. The resulting spectrum admits the structure of sheaf of $E_{\infty}$-ring spectra and one can 
ask whether this product structure is compatible with the Deligne-Beilinson product structure.
Indeed, this is the case.
\begin{lemma}
[Ring structure for differential cohomology with coefficients in $\pi_*(\op{KO})$]
We have a natural isomorphism of graded rings
$$
\widehat{H}^*(M;\QQ)\otimes \pi_*(\op{KO})\cong \widehat{H}^*(M;\QQ\otimes \pi_*(\op{KO}))\;,
$$
where the  ring structure on the left is induced by the Deligne-Beilinson cup product, and the ring 
structure on the right  is the canonical ring structure induced by the differential refinement of the
Eilenberg-MacLane  ring   spectrum $\mathscr{H}(\QQ\otimes \pi_*(\op{KO}))$ as discussed above.
\end{lemma}
\theproof
We first show that we have an isomorphism at the level of the underlying abelian groups. The refinement of 
$\QQ\otimes \pi_*(\op{KO})\cong \QQ[\alpha,\alpha^{-1}]$ can be presented by the sheaf of chain complexes
$$
{\rm diff}^n(\mathscr{H}\QQ[\alpha,\alpha^{-1}],i, \RR[\alpha,\alpha^{-1}])\simeq 
\mathscr{H}\Big({\rm cone}\big(\QQ[\alpha,\alpha^{-1}]\oplus \tau_{\geq n}\Omega^*[\alpha,\alpha^{-1}] \into \Omega^*[\alpha,\alpha^{-1}]\big)[-1]\Big)\;.
$$
By standard arguments in homological algebra (see e.g. \cite{GM}) one finds that the projection 
$$
\hspace{-2mm}
\xymatrix{
\pi:{\rm cone}\big(\QQ[\alpha,\alpha^{-1}]\oplus \tau_{\geq n}\Omega^*[\alpha,\alpha^{-1}]
\into \Omega^*[\alpha,\alpha^{-1}]\big)[-1]
\ar[rr]&&
{\rm cone}\big(\QQ[\alpha,\alpha^{-1}]\to \tau_{< n}\Omega^*[\alpha,\alpha^{-1}]\big)[-1]
}
$$ 
is an equivalence. Identifying the complex on the right degree-wise, we can organize this sheaf as follows
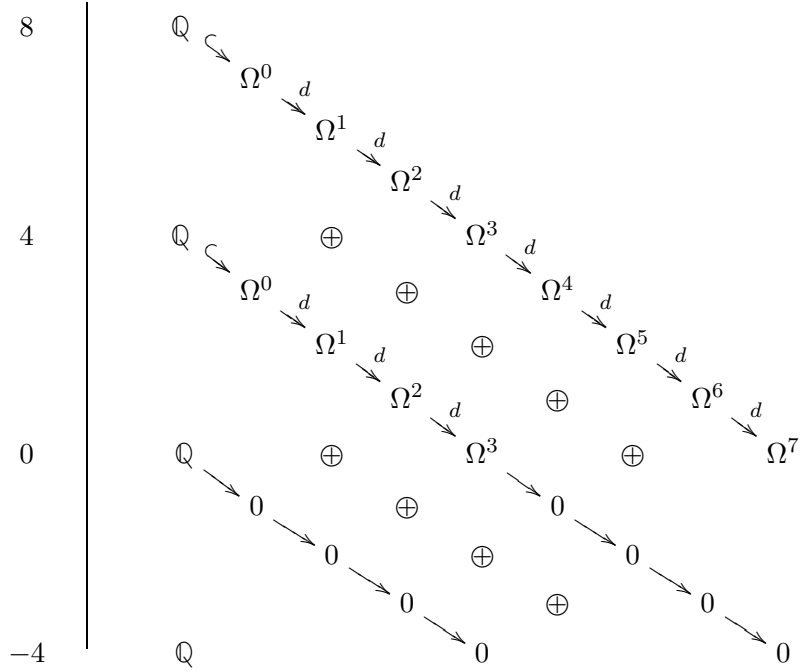
\vspace{-3mm}
$$
\xymatrix@C=8pt@R=2.5pt{
& \ar@{-}[ddddddddddddd]  & & & & & & & & & &  &
\\
8 & & & \QQ \; \ar@{^{(}->}[dr] & &  & &  & &  & &  &\; \; \;
\\
& & & & \Omega^0 \ar[dr]^-{d} & & & & & & & &
\\
& & & & & \Omega^1 \ar[dr]^-{d}& & & & & & &
\\
& & & & & & \Omega^2 \ar[dr]^-{d}& & & & & &
\\
4 & & & \QQ \; \ar@{^{(}->}[dr]& & \oplus & &   \Omega^3 \ar[dr]^-{d} & & & & & \; \; \; 
\\
& & & & \Omega^0 \ar[dr]^-{d}  & & \oplus & &  \Omega^4 \ar[dr]^-{d} & & & &
\\
& & & & & \Omega^1 \ar[dr]^-{d}& & \oplus & & \Omega^5 \ar[dr]^-{d}  & & &
\\
& & & & & &  \Omega^2 \ar[dr]^-{d}& & \oplus & & \Omega^6 \ar[dr]^-{d} & &
\\
0 & & &  \QQ \ar[dr] & & \oplus & & \Omega^3 \ar[dr]& & \oplus & & \Omega^7 & &
\\
&& &    &   0  \ar[dr] & & \oplus  & &   0 \ar[dr] &  & & & 
\; \; \;  
\\
& & & & &   0   \ar[dr] & & \oplus & &  0  \ar[dr] & & &
\; \; \; 
\\
& & & & & &  0   \ar[dr]  & & \oplus & &  0  \ar[dr] & &
\; \; \;
\\
{-4}  & & & \QQ& & &   &    0 & & & &  0   &
\; \; \;
\\
& & & & & & & & & & &  &
}
$$

\vspace{-3mm}
\noindent
where the diagonal arrows represent the action of the differential. This, in turn is immediately identified with the 
direct sum of the Deligne complexes ${\mathcal D}(4n,{\QQ})={\rm cone}(\QQ\into \tau_{< 4n}\Omega^*)[-1]$. 
Passing to cohomology gives the identification at the level of the underlying graded abelian groups. To see that this isomorphism preserves the ring structure observe that, since the isomorphism came from an identification at the 
level of sheaves of complexes, it is enough to note that there is (up to equivalence) a unique graded commutative 
ring structure which refines the underlying ring structures. Thus the map
$$
\xymatrix{
\cup:\bigoplus_{n}{\mathcal D}(4n,\QQ)\otimes \bigoplus_{m}{\mathcal D}(4m,\QQ)\cong
\bigoplus_{n,m}{\mathcal D}(4n,\QQ)\otimes {\mathcal D}(4m,\QQ)
\ar[rr]^-{\cup_{\rm DB}}&& \bigoplus_{n,m}{\mathcal D}(4n+4m,\QQ)}
$$
must present the canonical product structure on $\widehat{H}^*(M;\QQ\otimes \pi_*(\op{KO}))$. 
\endofproof

\medskip
Recall the definition of the \emph{differential} Pontrjagin character via the Deligne-Beilinson cup product 
(i.e. Definition \ref{Def-DiffPontChar}). We have not yet related this Pontrjagin character to a natural 
transformation of differential cohomology theories 
$\widehat{\op{KO}}(-)\to \widehat{H}^*(-;\QQ\otimes \pi_*(\op{KO}))$. In fact, we have the 
following characterization.

\begin{proposition}
[The differential Pontrjagin character via sheaves of spectra]\label{phchrfuni}
The Pontrjagin character ${\rm Ph}:\op{KO}\to \mathscr{H}(\pi_*(\op{KO})\otimes \QQ)$ admits a differential 
refinement at the level of sheaves of spectra
\(\label{dfphmp}
\widehat{\rm Ph}:\widehat{\op{KO}}
\xymatrix{\ar[r]&}
\widehat{\mathscr{H}}(\pi_*(\op{KO})\otimes \QQ)={\rm diff}(\mathscr{H}(\QQ\otimes \pi_*(\op{KO})),i,\RR\otimes \pi_*(\op{KO}))\;,
\)
such that the underlying natural transformation $\widehat{{\rm Ph}}:\widehat{\op{KO}}(-)\to \widehat{H}^*(-;\pi_*(\op{KO})\otimes \QQ)$ is uniquely characterized by the following properties
\begin{enumerate}[{\bf (i)}]
\item \emph{(Curvature)} $\mathcal{R}(\widehat{\rm Ph}\circ \widehat{{\rm cyc}}(V,\nabla))
=\widehat{\rm Ph}(\mathcal{F}_{\nabla})$. 
\item  \emph{(Topological realization)} $\mathcal{I}(\widehat{\rm Ph}\circ \widehat{{\rm cyc}}(V,\nabla))
={\rm Ph}(V)$.
\end{enumerate}
\end{proposition}

\newpage 

\theproof
We first show that there exists a refinement. From the definition of $\widehat{\op{KO}}$ and the fact that the
Pontrjagin character form takes rational periods on cycles (i.e., it represents an element in the image of 
$H^*(M;\mathbb{Q})\into H^*(M;\mathbb{R})$ term-wise), we have the homotopy commutative diagram
\(
\label{cmdiagph}
\xymatrix{
& \mathscr{H}(\tau_{\geq 0}\Omega^*(-;\pi_*(\op{KO})))\ar[d]\ar@{=}[r] 
& 
\mathscr{H}(\tau_{\geq 0}\Omega^*(-;\pi_*(\op{KO})))\ar[d]
\\
\mathbf{K}{\op{O}}_{\nabla}\ar[ru]^-{{\rm Ph}}\ar[r] \ar[dr]_-{\eta_{\mathbf{K}{\op{O}}_{\nabla}}} 
& \mathscr{H}(\Omega^*(-;\pi_*(\op{KO})))\ar@{=}[r] &
\mathscr{H}(\Omega^*(-;\pi_*(\op{KO})))
\\
& \op{KO}\ar[r]^-{\widetilde{\rm Ph}}\ar[u]_-{i \widetilde{\rm Ph}\simeq \Pi({\rm Ph})} &
\mathscr{H}(\pi_*(\op{KO})\otimes \QQ)\ar@{^{(}->}[u]_-{i}
\;.}
\)
By the universal property, this induces a morphism of differential ring spectra 
$$
\xymatrix{
\mathbf{K}{\op{O}}_{\nabla}\ar[r] & {\rm diff}\big(\op{KO},{\rm Ph},\pi_*(\op{KO})\big)
\ar[r]^-{\widehat{\rm Ph}}&
{\rm diff}\big(\mathscr{H}(\pi_*(\op{KO})\otimes \RR),{\rm id},\pi_*(\op{KO})\otimes \RR\big)}.
$$ 
By the definition of the Hopkins-Singer $\widehat{\op{KO}}$, the map $\widehat{{\rm Ph}}\circ \widehat{{\rm cyc}}$ satisfies the desired properties. We need only prove uniqueness of the underlying natural transformation. Arguing as in the proof of Proposition \ref{hclvcnko} proves the claim. 
\endofproof

\begin{corollary}
[Uniqueness of the differential Pontrjagin character]\label{dfpchcomp}
The differential Pontrjagin character given in Definition \ref{Def-DiffPontChar} is compatible with the unique transformation $\widehat{{\rm Ph}}:\widehat{\op{KO}}(-)\to \widehat{H}^*(-;\pi_*(\op{KO})\otimes \QQ)$ 
satisfying the properties of Proposition \ref{phchrfuni}. More precisely, Definition \ref{Def-DiffPontChar}
is given by precomposing $\widehat{\rm Ph}$ with the cycle map 
$$
\widehat{{\rm cyc}}:\pi_0{\rm Iso}({\rm Vect}_{\nabla})\longrightarrow \widehat{\op{KO}}(-)\;.
$$
\end{corollary}
\theproof
We only need to verify that the differential Pontrjagin character of Definition \ref{Def-DiffPontChar} 
satisfies the conditions of Proposition \ref{phchrfuni}. This, however, is clear since the Deligne-Beilinson 
cup product refines both the wedge product of forms and the cup product.
\endofproof

Given the compatibility provided by Corollary \ref{dfpchcomp}, we will henceforth not distinguish between 
various manifestations of the differential Pontrjagin character. We will only speak of \emph{the} unique 
differential Pontrjagin character. The differential Pontrjagin character satisfies similar properties as its 
topological counterpart (from Remark \ref{RemPropPont}).
\begin{proposition}
[Properties of  the differential Pontrjagin character]
\label{Prop-propPont}
Let $\widehat{{\rm cyc}}:{\rm Vect}_{\nabla}\to \widehat{\rm KO}(-)$ be the unique differential refinement of 
the cycle map which sends a vector bundle $V$ with connection $\nabla$ to its associated differential 
cohomology class. The differential Pontrjagin character $\widehat{\rm Ph}$ satisfies the following properties:
\begin{enumerate}[{\bf (i)}]
\item  \emph{(Naturality)} $\widehat{{\rm Ph}}$ is a natural transformation of differential cohomology theories.
\item   \emph{(Compatibility with the differential Chern character)} It is given by the composition
$$
\widehat{\rm Ph}:\widehat{\op{KO}}(-)\xrightarrow{\;\otimes \CC\;} \widehat{\op{K}}(-)
\xrightarrow{\;\widehat{\rm Ch}\;} \widehat{H}^*(-;\QQ[u,u^{-1}])\;.
$$
\item    \emph{(Correspondence)} Under the canonical maps $\mathcal{R}:\widehat{\op{KO}}(-)\to \Omega^*(-;\pi_*(\op{KO}))$ 
and $\mathcal{I}:\widehat{\op{KO}}(-)\to \op{KO}(-)$, we have 
$$
\mathcal{I}(\widehat{\rm Ph}\circ \widehat{{\rm cyc}}(V,\nabla))={\rm Ph}(V)
\qquad \text{and} \qquad 
\mathcal{R}(\widehat{\rm Ph}\circ \widehat{{\rm cyc}}(V,\nabla))={\rm Ph}(\mathcal{F}_{\nabla})\;.
$$
\item  \emph{(Whitney sum and tensor product)} It satisfies 
$$
\widehat{\rm Ph}((E,\nabla)\otimes (E^{\prime},\nabla^{\prime}))=\widehat{\rm Ph}(E,\nabla)\cup_{\rm DB} \widehat{\rm Ph}(E^{\prime},\nabla^{\prime})
$$
and
$$
\widehat{\rm Ph}((E,\nabla)\oplus (E^{\prime},\nabla^{\prime}))=\widehat{\rm Ph}((E,\nabla))+ \widehat{\rm Ph}(E^{\prime},\nabla^{\prime})\;.
$$
\end{enumerate}
\end{proposition}
\theproof
{\bf (i)} This follows immediately from the fact that $\widehat{\rm Ph}$ comes from a morphism of sheaves of spectra.

\noindent  {\bf (ii)} This follows immediately from the uniqueness clause of Proposition \ref{phchrfuni} along with the usual properties of the Pontrjagin character.  

\noindent {\bf (iii)} This  is a restatement of part of Proposition \ref{phchrfuni}.

\noindent  {\bf (iv)} This follows from uniqueness, the fact that the Deligne-Beilinson cup product refines the wedge product and the cup product, and the corresponding identities which are satisfied by the Pontrjagin character and the
Pontrjagin character form.
\endofproof

\section{Further  constructions in differential KO-theory}
\label{Ch-further}

\subsection{Differential $\op{KO}$-theory with compact support}
\label{Sec-compact}

In order to discuss orientation and pushforward, we will need to introduce differential $\op{KO}$-theory with 
compact support. Let $\widetilde{S}^n$ be the smooth sheaf in $\sh_{\infty}(\mathscr{M}{\rm an})$ obtained 
by quotienting the smooth $n$-dimensional disk $D^n$ by its boundary (see \cite{GS-cob} for a detailed discussion 
on this approach). The result is (of course) not a smooth manifold, but is a smooth sheaf. Similarly, given any vector 
bundle $V\to M$, we can form the Thom space ${\rm Th}(V)\to M$ as a quotient of the unit disk bundle by the sphere bundle, taken in the category of smooth sheaves.
\footnote{Note that one can also define more general Thom stacks as in \cite{GS-cob}.}
Since $\widehat{\op{KO}}$ is represented by a sheaf of spectra, we can make sense of the differential 
$\op{KO}$-theory of any smooth stack. More precisely, if $X$ is a smooth stack, we can define
$$
\widehat{\op{KO}}^n(X):=\pi_0\map(X_+,\widehat{\op{KO}}_n)\;,
$$
where $+$ is the functor which formally adjoins a basepoint and $\widehat{\op{KO}}_n:={\rm diff}(\Sigma^{-n}\op{KO},\widetilde{{\rm Ph}},\pi_{*+n}(\op{KO}))$ 
\footnote{We will henceforth always refer to this spectrum by the less burdensome notation 
$\widehat{\op{KO}}_n$.} We make the following natural definition analogous to the topological setting.

\begin{definition}[Differential KO-theory with vertical supports]
We define differential $\op{KO}$-theory with vertical compact support as 
$$
\widehat{\op{KO}}_{c}(V):=\widehat{\op{KO}}({\rm Th}(V))\;.
$$
\end{definition}

From the definition, it is clear that we can make the same definition for any sheaf of spectra. The sheaf of differential 
$n$-forms $\Omega^n$ is canonically pointed via the zero form. It admits a strict abelian group structure via the 
addition of differential $n$-forms and hence gives rise to a sheaf of spectra. The next example shows how our 
definition can be used to recover forms with compact support. 

\begin{example}[Differential forms with vertical compacts support on the Thom space]
\label{Ex-relform}
Let $V\to M$ be a real vector bundle on $M$. Then $\Omega^n({\rm Th}(V))$ is the group of differential 
forms with vertical compact support on $V$. Indeed, since $\Omega^n$ is zero-truncated, we have 
$$
\Omega^n({\rm Th}(V)):=\pi_0\map({\rm Th}(V),\Omega^n)=\hom(D(V)/S(V),\Omega^*)=
\left\{\vcenter{\xymatrix{P \ar[d] \ar[r] & \Omega^n(D(V))\ar[d]^-{r}
\\
\ast  \ar[r]^-{0} &\Omega^n(S(V))
}}\right\},
$$
where the diagram is a pullback and $r$ denotes the restriction to the boundary. \footnote{Note that for the set of \emph{pointed} maps $\hom(\ast,\Omega^*)\cong \ast$.} The pullback $P$ is immediately identified with the set of forms on the disk bundle which restrict to zero on the boundary.  Note that this can be rephrased as being the relative forms (relative to the boundary). 
\end{example}

The above definition of the smooth Thom space  is the analogue of the classical Thom space construction in the 
topological setting. In fact, topological $\op{KO}$-theory with compact support can be defined analogously via
$$
\op{KO}_c(V):=\op{KO}({\rm Th}(V))\cong \op{KO}(D(V)/S(V))\cong 
\pi_0\map(D(V)/S(V),\Omega^{\infty}\op{KO})\;.
$$
Fix an equivalence of infinite loop spaces $\Omega^{\infty}\op{KO}\simeq \op{BO}\times \ZZ$, where the latter 
is an $H$-space with respect to the additive structure induced by Whitney sum of vector bundles. Then for a pair 
of vector bundles $E,F\to V$ we have the ``virtual difference map" $E \ominus F:V\to (\op{BO}\times \ZZ)_{\oplus}$. Then this map descends to a well-defined map on the Thom space if and only if the map $E\ominus F$ is nullhomotopic 
on the boundary of the unit disk bundle $D(V)$. In this case, the map $E$ must be homotopic to $F$ on the boundary 
and there is an isomorphism of corresponding bundles $f:E\cong F$ on $\partial D(V)$.  Thus we recover the more classical definition of the group $\op{KO}_c(V)$, as the group of equivalence classes of triples $[E,F,f]$, where $E,F\to V$ are vector bundles and $f:E\to F$ is a bundle isomorphism outside some compact subset of $V$. The topological case is compatible with the geometric description in the following way.

\begin{proposition}[${\rm KO}_c$-theory of topological vs. smooth bundles]\label{prop-top-vs-sm}
Let $\pi:V\to M$ be a locally trivial, real vector bundle on a smooth compact manifold $M$. Let $\pi^{\prime}:V^{\prime}\to M^{\prime}$ denote the underlying topological vector bundle. Let $\op{KO}$ be a spectrum representing $\op{KO}$-theory and consider the corresponding geometrically discrete sheaf of spectra $\delta(\op{KO})$. Then we have an isomorphism 
$$
\op{KO}_c(V^{\prime})\cong \delta(\op{KO})_c(V)\;.
$$
\end{proposition}
\theproof
By the adjunction $\vert \cdot \vert \circ  \Pi \dashv \delta\circ {\rm sing}$, we have an equivalence 
$$
\map(\vert \Pi{\rm Th}(V)\vert, \Omega^{\infty}\op{KO})\simeq \map({\rm Th}(V), \Omega^{\infty}\delta({\rm sing}(\op{KO})))\;,
$$
so it suffices to prove that we have a weak equivalence $\vert \Pi{\rm Th}(V)\vert\simeq {\rm Th}(V^{\prime})$. 
Since both $\Pi$ and $\vert \cdot \vert$ are left adjoints, they preserve quotients. Moreover, for a smooth manifold 
$M$, $\vert \Pi(M) \vert$ can be identified with the topological realization of the smooth singular nerve of $M$
and it is well-known that the smooth singular nerve is equivalent (as a Kan-complex) to the topological nerve. 
These observations yield a commutative  diagram
$$
\xymatrix{
S(V^{\prime})\ar@{^{(}->}[r]& D(V^{\prime}) \ar[r] & {\rm Th}(V^{\prime})
\\
\vert {\rm sing}(S(V^{\prime}))\vert \ar@{^{(}->}[r]\ar[u]^-{\simeq}  & \vert {\rm sing}(D(V^{\prime}))\vert \ar[r]\ar[u]^-{\simeq}  & \vert {\rm sing}({\rm Th}(V^{\prime}))\vert \ar[u]^-{\simeq} 
\\
\vert \Pi(S(V))\vert \ar@{^{(}->}[r]\ar[u]^-{\simeq} & \vert \Pi(D(V))\vert \ar[r]\ar[u]^-{\simeq} & 
\vert \Pi({\rm Th}(V))\vert\simeq \vert \Pi(D(V))\vert /\vert \Pi(S(V))\vert \ar[u]^-{\simeq}
\;,}
$$
where the top vertical maps are induced by the counit of the equivalence $\vert \cdot \vert \dashv {\rm sing}$, the 
bottom maps are induced by the natural inclusion of the smooth singular nerve and the rows are cofiber sequences. 
\footnote{Note that for the equivalence between the second and third rows none of the 
topological realizations are needed, but they are required for matching the first row.}  
This gives the desired equivalence. 
\endofproof

In the subsequent section we will construct an explicit differential refinement of the Thom class in $\op{KO}$-theory. Not surprisingly, for a real vector bundle $V\to M$, this refinement is related to the geometry of the spinor bundle $\mathscr{S}(V)\to M$ and the explicit class in $\widehat{\op{KO}}$-theory will arise via the cycle map 
$\widehat{\rm cyc}:\op{{\bf K}O}_{\nabla}\to \widehat{\op{KO}}$. We will need the following explicit description of $\op{{\bf K}O}_{\nabla,c}(V):=\op{{\bf K}O}_{\nabla}({\rm Th}(V))$. 

\begin{definition}[Equivalence of real vector bundles with connections]
Let $V\to M$ be a locally trivial real vector bundle over a compact manifold $M$. Consider the collection of pairs 
of vector bundles with connection $(E,\nabla_{E}),(F,\nabla_{F})$ on $V$ along with isomorphisms of bundles \emph{with connection} $\sigma:(E,\nabla_{E})\to (F,\nabla_{F})$, defined on the boundary of the unit disk 
bundle $D(V)$.

\item {\bf (i)}  We say that two such triples $((E,\nabla_{E}),(F,\nabla_{F}),\sigma)$ and $((E^{\prime},\nabla_{E^{\prime}}),(F^{\prime},\nabla_{F^{\prime}}),\sigma^{\prime})$ are \emph{equivalent}
if there exists two other bundles $(G,{\nabla_{G}})$ and $(H,{\nabla_{H}})$ and isomorphisms
$$
(E,\nabla_{E})\oplus (G,{\nabla_{G}})\cong (E^{\prime},\nabla_{E^{\prime}})\oplus (H,{\nabla_{H}})
\qquad
\text{and}
\qquad
(F,\nabla_{F})\oplus (G,{\nabla_{G}})\cong  (F^{\prime},\nabla_{F^{\prime}})\oplus (H,{\nabla_{H}})
$$
which are compatible with the isomorphisms $\sigma$ and $\sigma^{\prime}$ on the boundary. 

\item {\bf (ii)}We denote the set of equivalence classes of such pairs as $L_1(D(V),S(V))$ and write its elements as 
$$[(E,\nabla_{E}),(F,F_{\nabla}),\sigma]\in L_1(D(V),S(V))\;.$$
\end{definition}

We provide a model for differential real K-theory with compact supports using an analogue of 
Atiyah's approach to K-theory found in \cite[p. 88]{At-K}. Hence this can be viewed as 
a differential refinement of the real counterpart of Atiyah's treatment.  

\begin{proposition}
[Differential KO-theory with compact supports as a refinement of the Atiyah model]
\label{expdkonab}
We have an isomorphism 
$$
L_1(D(V),S(V))\cong \op{{\bf K}O}_{\nabla,c}(V)\;.
$$
\end{proposition}
\theproof
By definition and properties of the mapping space functor, we have  isomorphisms
\bea
\op{{\bf K}O}_{\nabla,c}(V):=\op{{\bf K}O}_{\nabla}({\rm Th}(V))
&\cong& 
\op{{\bf K}O}_{\nabla}(D(V)/S(V)) 
\\
&\cong&
\pi_0\map(D(V)/S(V),\Omega^{\infty}\op{{\bf K}O}_{\nabla})
\\
&\cong & \pi_0(P(V))
\eea
where $P(V)$ is defined as the $\infty$-pullback 
$$
\left\{\vcenter{\xymatrix{
P(V) \ar[d] \ar[rr]^-{j} && \map(D(V),\Omega^{\infty}\op{{\bf K}O}_{\nabla})\ar[d]^-{i^*}
\\
\ast  \ar[rr]^-{0} && \map(S(V),\Omega^{\infty}\op{{\bf K}O}_{\nabla})
}}\right\}\;.
$$
Computing the homotopy pullback via the local formula, we find that an element in $P(V)$ is a vertex in 
$\map(D(V),\Omega^{\infty}\op{{\bf K}O}_{\nabla})$ whose restriction to $S(V)$ is nullhomotopic 
through an edge in $\map(S(V),\Omega^{\infty}\op{{\bf K}O}_{\nabla})$. Define the map 
$L_1(D(V),S(V))\to \pi_0(P)$ as follows.  With $\Omega^{\infty}\op{{\bf K}O}_{\nabla}$ 
representing the Grothendieck group completion of vector bundles with connection (see the Appendix), 
for each triple $[(E,\nabla_{E}),(F,F_{\nabla}),\sigma]$ choose a representative of 
$[(E,\nabla_{E})]-[(F,F_{\nabla})]$ in $\map(D(V),\Omega^{\infty}\op{{\bf K}O}_{\nabla})$. 
Since $\sigma$ defines an isomorphism on the boundary, it follows that, after composing with the
restriction map $i^*$, this representative must be nullhomotopic,
\footnote{The nullhomotopy depends on $\sigma$ uniquely up to a contractible space. This follows from explicit identification of the edges
in the mapping space. Such identifications are routine but tedious, so we will not spell out the details.}
giving rise to an element in $P$. This is not a well-defined map; however, two such choices of 
representatives for $[(E,\nabla_{E})]-[(F,F_{\nabla})]$ are homotopic in 
$\map(D(V),\Omega^{\infty}\op{{\bf K}O}_{\nabla})$. Let $h$ denote such a homotopy between 
two representatives. Then, if $\overline{\sigma}$ is the homotopy induced by $\sigma$, the induced 
nullhomotopy $\overline{\sigma}\circ i^*(h)$ gives rise to the same class in $\pi_0(P)$. Thus, we 
have a well defined map $f:L_1(D(V),S(V))\to \pi_0(P)$. Clearly this map is surjective. It is also 
injective since if $f([(E,\nabla_{E}),(F,{\nabla}_F),\sigma])=0$, then $[(E,\nabla_E)]-[(F,\nabla_F)]=0$ and the nullhomotopy on the boundary must be in the class of ${\rm id}:(E,\nabla_E)\vert_{S(V)}= (F,\nabla_F)\vert_{S(V)}$. It follows at once that $[(E,\nabla_E),(F,\nabla_F),\sigma] =0$. 
\endofproof

The smash product 
$M_+\wedge \widetilde{S}^n$ can be identified with the quotient of $M\times D^n$ by the canonical 
map $M\times \partial D^n\into M\times D^n$; see \cite{GS-cob} for a detailed discussion.
Via Example \ref{Ex-relform}, one finds that the group 
$\Omega^k(M_+\wedge \widetilde{S}^n)$ is precisely the group of forms on $M\times D^n$ with vertical 
support in the interior of $D^n$. These forms can be integrated over $D^n$, giving rise to a map 
$$
\int_{D^n} :\Omega^k(M_+\wedge \widetilde{S}^n)\longrightarrow \Omega^{k-n}(M)
$$
for each $k$. Extending by 4-periodicity, this gives rise to an integration map 
$$
\int:\Omega^*(M_+\wedge \widetilde{S}^n;\pi_*(\op{KO}))\longrightarrow \Omega^{*-m}(M;\pi_*(\op{KO}))\;.
$$
Since the underlying topological theory $\op{KO}$ is homotopy invariant, it follows that we have an equivalence 
$\delta(\op{KO})^*(M\wedge \widetilde{S}^n)\simeq \op{KO}^*(\Sigma^nM)$, where $\Sigma^n$ 
is the usual suspension. These observations combine to give the following.
\begin{proposition}
[Differential desuspension map]\label{desprdf}
There is a desuspension map 
$$
\hat{\sigma}_n^{-1}:\widehat{\op{KO}}^*(M\wedge \widetilde{S}^n)\longrightarrow \widehat{\op{KO}}^{*-n}(M)
$$
in differential $\op{KO}$-theory, such that for each $x\in \widehat{\op{KO}}(M\wedge \widetilde{S}^n)$ 
we have 
\begin{enumerate}[{\bf (i)}]
\item \emph{Compatibility with curvature)} $\int\mathcal{R}(x)=\mathcal{R}(\hat{\sigma}^{-1}(x))$,
\item \emph{Compatibility with topological realization)} $\sigma_n^{-1}\mathcal{I}(x)=\mathcal{I}(\hat{\sigma}_n^{-1}(x))$, where $\sigma_n^{-1}: \op{KO}^*(\Sigma^n M)\to \op{KO}^{*-n}(M)$ is the usual 
desuspension map in the underlying topological KO-theory. 
\end{enumerate}
\end{proposition}
\theproof
Since the elements in $\Omega^*(-\wedge \widetilde{S}^n)$ have vertical support in the interior of the 
$n$-disk $D^n$, the fiber integration of forms gives rise to a morphism of sheaves of complexes
$$
\xymatrix{
\int_{D^n}:\Omega^*(-\wedge \widetilde{S}^n)\ar[r]&\Omega^{*-n}(-)
}.
$$
Consequently, this induces a morphism of sheaves of spectra
$$
\xymatrix{
\int:\map\big(-\wedge \widetilde{S}^n,\mathscr{H}(\Omega^*(-;\pi_*(\op{KO})))\big)\ar[r]&
\mathscr{H}\big(\Omega^{*-n}(-;\pi_*(\op{KO}))\big)\;.
}
$$
Let $\widetilde{\mathcal{I}}:\widetilde{S}^n\to S^n$ be the canonical map induced by topological realization, where $S^n$ is the 
singular nerve of the topological $n$-sphere. For any homotopy invariant spectrum $E$, the map $\widetilde{\mathcal{I}}$ induces 
an equivalence $\widetilde{\mathcal{I}}^*:\map(-\wedge S^n,E)\xrightarrow{\simeq} \map(-\wedge \widetilde{S}^n,E)$. This applies for both $\op{KO}$ 
and $\mathscr{H}(\pi_*(\op{KO})\otimes \RR)$. Fix a homotopy inverse $(\widetilde{\mathcal{I}}^*)^{-1}$ of the above equivalence. Define the map
$$
\xymatrix{
\tilde{\sigma}_n:\map(-\wedge \widetilde{S}^n,\op{KO})
\ar[rr]^-{(\widetilde{\mathcal{I}}^*)^{-1}} &&
\map(-\wedge S^n,\op{KO}) \; \overset{\sigma_n} \simeq \; \Sigma^{-n} \delta\big(\op{KO}\big)}.
$$
In a similar way, define $\tilde{\rho}_n$ as the composite 
$$
\hspace{-1mm}
\xymatrix{
\tilde{\rho}_n:\map(-\wedge \widetilde{S}^n,\mathscr{H}(\pi_*(\op{KO})\otimes \RR))
\ar[rr]^-{(\widetilde{\mathcal{I}}^*)^{-1}} &&
\map(-\wedge S^n,\mathscr{H}(\pi_*(\op{KO})\otimes \RR)) \; \overset{\rho_n} \simeq \; \Sigma^{-n} \delta\big(\mathscr{H}(\pi_*(\op{KO})\otimes \RR)\big)}\!,
$$
where $\rho_n$ is the equivalence induced by the adjoints of the maps
$\Sigma^n \mathscr{H}(\pi_*(\op{KO})\otimes \RR)_k\to \mathscr{H}(\pi_*(\op{KO})\otimes \RR)_{k+n}$ and $\mathscr{H}(\pi_*(\op{KO})\otimes \RR)$ is some fixed $\Omega$-spectrum model. Then we have a homotopy commutative diagram of spans
$$
\xymatrix{
\map\big(-\wedge \widetilde{S}^n,\mathscr{H}(\tau_{\leq 0}\Omega^*(-;\pi_*(\op{KO}))\big)\ar@/_7pc/[dd]
\ar[rr]^-{\int}_{\ }="x"\ar[d] &&
\mathscr{H}\big(\tau_{\leq 0}\Omega^{*}(-;\pi_*(\op{KO})[-n])\big)\ar[d]\ar@/^6pc/[dd]
\\
\map\big(-\wedge \widetilde{S}^n,\mathscr{H}(\Omega^*(-;\pi_*(\op{KO})))\big)\ar@/^1pc/[d]^-{j^{-1}}
\ar[rr]^-{\int}="w"_{\ }="v" && 
\Sigma^{-n}\delta\big(\mathscr{H}(\Omega^*(-;\pi_*(\op{KO})))\big)\ar@/_1pc/[d]_{j^{-1}}
\\
\map(-\wedge \widetilde{S}^n,\mathscr{H}(\pi_*(\op{KO})\otimes \RR)\ar[rr]_-{\ }="t"^-{\tilde{\rho}_n}="u"\ar[u]^-{j} && 
\Sigma^{-n}\delta\big(\mathscr{H}(\pi_*(\op{KO})\otimes \RR)\big)\ar[u]_-{j}
\\
\map(-\wedge \widetilde{S}^n,\op{KO})\ar[rr]^-{\tilde{\sigma}_n}="s"\ar[u]^-{\widetilde{{\rm Ph}}} && 
\Sigma^{-n}\delta\big(\op{KO}\big)\;.  \ar[u]^-{\widetilde{\rm Ph}}
\ar@{=>}^-{{\rm id}} "s"-<1cm,0cm>;"t"+<1cm,-.2cm>
\ar@{=>}^-{\nu} "u"-<1cm,0cm>;"v"+<1cm,-.2cm>
\ar@{=>}^-{\rm id} "w"-<1cm,0cm>;"x"+<1cm,-.2cm>
}
$$
To see this, note that we can take as an explicit presentation of $\op{KO}$ and $\mathscr{H}(\pi_*(\op{KO})\otimes \RR)$ by $\Omega$-spectra as the topological realization of the sheaves of spectra $\op{\mathbf{K}O}_{\nabla}$ and $\mathscr{H}(\tau_{\geq 0}\Omega^*(-;\pi_*(\op{KO}))$. With this presentation, $\widetilde{{\rm Ph}}$ is the topological realization of the Pontrjagin character form and is a morphism of $\Omega$-spectra. Then, by definition of such a map, the bottom square commutes strictly. The top square commutes by definition of the map $\int$. 

It remains to show that the middle square commutes up to the choice of homotopy $\nu$. To this end, observe that $j$ is induced by the inclusion $i: \pi_*(\op{KO})\otimes \RR\into \Omega^*(-;\pi_*(\op{KO}))$. Integration over the fiber induces the inverse for the Thom isomorphism in de Rham cohomology, and similarly $\tilde{\rho}_n$ induces an inverse in real cohomology. A choice of Thom class $\nu$ in real cohomology gives a homotopy inverse of $\tilde{\rho}_n$ and $i$ is multiplicative and compatible with Thom classes. It follows that the middle square commutes up to the homotopy defined by $\nu$. The curved maps $\mathscr{H}(\tau_{\leq 0}\Omega^*(-;\pi_*(\op{KO}))\big)\to  \mathscr{H}(\pi_*(\op{KO})\otimes \RR)$ in the diagram are obtained by inverting the equivalence $j$. It follows that the induced map on pullbacks gives the desired map and it is clear from the construction that this map satisfies the two properties in the statement.
\endofproof

\subsection{Orientation and genus in the differential case}
\label{Sec-or-gen}

In Section \ref{Sec-Pont-class}, we defined the differential $\hat{A}$-genus as 
a formal power series expansion in the differential Pontrjagin classes (see Definition \ref{Def-diffA}). We 
would like to connect this genus with a Thom isomorphism theorem for $\widehat{\op{KO}}$. We considered 
orientation and genera associated with the underlying topological theory in Section \ref{Sec-orgen}. 
A KO-orientation is a Spin structure on $M$ by \cite{ABS}, while a differential $\widehat{\op{KO}}$-orientation 
is a Spin structure together with a Riemannian metric (see \cite{FMS}
\footnote{There the cohomology theory was denoted ${\op{KO}^\vee}$.}). The latter gives rise to a preferred Levi-Civita connection, which in turn lifts to a preferred Spin connection.

\medskip
As in Bunke's description in the general case \cite{Bun}, 
we offer the following definition.
\footnote{Note that for $\widehat{\op{KO}}$ there is no distinction between orthogonal 
and real bundles, unlike what we have seen for $\mathbf{K}{\rm O}$ and $\mathbf{K}{\rm O}_\nabla$.}
\begin{definition}[Differential KO-theory Thom class]
Let $V\to M$ be an $n$-dimensional real vector bundle on $M$. A \emph{differential} Thom class $\hat{\nu}$ in $\widehat{\op{KO}}$-theory is a differential cohomology class $\hat{\nu}\in \widehat{\op{KO}}(V)$ with $\mathcal{I}(\hat{\nu})=\nu$ a Thom class in ${\op{KO}}$-theory, where $\mathcal{I}$ is topological realization.
\end{definition}

Let $\pi:V\to M$ be a real vector bundle of rank $n\equiv 0 \mod 4$. Recall from Example \ref{Ex-LMbasic}
in Section \ref{Sec-orgen} 
the class in $\op{KO}_c(V)$ given by the difference element $\mathbf{S}(V):=[\pi^*\mathscr{S}(V)^+, \pi^*\mathscr{S}(V)^-,\mu]$ determined by the isomorphism 
\(\label{isosppmb}
\mu:\pi^*\S^+(V)\overset{\cong}{\longrightarrow} \pi^*\S^-(V)\;,
\)
where $\mu$ is induced fiberwise by Clifford multiplication by a vector $ e \in \RR^{4k}$. 
For $n\equiv 0 \mod 8$, ${\bf S}(V)$ is the Thom class in $\op{KO}$-theory, giving rise to the Thom isomorphism 
$$
{\bf S}(V)\cup \pi^*(-):\op{KO}(M)\overset{\cong}{\longrightarrow} \op{KO}_c(V)\;.
$$
For $n\equiv 4 \mod 8$ this is \emph{not} an isomorphism; nevertheless it gives rise to a map of the above form. 
A more homotopy theoretic approach goes as follows (see \cite{LM}). 
Let $p:E\to {\rm BSpin}(n)$ be the canonical bundle 
with Thom space ${\rm MSpin}_n$. Let $\mathbf{S}(E)\in \op{KO}_c(E):=\widetilde{\op{KO}}({\rm MSpin}_n)$ 
be the difference class defined above. Then for a vector bundle $\pi:V\to M$, we have a corresponding classifying 
map $g_\pi:M\to {\rm BSpin}(n)$ and $V$ can be identified with the pullback of the universal bundle 
$E\to {\rm BSpin}(n)$. This gives rise to a map of Thom spaces $f_{\pi}:{\rm Th}(V)\to {\rm MSpin}_n$ 
and we have 
$$
\mathbf{S}(V)=f_{\pi}^*\mathbf{S}(E)
$$
in $\op{KO}_c(V):=\widetilde{\op{KO}}({\rm Th}(V))$. If we take $V$ to be the normal bundle of an embedding $i:M\into \RR^{n+8k}$, for some sufficiently large $k$, precomposition with the Pontrjagin-Thom collapse map $S^{n+8k}\to {\rm Th}(V)$ gives rise to the $\hat{A}$-genus
$$
\hat{A}:\Omega^{\rm Spin}_n\longrightarrow 
\widetilde{\op{KO}}^0(S^{n+8k})\cong \widetilde{\op{KO}}^{-n}({\rm pt})\;,
$$
where $\Omega^{\rm Spin}_n$ is the Spin cobordism group in degree $n$.

\medskip
We now discuss the differential refinement. Let $\nabla$ be a connection on $V$ and let $\widetilde{\nabla}$ 
be the corresponding lift to a Spin connection, acting on sections of the Spin bundle 
$\mathscr{S}(V)\to M$. Since the splitting $\mathscr{S}(V)\cong \mathscr{S}(V)^+\oplus \mathscr{S}(V)^-$ 
came from the splitting of the spinor representation $\Delta_n$ into irreducibles $\Delta_n^\pm$, it
follows that $\widetilde{\nabla}$ preserves the $\ZZ/2$-grading. Consequently,  $\widetilde{\nabla}$ splits 
as a direct sum $\widetilde{\nabla}^+\oplus \widetilde{\nabla}^-$, with $\widetilde{\nabla}^{\pm}$ defined 
on the corresponding bundles $\mathscr{S}(V)^{\pm}\to M$. 
Now the map $\mu$ is not just an isomorphism of bundles, but we have the following.

\begin{lemma}[Clifford mutiplication respects connections]\label{lemm-cliff}
The Clifford multiplication map $\mu$ is 
an isomorphism of bundles \emph{with connection}. 
\end{lemma}
\theproof
Write 
a section $\psi$ of $\mathscr{S}(V)=\mathscr{S}(V)^+\oplus \mathscr{S}(V)^-$ as $(\psi^+,\psi^-)$. 
Consider the section $s:V\to \pi^*\mathscr{S}(V)$ defined by $s(v)=v$. Then, on the boundary of the unit disk 
$D(V)$, we have $s(v)\cdot s(v)=v^2=-1$ for all $v\in \partial D(V)$. Then the Leibniz rule with respect to 
the Clifford product implies that 
$$
0=\pi^*\widetilde{\nabla}(s\cdot s)=2(\pi^*\widetilde{\nabla}(s)\cdot s)\;,
$$
so that $s$ is parallel with respect to the pullback connection $\pi^*\widetilde{\nabla}$. From this, we have
\begin{eqnarray*} 
\mu((\pi^*\widetilde{\nabla}^+)(\pi^*\psi^+))&:=&
s\cdot ((\pi^*\widetilde{\nabla})(\pi^*\psi^+,0)) 
\\
&=& \pi^*\widetilde{\nabla}((0,s\cdot \pi^*\psi^+))
\\
&=&
\pi^*\widetilde{\nabla}^-(\mu(\pi^*\psi^+))\;,
\end{eqnarray*}
so that $\mu$ is indeed an isomorphism of bundles with connections.
\endofproof

From Lemma \ref{lemm-cliff}, we can consider the class in ${\bf K}{\rm O}_{\nabla}(M)$ 
given by the virtual difference (see Proposition \ref{expdkonab})
\(
\label{SVnablatilde}
\mathbf{S}(V)_{\widetilde{\nabla}}:=
\big[(\pi^*\mathscr{S}(V)^+,\pi^*\widetilde{\nabla}^+),
(\pi^*\mathscr{S}(V)^-,\pi^*\widetilde{\nabla}^-),\mu\big]\;.
\)
This gives a natural candidate for the refinement of the Thom class $\mathbf{S}(V)$, namely 
$$
\widehat{\mathbf{S}(V)}:=\widehat{{\rm cyc}}
\big([\pi^*(\mathscr{S}(V)^+,\widetilde{\nabla}^+),\pi^*(\mathscr{S}(V)^-,\widetilde{\nabla}^-),\mu]\big)\;,
$$
where $\widehat{\rm cyc}:\mathbf{K}{\rm O}_{\nabla}({\rm Th}(V))\to \widehat{\op{KO}}({\rm Th}(V))=:\widehat{\op{KO}}_c(V)$ is the cycle map. 
\footnote{Classical KO-theory of Thom spaces is studied recently in \cite{Do}.}
Recall from Lemma \ref{Lemm-a-roof} (stated more generally for the differential version) 
that $\widehat{A}(\mathcal{F}_{\nabla})$ is a unit in the graded ring of differential forms.
The following is known, although the precise statement usually does not take this form
(see, e.g. \cite[V, Theorem 4.14]{Ka} for the topological part). 
%

\begin{proposition}
[Formula for $\widehat{A}$-genus in terms of Pontrjagin characters] \label{Prop-formul}
Let $\pi: V \to M$ be a real vector bundle of rank $n\equiv 0 \mod 8$ with connection $\nabla$
and suppose $V$ admits a Spin structure with a lift $\widetilde{\nabla}$. 
Let $\hat{A}(\mathcal{F}_\nabla)$ denote the $A$-roof form, defined via Chern-Weil theory. 
\item {\bf (i)} At the level of the underlying de Rham cohomology classes, we have 
$$
\left[\int_{V/M}{\rm Ph}\big(\mathcal{F}_{\pi^*\widetilde{\nabla}^+}\big)
-{\rm Ph}\big(\mathcal{F}_{\pi^*\widetilde{\nabla}^-}\big)\right]=
\big[\big(\hat{A}(\mathcal{F}_{\nabla})\big)^{-1}\big]\;.
$$
\item {\bf (ii)} Hence, at the level of differential forms with vertical compact support, 
there is a unique element 
$\eta\in \Omega_c^{-1}(V;\pi_*(\op{KO}))/{\rm im}(d)$, such that 
$$
\int_{V/M}{\rm Ph}\big(\mathcal{F}_{\pi^*\widetilde{\nabla}^+}\big)-
{\rm Ph}\big(\mathcal{F}_{\pi^*\widetilde{\nabla}^-}\big)+d(\eta)=
\big(\hat{A}(\mathcal{F}_{\nabla})\big)^{-1}\;.
$$
\end{proposition}
\theproof
This follows immediately from the topological formula
$$
\pi_{V/M}{\rm Ph}\big([\pi^*\mathscr{S}(V)^+,\pi^*\mathscr{S}(V)^-,\mu]\big)
=\hat{A}(TM)^{-1}\;,
$$
where $\pi_{V/M}$ denotes the pushforward in de Rham cohomology, along with Proposition \ref{Prop-propPont}.
\endofproof

Note that the above statement could alternatively have been deduced from 
the well-known relation between the 
$\widehat{A}$-genus and the Todd genus (see \cite{cap} for such explicit 
combinatorial relations). 

\begin{remark} [Relation between Thom classes in versions of differential KO-theory]
Let $V\to M$ be a real vector bundle admitting Spin structure, equipped with connection $\nabla$. 
By definition,
$$
\mathcal{R}(\widehat{\mathbf{S}(V)})=
{\rm Ph}\big(\mathcal{F}_{\pi^*\widetilde{\nabla}^+}-\mathcal{F}_{\pi^*\widetilde{\nabla}^-}\big)\;.
$$
It then follows immediately from part ${\bf (ii)}$ of Proposition \ref{Prop-formul}
that the class $\widehat{\mathbf{S}(V)}$ is a differential refinement of $\mathbf{S}(V)$ which satisfies
$$
\int_{V/M}\mathcal{R}(\widehat{\mathbf{S}(V)}) = \hat{A}(\mathcal{F}_{\nabla})^{-1}+d\eta\;,
$$
where $\eta$ is a formal power series of forms $\eta=\eta_1+\eta_3+\hdots$ each with vertical compact support. 
\end{remark}

It turns out, as  we will see in Proposition \ref{Prop-uni-spin} below, that up to the appropriate notion of equivalence of Thom classes, this is the unique refinement which is \emph{natural} under pullbacks of bundles admitting a Spin structure.  
\begin{definition}[Equivalence of differential Thom classes]
Two differential Thom classes $\hat{\nu}$ and $\hat{\nu}^{\prime}$ on a bundle $V\to M$ are said to be \emph{homotopic} if there is a differential Thom class $h$ on the pullback bundle 
\(\label{pull-dfth}
\xymatrix@R=1.5em{
V\times I\ar[rr] \ar[d]&& V\ar[d]
\\
M\times I\ar[rr]^-{\rm pr} && M
}
\)
whose restriction to the endpoints of the interval $I$ are $\hat{\nu}$ and $\hat{\nu}^{\prime}$, 
respectively, 
and such that 
$$
\int_{V\times I/M\times I}\mathcal{R}(h)={\rm pr}^*\Big(\int_{V/M}\mathcal{R}(\hat{\nu})\Big)\;.
$$
\end{definition}
For example, in our case, $\mathcal{R}(\widehat{\mathbf{S}(V)})$ is the inverse of the square of the 
$\widehat{A}$-genus. 
If we do not impose naturality with respect to pullback, we can have many Thom classes, as we now illustrate. 
\begin{proposition}[Set of differential Thom classes]
\label{unithclko}
For $\pi:V\to M$ a vector bundle of rank $n\equiv 0 \mod 8$,  the set of homotopy classes of differential Thom classes refining a given Thom class $\nu$ and genus form are parametrized by the group
$$
\mathscr{O}{\rm r}(\pi)=
\frac{H^{-1}(M;\pi_*(\op{KO})\otimes \RR)}{\hat{A}(V)\cup {\rm Ph}({\rm KO}^{-1}(M))}\;.
$$
\end{proposition}
\theproof
Let $\hat{\nu}$ and $\hat{\nu}^{\prime}$ be two differential Thom classes refining a fixed Thom class $\nu$ which satisfy
$$
\int_{V/M}\mathcal{R}(\hat{\nu})=\int_{V/M}\mathcal{R}(\hat{\nu}^{\prime}).
$$
Then there is a fiberwise compactly supported form $\eta\in \Omega_c^{-1}(V)$, 
unique up to elements in the image of ${\rm Ph}:\op{KO}_c^{-1}(V)\to H^{-1}(V;\pi_*(\op{KO})\otimes \RR)$,
such that
$$
\hat{\nu}-\hat{\nu}^{\prime}=a(\eta)\;.
$$
Since $\mathcal{I}(\hat{\nu})=\mathcal{I}(\hat{\nu}^{\prime})=\nu$ it follows that $d\int_{V/M}\eta=0$. By the Riemann-Roch theorem for $\op{KO}$ \cite{AH-diff}, we have 
$$
\int_{V/M}{\rm Ph}(\op{KO}_c^{-1}(V))= \hat{A}(M)\cup {\rm Ph}(\op{KO}^{-1}(M))\;.
$$
Now define the difference class 
\(
\label{diff-class-form}
\delta(\hat{\nu},\hat{\nu}^{\prime}):=\left[\int_{V/M}\eta\right]\in \frac{H^{-1}(M;\pi_*(\op{KO})\otimes \RR)}{\hat{A}(V)\cup {\rm Ph}({\rm KO}^{-1}(M))}\;.
\)
We claim that $\delta=0$ if and only if $\hat{\nu}$ and $\hat{\nu}^{\prime}$ are homotopic, which implies the proposition. If $\delta=0$, then there is $y$ such that $\int_{V/M}\eta= \hat{A}(V)\wedge \widetilde{\rm Ph}(y)$, where 
$$
\widetilde{\rm Ph}:\op{KO}\longrightarrow \mathscr{H}\big(\Omega^*(-;\pi_*(\op{KO}))\big)
$$
is the canonical map in the Hopkins-Singer differential $\op{KO}$-theory. Thus, again by Riemann-Roch,
we can choose a closed form with vertical compact support $\beta$ such that 
$\beta= \widetilde{\rm Ph}(\nu\cup \pi^*y)$ and 
$$
\int_{V/M}\beta=\int_{V/M}\eta\;.
$$
Next consider the homotopy
$
h:=\hat{\nu}+a(t(\eta-\beta))\;.
$
Then $h_0:=h|_{t=0}=\hat{\nu}$ and 
$h_1:=h|_{t=1}=\hat{\nu}+a(\eta-\beta)=\hat{\nu}+a(\eta)=\hat{\nu}^{\prime}$. 
Using the Leibniz rule and the relation $\mathcal{R}a=d$ from the diamond diagram \eqref{kodfdiam}, 
we also have 
$
\mathcal{R}(h)=\mathcal{R}(\hat{\nu})+dt\wedge (\eta-\beta)+td\eta
$,
whence
$$
\int_{V/M}\mathcal{R}(h)={\rm pr}^*\int_{V/M}\mathcal{R}(\hat{\nu})
$$
and $\hat{\nu}$ is homotopic to $\hat{\nu}^{\prime}$ by definition. Conversely, if $\hat{\nu}$ and $\hat{\nu}^{\prime}$ are connected by a homotopy $h$ then, by the homotopy formula, we have
$$
\hat{\nu}-\hat{\nu}^{\prime}=a\Big(\int_{I\times V/V}\mathcal{R}(h)\Big)\;.
$$
Set $\eta:=\int_{I\times V/V}\mathcal{R}(h)$. By the pullback diagram \eqref{pull-dfth} and by virtue of the fact that $\hat{\nu}$ and $\hat{\nu}^{\prime}$ are homotopic, we have
$$\int_{V/M}\eta:=\int_{V/M}\int_{V\times I/V}\mathcal{R}(h)=\int_{M\times I/M}\int_{V\times I/M\times I}\mathcal{R}(h)=\int_{I\times M/M}{\rm pr}^*\Big(\int_{V/M}\mathcal{R}(\hat{\nu})
\Big)=0\;.$$
Hence $\delta=0$. 
\endofproof

Since the set of differential refinements of the Thom classes (which are compatible with the $\hat{A}$-genus form) are parametrized by a subgroup of $H^{-1}(M;\pi_*(\op{KO})\otimes \RR)$, demanding that the refinement is natural with respect to pullback of bundles admitting Spin structure allows us to use our standard trick for proving uniqueness. More precisely, arguing as in Proposition \ref{dfcycmp} gives the following.

\begin{proposition}[Uniqueness of differential refinement of the Spin KO Thom class via naturality]
\label{Prop-uni-spin}
Up to homotopy of differential Thom classes, the class $\widehat{{\bf S}(V)}$ is the unique class refining $\mathbf{S}(V)$ which is natural with respect to pullback of vector bundles of dimension $n\equiv 0 \mod 8$, and which admit Spin structure. 
\end{proposition}

Summarizing our observations, we have established the following theorem.
\begin{theorem}
[Characterization of differential KO Thom class]\label{thomrfko}
Let $V\to M$ be a real vector bundle of rank $n\equiv 0 \mod 8$, equipped with connection $\nabla$ and let $\mathbf{S}(V)_{\widetilde{\nabla}}:=[(\pi^*\S^+(V),\widetilde{\nabla}^+),(\pi^*\S^+(V),\widetilde{\nabla}^-);\mu]$ 
be the induced class from \eqref{SVnablatilde} in ${\rm Gr}(\pi_0{\rm Vect}_{\nabla}(V))$. Then the $\widehat{\op{KO}}_c$-class 
$$
\widehat{\mathbf{S}(V)}=\widehat{\op{cyc}}\big([(\pi^*\S^+(V),\widetilde{\nabla}^+),(\pi^*\S^+(V),\widetilde{\nabla}^-);\mu]\big)+a(\tilde{\eta})\in \widehat{\op{KO}}_c(V)
$$
is a differential refinement of the class $\mathbf{S}(V)=[\pi^*\S^+(V),\pi^*\S^-(V);\mu]$ satisfying the following:
\begin{enumerate}[{\bf (i)}]
\item \emph{(Compatibility with the topological Thom class)} $\mathcal{I}(\mathbf{S}(V)_{\widetilde{\nabla}})=\mathbf{S}(V)$.
\item \emph{(Compatibility with the classical $\hat{A}$-genus)} 
$$
\int_{V/M}\mathcal{R}(\mathbf{S}(V)_{\widetilde{\nabla}})=(\hat{A}(\nabla))^{-1} +d\eta\;,
$$
where $\eta:=\int_{V/M}\tilde{\eta}$
\item  \emph{(Thom isomorphism)}  The class $\mathbf{S}(V)$ is a Thom class, giving rise to an isomorphism 
$$
\mathbf{S}(V)\cdot \pi^*(-):\op{KO}(M)\longrightarrow \op{KO}_c(V)
$$
and $\widehat{{\bf S}(V)}$ is the unique differential refinement which is natural with respect to pullback of vector bundles admitting Spin structure, equipped with connection.
\end{enumerate}
\end{theorem}

We now see how this class can be put to use. 

\begin{proposition}
[The differential KO Thom injection formula]\label{thinjdfko}
The differential Thom class $\widehat{\mathbf{S}(V)}$ gives an injective 
(but not surjective) morphism
$$
\hat{\Phi}:=\widehat{\mathbf{S}(V)}\cup \pi^*(-):\widehat{\op{KO}}^*(M)
\longrightarrow \widehat{\op{KO}}_c^{*+n}(V)
$$
of $\widehat{\op{KO}}$-modules.
\end{proposition}
\theproof
As in the complex case (see \cite{Bun}), fix a differential cohomology class 
$\hat{x}\in \widehat{\op{KO}}$ with $\widehat{\mathbf{S}(V)}\cup \pi^*(\hat{x})=0$. Then the following holds
$$
0=\int_{V/M}\mathcal{R}(\widehat{\mathbf{S}(V)}\cup \pi^*(x))= \int_{V/M} \mathcal{R}(\widehat{\mathbf{S}(V)})\wedge \pi^*(\mathcal{R}(\hat{x}))=0\;.
$$
Furthermore, by Proposition \ref{thomrfko}, we have
\bea 
0=\int_{V/M}\mathcal{R}(\widehat{\mathbf{S}(V)})\wedge \pi^*\mathcal{R}(\hat{x}) &=& \Big(\int_{V/M}\mathcal{R}(\widehat{\mathbf{S}(V)})\Big)\wedge \mathcal{R}(\hat{x})
\\
&=&  \hat{A}(\mathcal{F}_{\widetilde{\nabla}})\wedge \mathcal{R}(\hat{x})+d\eta\wedge \mathcal{R}(\hat{x})
\\
&=& 1\wedge \mathcal{R}(\hat{x})+ \text{higher order terms} \;,
\eea
so that $\mathcal{R}(\hat{x})=0$. Hence $\hat{x}\in \op{KO}^*(M;U(1))$, the  flat KO-theory. Now the 
differential Thom isomorphism reduces to the Thom isomorphism in $\op{KO}^*(M;U(1))$ when restricted to 
flat classes. To see this, note that the differential Thom isomorphism is compatible with the Thom isomorphism 
in $\op{KO}^{*-1}(-;\RR)$. By the compatibility $\beta=\mathcal{I}j$, uniqueness of the Thom isomorphism 
in $\op{KO}(-;U(1))$, and the Five Lemma applied to diagram
$$
\xymatrix{
\op{KO}_c^{*-1}(V;\ZZ)\ar[r]\ar[d]^-{\cong} & \op{KO}_c^{*-1}(V;\RR)\ar[r]^-{\beta}\ar[d]^-{
\cong} & \op{KO}_c^{*-1}(V;U(1)) \ar[r]\ar[d]^-{\hat{\Phi}} & \op{KO}_c^{*}(V;\ZZ)\ar[d]^-{\cong}
\\
\op{KO}^{*-1}(M;\ZZ)\ar[r] & \op{KO}^{*-1}(M;\RR)\ar[r]^-{\beta} & \op{KO}^{*-1}(M;U(1)) \ar[r] & \op{KO}^{*}(M;\ZZ)
}
$$ 
it follows that the restriction of $\widehat{\Phi}$ to $\op{KO}(-;U(1))$ 
is indeed the Thom isomorphism. Hence $\hat{x}=0$ and $\hat{\Phi}$ is injective. 
\endofproof

In general, differential Thom maps are not isomorphisms, but are only injective maps.
We illustrate this with the following example (see also \cite{GS-cob}). 

\begin{example}[Thom map for ordinary differential cohomology]
Consider  ordinary differential cohomology $\widehat{H}(M;\mathbb{Z})$. Let $V\to M$ be an oriented, rank $n$, real vector bundle over $M$. Then we can find a compactly supported differential form $u\in \Omega_{\rm cl}^n(V)$ (fiberwise volume form) such that $\int_{V/M}u=1$. A choice of such a form gives a geometric representative for the Thom class. The wedge product with this form gives rise to a map 
\(\label{thfrmniso}
u\wedge \pi^*(-):\Omega_{\rm cl}^*(M)\longrightarrow \Omega_{\rm cl}^{n+*}(V)\;,
\)
inducing the Thom isomorphism in de Rham cohomology 
$[u]\wedge \pi^*(-):H_{\rm dR}^*(M)\to H^{*+n}_{\rm dR}(V)\;.$ If $u$ is normalized so 
that it has integral periods along cycles then $u$ gives rise to a differential cohomology class 
$\hat{u}\in \widehat{H}^n(V;\mathbb{Z})$ and the underlying class $u\in  H^n(V;\mathbb{Z})$ 
is a Thom class in integral cohomology. Arguing as in the proof in Proposition \ref{thinjdfko}, 
we see that this gives rise to an \emph{injection}
$$
\hat{u}\cup\pi^*(-):\widehat{H}^*(M;\mathbb{Z})\longrightarrow \widehat{H}^{*+n}(V;\mathbb{Z})
$$
but not a surjection. Indeed, if this map were \emph{surjective} then, in particular, we would have a surjection 
at the level of forms given by \eqref{thfrmniso}. But this map cannot be surjective for any choice of $u$. Indeed, 
let $d\alpha$ be any compactly supported exact $n$-form. Then for any form 
$\omega\in \Omega_{\rm cl}^*(M)$, Stokes theorem implies
\begin{align*}
\int_{V/M} \big(d\alpha\wedge \pi^*(\omega)+u\wedge \pi^*(\omega)\big)&=
0+\int_{V/M}u\wedge \pi^*(\omega)
\\
&=\Big(\int_{V/M}u\Big)\wedge \omega
\\
& = \omega\;.
\end{align*}
Hence integration over the fiber defines a left inverse to $u\wedge \pi^*(-)$, but has nontrivial kernel. 
\end{example}

With $\hat{\Phi}:\widehat{\op{KO}}^*(M)\to \widehat{\op{KO}}^{*+n}_c(V)$ the Thom injection defined 
by Proposition \ref{thinjdfko}, 
we can  consider the Thom injection $\hat{\Phi}_{H}$ in ordinary differential cohomology ${\widehat{H}}$. This map does not admit an inverse. However, there is a map $\Psi_{V/M}:\widehat{H}^k_c(V;\ZZ)$ which satisfies some of the familiar properties of underlying cohomological inverse. In fact, we have the following proposition.

\begin{proposition}[Pushforward in differential cohomology with compact support]\label{invdfth-iso}
Let $\pi:V\to M$ be a real vector bundle of rank $n$ and let $\hat{\nu}$ be a differential refinement of the Thom class in ordinary cohomology. Then there is a map $\Psi_{V/M}:\widehat{H}_c^k(V;\ZZ)\to \widehat{H}^{k-n}(M;\ZZ)$ which satisfies the following properties:
\begin{enumerate}[{\bf (i)}]
\item \emph{(Compatibility with curvature)}$\mathcal{R}(\Psi_{V/M}(\hat{x}))=\int_{V/M}\mathcal{R}(\hat{x})$.
\item \emph{(Compatibility with topological realization)} 
$\mathcal{I}(\Psi_{V/M}(\hat{x}))=\Phi_{H}(\mathcal{I}(\hat{x}))$.
\item \emph{(Projection formula)} We have
$$
\Psi_{V/M}(\pi^*(\hat{x})\cup_{\rm DB}\hat{y})=\hat{x}\cup_{\rm DB}\Psi_{V/M}(\hat{y})+a(h)\;,
$$
where $h\in H^{\ast-1}(M;\RR)$.
\item \emph{(Form representative)} $\hat{\Phi}_{H}\Psi_{V/M}$ is in the image of $a:\Omega^{k-1}(V)/{\rm im}(d)\to \widehat{H}^k(V;\ZZ)$.
\end{enumerate}
\end{proposition}
\theproof
In \cite{GT1},  a fiber integration in Deligne cohomology is constructed. Using Deligne cohomology as a model 
for $\widehat{H}^*(M;\ZZ)$, we can take $\Psi_{V/M}$ to be this map. Properties {\bf (i)}  and {\bf (ii)}   follow as part 
of the construction in \cite{GT1}. For property {\bf (iii)} observe that both the pushforward in integral cohomology 
and the fiber integratiojn of forms satisfy the projection formula. It follows at once that the difference 
$\Psi_{V/M}(\pi^*(\hat{x})\cup_{\rm DB}\hat{y})-\hat{x}\cup_{\rm DB}\Psi_{V/M}(\hat{y})$
is in the kernel of both $\mathcal{R}$ and $\mathcal{I}$, hence in the image of $a$. Property {\bf (iv)}  
follows from the fact that $\hat{\Phi}_{H}\Psi_{V/M}-{\rm id}$ is in the kernel of $\mathcal{I}$. 
\endofproof

We will now study the relation between these maps and 
the corresponding ones in differential KO-theory, via the refined $A$-genus.
The following refines relation \eqref{arfcmpc}.

\begin{corollary}
[Differential A-genus via Thom injection]\label{cor-invTh}
Let $V\to M$ be a real vector bundle of rank $n\equiv 0 \mod 8$. Assume $V$ admits Spin structure and fix a 
connection $\nabla$ on $V$. Let $\Psi_{V/M}$ be the integration map defined in Proposition \ref{invdfth-iso}. 
Then there is a unique secondary characteristic form $\eta(V,\nabla)\in \Omega^{-1}(M;\pi_*(\op{KO}))$ such that 
$$
a(\eta(V,\nabla))+ \big(\hat{\mathbb{A}}(V,\nabla)\big)^{-1}=
\Psi_{V/M}\widehat{{\rm Ph}}\big(\hat{\Phi}(1)\big)\;,
$$
where $\hat{\mathbb{A}}$ is the differential $A$-genus (see  Definition \ref{Def-diffA}).
\end{corollary}
\theproof
From Proposition \ref{unithclko} and Proposition \ref{phchrfuni}, we read off 
\bea
\int_{V/M}\mathcal{R}(\widehat{{\rm Ph}}(\hat{\Phi}(1)))
&=& \int_{V/M}\mathcal{R}(\widehat{{\rm Ph}}(\widehat{\mathbf{S}(V)}))
\\
&=& \int_{V/M}{\rm Ph}\big(\mathcal{F}_{\widetilde{\nabla}^+}\big)
-{\rm Ph}\big(\mathcal{F}_{\widetilde{\nabla}^-}\big)
\\
&=& \hat{A}(\mathcal{F}_{\nabla})^{-1}+d\eta\;.
\eea
Since $\mathcal{I}(\widehat{{\rm Ph}}(\Phi(1)))={\rm Ph}(\mathbf{S}(V))$, 
Proposition \ref{invdfth-iso} implies that (after modifying $\eta$ by some form in $H^{-1}(M;\RR)$), the differential cohomology class
$$
\Psi_{V/M}\widehat{{\rm Ph}}\big(\widehat{\mathbf{S}(V)}\big)-a(\eta)
$$
refines both $\hat{A}(\mathcal{F}_{\nabla})^{-1}$ and $\mathbf{S}(V)$. By our now standard argument for proving
uniqueness of characteristic forms, it follows from naturality of the class $\widehat{\mathbf{S}(V)}$, that there is 
a unique choice of $\eta(V,\nabla)$ which is natural with respect to pullback such that 
$$
\Psi_{V/M}\widehat{{\rm Ph}}\big(\widehat{\mathbf{S}(V)}\big)-a(\eta(V,\nabla))=\hat{\mathbb{A}}(V,\nabla)\;.
$$

\vspace{-7mm}
\endofproof


\subsection{Pushforward and Riemann-Roch in $\widehat{\rm KO}$-theory}
\label{Sec-push}

We now construct a pushforward in differential $\op{KO}$ theory, and use Corollary \ref{cor-invTh} 
to prove the differential Riemann-Roch theorem for $\widehat{\op{KO}}$. 
For smooth maps in the topological case, a general axiomatic description via the Pontrjagin-Thom 
construction is given in \cite{CK}. Recall also the recollection of the Riemann-Roch in classical 
KO-theory in Section \ref{Sec-orgen}.

\medskip
The refined desuspension $\hat{\sigma}$ and the collapse map $c$, 
together with the use of Thom spaces, will allow us to offer  the 
following. 

\begin{definition}[Pushforward in differential KO-theory]\label{def-pushfor}
Let $V\to M$ be a real vector bundle of rank $n\equiv 0$ mod 8, admitting a Spin structure and 
let $\hat{\Phi}$  the Thom injection in differential 
$\op{KO}$-theory $\hat{\Phi}(x):=\widehat{\mathbf{S}(V)}\cdot \pi^*(-)$ (see Proposition \ref{thinjdfko}).
Let $f:W\to M$ be a proper map. An embedding of $W\into \RR^n$ and $f$ gives rise to 
an embedding $i:W\into M\times \RR^n$. Let $\mathcal{N}$ denote the normal 
bundle of the embedding. We define the \emph{pushforward} of $f$ as the composition 
\vspace{-2mm}
$$
\xymatrix{
f_!:\widehat{\op{KO}}^*(W)\ar[r]^-{\widehat{\Phi}} & \widehat{\op{KO}}^{*+n-k}({\rm Th}(\mathcal{N}))
\ar[r]^-{c(i)^*} & \widehat{\op{KO}}^{*+n-k}(M\wedge \widetilde{S}^n)\ar[r]^-{\hat{\sigma}_n^{-1}} &
\widehat{\op{KO}}^{*-k}(M)},
$$
which depends on the embedding $i$. 
\end{definition}

\medskip
Recall the Nash Embedding Theorem \cite{Na} that every Riemannian manifold $M$ equipped 
with Riemannian metric $g$ can be isometrically embedded into some Euclidean space. 
We will use this in order to relate  $\mathbf{K}{\rm O}_{\nabla}(M)$
classes associated to the tangent bundle of such a manifold with those associated to the normal bundle, 
both taken with appropriate corresponding connections. 

\begin{lemma}[Splitting of natural  $\mathbf{K}{\rm O}_{\nabla}(M)$ classes via embedding into Euclidean space]\label{splitemb}
Let $(M,g)$ be a Riemannian manifold and let $i:(M,g)\into \RR^n$ be an isometric embedding, 
with $\RR^n$ is equipped with the standard metric. The normal bundle
$\mathcal{N}$ of the embedding inherits a metric given by
 $g^{\prime}(\phi(x),\psi(x)):=\langle \phi(x),\psi(x)\rangle$, where 
 $\phi,\psi:M\to \mathcal{N}$ are sections. We denote the corresponding
  Levi-Civita connection by $\nabla^{\perp}_{g}$. 
Then, under the cycle map 
$$\widehat{{\rm cyc}}:\mathbf{K}{\rm O}_{\nabla}(M)\to \widehat{\rm K}{\rm O}(M),$$
in Proposition \ref{csconvb}, we have
$$
\widehat{{\rm cyc}}([TM,\nabla_{g}])+\widehat{{\rm cyc}}([\mathcal{N},\nabla^{\perp}_{g}])=a(\xi)
$$
for some odd differential form $\xi\in \Omega^{-1}(M;\pi_*({\rm KO}))$.
\end{lemma}
\theproof
The Whitney sum decomposition ${\bf n}=TM\oplus {\cal N}$ gives the identification of underlying ${\rm KO}$-classes $[TM]+[{\cal N}]=0$. By the commutativity $\mathscr{F}={\cal I}\circ \widehat{{\rm cyc}}$ in Proposition \ref{csconvb}, it follows that 
$${\cal I}\Big(\widehat{{\rm cyc}}([TM,\nabla_{g}])+\widehat{{\rm cyc}}([\mathcal{N},\nabla^{\perp}_{g}])\Big)=[TM]-[{\cal N}]=0\;.
$$
Hence, by the exact sequence 
\begin{equation}\Omega^{-1}(M;\pi_*({\rm KO}))\overset{a}{\longrightarrow} 
\widehat{{\rm KO}}(M)\overset{\cal I}{\longrightarrow} {\rm KO}(M) \longrightarrow 0 \;,
 \label{kodifexact}
\end{equation}
the claim follows. 
\endofproof

\begin{remark}
In fact, the Lemma \ref{splitemb} can be improved. We have the relation
$$\widehat{{\rm cyc}}([TM,\nabla_{g}])=-\widehat{{\rm cyc}}([\mathcal{N},\nabla^{\perp}_{g}])$$
in $\widehat{\op{KO}}$, i.e., we can take $\xi=0$. Indeed, one can model $\widehat{\op{KO}}$ via structured vector bundles (see \cite{SSu1}). Since the curvature of $d$ vanishes, \cite[Lemma 1.16]{SSu1} implies that the Chern simons form $CS(d,\nabla_g\oplus \nabla_g^{\perp})$ vanishes. Hence $d$ and $\nabla_g\oplus \nabla_g^{\perp}$ give rise to the same class in $\widehat{\op{KO}}$. We will not provide more details, since Lemma   \ref{splitemb} will be sufficient for our needs. 
\end{remark}

\begin{remark}
The relation $[TM,\nabla_g]\oplus [{\cal N},\nabla_g^{\perp}]=0$ does not hold in ${\bf K}{\rm O}_{\nabla}(M)$. Indeed, the direct sum of projected connections is not the Levi-Civita connection on $\RR^n$. One can see this, for example, by observing that parallel transport in the direct sum connection preserves tangential components.  On the other hand, this clearly fails for the connection $d$ on $\RR^n$. 
\end{remark}

More generally, suppose we have a smooth map $f: (W,g)\to (M,h)$  between Riemannian manifolds. 
Equip $W$ with the metric $g+f^*h$. Let $i:W\into \RR^n$ be an isometric embedding of $(W,g)$ into 
$\RR^n$. Take the metric on the product $\RR^n\times M$ to be usual one defined by 
$k((x,y),(x^{\prime},y^{\prime})):=\langle x, y\rangle+h(x^{\prime},y^{\prime})$ for $x,y$ in $(T\RR^n)_p$ and $x^{\prime},y^{\prime}$ sections of $TM_q$. Then the map $x\mapsto (i(x),f(x))$ defines an isometric embedding. Then we have the following.

\begin{proposition}[Splitting of natural  $\widehat{{\rm KO}}(M)$ classes for maps between Riemannian manifolds]\label{Prop-split-MW}
For smooth map $f: (W,g)\to (M,h)$  between Riemannian manifolds, we have
an identification in $\widehat{{\rm K O}}$
$$
\widehat{{\rm cyc}}([TW,\nabla_{g+f^*h}])+\widehat{{\rm cyc}}([\mathcal{N},\nabla^{\perp}_{g}])=f^*\widehat{{\rm cyc}}([TM,\nabla_{h}])+a(\xi)\;,
$$ 
for some $\xi\in \Omega^{-1}(M;\pi_*({\rm KO}))$.
\end{proposition} 
\theproof
Let ${\cal N}$ denote the normal bundle of the embedding $j=(i,f)$. Then $TW\oplus {\cal N}\cong {\bf n}\oplus f^*TM$. Hence, in ${\rm KO}(M)$, we have $[TW]+[{\cal N}]=f^*[TM]$. Then by the commutativity $\mathscr{F}={\cal I}\circ \widehat{{\rm cyc}}$, it follows that  
$${\cal I}\Big(\widehat{{\rm cyc}}([TW,\nabla_{g+f^*h}])+\widehat{{\rm cyc}}([\mathcal{N},\nabla^{\perp}_{g}])-f^*\widehat{{\rm cyc}}([TM,\nabla_{h}])\Big)=0$$
and the claim follows again by the exact sequence \eqref{kodifexact}.
\endofproof

We now utilize this, Corollary \ref{cor-invTh}, and integration over the fiber to get the main result
in this section.  General fiber integration in generalized  (Eilenberg-Steenrod) cohomology is 
discussed in \cite[Section 3.4]{HS}\cite[around Prop. 2.1]{Fr-CS2}. Fiber integration of cocycles 
in Deligne cohomology is discussed explicitly in detail in \cite{GT1}\cite{GT2}\cite{DL}\cite{Lip}. 
The refinement to smooth higher moduli stacks is described in \cite{FSS1}. 
As in the classical case for the Riemann-Roch theorem (see \cite{Hir} and Section \ref{Sec-orgen}), 
we consider the following diagram involving pushforwards
$$
\xymatrix{
{\scriptstyle (E, \nabla_E) \; \in} \; \op{KO}_\nabla(W) \ar[dr]^{\widehat{\rm cyc}} &&&\\
{\scriptstyle \widehat{\rm cyc}(E, \nabla_E) \; \in } \hspace{-3.2cm} & 
\widehat{\op{KO}}(W)\ar[rr]^-{f_!}\ar[d]_-{\widehat{{\rm Ph}}} && 
\widehat{\op{KO}}(M)\ar[d]^-{\widehat{{\rm Ph}}} &
\\ 
{\scriptstyle\hat{\mathbb{A}}(TW, \nabla_{g+f^*h}) \; \in}  \hspace{-2.5cm} 
& \widehat{H}(W;\mathbb{Q})\ar[rr]^-{f_*} && \widehat{H}(M;\mathbb{Q}) 
& \hspace{-1cm} {\scriptstyle\ni \; \hat{\mathbb{A}}(TM, \nabla_h)}.
}
$$
The differential Riemann-Roch theorem gives the correction factor for the non-commutativity of this diagram.

\begin{theorem} [Riemann-Roch for differential  KO-theory]
Let $f:(W,g)\to (M,h)$ be a a proper smooth map between Riemannian manifolds. Let $\nabla_{h}$ denote the 
Levi-Civita connection on $M$ and equip $TW\to W$ with the Levi-Civita connection associated to the metric 
$g+f^*h$ . Let $E\to W$ be a real vector bundle with connection $\nabla_{E}$. Then we have the formula 
\(
\label{rrthmko}
\widehat{{\rm Ph}}(f_!(\widehat{{\rm cyc}}(E,\nabla_{E}))) \cup_{\rm DB}  \hat{\mathbb{A}}(TM,\nabla_{h})
=f_*(\widehat{{\rm Ph}}(\widehat{\rm cyc}(E,\nabla_{E}))\cup_{\rm DB} \widehat{\mathbb{A}}(TW,\nabla_{g+f^*h}))+a(\eta)\;,
\)
with $\eta$ some differential form in $\Omega^{-1}(M;\pi_*(\op{KO}))$. 
\end{theorem}
\theproof
Fix an embedding $W\into \RR^n$ and consider the induced embedding $W\into \RR^n\times M$ as discussed above. 
From Definition \ref{def-pushfor} of the pushforward and Corollary \ref{cor-invTh} for the differential genus
$\hat{\mathbb{A}}$, we immediately have the formula
\begin{align}\label{rreqcool}
f_*(\widehat{{\rm Ph}}((E,\nabla_{E})\cup \hat{\mathbb{A}}(\mathcal{N},\nabla^{\perp}_{g})) &=
f_*(\widehat{{\rm Ph}}(E,\nabla_{E})\cup_{\rm DB} (\Psi_{\mathcal{N}/M}\widehat{{\rm Ph}}(\widehat{\Phi}(1))+a(\eta)))
\\
&=f_*(\widehat{{\rm Ph}}(E,\nabla_{E})\cup_{\rm DB} \Psi_{\mathcal{N}/M}\widehat{{\rm Ph}}(\widehat{\Phi}(1)))+ f_*(\widehat{{\rm Ph}}(E,\nabla_{E})\cup_{\rm DB}a(\eta))\;. \nonumber
\end{align}
The image of $a$ is an ideal in the ring $\widehat{H}^*(M;\ZZ)$ and there is $\tilde{\eta}$ such that $a(\tilde{\eta})=f_*(\widehat{{\rm Ph}}(E,\nabla_{E})\cup_{\rm DB}a(\eta))$. Using this and properties 
{\bf (iii)}(Projection formula) and {\bf (iv)}(Form representative) of 
Proposition \ref{invdfth-iso} applied to the normal bundle $\mathcal{N}$ of the embedding $W\into \RR^n\times M$, 
we have
\begin{align*}
\hspace{-0mm}
f_*(\widehat{{\rm Ph}}(E,\nabla_{E})\cup_{\rm DB} \Psi_{\mathcal{N}/M}\widehat{{\rm Ph}}(\widehat{\Phi}(1))) 
&=
f_* \Psi_{\mathcal{N}/M}(\pi^*\widehat{{\rm Ph}}(E,\nabla_{E})\cup_{\rm DB}\widehat{{\rm Ph}}(\widehat{\Phi}(1)))+a(h)
\\
&= \hat{\rho}_n^{-1}c(i)^*\hat{\Phi}_{H}\Psi_{\mathcal{N}/M}(\pi^*\widehat{{\rm Ph}}(E,\nabla_{E})\cup_{\rm DB}\widehat{{\rm Ph}}(\widehat{\Phi}(1)))+a(h^{\prime})+a(h)
\\
&= \widehat{{\rm Ph}}(\hat{\sigma}_n^{-1} c(i)^*(\pi^*\widehat{\rm cyc}(E,\nabla_{E})\cdot \widehat{\Phi}(1)))+a(h+h^{\prime})
\\
&= \widehat{{\rm Ph}}(f_!(\widehat{\rm cyc}(E,\nabla_{E})))+a(h+h^{\prime})\;,
\end{align*}
for some $h,h^{\prime}\in H^{-1}(M;\pi_*(\op{KO}))$. Here we have used the Definition of the 
Thom injection  $\hat{\Phi}(\hat{x})=\pi^*(\hat{x}) \cdot \hat{\Phi}(1)$   from Proposition \ref{thinjdfko}.
Furthermore, in equalities one through four we have used, respectively, the projection formula, 
the definition of $f_*$, the naturality of $\widehat{\rm Ph}$, and the definition of $f_!$.
Combining this with equation \eqref{rreqcool}, using the relation in Proposition \ref{Prop-split-MW} (recalling also that the image of $a$ is an ideal),
and the fact that $\hat{\mathbb{A}}$ is natural with respect to pullback and multiplicative 
with respect to Whitney sums, we have 
$$
f_*(\widehat{{\rm Ph}}(E,\nabla_{E})\cup_{\rm DB} f^*\hat{\mathbb{A}}(TM,\nabla^{\perp}_{g})^{-1}\cup_{\rm DB}\hat{\mathbb{A}}(TW,\nabla^{\perp}_{g+f^*h}))+a(\tilde{\xi})=\widehat{{\rm Ph}}(f_!(\widehat{\rm cyc}(E,\nabla_{E})))+a(h+h^{\prime})+a(\tilde{\eta})\;,
$$
for some form $\tilde \xi\in \Omega^{-1}(M;\pi_*({\rm KO}))$. Finally, using the projection formula for $f_*$ and absorbing forms together into a single $\tilde{\eta}$, we have the formula
$$
f_*(\widehat{{\rm Ph}}((E,\nabla_{E})\cup_{\rm DB}\hat{\mathbb{A}}(TW,\nabla^{\perp}_{g+f^*h}))\cup_{\rm DB}\hat{\mathbb{A}}(TM,\nabla^{\perp}_{g}))^{-1}=\widehat{{\rm Ph}}(f_!(\widehat{\rm cyc}(E,\nabla_{E})))+a(\tilde{\eta})\;.
$$
Multiplying both sides by $\hat{\mathbb{A}}(TM,\nabla^{\perp}_{g})$, using that the image of $a$ is an ideal,
and absorbing forms into $\tilde{\eta}$ gives the result.
\endofproof

The differential form $\eta$ appearing in the Riemann-Roch theorem seems to be related to the 
$\eta$-form of Bismut-Cheeger \cite{BCh}. More precisely, we observe that $\tilde{\eta}=\int_{\mathcal{N}/M}\eta$ 
for some form $\eta$ with fiberwise compact support. Then pushing forward to the point $M=\ast$ and 
applying the curvature map $\mathcal{R}$ to formula \eqref{rrthmko} gives  
$$
{\rm Index(D)}=\int_{M}{\rm Ph}(\mathcal{F}_{\nabla_{E}})\wedge 
\hat{A}(\mathcal{F}_{\nabla_{g+f^*h}})+\int_{M}d\eta \;.
$$
If $M$ has nonempty boundary then, by the Atiyah-Patodi-Singer index theorem, it follows
that $\int_{\partial M}\eta$ is the $\eta$-invariant of the twisted Dirac operator $D$. 
This subject deserves further development elsewhere.

\section{ The Atiyah-Hirzebruch spectral sequence} 
\label{Ch-AHSS}

\subsection{The AHSS for KO-theory} 
\label{Sec-AHSS-KObare}

We start by recalling the Atiyah-Hirzebruch spectral sequence (AHSS) for topological KO-theory, 
summarizing the main results in the literature, which will serve as a basis for the differential version in the 
next section. In addition to providing one place for material scattered in the literature, 
we also extend existing results in the case of KO-theory. 

\medskip
One major difference between the real and the complex theory consists of the following
observations which make the former much more delicate to deal with. 
In an AHSS for a CW-complex or a spectrum which is  
bounded below, the images of all differentials $d_r$ with $r \geq 2$ are torsion subgroups.
One consequence is that if $E$ has torsion-free coefficient group $\pi_*(E)$ and 
$R$ is torsion-free ordinary homology, then all differentials $d_r$ with $r \geq 2$ vanish;
see  \cite[ Bemerkung 14.18]{Dold}\cite{Arl}. 
That is, for torsion-free $H_*(X)$, the AHSS for cohomology collapses at the $E_2$-page.
In particular, ${\rm K}^*(X)$ is torsion-free. This is in sharp contrast with real theory ${\rm KO}(X)$
where the coefficient group $\pi_*({\rm KO})$ already has torsion and so the above statement will
not hold. On the other hand, unlike the case of ${\op K}$-theory, 
requiring ${\rm KO}(X)$ of a space $X$ to be torsion-free is severely restrictive. 
It is the presence of torsion (inevitable for coefficients and almost inevitable for a vast class of spaces)
which prevents one from making naive statements
about direct analogies between the two theories.

\medskip
If $X$ is any space with $p$-skeleton $X_p$, there is an Atiyah-Hirzebruch (AH) 
filtration on $\op{KO}(X)$ obtained by letting $\op{KO}^*(X)_p=\bigcap_f \op{Ker}(f^*)$, 
where $f$ runs over all maps $f: Y \to X$, $\dim (Y)<p$, with $Y$ a finite complex \cite{AH}. 
The AHSS for $\op{KO}^*(X)$ is the spectral sequence $\{E_r^{p, q}, d_r\}$ with
(see \cite[pp. 336-341]{Sw}) 
$$
E_2^{p, q} \cong  H^p(X; \op{KO}^q)
\qquad \text{and} \qquad 
E_\infty^{p, q} \cong  F^{p, q}/F^{p+1, q-1}=\op{KO}^{p+q}(X)_p/ \op{KO}^{p+q}(X)_{p+1},
$$
where 
$
F^{p, q}=\op{Ker}[i^*_0: \op{KO}^{p+q}(X) \to \op{KO}^{p+q}(X_{p-1})]
$.
If $x \in \op{KO}^0(X)$ has filtration $p$, then 
$x$ determines a coset $[x]\subset H^p(X; \op{KO}^{-p}({\rm pt}))$. 
The differentials $d_r^{p, q}: E_r^{p, q} \to E_r^{p+r, q-r+1}$
satisfy $d_r^{p, q}=\Omega d_r^{p+1, q}$
and $d_r^{p, q}=d_r^{p, q+8}$. Furthermore, 
$d_2^{8t, 0}$,
$d_2^{8t, -1}$ and 
$d_3^{8t, -2}$ are induced by the cohomology 
operations defined by the $k$-invariants 
$$
k^{8t, 8t} \in H^{8t+2}(K(\Z, 8t), \Z/2)\;,
\qquad 
k^{8t, 8t} \in H^{8t+2}(K(\Z/2, 8t), \Z/2)\;,
\qquad 
k^{8t, 8t+1} \in H^{8t+3}(K(\Z/2, 8t), \Z)\;,
$$
respectively \cite[Theorem 3.4]{Ma}. 
In the $E_2$-term  
the nonzero differentials may then be described as $Sq^2 \otimes \cdot \eta$, where
$\cdot \eta$ denotes multiplication by $\eta\in \op{KO}_1$, precomposed with reduction
modulo 2 if necessary. Explicitly, \cite[equations (1.3)]{Fu1} \cite[Theorem 4.2]{Tho}
\bea
d_2^{p, -8t}&=&Sq^2\rho_2 : H^p(X; \Z) \longrightarrow H^{p+2}(X; \Z/2)\;,
\\
d_2^{p, -8t-1}&=&Sq^2 : H^p(X; \Z/2) \longrightarrow H^{p+2}(X; \Z/2)\;,
\\
d_3^{p, -8t}&=&\beta_2 Sq^2 : H^p(X; \Z/2) \longrightarrow H^{p+3}(X; \Z)\;,
\eea
where $\beta_2$ is the Bockstein operator associated with the exact 
coefficient sequence
$
0 \to \Z \xrightarrow{\times 2}  \Z \xrightarrow{\rho_2} \Z/2 \to 0
$, with $\rho_2$ being mod 2 reduction.
%
Adding  integral operations to this and retaining the origin of the coefficients,
the differentials are given as follows (see \cite{ABP}[proof of Lemma 5.6])
\footnote{ Attributed to \cite{Bott} there.} 
\begin{eqnarray}
\label{d2d5-class}
Sq^2 \rho_2: H^p(X; \op{KO}^0({\rm pt})) &\longrightarrow & H^{p+2}(X; \op{KO}^{-1}({\rm pt}))\;,
\\ \nonumber
Sq^2 : H^p(X; \op{KO}^{-1}({\rm pt})) &\longrightarrow & H^{p+2}(X; \op{KO}^{-2}({\rm pt}))\;,
\\ \nonumber
\beta_2 Sq^2 : H^p(X; \op{KO}^{-2}({\rm pt})) &\longrightarrow & H^{p+3}(X; \op{KO}^{-4}({\rm pt}))\;,
\\ \nonumber
\beta_2 Sq^4 \rho_2: H^p(X; \op{KO}^0({\rm pt})) &\longrightarrow & H^{p+5}(X; \op{KO}^{-1}({\rm pt}))\;.
\end{eqnarray}
Due to the relation $\rho_2 \beta_2 Sq^4 \rho_2=Sq^1 Sq^4 \rho_2=
Sq^4 Sq^1 \rho_2 + Sq^2 Sq^3 \rho_2$, the fourth operation reduces via
$\rho_2 \beta_2 Sq^4 \rho_2=Sq^2 Sq^3 \rho_2$. 
%
Indeed, in \cite{Gh} the latter is related to the AHSS differential for KO-theory with 
$\Z/2$-coefficients $\op{KO}(-;\Z/2)$ acting as $d_5=Sq^2 \beta_2 w_2$ on the 
Thom class.  Mod-$2^n$ real K-theory $\op{KO}(-;\ZZ/2^n)$ has also been investigated in 
\cite{Za}. Higher differentials have been considered for certain complexes e.g. in 
\cite{KH}\cite{KH2}\cite{KKO}\cite{KKO1}\cite{KO}, where the 
differentials appear in degrees 2 mod 8.

\medskip
In \cite{Buc1}\cite{Buc2} it was shown that for complex $K$-theory,  the general formula
for the differentials is given as:
\(
\label{Buch-diff}
d
_{2r(p-1)+1}(x)_p=\epsilon_r(\beta_p P^r)(x/p^{r-1})\cdot v \;,
\)
where $d(-)_p$ denotes the $p$-component of the differential, $\epsilon_r\neq 0 \mod p$, 
$v\in K^{-2r(p-1)}(\ast)$ a power of the Bott generator $u$, and $\beta P^r$ is an integral 
Steenrod operation given by the mod $p$ Bockstein $\beta_p$ acting on the Steenrod reduced power 
operation $P^r$ of degree $r$ at the prime $p$. We can use this and complexification 
to identify a large portion of the differentials for the real counterpart for classes $x$ in KO-theory 
of a finite CW-complex $X$ which are divisible by certain powers of $p$. Recall 
(see \eqref{Rat-coeff})  that the rationalization of
KO-theory is 4-periodic generated by two elements $\alpha$ and $\beta$ of degrees 4 and 8,
respectively.

\begin{proposition}[Differentials in the KO-theory AHSS]\label{diff-Ko-ahss} 
The differentials for KO-theory are given according to the prime as

\item [{\bf (i)}] For $p$ odd, we have the formula 
$$
d^{s,4t}_{4r(p-1)+1}(x)_p=\epsilon_{2r}(\beta_p P^{2r})(x/p^{2r-1})\cdot \tilde{\alpha}\;.
$$

\item [{\bf (ii)}] For $p=2$, we also have the formula 
$$
d^{s,8t}_{8r+1}(x)_2=\epsilon_{4r}(\beta_2 Sq^{8r})(x/2^{4r-1})\cdot \tilde{\alpha}\;,
$$
where $\tilde{\alpha}$, a power of $\a$, is the generator in $\pi_{-4t}(\op{KO})$ appearing in 
$E^{s,t}=H^s(M; \pi_{-4t}(\op{KO})$ (see Theorem \ref{thm-ahssDiff} below); explicitly, under complexification 
$$
c (\tilde{\alpha})=\left\{ 
\begin{array}{ll}
2u^{2t} & \text{for $t$ odd}, \\
u^{2t} & \text{for $t$ even}.
\end{array}
\right.
$$

\end{proposition}

\theproof
The complexification map $c:\op{KO}\to \op{K}$
induces a morphism of spectral sequences. The Bott sequence  (see Section \ref{Sec-cplx})  
identifies this map on coefficients with  an isomorphism 
$\op{KO}^{8t}({\rm pt})\overset{c}{\cong} \op{K}^{8t}({\rm pt})$ 
and the map 
$$
\Z\cong \op{KO}^{4t}({\rm pt})\overset{\times 2}{\longrightarrow} \op{K}^{4t}({\rm pt})\cong \Z\;.
$$
with $t\not\equiv 0$ mod 2. The induced map on the corresponding AHSS's gives rise to the commutative diagrams 
$$
\xymatrix{
H^s(M;\Z)_p\ar@{=}[r]\ar[d]_-{d^{s,8t\ \prime}_{8r(p-1)+1}} &  H^s(M;\Z)_p\ar[d]^-{d^{s,8t}_{8r(p-1)+1}}
\\
H^{s+8r(p-1)+1}(M;\Z)_p \ar@{=}[r] &  H^{s+8r(p-1)+1}(M;\Z)_p
}
$$
and 
$$
\xymatrix{
H^s(M;\Z)_p\ar[r]^-{\times 2}\ar[d]_-{d^{s,4+8t\ \prime}_{4r(p-1)+1}} &  H^s(M;\Z)_p\ar[d]^-{d^{s,4+8t}_{4r(p-1)+1}}
\\
H^{s+4r(p-1)+1}(M;\Z)_p\ar[r]^-{\times 2} &  H^{s+4r(p-1)+1}(M;\Z)_p
\;,}
$$
where $H^*(M;\Z)_p$ denotes the $p$-primary part. Away from $p=2$, this gives the desired identification
via \eqref{Buch-diff}. For $p=2$, the identification follows from the first diagram. 
\endofproof

We will use this proposition in the next section to provide conditions on classes in  the image of the 
Pontrjagin character map, using the differential AHSS we develop
there. 

\subsection{AHSS for differential KO-theory}
\label{Sec-AHSS-KOhat}

The AH filtration on the Hopkins-Singer type differential theory $\widehat{\op{KO}}(X)$ is induced from 
the filtration on the underlying topological theory $\op{KO}(X)$ via the natural topological realization map
$\mathcal{I}:\widehat{\op{KO}}(X) \to {\op{KO}}(X)$ (see the statement after Lemma B.24 in \cite{Fr}). 
In \cite{GS3}, we constructed a general AHSS which works for any differential cohomology theory. We briefly 
summarize this construction. Let $M$ be a compact manifold and consider the simplicial object given 
by resolving $M$ in $\sh_{\infty}(\mathscr{M}{\rm an})$, by the nerve of a good open cover 
$\{U_{\alpha}\}$. This simplicial object has a natural filtration through its skeleta and induces a filtration 
on a differential cohomology group $\mathscr{E}^*(M)$ via
$$
\mathscr{E}^*(M)_s:=\ker\big(i^*:\mathscr{E}^*(M)\to \mathscr{E}(X_s)\big)\;,
$$
with $X_s:={\rm sk}_s(C(\{U_{\alpha}\}))$ and $i^*$ induced by the canonical map 
$i: {\rm sk}_s(C(\{U_{\alpha}\})\to M$. Through the Nerve Theorem one finds that this filtration 
is compatible with the canonical map $\mathcal{I}$ and agrees with the one discussed in \cite{Fr}. 

\medskip
Recall (see  expression\eqref{Rat-coeff}) that KO-theory is rationally 4-periodic with  
$\pi_* {\rm KO} \otimes \QQ=\QQ [\alpha, \alpha^{-1}]$.
Using the general methods of \cite{GS3}\cite{GS5} specialized to $\op{KO}$, we can immediately 
establish the following.

\begin{theorem} [AHSS for differential KO-theory]\label{thm-ahssDiff}
Let $C(\{U_{\alpha}\})$ be the {\v C}ech nerve of a good open cover of a compact manifold $M$. Filtering 
the realization of this simplicial object by skeleta gives rise to a half-plane spectral sequence 
for $\widehat{\op{KO}}$ taking the following form 
\footnote{Note that $\widehat{\rm KO}_0={\rm diff}({\rm KO}, \widetilde{\rm Ph}, \RR[\a, \a^{-1}])$ is
the sheaf of spectra associated to differential KO-theory in degree zero.}
\(
E^{s,t}=
\left\{ \begin{array}{cc}
\Omega^0_{\ZZ,{\rm cl}}(M;\pi_*(\op{KO}))& s=t=0
\\
H^{s}(M; \op{KO}^{t}({\rm pt}; U(1))) & t<0 
\\
H^{s}(M;\op{KO}^{t}({\rm pt})) & t>0
\end{array}\right\} \; \Longrightarrow  \;  \big(\widehat{\op{KO}_0}\big)^{s, t}(M)\;,
\)
where the elements in the group $\Omega^0_{\ZZ,{\rm cl}}(M;\pi_*(\op{KO}))$ are formal combinations 
$$\omega=\omega_0+\omega_4+\omega_8+\hdots $$
with $\omega_{4k}$ a differential form of degree $4k$ on $M$ and the zero form $\omega_0$
is a constant $\ZZ$-valued function.
\end{theorem}

The following was also shown in  \cite{GS3} for the case of complex K-theory, but the argument 
works in general for any sheaf of spectra, including KO-theory. Nevertheless, we provide an argument 
here in order to be self-contained. Recall the Pontrjagin character from Section \ref{Sec-PontCh}. 

\begin{proposition}[The $E_\infty$-term in the differential KO-theory AHSS] \label{prop-E00SS}
We have 
\bea
E_{\infty}^{0,0}&=&\widehat{\op{KO}}(M)/\op{KO}^{-1}(M;U(1))
\\
&\cong& {\rm Im}({\rm Ph})\;,
\eea
where ${\rm Ph}$ is the Pontrjagin character form ${\rm Ph}:\pi_0\big({\rm Iso}\big({\rm Vect}_{\nabla}(M)\big)\big)\to \bigoplus_k\Omega^{4k}(M)$.
\end{proposition}
\theproof
First recall that, from the construction of the Hopkins-Singer differential $\op{KO}$, the composition 
$\widehat{\rm cyc}:\mathbf{K}{\rm O}_{\nabla}(M)\to \widehat{\op{KO}}(M)\xrightarrow{\mathcal{R}}\bigoplus_k\Omega^{4k}(M)$ is precisely the Pontrjagin character form map ${\rm Ph}$. Thus, the image of 
$\mathcal{R}$ is precisely the image of ${\rm Ph}$. Identifying $\op{KO}^{-1}(M;U(1))$ with the flat part of 
$\widehat{\op{KO}}(M)$, it then follows immediately from the differential cohomology diamond diagram 
\eqref{kodfdiam} that 
\(
\label{grqdfkoph}\widehat{\op{KO}}(M)/\op{KO}^{-1}(M;U(1))\cong {\rm Im}({\rm Ph})\;.
\)
Hence, we need only show that the left hand side is indeed $E_{\infty}^{0,0}$. By definition of the filtration and 
$\widehat{\op{KO}}$, the kernel
$$
\ker\Big(i^*:\widehat{\op{KO}}(M)\longrightarrow
\widehat{\op{KO}}\big({\rm sk}_0(C(\{U_{\alpha}\})\big)\cong 
\ZZ\oplus \prod_{\alpha}\Omega_{\rm cl}^0\big(U_{\alpha};\pi_*(\op{KO})\big)\Big)
$$
can be identified with those classes in $\widehat{\op{KO}}(M)$ whose local curvature forms vanish on each patch 
of the good open cover $\{U_{\alpha}\}$. But this is precisely the kernel of $\mathcal{R}$, which by exactness 
of the diamond is the image of $\op{KO}^{-1}(M;U(1))$ in $\widehat{\op{KO}}(M)$. The elements in 
$E^{0,0}_{\infty}$ are, by definition, those elements which converge to the graded quotient 
$\widehat{\op{KO}}(M)_0$ which we have shown is precisely the left hand side of 
the isomorphism \eqref{grqdfkoph}.
\endofproof

Next, for $-t\leq s$ and $s\neq 0$ we have an isomorphism
\bea
E_{\infty}^{s,t} &=&\widehat{\op{KO}}_{\rm flat}(M)_s/\widehat{\op{KO}}_{\rm flat}(M)_{s+1}
\\
&\cong& \op{KO}^{-1}(M;U(1))_{s}/\op{KO}^{-1}(M;U(1))_{s+1}\;,
\eea
where the right hand side is the usual AH filtration on for the theory $\op{KO}^{-1}(-;U(1))$. This is more or less implicit in \cite{GS3}, but we include a proof of this fact here.

\begin{proposition}
[AHSS for flat KO-theory vs. that of KO-theory with U(1)-coefficients]
For the flat theory $\widehat{\op{KO}}^*_{\rm flat}$, the refined AHSS of \cite{GS3} is isomorphic 
to usual AHSS for $\op{KO}^{*-1}(-;U(1))$.
\end{proposition}
\theproof
We will make use of the functors appearing in the quadruple adjunction \eqref{quad-adj}.
The theory $\widehat{\op{KO}}_{\rm flat}$ is representable by the constant sheaf of 
spectra $\delta\big(\op{KO}_{U(1)}\big)$. Let $C(\{U_{\alpha}\})$ be the realization 
of the {\v C}ech nerve of a good open cover $\{U_{\alpha}\}$ of $M$. Then, by the cohesive adjunction,
$$
\map\big(C(\{U_{\alpha}\}),\delta\big(\op{KO}_{U(1)}\big)\big)\simeq 
\map\big(\Pi(C(\{U_{\alpha}\})),\op{KO}_{U(1)}\big)\;.
$$
Now, since the cover $\{U_{\alpha}\}$ is good, it follows from the Nerve Theorem that 
$\vert \Pi(C(\{U_{\alpha}\}))\vert \simeq X\simeq M$, where $X$ is a finite CW-complex 
with $s$-skeleton $X_s$ given by attaching a copy of $\Delta^s$ to $X_{s-1}$ via the 
combinatorics of the cover. It follows at once that $\Pi$ preserves the respective filtrations 
and induces an isomorphism on corresponding exact couples.
\endofproof

\begin{remark}[The $E_2$-page]
From the above results, it follows that the $E_2$-page of the spectral sequence takes the form 
depicted in Figure 1.
Note that some of the entries (in particular, those with $U(1)$-coefficients) are shifted in degree 
as compared with the AHSS for $K{\rm O}$. This is a typical effect which arises from the identification 
$\widehat{\op{K}}^*_{\rm flat}(-)\cong \op{KO}^{*-1}(-;U(1))$ (see \cite{GS3}\cite{GS5}  
for more detail).
\end{remark} 

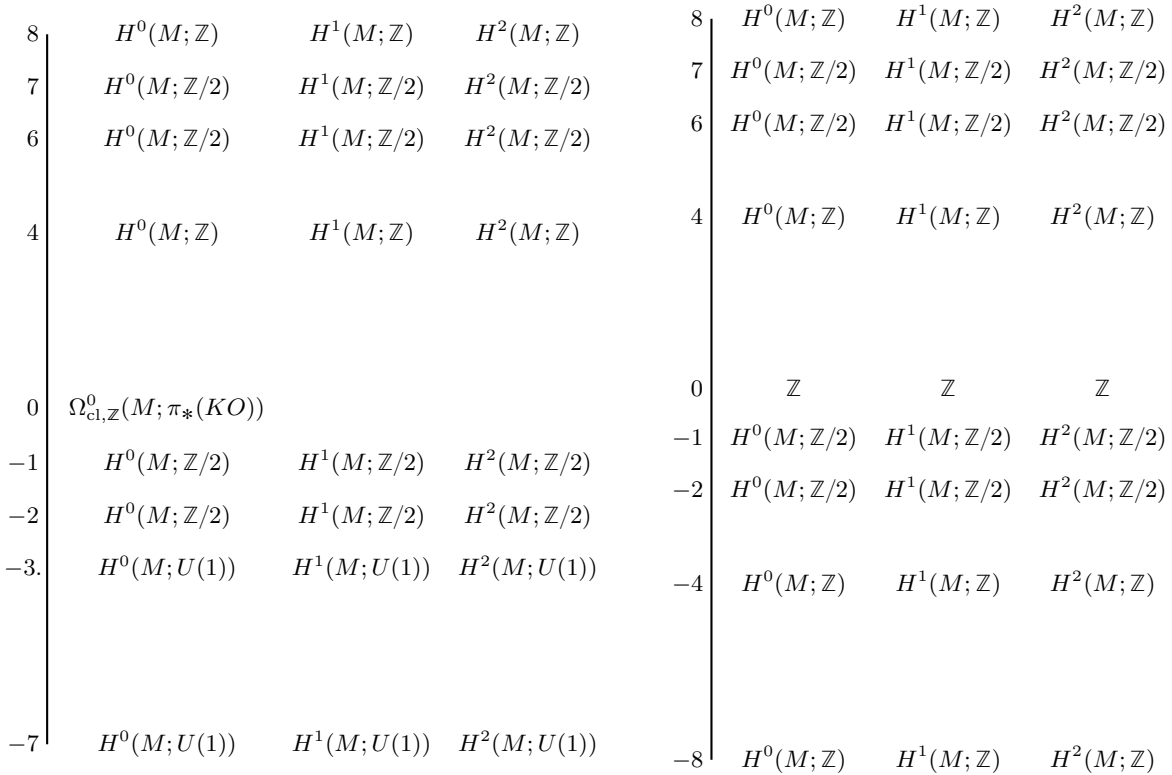
\begin{figure}[h!]\label{theE2sps}
\begin{minipage}{.5\textwidth}
\centering
\footnotesize
\begin{tikzpicture}
\matrix (m) [matrix of math nodes,
nodes in empty cells,nodes={minimum width=3ex,
minimum height=3ex,outer sep=0pt},
column sep=1ex,row sep=1ex]{
         \phantom{-}8    &  H^0(M;\ZZ) &  H^1(M;\ZZ)  & H^2(M;\ZZ)   \\
          \phantom{-} 7  & H^0(M;\ZZ/2)& H^1(M;\ZZ/2) &  H^2(M;\ZZ/2) \\
           \phantom{-}6  & H^0(M;\ZZ/2)& H^1(M;\ZZ/2) &  H^2(M;\ZZ/2) \\
             & & & \\
        \phantom{-} 4    &  H^0(M;\ZZ) &  H^1(M;\ZZ)  & H^2(M;\ZZ)   \\
             & & &  \\
             & & & \\
             & & &  \\
          \phantom{-}0   &  \Omega^0_{\rm cl,\ZZ}(M;\pi_*(KO))  & & \\
      -1    & H^0(M;\ZZ/2)& H^1(M;\ZZ/2) &  H^2(M;\ZZ/2)\\
       -2   & H^0(M;\ZZ/2)& H^1(M;\ZZ/2) &  H^2(M;\ZZ/2) \\
      -3.   & H^0(M;U(1)) &  H^1(M;U(1))  & H^2(M;U(1))  \\
             & & &\\
             & & &  \\
             & & & \\
        -7  & H^0(M;U(1)) & H^1(M;U(1))  &  H^2(M;U(1))  \\
\quad\strut &   \strut \\ };
\draw[thick] (m-1-1.east) -- (m-16-1.east) ;
\end{tikzpicture}
\end{minipage}%
\begin{minipage}{.5\textwidth}
\centering
\footnotesize
\begin{tikzpicture}
\matrix (m) [matrix of math nodes,
nodes in empty cells,nodes={minimum width=3ex,
minimum height=3ex,outer sep=0pt},
column sep=1ex,row sep=1ex]{
        \phantom{-} 8   &  H^0(M;\ZZ) &  H^1(M;\ZZ)  & H^2(M;\ZZ)   \\
           \phantom{-}7 & H^0(M;\ZZ/2)& H^1(M;\ZZ/2) &  H^2(M;\ZZ/2) \\
          \phantom{-} 6 & H^0(M;\ZZ/2)& H^1(M;\ZZ/2) &  H^2(M;\ZZ/2) \\
            & & &\\
         \phantom{-}4   &  H^0(M;\ZZ) &  H^1(M;\ZZ)  & H^2(M;\ZZ)   \\
            & & & \\
            &  & &\\
            &  & & \\
 \phantom{-}0  & \ZZ & \ZZ & \ZZ  \\
       -1  & H^0(M;\ZZ/2)& H^1(M;\ZZ/2) &  H^2(M;\ZZ/2)\\
       -2  & H^0(M;\ZZ/2)& H^1(M;\ZZ/2) &  H^2(M;\ZZ/2) \\
            & & &  \\
        -4 &H^0(M;\ZZ) &  H^1(M;\ZZ)  & H^2(M;\ZZ) \\
            & & & \\
            & & & \\
            & & &  \\
        -8 &H^0(M;\ZZ) &  H^1(M;\ZZ)  & H^2(M;\ZZ) \\        
\quad\strut &   \strut \\ };
\draw[thick] (m-1-1.east) -- (m-17-1.east) ;
\end{tikzpicture}
\end{minipage}
\vspace{-10mm}
\caption{The $E_2$-page for $\widehat{\rm KO}$ on the left vs. the $E_2$-page for $\op{KO}$ on the right.}
\end{figure}

\medskip
The spectral sequence mixes the topological and geometric data in a nontrivial way. One sees, for example, 
that the differential forms which survive to the $E_{\infty}$-page are constrained by the topological data.
In fact, even the torsion groups in lower degrees have an effect, as we shall see. 

\subsection{Identifying the differentials}
\label{Sec-Ident}

In this section we identify the differentials, which arise either as differential refinements of the
underlying topological differentials or as novel ones which mix the geometric and topological 
data. 

\begin{remark}[Observations on the differentials and corresponding maps]
We make a few preliminary observations.
\vspace{-1mm}
\begin{enumerate}[{\bf (i)}]
\item  The differentials contained only in the 1st quadrant of the 
spectral sequence reduce to the underlying topological differentials. 
\item As noted in \cite{GS3}, 
the canonical map $\mathcal{I}:\widehat{\op{KO}}(-)\to \op{KO}(-)$ induces an isomorphism on 
this first quadrant. 
\item From the hexagon diagram \eqref{kodfdiam}, the map ${\cal I}$ reduces to the Bockstein map $\beta:\op{KO}^{-1}(-;U(1))\to \op{KO}(-)$, associated to the fiber sequence
$$
\op{KO}\longrightarrow \op{KO}_{\mathbb{R}}\longrightarrow \op{KO}_{U(1)},
$$
on the flat part of $\widehat{\op{KO}}(-)$ and hence $\beta$ induces a morphism on 4th quadrants. 
\item The map $j$ in the hexagon diagram \eqref{kodfdiam}
induces an isomorphism on 4th quadrants. \end{enumerate}
\end{remark}

With these morphisms in hand, for the entries of  the form $H^i(M; \Z/2)$
we have the following identifications.

\begin{proposition}[The second differential for $\widehat{\op{KO}}$-theory AHSS at the $\ZZ/2$-entries] \label{d2difko}
The differentials of the form
$$
d_{2}:H^p(M;\ZZ/2)\longrightarrow H^{p+2}(M;\ZZ/2)
$$
appearing in the fourth quadrant of Figure 1
agrees with the topological case and are given by
$$
Sq^2:H^p(M;\ZZ/2)\longrightarrow H^{p+2}(M;\ZZ/2)\;.
$$
\end{proposition}
\theproof
Consider the morphism of spectral sequences arising from the canonical topological realization 
map $\mathcal{I}:\widehat{\op{KO}}(-)\to \op{KO}(-)$. This gives rise to a commutative diagram 
\(\label{impcmsq}
\xymatrix{
H^p(M;\ZZ/2)\ar[r]^-{\mathcal{I}'}\ar[d]^-{d_2} & H^p(M;\ZZ/2)\ar[d]^-{Sq^2}
\\
H^{p+2}(M;\ZZ/2)\ar[r]^-{\mathcal{I}'}  & H^{p+2}(M;\ZZ/2)
\;,
}
\)
where $\mathcal{I}'$ is induced by the corresponding map on coefficients. The coefficients on the left hand side
of \eqref{impcmsq} were determined from computing $\op{KO}({\rm pt};U(1))$ via the Bockstein sequence
$$
\xymatrix{
\op{KO}^*(-)\ar[r]&
\op{KO}^*(-;\RR)\ar[r]&
\op{KO}^*(;U(1))
\ar[r]^-{\beta_{U(1)}} & 
\op{KO}^{*+1}(-)
}.
$$
In particular, $\beta_{U(1)}$ gives rise to an isomorphism 
$\beta_{U(1)}:\op{KO}^{-2}({\rm pt};U(1))\cong \op{KO}^{-1}({\rm pt})\cong \ZZ/2$ and 
$\beta_{U(1)}:\op{KO}^{-3}({\rm pt};U(1))\cong  \ZZ/2$. \footnote{Recall again that we have
degree shifts up by 1 for the groups with  $U(1)$-coefficients.} Given the compatibility of
$\beta_{U(1)}$ with $\mathcal{I}'$ in the diamond diagram \eqref{kodfdiam}, i.e. $\mathcal{I}'j=\beta_{U(1)}$,
it follows that the top horizontal arrow in \eqref{impcmsq} is an isomorphism and hence $d_2=Sq^2$. 
\endofproof

\begin{proposition}
[Formula for the second differential in 4th quadrant of AHSS for $\widehat{\rm KO}$-theory]\label{prop-formula}
The differential 
$$
d_2:H^p(M;\ZZ/2) \xymatrix{\ar[r]&} H^{p+2}(M;U(1))\;,
$$
is given by $jSq^2$, with $j$ being the canonical inclusion 
$$
j:H^{p+2}(M;\ZZ/2) \xymatrix{\ar@{^{(}->}[r]&} H^{p+2}(M;U(1))
$$
induced by the inclusion of $\ZZ/2$ into $U(1)$ as the square roots of unity.
\end{proposition}
\theproof
The second differential $d_2$ defines a degree 2 increasing cohomology operation. From the exact sequence 
$\ZZ/2\into U(1)\xrightarrow{\times 2} U(1) $ and divisibility of $U(1)$, we get an exact sequence
(i.e., the Ext term vanishes) 
$$
0 \longrightarrow H^{p+2}(M;\ZZ/2)\xrightarrow{\;\;j\;\;} H^{p+2}(M;U(1))
\longrightarrow H^{p+2}(M;U(1))\longrightarrow 0\;.
$$
Hence the 2-torsion subgroup of $H^{p+2}(M;U(1))$ is precisely the image of $H^{p+2}(M;\ZZ/2)$ under 
$j$. As $H^{p+2}(M;\ZZ/2)$ is a $\ZZ/2$-vector space, the image 
of $d_2$ is necessarily 2-torsion. Then, by universal considerations, it follows that $d_2=j\phi$ for some operation $\phi:H^p(-;\ZZ/2)\to H^{p+2}(-;\ZZ/2)$. Since $H^*(K(\ZZ/2,n);\ZZ/2)$ is a polynomial ring generated by Steenrod squares, it follows that the only degree 2 increasing operations are $(Sq^1)^2$ and $Sq^2$, but the former vanishes by the Adem relations.
\endofproof

On the $E_{4k}$-page we have  differentials that involve geometric data and take the form 
\begin{equation}\label{geodf}
d_{4k}:\Omega^*_{\rm cl,\ZZ}(M;\pi_*(\op{KO})) \longrightarrow  {\rm ker}(d_{4k-1})/{\rm im}(d_{4k-1})\;,
\end{equation}
and it is interesting to identify what condition the cokernel forces on $d_{4k}(\omega)$ so that it vanishes. The next theorem provides a powerful characterization of these differentials, which witness an interesting 
geometry-topology interplay. To that end, we will need the following.

\begin{lemma}[Vanishing of $d_{4k}$ on exact forms]\label{phexfm}
Any element in $\Omega^0_{\rm \ZZ,{\rm cl}}(M;\pi_*(\op{KO}))$ of the form
$$\omega=\omega_0+d\eta_3+d\eta_7+\hdots \;,$$
with $\eta_{4k-1}$ a differential form of degree $4k-1$, lifts through the curvature map $\mathcal{R}:\widehat{\op{KO}}(M)\to \Omega^0_{\ZZ,{\rm cl}}(M;\pi_*(\op{KO}))$. Hence, all the differentials $d_{4k}:E^{0,0}_{4k}\to E^{4k,4k-1}_{4k}$ vanish on $\omega$.
\end{lemma}
\theproof
The commutativity of the diamond diagram \eqref{kodfdiam} implies 
$$
\mathcal{R}a(\eta_{3}+\eta_7+\hdots )=d\eta_3+d\eta_7+\hdots=\omega-\omega_0\;.
$$
Let $1\in \op{KO}(M)$ be the unit. Then $\omega_0\cdot 1\in \widehat{\op{KO}}(M)$ and $\omega_0\cdot 1+a(\eta_{3}+\eta_7+\hdots )$ defines such a lift through the curvature. See Proposition \ref{prop-E00SS}.
\endofproof

The differentials \eqref{geodf} are very rich and there seems to be no general formula which captures all degrees. We will compute this differential in low degrees. We also have the following theorem, which computes the first nonvanishing $d_{4k}$. 

\begin{theorem}[The differential $d_{4k}$ in AHSS for $\widehat{\rm KO}$-theory]\label{thm-d4d8}

\item [{\bf (i)}] For a closed graded differential form $\omega:=\omega_{8k}+\omega_{8k+4}+\hdots \in 
\Omega^0_{{\rm cl}, \ZZ}(M;\pi_*(\op{KO}))$, the differential $d_{8k}(\omega)$ is the image of the class 
$$
[\omega_{8k}]  \!\!\! \mod \ZZ  \;\; \in H^{8k}(M;U(1))
$$ 
in the quotient $\ker(d_{8k-1})/{\rm im}(d_{8k-1})$. 

\item [{\bf (ii)}] For a closed graded differential form $\omega=\omega_{8k+4}+\omega_{16k}+\hdots $, the differential $d_{8k+4}(\omega)$ is given by taking the image of the class 
$$
\tfrac{1}{2} [\omega_{8k+4}] \!\!\! \mod \ZZ  \;\; \in H^{4k}(M;U(1))
$$ 
in the quotient $\ker(d_{8k+3})/{\rm im}(d_{8k+3})$.
\end{theorem}
\theproof
{\bf (i)} For each $n\in \NN$, there is a smooth manifold $N$ and an $r=\max\{n,{\rm dim}(M)\}+1$ connected map $N\to \prod_{l\geq k}K(\RR,4l)$ factorizing $\widetilde{\omega}:M\to \prod_{l\geq k}K(\RR,4l)$, where $\widetilde{\omega}$ is a map representing the class $[\omega]$. By Lemma \ref{phexfm}, $d_{8k}$ factors through the de Rham cohomology groups. By naturality of the differentials, it is sufficient to determine $d_{8k}$ on $N$, which we can assume has the homotopy type of $\prod_{m\geq l\geq k}K(\RR,4l)$, for some large $m$. The K\"unneth formula with coefficients in both $\RR$ and $\ZZ$, along with the exponential sequence 
$$
\hdots \longrightarrow H^*(N;\ZZ)\longrightarrow H^*(N;\RR) \longrightarrow 
H^*(N;U(1))\overset{\beta}{\longrightarrow} H^*(N;\ZZ)
\longrightarrow \hdots,
$$
and degree considerations imply that 
$H^{4k}(N;U(1))\cong U(1)$, and $d_{8k}$ is a natural operation corresponding to the homotopy class of a map $K(\RR,4k)\to K(U(1),4k)$. Let ${\rm exp}(2\pi i\times ):K(\RR,4k)\to K(U(1),4k)$ be the map induced by the  exponential map. Then every class in $H^{8k}(K(\RR,4k),U(1))\cong U(1)$ is of the form ${\rm exp}(\lambda 2\pi i\times )$, for some $\lambda \in \RR$. Hence,
$$
d_{4k}={\rm exp}(\lambda 2\pi i\times )
$$
for some $\lambda$. To see that $\lambda=1$ for $4k\equiv 0 \mod 8$, take $M=S^{8k}$. The canonical bundle $\beta\to S^8$ satisfies $\langle {\rm Ph}(\beta),[S^8]\rangle=1$. Fix a connection $\nabla$ on $\beta$ and let ${\rm Ph}(\nabla)$ be the corresponding Chern-Weil form. Then 
$$
1=d_{8}({\rm Ph}(\nabla)_{8})={\rm exp}(\lambda 2\pi i\langle {\rm Ph}(\beta),[S^8]\rangle) 
\; \Rightarrow \; \lambda\in \ZZ\,.
$$
Since $\beta$ generates ${\rm KO}(S^8)$, we deduce $\lambda=1$. The same argument applied to the bundles $\beta^{\boxtimes k}\to \prod_{i=1}^kS^{8}$ gives $\lambda=1$ for $d_{8k}$, since $\beta^k$ generates ${\rm KO}(S^{8k})$. 

For part {\bf (ii)}, we use $M=S^{8k+4}$. When $k=0$, the quaternionic Hopf bundle $\alpha\to S^4$ generates ${\rm KO}(S^4)$ and $\langle {\rm Ph}(\alpha),[S^4]\rangle=2$. Again, fixing a connection on $\alpha$ and using that $\alpha$ is a generator, we find that $2\lambda=1 \Rightarrow \lambda=\frac{1}{2}$. Arguing analogously for the bundles $\alpha\boxtimes(\beta^{\boxtimes k})\to S^4\times \prod_{i=1}^kS^8$ gives $\lambda=\frac{1}{2}$ in the general case for $d_{8k+4}$.
\endofproof

\begin{corollary}
Let $M$ be any closed smooth manifold. Let $\omega$ be a differential form of degree $\vert \omega\vert\equiv 0 \mod 4$ and assume $M$ is $\vert \omega\vert-1$ connected, the following condition is necessary (but not sufficient) for $\omega$ to be a Pontrjagin character form.
\begin{enumerate}
\item  If $\vert \omega\vert\equiv 0\mod 8$  then $\omega$ is integral.
\item If $\vert \omega \vert\equiv 4\mod 8$, then $\omega$ is integral and divisible by 2.
\end{enumerate}
\end{corollary}
\theproof
The Simons-Sullivan description of K-theory via structured vector bundles applies in the ${\rm KO}$-setting. By construction, the Pontrjagin character factors through the cycle map $\widehat{\rm cyc}:{\bf K}{\rm O}_{\nabla}(M)\to \widehat{{\rm KO}}(M)$. It follows that any form lifting to $\widehat{{\rm KO}}(M)$ also lifts to ${\bf K}{\rm O}_{\nabla}(M)$ as the pontrjagin character of some (virtual) vector bundle with connection. By Theorem \ref{thm-d4d8}, if $\omega$ lifts to $\widehat{\rm KO}(M)$, and hence ${\bf K}{\rm O}_{\nabla}(M)$, all differentials must vanish on $\omega$. In particular, $d_{\vert \omega\vert}$ must vanish. The claim then follows from Theorem \ref{thm-d4d8}. 
\endofproof


The following example reveals the intricate interplay between the topology and geometry of manifolds. 

\begin{example}[Vector bundles on $\CC P^2$]
Consider the complex projective plane $\CC P^2$. The relevant terms in the fourth 
quadrant (on the $E_2$ page) of the spectral sequence for ${\rm KO}$ are given as follows.
\begin{center}
{\footnotesize
\begin{tikzpicture}
\matrix (m) [matrix of math nodes,
nodes in empty cells,nodes={minimum width=3ex,
minimum height=3ex,outer sep=0pt},
column sep=1ex,row sep=1ex]{
          \phantom{-}0   &  \Omega^0_{\rm cl,\ZZ}(\CC P^2;\pi_*(KO))  & & & & \\
      -1    & \ZZ/2& 0 &  \ZZ/2 & 0 & \ZZ/2\\
       -2   & \ZZ/2& 0 &  \ZZ/2 & 0 &\ZZ/2 \\
      -3  & U(1) & 0  & U(1) & 0 &U(1)  \\
\quad\strut &   \strut \\ };
\draw[thick] (m-1-1.east) -- (m-5-1.east);
\end{tikzpicture}}
\end{center}

\vspace{-4mm} 
\noindent By Proposition \ref{d2difko}, $d_2$ is given by $Sq^2$ on groups with 
target $\ZZ/2$ and $jSq^2$ on those with target $U(1)$. Now $Sq^2$ maps the 
generator $x\in H^2(\CC P^2;\ZZ/2)\cong \Z/2$ to $x^2$. The generator 
$x=w_2({\cal L})$, where ${\cal L}\to \CC P^2$ is the pullback of the canonical
 line bundle and $w_2$ denotes the second Stiefel-Whitney class. From the Wu formula 
 and Proposition \ref{thm-d4d8}, we see that 4-form $\omega$ lifts to $\widehat{{\rm KO}}(M)$ provided
$$
{\rm exp}\Big(\pi i \int_{\CC P^2}\omega \Big)= j_2\left\langle Sq^2(w_2({\cal L})),[\CC P^2]\right\rangle=j_2\left\langle w_2(\CC P^2)\cup w_2({\cal L}),[\CC P^2]\right\rangle.
$$
Applying this to the $\hat A$ genus, we have 
$$
{\rm exp}\Big(\pi i\int_{\CC P^2}\frac{p_1(\CC P^2)}{24}\Big)={\rm exp}\big(\tfrac{3}{24}\pi i\big)\not\in \ZZ/2\subset U(1)\,.
$$
Hence, $\hat A$ is not the Pontrjagin character of a vector bundle. Of course this is expected, since $\CC P^2$
 is not spin.

\end{example}

\begin{remark}[Pontrjagin character (forms)]
{\bf (i)} The forms $\omega$ which lift to differential $\op{KO}$-theory are precisely those forms which 
are in the image of the Pontrjagin character, by Proposition \ref{prop-E00SS}. Thus, we see from 
Theorem \ref{thm-d4d8} that the differentials in the 4th quadrant of the spectral sequence completely determine the image of the Pontrjagin character. 
\item {\bf (ii)} 
In the case 
of complex $K$-theory, it is known that the structure of the classical AHSS is intimately related to the values 
of the Chern character (see \cite{AH}\cite{AH2}\cite{Buc1}\cite{Buc2}). 
In the real topological case, such a relationship with the Pontrjagin character is expected but does not seem to 
have been investigated before. We find it striking the our construction of the AHSS for differential KO-theory 
makes  the relationship between the Pontrjagin character form and the differentials in the AHSS 
completely manifest.
\end{remark}

\subsection{Consequences for the Pontrjagin character}
\label{Sec-IntPont}

We now turn to the problem of characterizing the image of the Pontrjagin character map. Before presenting 
the  main theorem in this section, we will need to refine the result of Buchstaber \cite{Buc1}\cite{Buc2} 
to the differential setting. We begin by observing that if the cohomology groups of a space $X$ are finitely 
generated, the torsion part in $H^*(X;\ZZ)$ is precisely the image of the Bockstein map 
$\beta_{U(1)}:H^*(X;U(1))\to H^*(X;\ZZ)$ associated to the exponential sequence. Then, for every 
finite-dimensional smooth manifold $M$, it follows  that we have a short exact sequence
$$
0\longrightarrow \frac{H^*(M;\RR)}{H^*(M;\ZZ)_{\rm free}}\longrightarrow 
H^*(M;U(1))\xrightarrow{\beta_{U(1)}} H^{*+1}(M;\ZZ)_{\rm tor}\longrightarrow 0\;.
$$
Since  $T^{b_*}:=\frac{H^*(M;\RR)}{H^*(M;\ZZ)_{\rm free}}$ is divisible as an abelian group it follows 
that ${\rm Ext}^1(H^{*+1}(M;\ZZ)_{\rm tor},T^{b_*})=0$ and the sequence splits (noncanonically). 
This gives an identification for each $k$
\(\label{splu1zr}
H^k(M;U(1)){\cong} H^{k+1}(M;\ZZ)_{\rm tor}\oplus T^{b_k}\;,
\)
where $T^{b_k}$ denotes the torus of dimension the $k$-th Betti number $b_k:={\rm dim}\big(H^k(M;\RR))$. 
Using this decomposition, we can study the \emph{$p$-primary part of the differentials}
in the AHSS by restricting to  the differentials factoring through the subgroup 
$(H^{k+1}(M;\ZZ)_{\rm tor})_p\oplus T_p^{b_k}$, where $T^{b_k}_p$ denotes 
the $p$-torsion subgroup of the torus. The full differential will be a product of these 
components (in $H^k(M;U(1))$) as $p$-varies. 
We consider maps 
of the form 
\(
\label{def-jp}
j_{p^m}:H^k(M;\ZZ/p^m)\longrightarrow H^k(M;U(1))
\)
induced by a representation of the $p$-primary group $\ZZ/p^m$. Realizing $\ZZ/p^m$ as the 
$p^m$-roots of unity in $U(1)$, we can regard $j_{p^m}$ simply as the inclusion. 

\begin{remark}[Image of $d_{4k}$]
Note that the image of the differentials $d_{4k}:E^{0,0}_{4k}\to E^{4k,-4k+1}_{4k}$ necessarily lands 
in the toric part of the decomposition $H^k(M;U(1)){\cong} H^{k+1}(M;\ZZ)_{\rm tor}\oplus T^{b_k}$. 
Indeed, for a differential form $\omega_{4k}$ of degree $4k$, $d_{4k}(\omega_{4k})=[\omega_{4k}]\mod \ZZ$ 
lies in the image of the exponential map $e:H^{4k}(M;\RR)\to H^{4k}(M;U(1))$, hence the kernel of $\beta_{U(1)}$. 
\end{remark}




With these definitions, and with the same definition for the generators $\tilde{\a}$ as in 
Proposition \ref{diff-Ko-ahss} above,  we have the following.

\begin{proposition} [Odd-primary part of the differential]\label{prop-oddprim}
Let $M$ be a compact manifold. We have the following formulas for the $p$-primary 
component of the differentials in $\widehat{\op{KO}}$ for $s,t\neq 0$

\item [{\bf (i)}] For $p$ odd:
$$
d^{s,-4t+1}_{4r(p-1)+1}(x)_p=\epsilon_{2r}\cdot j_p\big( P^{2r}(\beta_{U(1)}(x/p^{2r-1}))\cdot \tilde{\alpha} \big).
$$

\item [{\bf (ii)}] For $p=2$:
$$
d^{s,-8t+1}_{8r+1}(x)_2= \epsilon_{4r}\cdot j_2\big( Sq^{8r}(\beta_{U(1)}(x/2^{4r-1}))\cdot \tilde{\alpha} \big),
$$
where $j_p$ is given by  \eqref{def-jp} and $\epsilon_{k}\neq 0$ mod $p$ for $k=2r$, $4r$.
\end{proposition}
\theproof
%
The morphism of spectral sequences induced by the Bockstein map associated to the fiber sequence 
$$
\op{KO}\longrightarrow \op{KO}_{\RR}\longrightarrow \op{KO}_{U(1)}
$$
is the zero map at the $E_n^{s,4t-1}$ entries. In particular, starting with the $E_2$-page we have a well-defined 
boundary map at these entries. By the argument in \cite[Proposition 17]{GS1}, we can identify the boundary map 
with the Bockstein morphism $\beta_{U(1)}$ on ordinary cohomology associated to the exact sequence 
$\ZZ\to \RR\to U(1)$. Then, by Buchstaber's  result \cite{Buc1}\cite{Buc2}, the fact that the boundary map 
commutes with the differentials, and the compatibility $\beta_p=\beta_{U(1)}{j}$, we immediately identify 
the $p$-primary part
$$
d^{s,-4t+1}_{4r(p-1)+1}(x)_p= 
\epsilon_{2r}\cdot j_p(\beta_{p} P^{2r}\beta_{U(1)})(x/p^{2r-1})\cdot \tilde{\alpha} \;.
$$

\vspace{-8mm}
\endofproof

Proposition \ref{prop-oddprim} allows us to establish the following result.

\begin{theorem}[Congruence for the Pontrjagin character]\label{intphnum}
Let $M$ be a compact manifold.  
Define the function $\phi:\RR\to \ZZ$ by 
$$\phi(x)=\left\{\begin{array}{cc}
[x], & x\not\in \ZZ
\\
x-1, & x\in \ZZ\end{array}\right.,
$$
where $[x]$ denote the greatest integer $<x$, for $x$ a real number. Consider the number
$$
n_k=\prod_{p\geq 2}p^{\phi\big(\tfrac{k}{p}\big)}.
$$ 
Then a necessary condition for a differential $4k$-form $\omega$ to be the $4k$-component of a Pontrjagin form is that $n_k\cdot \omega$ is integral, modulo denomiators involving a power of 2.  
\end{theorem}
\theproof
If the form $\omega$ lifts to differential $\op{KO}$, then it must
be a permanent cycle. By Theorem \ref{thm-d4d8}, we must have
$$
d_{4k}(\omega_{4k})=
\tfrac{1}{2^{k\;{\rm mod}\; 2}}[\omega_{4k}]\mod \ZZ
$$ 
is in the image of $d_{4r(p-1)+1}$ for all $4r(p-1)+1<4k$, i.e., $r<\frac{4k-1}{4(p-1)}$. As in \cite{Buc2}, 
one can show that there exists a smooth manifold $M$ with the following properties 
\footnote{This is achieved by approximating the 
Eilenberg-MacLane complexes 
by smooth manifolds via a sequence of surgeries.} 
\begin{itemize}
\item Fix a prime $p\neq 2$. For any $r<\frac{4k-1}{4(p-1)}$, there is a class $\beta_{U(1)}(x)$ of degree $4k-4r(p-1)-1$, 
which is nonzero and divisible by $p^{2r-1}$, provided $\rho_p\beta_{U(1)}(x)$ lies in the nonvanishing regime 
of $P^{2r}$.  
\item For every such $\beta_{U(1)}(x)\neq 0$, the primary operations 
$j_p P^{2r}(\beta_{U(1)}(x)/p^{2r-1})\neq 0$.
\footnote{Note that we are taking $0$ to be the identity element in $H^*(M; U(1))$. We could have used
$1=\exp(2\pi i n)$ for the multiplicative identity, but that would lead to possible confusion upon 
asking for what it means for a differential to `vanish'.}
\item The operation  $\beta_{U(1)}j_p P^{2r}(\beta_{U(1)}(x)/p^{2r-1})=0$.
\end{itemize}
Taking degrees into consideration, we find that the regime of $r$'s for which ${P}^{2r}$ necessarily vanishes 
is precisely $r<[\frac{4k-1}{4p}]=[\frac{k}{n}-\frac{1}{4p}]=\phi\big(\frac{k}{n}\big)$. Then for all 
$r<{\rm min}\{\phi\big(\frac{k}{n}\big),[\frac{4k-1}{4(p-1)}]\}=\phi\big(\frac{k}{n}\big)$, we have that 
$j_p P^{2r}(\beta_{U(1)}(x)/p^{2r-1})\neq 0$ and the differential $d_{4k}:E^{0,0}_{4k}\to E^{4k,4k-1}_{4k}$ 
vanishes precisely if $d_{4k}(\omega_{4k})=j_pP^{2r}(\beta_{U(1)}(x)/p^{2r-1})$, for some $x$. This image is in the $p$-torsion subgroup of the torus $T^{b_{4k}}$, since 
\\
$\beta_{U(1)}j_p P^{2r}(\beta_{U(1)}(x)/p^{2r-1})=0$. 
If $r=1$,  then this implies that
$2^mp\cdot \omega_{4k}$ is integral, for some $2^m$ arising from the image of the differentials with source $H^*(M;\ZZ/2)$. 
Passing to the quotient, we find that, for $r=2$, $p^2\cdot \omega_{4k}$ is integral, modulo denominators involving powers of 2.  Iterating this process, we arrive at the integer $n=\phi\big(\frac{k}{n}\big)$. 
\endofproof

One would like to be able to construct interesting Pontrjagin character forms using  only the cohomology of the manifold $M$. The following gives such a construction.

\begin{proposition}[Denominators of Pontrjagin characters] \label{p-tordfko}
Fix an integer $\ell$ and a prime $p$. Let $r$ and $k$ be integers with $4k+r(p-1)=\ell$ and consider 
the following property of the pair $(r,k)$. There is a class $x_r\in H^{4k}(M;\ZZ)$ (which depends on $r$) 
such that,  for $p\neq 2$,
\begin{enumerate}
\item $x_r$ is $p$-torsion and is divisible by $p^{2r-1}$, but not $p^{2r}$,
\item $P^{2r}(\rho_p(x_r/p^{2r-1}))\neq 0$, 
\item $\beta_pP^{2r}(\rho_p(x_r/p^{2r-1}))=0$,
\end{enumerate}
with $\beta_p:H^*(-;\ZZ/p)\to H^{*+1}(-;\ZZ)$ the Bockstein associated to the sequence 
$\ZZ\xrightarrow{\times p} \ZZ \to \ZZ/p$. Call such pairs admissible. Then there exists a vector bundle $\xi_\ell$, depending on the sum $\ell=4k+4r(p-1)$, such that the $p$-primary part of the denominator of ${\rm Ph}_{\ell}(\xi)$ is precisely $p^{s}$, 
where $s$ is the number of admissible pairs. 
\end{proposition}
\theproof
Fix any admissible pair $(r,k)$. Since $x_r$ is torsion, it must be in the image of $\beta_{U(1)}$. Since $x_r$ is divisible by $p^{2r-1}$, we can form the class $x^r/p^{2r-1}$. Since $x_r$ is not divisible by $p^{2r}$, the mod $p$ reduction of this class does not vanish. Since $P^{2r}(\rho_p(x_r/p^{2r-1}))\neq 0$, it follows that $j_pP^{2r}(\rho_p(x_r/p^{2r-1}))\neq 0$, since $j_p$ is injective. Finally, since $\beta_pP^{2r}(\rho_p(x_r/p^{2r-1}))=0$ and $\beta_{U(1)}j_p=\beta_p$, it follows that $j_pP^{2r}(\rho_p(x_r/p^{2r-1}))$ is in the image of the exponential map $H^{\ell}(M;\RR)\to H^{\ell}(M;U(1))$. Choose a differential form $\omega_{\ell}$ whose de Rham class maps to $j_pP^{2r}(\rho_2(x_r/p^{2r-1}))$ under the exponential. Being in that image, $\omega_{\ell}$ is killed by 
$d_{\ell}:E^{0,0}_{\ell}\to E^{\ell,-\ell+1}_\ell$ and must be the degree $\ell$ component of a Pontrjagin character. Since $j_pP^{2r}(\rho_2(x_r/p^{2r-1}))$ is non-vanishing and lies in $p$-torsion, it follows that $[\omega_{\ell}]$ must have denominator $p$. Moreover, each distinct admissible pair contributes a copy of $p$, yielding the denominator $p^{s}$.
\endofproof

We now illustrate the utility of the spectral sequence via some concrete applications of Theorem 
\ref{thm-d4d8}. We begin with a few representative statements, and continue with more 
applications in subsequent sections.

\begin{corollary}[Classes contributing to denominators] 
Suppose the order of the $p$-primary component $H^{4k}(M;\ZZ)_p$ is $p^{2r}$ and suppose 
$P^{2r}:H^{4k}(M;\ZZ/p)\to H^{4k+4r(p-1)}(M;\ZZ/p)$ is nonvanishing for all classes in the 
image of the mod $p$ reduction and has image in the mod $p$ reduction. Then any class 
$x\in  H^{4k}(M;\ZZ)_p$ of order $p^{2r}$ contributes a copy of $p$ to  the product in 
the denominator of ${\rm Ph}_{4k+4r(p-1)}(\xi)$ for some vector bundle $\xi$.
\end{corollary}
\theproof
The element $p^{2r-1}x$ is $p$-torsion since the order is $p^{2r}$. Moreover, $\rho_p(x)\neq 0$ since otherwise, $x=py$ for some $y$ and $p^{2r-1}x=p^{2r}y=0$. Since $P^{2r}$ is nonvanishing on the image of $\rho_p$ it 
follows that $j_pP^{2r}(\rho_2(x))\neq 0$. Hence, the claim follows from Proposition \ref{p-tordfko}.
\endofproof

\begin{proposition}[Degree 8 components of the Pontrjagin character]\label{prop-deg8ph}
Let $(V,\nabla)$ be a real vector bundle on a compact (orientable) 8-manifold $M$. Suppose all elements in the image of 
$\beta_{U(1)}:H^{3}(M;U(1))\to H^4(M;\ZZ)$ are divisible by $2$. Assume also that $H^1(M;\ZZ/2)\cong H^2(M;\ZZ/2)\cong 0$  \footnote{This assumption is probably not necessary, as we expect the differentials $d_6$ and $d_7$ to be given by Steenrod squares of degree exceeding $1$ and $2$. However, since we have not computed these differentials, we will add this condition to ensure that they vanish.}. Then 
$$
\int_M{\rm Ph}(V,\nabla)_8=\frac{1}{12}\int_M p_1^2(V,\nabla)-2p_2(V,\nabla)\in \ZZ
$$
\end{proposition}
\theproof
Note that the differentials can be  labeled schematically as $d_{i}^{{\rm right}, {\rm down}}$,
whose effect on the grid  is to move  to the right by $i$ and downwards by $i-1$. 
From this we see that the only nonzero differentials which can have target in bidegree $(8,-7)$ are 
$d_5$, $d_6$ and $d_7$. 
The differentials $d^{2,-2}_6$ and $d^{1,-1}_7$ vanish since $H^1(M;\ZZ/2)\cong H^2(M;\ZZ/2)\cong 0$. By Proposition 
\ref{diff-Ko-ahss}, we have the formula
$$
d^{3,-3}_5=j_2 Sq^4 \rho_2\beta_{U(1)}: H^3(M; U(1)) \longrightarrow H^8(M; U(1))\;.
$$
Since $\beta_{U(1)}(x)$ is divisible by $2$ for all $x\in H^3(M;U(1))$, it follows that 
$\rho_2\beta_{U(1)}(x)=0$ and the differential vanishes. Thus, for a differential form 
$\omega_{8}$ to be the degree 8 component of a Ponrjagin character form, it is necessary 
and sufficient that $\omega_8$ must represent an integral class.
\endofproof

\medskip
We will consider consequences and applications of the above results in Section 
\ref{Sec-higher}.

\subsection{Computations of $\widehat{\rm KO}$-groups of spheres}
\label{Sec-comp} 

In this section we will put to work the machinery we developed in Section  \ref{Sec-AHSS-KOhat}
of this section. We will demonstrate how the AHSS for $\widehat{\rm KO}$-theory 
is a powerful tool for computations, as has been shown for other theories in
\cite{GS3}\cite{GS4}\cite{GS5}\cite{GS6}.  We do this for classical spaces with high symmetry, namely 
spheres, leaving a more comprehensive account to be developed more fully elsewhere.

\medskip
We first start by recalling basic computations in classical KO-theory in relation to spheres.
For $S^1$, $\op{KO}(S^1)=\ZZ\oplus \ZZ/2$ generated by the M\"obius bundle 
and the trivial bundle. This then implies that the reduced theory is given by 
$\widetilde{\rm KO}(S^1)=\ZZ/2$.  For higher spheres, i.e. $n\geq 1$, 
$\widetilde{\rm KO}(S^n)\cong \pi_n(\mathrm{BO})\cong \pi_{n-1}(\mathbf(O))$
is  isomorphic to  $\ZZ$, $\ZZ/2$, or  $0$ according to 
the value of $n$ mod 8. That is, one has the following pattern:

\begin{center} 
\begin{tabular}{|c||cccccccc|}
\hline
$n \mod 8$ & $1$ & $2$ & $3$ & $4$ & $5$ & $6$ & $7$ & $8$
\\
\hline
\hline
$\widetilde{\rm KO}(S^n)$ & $\ZZ/2$ & $\ZZ/2$ & $0$ & $\ZZ$ & $0$ & $0$ & $0$& $\ZZ$
\\
\hline
\end{tabular} 
\end{center}

\medskip
Let $X$ be a CW-complex of finite dimension.
Bott periodicity of the orthogonal group in terms of KO-theory
gives that the map 
$$
\op{KO}^*(X ) \otimes \op{KO}^*( S^{8n}) \xrightarrow{\;\cong\;} \op{KO}^*(X \times S^{8n})\;, 
$$
induced by the tensor product of bundles, is an isomorphism \cite{Bo} \cite{Real}. 
There is also an isomorphism $\wedge: \widetilde{\rm KO}(X) \otimes \widetilde{\rm KO}(S^8)
\cong  \widetilde{\rm KO}(X \wedge S^8)=\widetilde{\rm KO}(\Sigma^8X)$. That is, 
$\wedge \sigma: \widetilde{\rm KO}(X) \cong \widetilde{\rm KO}(\Sigma^8X)$ for a generator
$\sigma$ of $\widetilde{\rm KO}(S^8)$ (see \cite[Theorem 5.13]{MT}).
A proof at the level of spectra is given in \cite{Ki2}.

\begin{remark}[Complex vs. real K-theory of spheres] \label{rem-cplkr} Recall (see Section \ref{Sec-cplx} in Section \ref{Ch-diffKO})
the two main maps back  and forth between KO-theory and K-theory given as (see \cite{Str}) 
$$
\begin{tabular}{ccccc} 
$m_O: \op{KO} \longrightarrow \op{K}$ & \quad $\eta \mapsto 0$, & \quad $\alpha \mapsto 2 u^2$, & \quad $\beta \mapsto u^4$, &
\\
$f_U: \op{K} \longrightarrow \op{KO}$ & \quad  $1 \mapsto 2$,  & \quad $u \mapsto \eta^2$, & \quad  $u^2 \mapsto \alpha$, & \quad
$u^3 \mapsto 0$.
\end{tabular} 
$$
\begin{enumerate}[{\bf (i)}]
\item In one direction (see \cite[Ch. IV, Theorem 5.12]{MT})

\begin{itemize}
\item If $n\equiv 0$ mod 8 then $c: \widetilde{\rm KO}(S^n)\cong \widetilde{\rm K}(S^n)$.
\item If $n\equiv 1, 2$ mod 8 then $\wedge \eta: \widetilde{\rm KO}(S^{n-1})\cong \widetilde{\rm KO}(S^n)$
is an isomorphism. 
\item If $n\equiv 4$ mod 8 then $c: \widetilde{\rm KO}(S^n)\to \widetilde{\rm K}(S^n)$ is an isomorphism 
onto $2(\widetilde{\rm K}(S^n))$. 
\end{itemize}

\item In the other direction (\cite[Ch. IV, Theorem 6.1]{MT}):
For $n>0$, the map $r: \widetilde{\rm K}(S^n) \to \widetilde{\rm KO}(S^n)$ 
is 
\begin{itemize} 
\item an isomorphism onto $2(\widetilde{\rm KO}(S^n))$ if $n\equiv 0$ mod 8. 
\item an epimorphism if $n\equiv 2$ mod 8.
\item an isomorphism if $n\equiv 4$ mod 8.
\end{itemize}
\end{enumerate}
\end{remark}

We now consider differential refinements of the calculations of KO-theory of 
spheres. As we saw above in the topological case, the latter are closely related
to the coefficients of the theory. Indeed, something analogous will occur in the 
case of $\widehat{\rm KO}$-theory.

\medskip
\begin{proposition}[Differential KO-theory of $4n$-spheres]
Fix a Riemannian metric $g$ on the sphere $S^{4n}$, with $n\not\equiv 0 \mod 2$. Then we have 
an isomorphism 
\footnote{Note that we are still considering differential refinements of \emph{reduced} 
KO-theory, so that the expressions below will have an extra copy of $\ZZ$ from the topological side. 
Each factor of $\langle d {\rm vol}\rangle_\ZZ$ will be the harmonic form representative 
corresponding to that topological class.}
$$
\widehat{\op{KO}}(S^{4n})\overset{g}{\cong} \langle 2\cdot d{\rm vol} \rangle_{\mathbb{Z}}
\oplus \bigoplus_{k=1}^{n}d\Omega^{4k-1}\;,
$$
where $d{\rm vol}$ is a normalized harmonic volume form on $S^{4n}$ and 
$\langle d{\rm vol} \rangle_{\mathbb{Z}}$ denotes its $\mathbb{Z}$-linear span.
\end{proposition}
\theproof
We start with cohomology, where 
$$
H^k(S^{4n};A)\cong \left\{\begin{array}{cc}
A & k=0,4n,
\\
0 & \text{otherwise},
\end{array}\right.
$$
for any abelian group $A$. Then,  lacunary considerations show that the only relevant 
nonvanishing differentials are 
$$
d_{4n}:\Omega_{\rm cl,\ZZ}^0(S^{4n};\pi_*(\op{KO}))\longrightarrow H^{4n}(S^{4n};U(1))\cong U(1)\;.
$$
Recall that a choice of metric $g$ gives rise to a Hodge decomposition at the level of differential forms 
$$
\Omega^n(M)\cong \mathcal{H}_{\Delta}(M)\oplus d\Omega^{n-1}(M)\oplus d^{\dagger}\Omega^{n+1}(M)\;,
$$
where $\Delta$ denotes the Laplace-Beltrami operator, $d^{\dagger}$ is the adjoint of $d$ (with respect to $g$) 
and $\mathcal{H}_{\Delta}(M)$ denotes the subspace of harmonic forms on $M$. By Theorem \ref{thm-d4d8}, 
$d_{4n}(\omega)=\frac{1}{2}[\omega_{4n}]=0$ mod $\ZZ$, which vanishes precisely if $\frac{1}{2}\omega_{4n}$
has integral periods. Using the Hodge decomposition, we can identify all such $\omega_{4n}$'s as differential forms 
that read $n\cdot 2\cdot d{\rm vol}+d\xi$, with $n\in \ZZ$, $d{\rm vol}$ the normalized harmonic volume 
form on $M$, and $\xi=\xi_3+\xi_7+\hdots$ is an arbitrary differential form in $\Omega^{-1}(M;\pi_*(\op{KO}))$. This  immediately gives the desired identification (which depends on the metric $g$). 
\endofproof

In degrees $8k$, the identification of the differentials $d_{8n}(\omega)=[\omega_{8n}]\mod \ZZ$ along with a similar argument gives the following.

\begin{proposition}
[Differential KO-theory of $8n$-spheres]
Fix a Riemannian metric $g$ on a sphere $S^{8n}$. Then we have an induced isomorphism 
$$
\widehat{\op{KO}}(S^{8k})\overset{g}{\cong}  \langle d{\rm vol} \rangle_{\mathbb{Z}} \oplus \bigoplus_{k=1}^{2n}d\Omega^{4k-1}\;.
$$
\end{proposition}

We also have the following relationships, akin to the topological case. The maps $c$ and $\eta\wedge$ in Remark \ref{rem-cplkr} are extended to the differential case in the obvious way
as in the differential refinement of the Bott sequence in Section \ref{Sec-cplx-diff}. 

\begin{proposition}[$\widehat{\rm KO}$-theory vs. $\hat{\op K}$-theory of spheres]
Let $\widehat{K}$ denote complex differential $\op{K}$-theory. 
\item {\bf (i)} Completely analogously 
to the real case, we have the identification 
$$
\widehat{\op{K}}(S^{2n})\overset{g}{\cong} \langle d{\rm vol} \rangle_{\ZZ} \oplus \bigoplus_{k=1}^{n}d\Omega^{2k-1}\;.
$$
\item {\bf (ii)} We have the following identifications:
\begin{enumerate}[{\bf (a)}] 
\item In one direction

\begin{itemize}
\item If $n\equiv 0 \mod 8$, then $c:\widehat{\op{KO}}(S^n)\cong \widehat{\op{K}}(S^n)$.
\item If $n\equiv 1,2 \mod 8$ then $\wedge \eta :\widehat{\op{KO}}(S^{n-1})\cong \widehat{\op{KO}}(S^{n})$.
\item If $n\equiv 4\mod 8$ then $c:\widehat{\op{KO}}(S^n)\to \widehat{\op{K}}(S^n)$ is an isomorphism onto $2\widehat{\op{K}}(S^n)$.
\end{itemize}

\item  In the other direction: for $n>0$, the map $r:\widehat{\op{K}}(S^n)\to \widehat{\op{KO}}(S^n)$ is

\begin{itemize}
\item an isomorphism onto $2\widehat{\op{KO}}(S^n)$.
\item an epimorphism if $n\equiv 2 \mod 8$.
\item an isomorphism $n\equiv 4\mod 8$.
\end{itemize}
\end{enumerate}

\end{proposition}
\theproof
To prove {\bf (a)}, recall that the complexification map sends $c:\alpha\mapsto 2u^2$ and $c:\beta\mapsto u^4$, 
where $\alpha$ generated $\widetilde{\op{KO}}(S^4)$ and $\beta$ generates $\widetilde{\op{KO}}(S^8)$. 
The morphism of spectral sequences induced by complexification therefore gives rise to the commutative diagrams 
$$
\xymatrix{
\Omega^0_{\ZZ,{\rm cl}}(S^{8n+4};\pi_*(\op{KO}))\ar[r]^{\times 2}\ar[d]_{d_{8n+4}} & \Omega^0_{\ZZ,{\rm cl}}(S^{8n+4};\pi_*(\op{K}))\ar[d]^{d^{\prime}_{8n+4}}
\\
H^{8n+4}(S^{8n+4};U(1))\ar[r]^{\times 2} & H^{8n+4}(S^{8n+4};U(1))
}
\qquad \text{and} \qquad
\xymatrix{
\Omega^0_{\ZZ,{\rm cl}}(S^{8n};\pi_*(\op{KO}))\ar[r]\ar[d]_{d_{8n}} & \Omega^0_{\ZZ,{\rm cl}}(S^{8n};\pi_*(\op{K}))\ar[d]^{d^{\prime}_{8n}}
\\
H^{8n}(S^{8n};U(1))\ar@{=}[r] & H^{8n}(S^{8n};U(1))\;.
}
$$
From the identification of the differentials, we see that the form $d{\rm vol}$ is mapped to $2d{\rm vol}$ in dimensions $8n+4$ and to $d{\rm vol}$ in dimensions $8n$, where $d{\rm vol}$ is the normalized harmonic volume form (with respect to some fixed metric). Combined with the Hodge decomposition, this immediately gives rise to the first and third identifications. From the differential cohomology diamond diagram \eqref{kodfdiam} we calculate 
$$
\widehat{\op{KO}}(S^n)\cong \ZZ/2\oplus \bigoplus_{k=1}^{n}d\Omega^{4k-1}\;,
$$
for both $n\equiv 1,2\mod 8$. The map $\eta\wedge$ only acts on the $\ZZ/2$ factor and gives the desired identification. Similar considerations give the identifications in part {\bf (b)}.
\endofproof

Note that in the proof we could alternatively use the differential Bott sequence from Proposition \ref{prop-diffBott}.

\section{Applications}
\label{Ch-apps}

\subsection{Higher structures}
\label{Sec-higher}

In this section we will consider sample applications of results in previous sections, 
mainly Section \ref{Sec-AHSS-KOhat} in Section \ref{Ch-AHSS}, to further congruences
involving Pontrjagin classes. The latter appear essentially as generators of the corresponding 
cohomology rings of higher tangential structures (see \cite{SSS1}\cite{SSS2}\cite{SSS3}\cite{9}). 

\medskip
First we will consider some classical constructions.
Let $f: X \to Y$ be a continuous map between compact 
differentiable manifolds, such that $\dim(Y) - \dim (X)=0$ mod 8. Suppose that 
$\nu_f=f^*(TY) - TX$ is provided with a stable Spin structure, in the sense of \cite{AH-diff}. 
Then we have the relation \cite[Theorem V.4.19]{Ka}:
\(
\label{RR-KO}
{\rm Ph}(f_*^{KO}(x))=f_*^H(A(\nu_f))\cdot {\rm Ph}(x))\;,
\)
where $f_*^{KO}$ is the Gysin homomorphism in KO-theory. 
We recall the following integrality results depending on the dimension mod 4. 
The $\widehat{A}$-genus of a compact differentiable Spin manifold is an integer, by
the classical Borel-Hirzebruch theorem \cite{BH}
(see Corollary to \cite[Theorem 26.3.1]{Hir}). For $X$ be an oriented differentiable manifold
with $\dim X\equiv 4$ mod 8 and $w_2(X)=0$, the $\widehat{A}$-genus $A(TX)$ is an 
even integer  (see \cite[Theorem 26.3.2]{Hir}).

\medskip
In fact, the latter can be refined as follows \cite[Theorem V.4.24]{Ka}:
For each element $x$ of $\op{KO}(X)$ the value of $\widehat{A}(TX) {\rm Ph}(x)$ on the fundamental 
class of $X$ is an even integer. This is obtained by applying the relation \eqref{RR-KO} to 
the constant map from $X$ to $S^4$. Then the Pontrjagin character 
${\rm Ph}: \op{KO}(S^4) \to H^*(S^4; \QQ)$ isomorphically maps $\op{KO}(S^4)$ 
onto $2H^4(S^4; \ZZ)=2\ZZ \subset \QQ=H^4(S^4; \QQ)$, since
$\ch$ isomorphically maps $\op{K}(S^4)$ onto $H^4(S^4; \ZZ)$ (the Chern character is 
integral on spheres in general) and since the complexification homomorphism 
$\op{KO}(S^4) \cong \ZZ \to K(S^4)\cong \ZZ$ is multiplication by 2
(which can be deduced from the Bott sequence \ref{Bott-seq}).
On the other hand, $f_*^H$ may be factored into $b_*^H a_*^H$, where
$a: X \to {\rm pt}$ and $b: {\rm pt} \to S^4$. Since $b_*$ is an isomorphism 
from $H^0({\rm pt}; \ZZ)$ to $H^4(S^4; \ZZ)$, the expression $f_*^H({\rm Ph}(x) \widehat{A}(\nu_f))$
is actually the value of ${\rm Ph}(x) \widehat{A}(TX)$ on the fundamental class of $X$ when we 
identify $H^4(S^4; \ZZ)$ with $\ZZ$. Note that this has been used in 
the above works in characterizing the obstruction classes associated with the higher 
tangential structures.

\begin{remark}[Non-Spin and inverting $2$]\label{rem-invert}
The Spin condition may be avoided in the above 
integrality results by considering the theory $\op{KO}(X) \otimes \ZZ[\tfrac{1}{2}]$. Let $E$ be 
an oriented vector bundle of rank $n$ over a compact space $X$. Then 
there is an isomorphism $\op{KO}(X) \otimes \ZZ[\tfrac{1}{2}]\cong \op{KO}^n(X^E) \otimes \ZZ[\tfrac{1}{2}]$
as $\op{KO}(X)$-modules (stated in \cite[IV.8.13 \& Remark V.4.26]{Ka} for the complex case).
Up to a power of 2, the $\hat{A}$-genus of {\it any} oriented manifold, evaluated on the 
fundamental class, is an integer.  For instance, from $\hat{A}_1=\tfrac{1}{2^2\cdot 3}p_1$, 
one has that $p_1$ is divisibly by 3 for a 4-manifold,
which is compatible with the L-genus being integral for oriented manifolds, as 
$L_1=\tfrac{1}{3}p_1$. 
\end{remark}


We now illustrate the utility of Proposition \ref{prop-deg8ph} by some examples. Some of these examples
are known; nevertheless the following provide alternative proofs for these known facts, purely based on the 
spectral sequence. Recall that in the topological case, the differential $d_5$ is given by the operation 
$\beta_2Sq^4\rho_2$ (see the last expression in \eqref{d2d5-class}). 
In the differential case, the argument used in the proof of Proposition
\ref{prop-oddprim} applies (i.e. using the boundary homorphism). This gives us the following identification. 

\begin{lemma}[Identifying the differential $d_5$]\label{lemm-d5}
The differential $d^{s,-4t+1}_5$ in the spectral sequence for $\widehat{\op{KO}}$ is given by 
$$
d^{s,-4t+1}_5=j_2Sq^4 \rho_2 \beta_{U(1)}:H^{s}(M;U(1))\longrightarrow H^{s+5}(M;U(1))\;.
$$
\end{lemma}

We now make use of this in the following examples. 

\begin{example}[Integrality of the Fivebrane obstruction]\label{ex-5b}
The obstruction to having a Fivebrane structure is measured by the integral class $\tfrac{1}{6}p_2$ \cite{SSS2}\cite{SSS3}. 
In fact, one can deduce the integrality of $\tfrac{1}{6}p_2$ purely from the spectral sequence. Indeed, if $M$ admits 
Spin structure, it follows from the Wu formula that the degree four Wu class is given by $\nu_4=w_4$. Moreover, 
since the third integral Stiefel-Whitney class vanishes $W_3=0$, the mod 2 reduction of $\frac{1}{2}p_1$ coincides 
with $w_4$ (see \cite{memb}\cite{Wuc} for detailed  discussions). Since $\frac{1}{2}p_1=0$, it follows that $d_5=jSq^4\rho_2\beta_{U(1)}=0$ on $H^3(M;U(1))$. Thus, $d_8$ vanishes if and only if 
$$
{\rm Ph}_8=\tfrac{1}{12}\big(p_1^2-2p_2)=\tfrac{1}{6}p_2
$$
is integral.
\end{example}

Note that the above is related to the class $I_8:=\tfrac{1}{48}\big(p_2- Q_1^2 \big)$ highlighted in \cite{KSpin},
which can be written in terms of $Q_1$ and $Q_2$, the first and second Spin characteristic classes defined in 
\cite{Th3}; see also Section \ref{Sec-Classicalw} later. However, the class $I_8$ is not necessarily integral;
for example, $\HH P^2$ admits 
Spin structure and its total Pontrjagin class is given by
$$
p(\mathbb{H}P^2)=\frac{(1+x)^{6}}{1+4x}=1 + 2 x + 7 x^2 \;.
$$
So $p_1=2x$, $p_2=7x^2$ and $\tfrac{1}{48}\big(p_2-Q_1^2)=\frac{7}{48}x^2-\frac{1}{48}x^2=\frac{1}{8}x^2$ is not integral because $x^2$ generates degree $8$.

\begin{example}[Integrality of the Ninebrane obstruction]
For a Fivebrane manifold $M$, the class $\tfrac{1}{6}p_2$ vanishes. The next rational obstruction is associated with 
a Ninebrane structure in the Whitehead tower for $\op{BO}$ and is $\tfrac{1}{240}p_3$ \cite{9}. Let us see what integrality condition we can 
achieve via the spectral sequence. We observe that  the relevant differentials in this case are $d_5$, $d_6$, $d_7$, and $d_9$. 
Since $M$ has Spin structure, $d_5$ vanishes as in Example \ref{ex-5b}.  The differentials $d_6$ and $d_7$ vanish 
by degree considerations. 
The $p$-primary part of $d_9:E^{3,-3}\to E^{12,-11}$ is given in Proposition \ref{prop-oddprim}. 
Setting $8=4r(p-1)$, we see that for $p=2$, $r=2$ 
and for $p=3$, $r=1$. But both $Sq^8$ and $P^2_3$ vanish by degree considerations. The differential 
$d_{12}:E^{0,0}_{12}\to E^{-11,12}_{12}$ evaluated on a differential form representing ${\rm Ph}_{12}$ 
necessarily vanishes, as this class lifts to $\widehat{\op{KO}}$ (See Proposition \ref{prop-E00SS}). 
It then follows from Theorem \ref{thm-d4d8} that $\frac{1}{2}{\rm Ph}_{12}$ represents an integral 
class. Using $\tfrac{1}{2}p_1=0$ and $\tfrac{1}{6}p_2=0$, we have 
$$
\tfrac{1}{2}{\rm Ph}_{12}=\tfrac{1}{2^4\cdot 3^2\cdot 5}(p^3_1-3p_2p_1+3p_3)=
\tfrac{1}{2^4\cdot 3\cdot 5}p_3=\tfrac{1}{240}p_3
$$
is integral. 
\end{example}

The next example illustrates how the spectral sequence can be utilized to find differential forms which do 
not lift to $\widehat{\op{KO}}$-theory. We will concentrate on the use of odd primes. Even for  the prime $p=3$,
one needs sufficiently large degrees. 

\begin{example}[Forms not lifting to $\widehat{\rm KO}$-theory] 
Let $\omega_{16}$ be a differential form of degree $16$ with integral periods on a manifold $M$. 
Then, for any positive integer $m>1$, the differential form $v_{16}=\frac{1}{3m}\omega_{16}$ is 
not in the degree $16$-component of a Pontrjagin character form. Indeed, this immediately follows 
from the spectral sequence by applying Theorem \ref{thm-d4d8}. Since
$d_{16}(v_{16})=[v_{16}]\mod \ZZ$ is not 3-torsion, it does not lie in the image of 
$d_9:=jP^2_3\beta_{U(1)}:E^{7,-7}_{9}\to E_9^{16,-15}\subset H^{16}(M;U(1))/{\rm im}(d_8)$, which is $3$-torsion. 
\end{example}


The following example demonstrates an interplay between integral and mod 2 classes in degree four. 

\begin{example}[Quadratic forms, String structures and $d_5$]
Let $M$ be a smooth manifold for which the 5th integral Stiefel Whitney class $W_5=0$, but $w_4\neq 0$. 
Since $0=W_5=\beta_2(w_4)$, it follows that $w_4$ is in the image of the mod 2 reduction $\rho_2$. 
Let $\lambda$ denote such an integral lift and suppose that $\lambda\in H^4(M;\ZZ)$ has order at least $4$ 
and that we can choose $\lambda\in H^4(M;\ZZ)_2$ with $2\lambda\neq 0$. Now $\lambda \neq 2y$
for any class $y$, as otherwise its mod 2 reduction would vanish. Moreover, since it is torsion, 
we have $\lambda=\beta_{U(1)}\tilde{\lambda}$ for some $\tilde{\lambda}\in H^3(M;U(1))$. Such a 
$\tilde{\lambda}$ is  an instance of an \emph{integral Wu structure} in the sense of \cite{HS}. 
Then if $Sq^4(w_4)=w_4^2\neq 0$,
it follows from Lemma \ref{lemm-d5} that $d_5:E_4^{3,-3}\to E^{8,-7}_4$ is nonzero. Hence 
for a differential form $\omega_{8}$ to lift to the degree $8$ component of the Pontrjagin character
form, it is both necessary and sufficient that $[\omega_8]$ mod $\ZZ=0, \frac{1}{2}$. Equivalently,
for each cycle $c:\Delta^n\to M$, this is the requirement  
\(\label{anlineb}
{\rm exp}\Big(2\pi i\oint_c \omega_8\Big) =\pm 1\;.
\)
Now fix any class $x\in H^3(M;U(1))$. The combination $\lambda+2x$ satisfies 
$\rho_2\beta_{U(1)}(\tilde{\lambda}+2x)=\rho_2\beta_{U(1)}\tilde{\lambda}$ and hence
$$
d_5(\tilde{\lambda}+2x)=d_5(\tilde{\lambda})=Sq^4(w_4)=w_4^2\neq 0\;.
$$
If we restrict to the classes $x$ that lie in the intermediate Jacobian $H^3(M;\RR)/H^3(M;\ZZ)=T^{b_3}$ and 
we choose differential forms $(\omega_8)_x$ which smoothly depend on $x\in T^{b_3}$, the exponentials \eqref{anlineb} will give rise to a section of a real line bundle over $T^{b_3}$. Note that this enters into  
the description of the anomaly line bundle for the fivebrane partition function
(see \cite{Wi-5}\cite{HS}\cite{memb}\cite{Wuc}) but we will elaborate further on this aspect elsewhere. 
We will consider lifts of Stiefel-Whitney classes in Section \ref{Sec-Classicalw}.
\end{example}

%
%

One can come up with many such results, but we will content ourselves with the above as representative samples.

\newpage
\subsection{Adams operations in differential KO-theory}
\label{Sec-Adams}

We would like to provide refinements of the Adams operations from KO-theory to 
differential KO-theory  as  natural transformations of differential K-theory. In fact, we do so for both 
$\widehat{\rm KO}$ and for ${\bf K}{\rm O}_\nabla$, the former being analogous to the 
complex case  in \cite{Bu-Ad} (see also \cite{LiuM} from a geometric perspective), 
while the latter has not been considered before even in the complex case. 
The differential Adams operations have the expected geometric meaning. If a differential K-theory 
class is, e.g., given by a line bundle with connection, then $\widehat{\psi}^k$ assigns to it the class 
represented by the $k$th tensor power. In fact, in \cite{Bu-Ad} only the complex K-theory 
operations are considered. So before addressing the differential real case, 
we first recall the corresponding operations in the topological real case, 
i.e., for KO-theory.

\medskip
\noindent {\bf Adams operations in KO-theory.} 
These are discussed in detail, for instance, in \cite[Lecture 10, 11]{Mill} and with a view towards
applications in \cite{MPS}. 
For $E$ and $E'$ real vector bundles, there is a natural isomorphism 
$$
\Exterior ^k (E \oplus E')= \bigoplus_{i+j=k} \Exterior^i E \oplus \Exterior^j E'
$$
so that one defines $\Lambda_t(E)=\sum_{i\geq 0} t^i \Exterior^i(E)$. Consequently,  
the exponential law 
$
\Lambda_t (E \oplus E')=  \Lambda_t E \cdot \Lambda_t E'
$
holds and $\Lambda_t(E)$ induces a homomorphism of commutative monoids
$\Lambda_t(E): {\rm Vect}(X) \to 1 + t \op{KO}(X)[\![t]\!]$, where the target 
is the group of formal power series in $t$ with coefficients in $\op{KO}(X)$ and
constant term one. There is a group homomorphism 
$\lambda_t: \op{KO}(X) \to 1 + t \op{KO}(X)[\![t]\!]$
such that the following diagram commutes
$$
\xymatrix@R=1.6em{
{\rm Vect}(X)
\ar[rr]^-{\Lambda_t} \ar[d] 
&&
1 + t \op{KO}(X)[\![t]\!]\;.
\\
\op{KO}(X) \ar[urr]_-{\lambda_t}
}
$$
\begin{remark} [Generalities on Adams operations for KO-theory] 
We highlight the following facts:

\vspace{-2mm}
\begin{enumerate}[{\bf (i)}]
\item Unlike the complex case, however, the Adams operations for KO-theory are {\it not}
uniquely characterized by what happens for line bundles 
$\psi^k(L)=L^{\otimes k}$ and additivity $\psi^{kl}=\psi^k \psi^l$. 

\vspace{-2mm} 
\item Nevertheless, 
the operations can be defined as follows.  One observes the triviality of $L^{\otimes 2}$:
For a real line bundle $L$ we have $L^*\cong L$ so that $L^{\otimes 2} \cong L \otimes L^* 
\cong {\rm Hom}(L, L) \ni 1$, which is a section of the bundle so that $L^{\otimes 2}$ is trivial . 
So one can define ${\psi^k(E)}$ to be $E$ for $k$ odd and to be $\psi^0(E)$ when $k$ is even,
where $\psi^0(E)$ is the trivial bundle of dimension equal to the dimension of $E$ over the 
basepoint of $X$. 

\vspace{-2mm} 
\item Adams operations in KO-theory can be deduced from those of the complex theory 
(see \cite[Lecture 12]{Mill}). For $G$ a compact Lie group, 
consider the ring of virtual real representations  $RO(G)$. Complexification of real representations
gives a map $c: RO(G) \to R(G)$  to the ring of virtual complex representations  $R(G)$. 
Moreover, the inclusion $i: O(n) \into U(n)$ induces 
$i^*: R(U(n)) \to R(O(n))$. The real Adams operations come
from the following additive $RO$-sequences
$\{ \psi^k_{\RR, n} \in RO(O(n)), n \geq 0\}$:
\vspace{-2mm} 
$$
\xymatrix@R=.2em@C=-1pt{
&&&& RO(O(n)) 
\ar@{^{(}->}[ddd]^c 
& 
\ni  \psi^k_{\RR, n}:=s_k(\Lambda^1, \dots, \Lambda^n)
\ar@{|->}[ddddl]
\\
\\
\\
R(U(n)) \ar[rrrr]^{i^*} &&&& 
R(O(n)) 
\\
\psi^k_n \ar@{|->}[rrrr]&&&& \;\;\;  i^*\psi^k_n  \phantom{AA}
}
$$

\vspace{-2mm} 
\noindent As the generators of $R(O(n))$ are complexifications of some 
real representations and the complexification $c:RO({\rm O}(n)) \to 
R({\rm O}(n))$  is injective,  the real representation ring is 
isomorphic to the complex representation ring \cite{Min}. 

\vspace{-3mm} 
\item Adams operations on generators: The action of $\psi^k$ 
on $\widetilde{KO}^*(S^0)$ is given by (see e.g. \cite{Wat1}) 
$$
\psi^k(\eta)=k \eta\;, \qquad
\psi^k(\alpha)=k^2 \alpha\;, \qquad
\psi^k(\beta)=k^4 \beta\;.
$$
\item An explicit recursive formula is given by (see \cite{Rog}) 
\vspace{-2mm} 
$$
\psi^r(E)=\sum_{i=1}^{r-1} (-1)^i \Lambda^i(E) \otimes \psi^{r-i}(E) + (-1)^r \Lambda^r(E)\;.
$$
\end{enumerate}
\end{remark}

The Adams operations can be refined to the differential setting in a natural way, namely, by incorporating 
connections. Recall that if $E$ is a vector bundle with connection $\nabla$, then the $k$th exterior 
power bundle $\bigwedge^k(E)$ admits a connection $\nabla_k$ which operates on a smooth section 
$\sigma_1\wedge \sigma_2\wedge \hdots \wedge \sigma_k$ of $\bigwedge^k(E)$ by 
\(\label{def-conn-exterior}
\nabla_k(\sigma_1\wedge \sigma_2\wedge \hdots \wedge \sigma_k)=
\sum_i(-1)^{i-1}\sigma_1\wedge \sigma_2\wedge \hdots \nabla(\sigma_i) \hdots  \wedge \sigma_k\;.
\)
\begin{definition}[Exterior class in in ${\mathbf{K}{\rm O}}_{\nabla}(M)$]
We can form the class $\Lambda^{\nabla}_t(E,\nabla)$ in ${\mathbf{K}{\rm O}}_{\nabla}(M)$ by setting
$$
\Lambda^{\nabla}_t(E,\nabla):=\sum_{i\geq 0}t^i[\Exterior^i(E),\nabla_i]\;.
$$
\end{definition} 
\noindent This class has the following properties:
\begin{enumerate}[{\bf (i)}]
\item $\Lambda^{\nabla}_t(E,\nabla)$ recovers the topological class $\Lambda_t(E)$ under the forgetful map $\mathscr{F}:\mathbf{K}{\rm O}_{\nabla}(M)\to \op{KO}(M)$. 
\item The map $\Lambda^{\nabla}_t$ again induces a homomorphism of commutative monoids. 
\item  By the adjunction ${\rm Gr}\dashv i$ (see Appendix), this again induces a natural transformation $\lambda^{\nabla}_t$, filling the diagram 
\(
\label{triang-konab}
\xymatrix@R=1.6em{
{\rm Iso}({\rm Vect}_{\nabla}(M))
\ar[rr]^-{\Lambda^{\nabla}_t} \ar[d]_-{\eta} 
&&
1 + t\cdot \mathbf{K}{\rm O}_{\nabla}(M)[\![t]\!]
\\
\mathbf{K}{\rm O}_{\nabla}(M) \ar[urr]_{\lambda^{\nabla}_t}
}\;,
\)
where $\eta$ is the unit of the adjunction $\eta:\mathbb{1}\to i{\rm Gr}$.
\end{enumerate} 
With this setup we are now ready to provide the following. 

\begin{definition}[Adams operations in $\mathbf{K}{\rm O}_{\nabla}$]\label{adam-konab}
We define the Adams operation in $\mathbf{K}{\rm O}_{\nabla}$ by the recursive formula 
$$\psi^r(E,\nabla)=\sum_{i=1}^{r-1} (-1)^i [\Lambda^i(E),\nabla_i] \otimes \psi^{r-i}(E,\nabla) + (-1)^r [\Lambda^r(E),\nabla_r].
$$
Alternatively, using the formal power series $\lambda^{\nabla}_t$, we can define $\psi^r(E,\nabla)$ as the coefficients of the formal power series
$$
\sum_{r=0}^{\infty}\psi^r(E,\nabla)t^r:=\psi^0(E,\nabla)-t\frac{d}{dt}\log\big(\lambda^{\nabla}_{-t}(E,\nabla)\big)
$$
where we set $\psi^0(E,\nabla):=({\bf rk}(E),d)$ and $\log$ is interpreted via its formal power series expansion.
\end{definition}
These operations are well-defined classes only for the degree zero cohomology group 
$\mathbf{K}{\rm O}^0_{\nabla}(M)$ (which is also usually referred to simply as 
$\mathbf{K}{\rm O}_{\nabla}(M)$).  In higher degrees, one needs to take more care. 
Note that in the topological case, the Adams operations do not commute with the Bott periodicity 
isomorphisms. To extend the operation to the graded groups
$\op{KO}^n(X)=\widetilde{\op{KO}}(\Sigma^mX)$, where $n+m=8k$, one observes 
that the complex Adams operation maps the generator $u=(1-H)\in \widetilde{\op{K}}^2({\rm pt})$ to 
$$
\psi^r(u)=\psi^r(1-H)=1-H^r=1-(1-u)^r=1-(1-ru)=ru\;.
$$
Hence, we must localize by inverting $r$ and define $\psi^r$ on $\op{K}[\tfrac{1}{r}]$ by $(1/r^k)\psi^r$. Similar considerations hold for $\op{KO}$, allowing us to extend via Bott periodicity. We make the following definition.

\begin{definition}[Adams operations in localized KO-theory]
We define the Adams operation in localized $\op{KO}$-theory $\op{KO}^n(X)[\tfrac{1}{r}]$ by $(1/r^k)\psi^r$. These operations then induce ring spectrum maps 
$\psi^r: \op{KO}[\tfrac{1}{r}] \to \op{KO}[\tfrac{1}{r}]$ (see e.g. \cite{Rog}) such that the diagrams
\vspace{-3mm} 
\(\label{tpbtkoad}
\xymatrix{
\op{KO}^\ast(M)[\tfrac{1}{r}] \ar[d]^-{\beta}\ar[rr]^-{r^{-k}{\psi}^r} && 
\op{KO}^\ast(M)[\tfrac{1}{r}]\ar[d]^-{\beta}
\\
\op{KO}^{\ast+8k}(M)[\tfrac{1}{r}] \ar[rr]^-{{\psi}^r} &&
\op{KO}^{\ast+8k}(M)[\tfrac{1}{r}]
}
\)

\vspace{-1mm} 
\noindent commute, where $\beta$ is the Bott periodicity map for $\op{KO}$-theory, 
and satisfy the compatibility $\psi^r\circ \psi^s=\psi^{rs}$ (at the cohomological level). 
\end{definition}

\newpage

Unfortunately, in the case of $\mathbf{K}{\rm O}_{\nabla}$, one does not have any form of Bott periodicity 
available. Hence, extending the operations introduced in Definition \ref{adam-konab} in a similar way 
to $\mathbf{K}{\rm O}^n_{\nabla}(M)$ in not possible for $n>0$. For this reason, we will develop the
extended Adams operations in the Hopkins-Singer type theory 
$\widehat{\op{KO}}$, paralleling the construction in \cite{Bu-Ad} in the complex case. We then reconnect 
these operations with the differential Adams operations from Definition \ref{adam-konab} via the cycle map.

\medskip
We now consider Adams operations for differential forms.
These are related to cohomology operations in 4-periodic rational cohomology. 
More precisely, recalling that $4\beta=\a^2$ from the coefficients \eqref{eq-coeff},
we have the following.
\begin{definition}[Adams operations on 4-periodic rational cohomology]
We define natural operations 
$$
\psi^r_{H}:H^*(M;\mathbb{Q}[\alpha,\alpha^{-1}])\longrightarrow
H^*(M;\mathbb{Q}[\alpha,\alpha^{-1}])\;,
$$
such that $\psi^r_{H}(x)=x$ for $x\in H^*(M;\QQ)$ and $\psi^r_{H}(\alpha)=\tfrac{1}{r^2}\alpha$. 
\end{definition}

Classically one has the compatibility with the Pontrjagin character
\(
\label{class-pont-comp}
{\rm Ph}\circ \psi^r=\psi^r_{H}\circ {\rm Ph}\;.
\)
Indeed, the two Adams operations can be related via rationalization. 
To see this, let $M_p \subset \op{KO}(X) \otimes \QQ$ be the subspace on which 
$\psi^k$ acts by multiplication with $k^p$. Then there is a direct sum 
decomposition \cite[Theorem 16.1]{Bott}  
$\op{KO}^*(X) \otimes \QQ = \sum_{p=0}^\infty M_p$. Furthermore, 
$M_p=0$ if $p$ is odd. 

\medskip
The operations on 4-periodic rational cohomology carry over to natural operations 
on 4-periodic forms in a natural way. 
\begin{definition}[Adams operations on forms]
We define the natural operations 
$$
\psi^r_{\rm form}:\Omega^*(M;\mathbb{R}[\alpha,\alpha^{-1}])\longrightarrow
\Omega^*(M;\mathbb{R}[\alpha,\alpha^{-1}])
$$
via $\psi^r_{\rm form}(\omega)=\omega$ with $\omega\in \Omega^*(M)$ 
and $\psi^r_{\rm form}(\alpha)=\tfrac{1}{r^2}\alpha$. 
\end{definition}

We now refine the Adams operations to differential $\op{KO}$-theory $\widehat{\rm KO}$. 
In the complex case, the existence 
and uniqueness of these operations was established  in \cite[Theorem 3.1]{Bu-Ad}. With minor 
modifications, this construction carries over to the real case. We now show how to make these modifications.

\begin{theorem}[Adams operations in differential (localized) KO-theory]\label{adamdfko}
{\bf (i)} There exists a unique natural transformation 
$$
\widehat{\psi}^r:\widehat{\op{KO}}^*(M)[\tfrac{1}{r}]\; \longrightarrow \;
\widehat{\op{KO}}^*(M)[\tfrac{1}{r}]
$$
of set-valued functors on the category of compact manifolds such that 
\(\label{admrfn}
\mathcal{I}\circ \widehat{\psi}^r=\psi^r\circ \mathcal{I}
\; \qquad \text{and} \qquad 
\mathcal{R}\circ \widehat{\psi}^r=\psi_{\rm form}^r\circ \mathcal{R}\;.
\)
\item {\bf (ii)} For $-6\leq i \leq 0$, the diagrams 
\(\label{addgko}
\xymatrix{
\widehat{\op{KO}}^i(S^1\times M)\ar[rr]^-{\widehat{\psi}^r} \ar[d]_-{\int} && \widehat{\op{KO}}^i(S^1\times M)\ar[d]^-{\int}
\\
\widehat{\op{KO}}^{i-1}(M)\ar[rr]^-{\widehat{\psi}^r} && \widehat{\op{KO}}^{i-1}(M)
}
\)
commute, $\widehat{\psi}^r$ preserves the ring structure and satisfies $\widehat{\psi}^r\circ \widehat{\psi}^s=\widehat{\psi}^{rs}$.
\end{theorem}
\theproof
We will break the proof up into stages as in  the complex case \cite{BS} as follows:

\noindent Claim 1. In degree zero, there is a unique natural transformation of set-valued functors
\(
\label{rfnadmcl} \hat{\psi}^r:\widehat{\op{KO}}(M)\longrightarrow \widehat{\op{KO}}(M)
\)
satisfying \eqref{admrfn}. 
\\
\noindent Claim 2. The unique transformations \eqref{rfnadmcl} preserves the ring structure.
\\
\noindent Claim 3. Going up in degrees, there are unique natural transformations 
\(
\label{rfnadmcl2} 
\hat{\psi}^r:\widehat{\op{KO}}^i(M)\longrightarrow \widehat{\op{KO}}^i(M)\;,
\)
with $-6\leq i < 0$, satisfying \eqref{admrfn}, making the diagrams \eqref{addgko} commute and which preserve the ring structure.
The proof then concludes with an extension of the operations via Bott periodicity.
To define $\widehat{\psi}^r$ in all degrees, we note that the compatibility \eqref{admrfn} and diagram \eqref{tpbtkoad} forces us to extend via Bott periodicity as
$$
\xymatrix{
\widehat{\op{KO}}^*(M)[\tfrac{1}{r}]\ar[rr]^-{r^{-k}\widehat{\psi}^r}\ar[d] &&
\widehat{\op{KO}}^*(M)[\tfrac{1}{r}]\ar[d]
\\
\widehat{\op{KO}}^{*+8k}(M)\ar[rr]^-{\widehat{\psi}^r} &&
\widehat{\op{KO}}^{*+8k}(M)
\;.}
$$
%
%
%
%
\noindent {\it Proof of Claim 1.}
The space $\op{KO}_0$ has the homotopy type of $\mathbb{Z}\times \op{BO}$.  \cite[Proposition 2.3]{BS}, there is a sequence of pointed smooth manifold $N_i$, together with maps $n_i:N_i\to N_{i+1}$ and $x_i:N_i\to \mathbb{Z}\times \op{BO}$ satisfying the following:

\vspace{-2mm}
\begin{enumerate}
\item $N_i$ is homotopy equivalent to an $i$-dimensional CW-complex.
\item The map $n_i$ is $(i-1)$-connected.
\item $n_i:N_i\into N_{i+1}$ is an embedding of submanifolds.
\item The diagrams
$$
\xymatrix@R=1.3em{
N_i\ar[rr]^-{n_i} \ar[dr]_-{x_i} && N_{i+1}\ar[dl]^-{x_{i+1}}
\\
& \op{BO}\times \ZZ&
}
$$
commute
\item $x_i$ is $i$-connected.
\end{enumerate}

\vspace{-2mm}
\noindent Let $u:\mathbb{Z}\times \op{BO}\to \mathbb{Z}\times \op{BO}$ be the identity map representing the universal class $\op{KO}(\mathbb{Z} \times \op{BO})$. By \cite[Proposition 2.6]{BS}, we can choose a sequence $\hat{u}_i\in \widehat{\op{KO}}^0(N_i)$ such that $\mathcal{I}(\hat{u}_i)=x_i^*u$ and $n_i^*\hat{u}_{i+1}=\hat{u}_i$ for all $i\geq 0$. By \cite[Lemma 3.8]{BS}, we also have for $r>2i+2$, $H^{2j+1}(N_r;\RR)=0$ for all $j\leq i$. Then the requirements \eqref{admrfn} uniquely determine $\hat{\psi}^r(\hat{u}_i)$ up to elements in the image 
\(\label{ambchui}
a:\bigoplus_{2j+1\geq i}H^{2j+1}(N_i;\RR)\subset \Omega^*(N_i;\pi_*(\op{KO}))/{\rm im}(d) 
\longrightarrow
\widehat{\op{KO}}(N_i)\;.
\)
Fix $\hat{y}\in \widehat{\op{KO}}(M)$ and let $f:M\to N_i$ be such that $\mathcal{I}(\hat{y})=f^*x_i^*(u)$. Choose a form $\rho$ such that 
$\hat{y}=f^*\hat{u}_i+a(\rho)$. Any two such choices necessarily necessarily differ by an element in the image of the Pontrjagin character. By the compatibility
$$
\psi^r_{H}({\rm Ph}(x))={\rm Ph}(\psi^r(x))\;,
$$
it follows at once that $a(\rho)$ does not depend on the choice of $\rho$. Then define
\(\label{dfdfadop}\widehat{\psi}^r(\hat{y})=f^*\widehat{\psi}^r(\hat{u}_i)+a(\psi^r_{\rm form}\rho)\;.\)
Since the ambiguity in uniquely defining $\widehat{\psi}^r(\hat{u}_i)$ is measured by the image of  \eqref{ambchui}, it follows that $f^*\widehat{\psi}^r(\hat{u}_i)$ does not depend on the choice of refinement of $\psi^r(u_i)$. The assignment also does not depend on the choice of map $f:M\to N_i$, which can be shown using the homotopy formula as in \cite[Lemma 3.2]{BS}. It remains to show that the well-defined map 
$$
\widehat{\psi}^r:\widehat{\op{KO}}(M)\longrightarrow
\widehat{\op{KO}}(M)
$$
is natural in $M$ and is unique. With the replacements made above, this follows verbatim as in \cite[Lemma 3.3]{BS}.


\medskip
\noindent {\it Proof of Claim 2.}
For complex K-theory, this is \cite[Lemma 3.4]{Bu-Ad}. The proof there uses only the
formula \eqref{dfdfadop} and the method of proof in \cite[Theorem 3.6]{BS}, which holds 
for $E=\op{KO}$. 

\medskip
\noindent {\it Proof of Claim 3.}
Choose $\hat{e}\in \widehat{\rm KO}^{-1}(S^1)$ such that $\int \hat{e}=1$ and $R(\hat{e})$ is the unique 
normalized rotation invariant form. For $\hat{x}\in \widehat{\op{KO}}^{-1}(M)$, the commutative diagram \eqref{addgko} forces us to define
$
\widehat{\psi}^r(\hat{x}):=\int \widehat{\psi}^r(\hat{e}\times \hat{x})
$
on $\widehat{\op{KO}}^{-1}(M)$. Similarly, we define
\footnote{Note that one could alternatively formulate this in the reduced case using the differential 
desuspension map we introduced in Proposition \ref{desprdf}, but the above suffices for the current purpose.}
$$
\widehat{\psi}^r(\hat{x}):=\int\int \widehat{\psi}^r(\hat{e}\times \hat{e} \times \hat{x})
$$ 
on $\widehat{\op{KO}}^{-2}(M)$ and so on. The choice of $\hat{e}$ is unique up to elements of the form 
$a(c)$ with $c\in \oplus_k H^{4k}(S^1;\RR)\cong H^0(S^1; \RR)\cong \RR$. 
It follows from the commutativity of diagram \eqref{addgko} that the 
definitions do not depend on the choice of $\hat{e}$. Moreover, these transformations satisfy relations 
\eqref{admrfn}, $\widehat{\psi}^r\circ \widehat{\psi}^s=\widehat{\psi}^{rs}$, and are multiplicative 
by direct computation (as in the complex case). 
\endofproof

We now conclude this section by reconnecting with the geometrically defined Adams operations in
Definition \ref{adam-konab}. 
\begin{proposition}[Compatibility of differential Adams operations]
We have a commutative diagram (in degree zero)
$$
\xymatrix{
\mathbf{K}{\rm O}_{\nabla}(M)\ar[rr]^-{\widehat{\rm cyc}}\ar[d]_-{\psi^k_{\nabla}} && 
\widehat{\op{KO}}(M)\ar[d]^-{\hat{\psi}^k}
\\
\mathbf{K}{\rm O}_{\nabla}(M) \ar[rr]^-{\widehat{{\rm cyc}} }&&  \widehat{\op{KO}}(M)
\;.}
$$
\end{proposition}
\theproof
Recall from Definition \ref{def-cyc-dif} and Proposition \ref{dfcycmp} 
that the cycle map $\widehat{{\rm cyc}}$ is the unique transformation 
satisfying $\mathcal{R}(\widehat{{\rm cyc}}(V,\nabla))={\rm Ph}(\mathcal{F}_{\nabla})$ and $\mathcal{I}(\widehat{{\rm cyc}}(V,\nabla))=[V]$, where $[V]$ is the isomorphism class of the underlying vector bundle (regarded as a class in $\op{KO}$). Moreover, by Theorem \ref{adamdfko}, the operations $\hat{\psi}^k$ are the unique refinements such that 
$\mathcal{I}\circ \hat{\psi}^k=\psi^k\circ \mathcal{I}$ and $\mathcal{R}\circ \hat{\psi}^k=\psi_{\rm form}^k\circ \mathcal{R}$. The same argument as in the proof there shows that the composition $\hat{\psi}^k\circ \widehat{{\rm cyc}}$ is the unique natural transformation characterized by the two conditions:
\begin{itemize}
\item [{\bf (i)}]  \emph{(Compatibility with topological realization)}
$\mathcal{I}\circ \hat{\psi}^k\circ \widehat{{\rm cyc}}(V,\nabla)=\psi^k([V])$.
\item [{\bf (ii)}] \emph{(Compatibility with curvature and forms)} $\mathcal{R}\circ \hat{\psi}^k\circ \widehat{{\rm cyc}}(V,\nabla)=\psi^k_{\rm form}({\rm Ph}(\mathcal{F}_{\nabla}))$.
\end{itemize}
By commutativity of  diagram \eqref{triang-konab} and noting the expression \eqref{def-conn-exterior}
for $\nabla_k$, we have 
$$
\mathcal{R}(\widehat{{\rm cyc}}\circ \psi^k_{\nabla}[V,\nabla])=\mathcal{R}\big(\widehat{{\rm cyc}}\big[ \Exterior^k(V),\nabla_k\big]\big)={\rm Ph}( \Exterior^k(V),\nabla_k)=\psi^k_{\rm form}{\rm Ph}(\mathcal{F}_{\nabla})\;,
$$
where we have also used multiplicativity of Ph  with respect to tensor product (see Remark \ref{RemPropPont}) and compatibility
on forms in expression \eqref{class-pont-comp} above. 
Similarly, we also have 
$$\mathcal{I}(\widehat{{\rm cyc}}\circ \psi^k_{\nabla}([V,\nabla])=\mathcal{I}(\widehat{{\rm cyc}}\big[ \Exterior^k(V),\nabla_k\big]=[\Exterior^k(V)]=\psi^k([V])\;.$$
Therefore, the transformation $\widehat{{\rm cyc}}\circ \psi^k_{\nabla}$ satisfies both 
properties ${\bf (i)}$ and ${\bf (ii)}$ and  we establish that indeed
$\widehat{{\rm cyc}}\circ \psi^k_{\nabla}=\hat{\psi}^k\circ \widehat{{\rm cyc}}$.
\endofproof

\subsection{A differential Wu formula} 
\label{Sec-Classicalw}

\subsubsection{Classical Stiefel-Whitney classes and their integral lifts} 

The differentials in the AHSS, as we indicated
in Section \ref{Sec-AHSS-KObare} and Section \ref {Sec-AHSS-KOhat},
are Steenrod squares or power operations. At the level of Thom spaces
and focusing on the prime 2,  Steenrod squares act on the Thom class
as Stiefel-Whitney classes. Indeed, in \cite{GS2} differential refinements 
of Steenrod squares are developed, and then used in \cite{GS3}\cite{GS5}\cite{GS6}
in characterizing the differentials in the spectral sequence. 

\medskip
In this section we consider instead 
differential refinements of the corresponding Stiefel-Whitney classes. 
The results of this section are known but we have not found 
a single source that contains all of 
the ingredients. Furthermore, we believe it worthwhile to collect
the results and discussions that we need in order to make the paper as 
self-contained and the arguments as streamlined as possible.

\medskip
Recall that the Stiefel-Whitney classes  $w_i$ are defined as certain characteristic classes 
in the $\ZZ/2$ cohomology ring
$
H^*(B{\rm O};\ZZ/2)
$,
where ${\rm O}:=\varinjlim{\rm O}(n)$ is the stable orthogonal group. In fact, if we restrict 
to dimension $n$, the cohomology ring $H^*(B{\rm O}(n); \ZZ/2)$ is the polynomial algebra 
on $w_i$, $i\leq n$. In particular, the first Stiefel-Whitney class arises as the generator 
of the polynomial algebra
$$
H^{*}(\RR P^{\infty}; \ZZ/2)\cong \ZZ/2[w_1]\;.
$$
The inclusion of diagonal matrices 
$(\ZZ/2)^n={\rm O}(1)^n \subseteq {\rm O}(n)$ induces a map at the level of classifying spaces 
$\alpha: (\RR P^\infty)^n=(B\ZZ/2)^n \to B{\rm O}(n)$. The corresponding map on cohomology 
$\alpha^*: H^*(B{\rm O}(n); \ZZ/2) \to H^*((\RR P^\infty)^n; \ZZ/2)=\ZZ/2[x_1, \cdots, x_n]$, 
is injective. The image is the subring $\ZZ/2[\sigma_1, \cdots, \sigma_n] \subset 
\ZZ/2[x_1, \cdots, x_n]$ generated by the elementary symmetric polynomials
$\sigma_i(x_1, \cdots, x_n)$ of degree $i$. Then define 
$w_i\in H^i(B{\rm O}(n); \ZZ/2)$ to be the preimage of $\sigma_i \in H^i((\RR P^\infty)^n; \ZZ/2)$
under $\alpha^*$. This then gives (see \cite{Q}\cite{Hu}\cite{MiS})
$$
H^*(B{\rm O}(n); \ZZ/2) = \ZZ/2[w_1, \cdots, w_n]\;.
$$  
Two-torsion is special when considering structures
and classes associated with the orthogonal group ${\rm O}(n)$ and 
the Spin group ${\rm Spin}(n)$, as torsion elements in $H^*(B{\rm O}(n); \Z)$ 
and $H^*(B{\rm Spin}(n);\Z)$ are all of order 2 \cite{BH}\cite{Ko}. 
Generators for the rings $H^*(B{\rm O}; \Z)$ and $H^*(B{\rm SO}; \Z)$
are given in \cite{Th1}\cite{Th2}, 
while relations between the generators are given in \cite{Br}\cite{Fe}. 
Generators for the ring $H^*(B{\rm Spin}; \Z)$ are given in \cite{Th3}.

\medskip
The following is important for the differential refinement, since integral lifts of 
mod 2 characteristic classes serve as part of the geometric data, as explained 
in \cite{HS}. Hence we find it worthwhile to discuss integral lifts of Stiefel-Whitney 
classes as we hope several subtle points will be appreciated. 
We are seeking class(es) ${\mathcal C}_i$ which
enjoy a relation of the type
$$
w_i = {\mathcal C} \mod 2\;, 
$$
which necessarily requires ${\mathcal C}_i$ to be of the same degree
as $w_i$.

\medskip
We now consider  specific classes which appear as generators of the above rings. 
To define these, one uses the coefficient short exact sequence
$
0\to \ZZ\overset{\times 2}{\longrightarrow}\ZZ\to \ZZ/2 \to 0
$
which gives the long exact Bockstein sequence on cohomology groups
\(
\label{les}
\xymatrix{
\hdots \to H^{n-1}(BO; \ZZ/2)\ar[r]^-{\beta} & H^{n}(B{\rm O}; \ZZ)\ar[r]^{\times 2}&
H^{n}(B{\rm O}; \ZZ)\ar[r]& H^{n}(B{\rm O}; \ZZ/2) \ar[r] & \hdots}
\)
and one defines the $(i+1)$ \emph{integral} Stiefel-Whitney class as $W_{i+1}:=\beta_2(w_i)$. 
These classes are 2-torsion within integral cohomology and so they still carry with them
the information that they arise from mod cohomology classes. 
Note that the whole subset $\beta_2(H^*(B{\rm O}(n); \Z/2)) \subset H^*(B{\rm O}(n); \Z)$ 
consists entirely of 2-torsion elements, the most prominent elements of which are 
the integral Stiefel-Whitney classes. 

\begin{remark}[Properties of the integral Stiefel-Whiteny classes]
Note the following:

\item {\bf (i)} $W_{j}$  measures whether the $(j-1)$st Stiefel-Whitney class $w_{j-1}$ is 
the mod 2 reduction of an integral class. 

\item {\bf (ii)}  All the $W_j$'s vanish when pulled back to $B{U}(n)$ under the maps
$B{U}(n) \to B{O}(2n)$  (see, e.g., \cite[Appendix B]{LM}). 

\item {\bf (iii)} One can describe these classes via subalgebras of $H^*(B{\rm SO}(n); \Z)$ 
(see \cite{BCM}).
Indeed, within the integral cohomology rings $H^*(B{\rm SO}(2n); \Z)$ or
$H^*(B{\rm SO}(2n+1); \Z)$ one can consider a subalgebra of integral elements of order 2
as 
\(
P_{\Z/2}\left\{ W_3, W_5, \cdots, W_{2n-1} \right\}\;.
\label{PZ2}
\) 
Furthermore, in the latter ring, one has $W_{2n+1}=\chi$, the Euler characteristic. 
\end{remark}

\medskip
The mod 2 reduction for real vector bundles
is rather subtle. It is not quite the case that the mod 2 reduction of the 
generators of $H^*(B{\rm SO}; \Z)$, i.e. the Pontrjagin classes, coincide with
the generators of $H^*(B{\rm SO}; \Z_2)$, which are the Stiefel-Whitney classes. 
The reason boils down to the presence of 2-torsion, as indicated above. 
The Pontrjagin classes for oriented real bundles reduce to some combination
of the Stiefel-Whitney classes, but we are interested in the lift of `bare'
Stiefel-Whitney classes. 
Note that in $H^*(B{\rm SO}(n); \Z_2)$, we have $\beta (w_{2i})=w_{2i+1}$ 
so that there exits a unique class of order 2 in $H^*(B{\rm SO}(n); \Z)$ 
which reduces mod 2 to $w_{2i+1}$. Indeed, these are  the integral 
Stiefel-Whitney classes $W_{2i+1}$.

\begin{lemma}[Integral lifts of Stiefel-Whitney classes]\label{Lemm-intSW}
The even Stiefel-Whitney classes $w_{2i}$ have no integral lifts in $H^*(B{\rm SO}(n), \Z)$.
The odd Stiefel-Whitney classes $w_{2i+1}$ have unique lifts in $H^*(B{\rm SO}(n), \Z)$
given by the integral Stiefel-Whitney classes $W_{2i+1}$. 
\end{lemma}

\theproof
Denote by $\beta_2$ and $\overline{\beta}$ the connecting 
homomorphisms, i.e. the Bockstein homomorphisms, for the two coefficients 
sequences $\Z \overset{\times 2}\longrightarrow \Z \overset{\rho_2}\longrightarrow \Z/2$ 
and $\Z/2 \overset{\times 2}\longrightarrow \Z/4 \overset{\rho_2}\longrightarrow \Z/2$, 
respectively, where $\rho_2$ denotes mod 2 reduction. Note that $\rho_2$ is a ring homomorphism. 
The three operations are related by the following diagram 
\(
\xymatrix@R=1.5em{
& \Z \ar[dr]^{\rho_2} & \\
\Z/2 \ar[ur]^{\beta_2} \ar[rr]^{\overline \beta} && \Z/2\;.
}
\)
Applied to the even/odd classes gives
\(
\xymatrix@R=1.5em{
& W_{2i+1} \ar@{|->}[dr]^{\rho_2} & \\
w_{2i} \ar@{|->}[ur]^{\beta_2} \ar@{|->}[rr]^{\overline \beta} && w_{2i+1}
}
\qquad \text{and} \qquad
\xymatrix@R=1.5em{
& W_{2i+2}\ar@{|->}[dr]^{\rho_2} & \\
w_{2i+1} \ar@{|->}[ur]^{\beta_2} \ar@{|->}[rr]^{\overline \beta} && 0\;.
}
\)
At the level of expressions, in $H^*(B{\rm O}; \Z_2)$  one has $\overline{\beta}w_{2i+1}=w_1 w_{2i+1}$ and 
$\overline{\beta} w_{2i}= w_{2i+1}+ w_1 w_{2i}$ (see e.g. \cite[Lemma 3.17]{Ha}), 
for which the terms involving $w_1$ vanish in the oriented case. 
The structure of the rings show that only one generator 
occurs in an odd dimension. 
\endofproof

\medskip
If we are willing to venture away from the Lie groups  ${\rm O}(n)$ and ${\rm SO}(n)$
then there are more possibilities, particularly for the even case.
Note that $\rho_2\beta_2(w_{2i})=\rho_2 W_{2i+1}=w_{2i+1}$,
which is not zero, so that there is no lift for $w_{2i}$. However, for space/bundles
for which $w_{2i+1}$ happens to be zero then we do have a lift.  
Beyond real bundles, for us there are other two 
sources of mod 2 reductions of the $w_i$'s: The first is  
Chern classes $c_i$ of complex bundles  ${\rm Spin}^c$ are unique integral lifts of $w_{2i}$.
The second is Spin bundles, as Spin characteristic classes $Q_i$ (under some conditions \cite{Th3})
are unique integral lifts of $w_{4i}$. 
Alternative (more general) lifts are given by Hopkins and Singer \cite{HS} as the Spin Wu classes 
$v^{\rm Spin}\in H^*(B{\rm Spin};\ZZ)$. 

\medskip
Note that there are geometric approaches to Steifel-Whitney classes; for instance:

\begin{enumerate}[{\bf (i)}]

\item \emph{{\v C}ech description of SW classes:} 
In \cite{McL} an explicit {\v C}ech cocycle representing the $k$-th Stiefel-Whitney
class $w_k$ of a vector bundle is constructed, 
involving only the transition function of the bundle. 
For each point $x$ in $U_{i_0} \cap \cdots \cap U_{i_k}$, 
a $(k-1)$-cocycle on $SO(n)/SO(k-1)$ is constructed which 
depends continuously on $x$. 

\item  \emph{ Cheeger-Simons refinement of SW classes:}
Differential refinements of Steifel-Whitney classes have been 
studied in \cite{HZ} \cite{Zw}, in which 
a very interesting secondary theory is developed for them. 
(This seems to work in all parities). 
\end{enumerate}

\subsubsection{Differential refinements of Stiefel-Whitney classes and the differential Wu-formula}

We now turn our attention to the natural differential refinements of the Stiefel-Whitney classes. 
There are many ways one could interpret such refinements. One way is to simply ask for a differential 
characteristic class $\hat{x}$, such that $\rho_2\mathcal{I}(\hat{x})=w_i$. In the language of
smooth stacks, we can phrase this by asking for a homotopy commutative diagram 
\(\label{diff-stwhit}
\xymatrix{
\BB {\rm O}_{\nabla}\ar[rr]^-{\hat{c}_i}\ar[d]_{\mathcal{I}} &&
\BB^iU(1)_{\nabla}\ar[d]^-{\rho_2\mathcal{I}}
\\
\op{BO}\ar[rr]^-{w_i} && {\rm B}^i \ZZ/2\;.
}
\)
Such lifts are very interesting to study and do appear naturally in other contexts
(see, e.g., \cite{HS}). However, they are not the most natural 
definition as one needs extra structure, as highlighted above
(see also \cite{memb}\cite{Wuc}). Indeed, since $\beta_2w_{i}$ 
is nonvanishing in general, one cannot even expect an integral lift in general. Nevertheless, 
under certain conditions one does have such a lift. However, the correct homotopy theoretic picture must 
change from the simple commutative diagram \eqref{diff-stwhit}. For example, for 
$w_2$  in general one has $\beta_2w_2=w_3+w_1w_2\neq 0$, so that $w_2$ cannot be 
lifted through $\rho_2$. But, from the Pasting Law for homotopy pullbacks, the 
cover $\op{BSpin}^c\to \op{BO}$ fits into the double Cartesian square
$$
\xymatrix{
\op{BSpin}^c\ar[rr]^{c_1}\ar[d] && {\rm B}^2\ZZ \ar[rr]\ar[d]^-{\rho_2}&& \ast\ar[d] 
\\
\op{BO}\ar[rr]^-{w_2}\ar@/_1pc/[rr]_{W_3} && {\rm B}^2 \ZZ/2\ar[rr]^-{\beta_2} && {\rm B}^{3}\ZZ
\;,}
$$
where the map $c_1$ can be identified with the first Chern class of the canonical line bundle associated to the 
${\rm Spin}^c$ structure. This means that on ${\rm BSpin}^c$ we indeed have an integral refinement of 
$w_2$. Then it makes sense to extend to a diagram of the form \eqref{diff-stwhit}, 
with $\BB{\rm Spin}^c_{\nabla}$ replacing $\BB{\rm O}_{\nabla}$, i.e.,  differential 
refinements of $w_2$ for ${\rm Spin}^c$ bundles. 

\medskip
Instead of asking for lifts of the form \eqref{diff-stwhit}, the most natural definition is to ask for a differential 
characteristic form which associates to each orthogonal vector bundle with connection a higher \emph{real} 
line bundle with connection, with a higher ${\rm O}(1)=\ZZ/2$ structure. Since this group is discrete, every such line 
bundle (and even higher line bundle) is necessarily flat. Thus the moduli stack classifying such bundles is 
simply $\BB^i \ZZ/2\simeq {\rm B}^i{\rm \ZZ/2}$. We have the following.

\begin{proposition}[Correspondence of $\ZZ/2$-cohomology of classifying stack]\label{prop-corZ2}
We have an isomorphism
$$
H^*(\BB{\rm O}(n)_{\nabla};\ZZ/2)\cong 
H^*(\op{BO}(n);\ZZ/2)\cong \ZZ/2[w_1,w_2,\hdots,w_n]\;,
$$
which is given by precomposition with the unit of the adjunction $\mathcal{I}:\mathbb{1}\to \delta\Pi$. 
\end{proposition}
\theproof
By definition, the cohomology groups of the stack $\BB{\rm O}(n)_{\nabla}$ are given by 
$$
H^i(\BB{\rm O}(n)_{\nabla};\ZZ/2):=\pi_0\map(\BB{\rm O}(n)_{\nabla};\BB^i\ZZ/2)\;.
$$
Since $\ZZ/2$ is geometrically discrete, we have $\BB^i\ZZ/2\simeq {\rm B}^i\ZZ/2$ and by the 
cohesive adjunction $\Pi\dashv \delta$ (see the Appendix), we find that 
$$
\pi_0\map(\BB{\rm O}(n)_{\nabla},B^i\ZZ/2) \cong 
\pi_0\map(\op{BO}(n),{\rm B}^i\ZZ/2)
\cong H^i(\op{BO}(n);\ZZ/2)\;,
$$
where the isomorphism is induced by the unit $\mathcal{I}$. Since the ring structure is induced from 
multiplication in $\ZZ/2$, it is clear that this map preserves the ring structure. 
\endofproof

Proposition \ref{prop-corZ2} indicates that, in general, the only appropriate definition for the differential 
Stiefel-Whitney classes is the following.
\begin{definition}[Differential Stiefel-Whitney classes]\label{diff-SW}
\noindent {\bf (i)}. We define the \emph{$i$-th differential Stiefel-Whitney class} $\hat{w}_i$ as the 
pullback $\mathcal{I}^*w_i$, where $\mathcal{I}:\BB{\rm O}(n)_{\nabla}\to \op{BO}(n)$ is the canonical 
map induced by the unit of the adjunction $\mathcal{I}:\mathbb{1}\to \delta\Pi$. We define the
\emph{total differential Stiefel-Whitney class} as 
$$
\hat{w}:=\hat{w}_1+\hat{w}_2+\hat{w}_3+\hdots \;. 
$$
\noindent {\bf(ii)}. We define the \emph{$i+1$ differential integral Stiefel-Whitney class}
$\widehat{W}_{i+1}$ as $j\circ j_2(w_i)$, where 
$$
\xymatrix{
j\circ j_2:H^i(-;\ZZ/2) \;\ar@{^{(}->}[r]&  H^i(-;U(1)) \;\ar@{^{(}->}[r]& \widehat{H}^i(-;\ZZ)
}.
$$
We define the \emph{total differential integral Stiefel-Whitney class} as 
$$
\widehat{W}:=\widehat{W}_1+\widehat{W}_2+\widehat{W}_3+\hdots \;.
$$
\end{definition}

\begin{remark}[Lifts vs. refinements]
Note that the classes $\hat{w}_i$ should not be confused with differential refinements of $w_i$ in the sense 
that $\rho_2\mathcal{I}(\hat{w}_i)=w_{i}$. Indeed, from Lemma \ref{Lemm-intSW}, it is not even true in general that $w_i$ admits an integral lift in general. However, the differential \emph{integral} Stiefel-Whitney classes of Definition \ref{diff-SW} are compatible with the integral Stiefel-Whitney classes. The next proposition makes this more precise.
\end{remark}

\begin{proposition}[Compatibility of differential Stiefel-Whitney classes]
The class $\widehat{W}_i$ refines $W_i$ in the sense that $\mathcal{I}(\widehat{W}_i)=W_i$.
\end{proposition} 
\theproof
Recall that from the differential cohomology diamond diagram \eqref{kodfdiam}, we have 
$\mathcal{I}j=\beta_{U(1)}$. Moreover, the Bockstein $\beta_{U(1)}$ is compatible with 
the Bockstein $\beta_2:H^*(-;\ZZ/2)\to H^{*+1}(-;\ZZ)$ in the sense that $\beta_{U(1)}j_2=\beta_2$. 
Then we calculate
$$
\mathcal{I}(\widehat{W}_i)=\mathcal{I} j j_2(w_i)=\beta_{U(1)}j_2(w_i)=\beta_2(w_i)=W_i\;.
$$

\vspace{-5mm}
\endofproof

With the definition of the differential Stiefel-Whitney classes in hand, it is now straightforward to prove the following.

\begin{theorem}[Differential Wu formula]
Let $V\to M$ be a real vector bundle on $M$ and let $\hat{\nu}$ be a differential Thom class in $\widehat{H}_c^*(V;\ZZ)=\widehat{H}^*({\rm Th}(V);\ZZ)$. Then we have the formula 
$$
\Psi_{V/M}\widehat{Sq}(\hat{\nu})=\hat{w}\;,
$$
where $\widehat{Sq}$ is the total differential Steenrod square as defined in \cite{GS2} and $\Psi_{V/M}$ is the pushforward in differential cohomology (see Proposition \ref{invdfth-iso}).
\end{theorem}
\theproof
From \cite{GS2}, we have the formula $\widehat{Sq}=j \circ j_2 \circ Sq \circ \rho_2 \circ  \mathcal{I}$.
Let $\psi_{V/M}$ be the Thom isomorphism in $U(1)$-cohomology and let $\overline{\psi}_{V/M}$ be the Thom isomorphism in $\ZZ/2$-cohomology. Since the differential Thom class refines the Thom class in integral cohomology,
the classical Wu  formula and compatibility of the Thom isomorphism in differential cohomology  
$U(1)$-cohomology (see Proposition \ref{thinjdfko}) give
\bea
\Psi_{V/M}\widehat{Sq}(\hat{\nu}) &=&  
\Psi_{V/M}\circ j\circ j_2\circ Sq\circ \rho_2\circ \mathcal{I}(\hat{\nu})
\\
&=& j\circ \psi_{V/M}\circ j_2\circ Sq(\nu)
\\
&=& j\circ j_2\circ \overline{\psi}_{V/M}Sq(\nu)
\\
&=:& \hat{w}\;.
\eea
The equality between the second and third line comes from the compatibility of the Thom isomorphism in $U(1)$-cohomology with $\ZZ/2$-cohomology.
\endofproof

\section{Appendix: Higher categorical framework}
\label{Ch-cat}

\subsection{Smooth stacks  and sheaves of spectra}

Let $\mathscr{M}{\rm an}$ be the category of smooth manifolds. This category becomes a site with the 
Grothendieck topology given by good open covers $\{U_{\alpha}\}$ of smooth manifolds (i.e. all finite 
intersections are contractible).  

\begin{definition}[Prestacks]
A prestack is a functor 
$$
X:\mathscr{M}{\rm an}^{\rm op}\longrightarrow \infty\mathscr{G}{\rm pd}.
$$
\end{definition}
The $\infty$-category of presheaves $\mathscr{P}\mathscr{S}{\rm h}_{\infty}(\mathscr{M}{\rm an})$ 
is defined as the $\infty$-category of functors ${\rm Fun}(\mathscr{M}{\rm an}^{\rm op},\infty\mathscr{G}{\rm pd})$. 
A prestack $X$ is called a \emph{stack}, if it satisfies descent along covers. That is, if 
\(\label{descfstak}
X(M)\simeq {\rm lim}\left\{\vcenter{
\xymatrix{
\prod_{\alpha}X(U_{\alpha})\ar@<.05cm>[r]\ar@<-.05cm>[r] & 
\prod_{\alpha\beta}X(U_{\alpha\beta})\ar@<.1cm>[r]\ar[r]\ar@<-.1cm>[r] & \hdots 
}}
\right\}\;.
\)
The full sub $\infty$-category on such objects is the $\infty$-category of \emph{smooth stacks} and 
is denoted $\mathscr{S}{\rm h}_{\infty}(\mathscr{M}{\rm an})$. According to 
\cite[Proposition 5.5.4.15]{Lur}, the canonical inclusion $i:\mathscr{S}{\rm h}_{\infty}(\mathscr{M}{\rm an})\into \mathscr{P}\mathscr{S}{\rm h}_{\infty}(\mathscr{M}{\rm an})$ 
admits a left adjoint 
$$
L:\mathscr{P}\mathscr{S}{\rm h}_{\infty}(\mathscr{M}{\rm an})\longrightarrow
\mathscr{S}{\rm h}_{\infty}(\mathscr{M}{\rm an})\;,
$$
which localizes at the strongly saturated class of morphisms generated by the covering maps 
$\coprod_{\alpha}U_{\alpha}\to M$. This functor is called the \emph{stackification} functor. 
The general theory of such localizations of presheaves is discussed in detail in \cite[Chapter 6]{Lur}. 
In particular, $L$ is also exact (preserves finite limits).

%


%

\begin{definition}[Presheaves of spectra]
A presheaf of spectra is a functor 
$$
X:\mathscr{M}{\rm an}^{\rm op}\longrightarrow \mathscr{S}{\rm p}.
$$
\end{definition}
The $\infty$-category of presheaves $\mathscr{P}\mathscr{S}{\rm h}_{\infty}(\mathscr{M}{\rm an};\mathscr{S}{\rm p})$ 
is defined as the $\infty$-category of functors ${\rm Fun}(\mathscr{M}{\rm an}^{\rm op},\mathscr{S}{\rm p})$. In the same way that a smooth prestack is a stack if it satisfies descent along covers, so a presheaf of spectra is called a \emph{sheaf} of spectra if it satisfies descent on such covers (i.e., we again have the limit diagram \eqref{descfstak}, but with $X$ taking values in spectra). Again, the canonical inclusion $i:\mathscr{S}{\rm h}_{\infty}(\mathscr{M}{\rm an};\mathscr{S}{\rm p})\into \mathscr{P}\mathscr{S}{\rm h}_{\infty}(\mathscr{M}{\rm an};\mathscr{S}{\rm p})$ 
admits a left adjoint 
$$
L:\mathscr{P}\mathscr{S}{\rm h}_{\infty}(\mathscr{M}{\rm an};\mathscr{S}{\rm p})\longrightarrow
\mathscr{S}{\rm h}_{\infty}(\mathscr{M}{\rm an};\mathscr{S}{\rm p})\;,
$$
which localizes at the covering maps 
$\coprod_{\alpha}U_{\alpha}\to M$. In the case of sheaves of spectra, this functor is called 
the \emph{sheafification} functor.

\medskip
The $\infty$-category of spectra $\mathscr{S}{\rm p}$ admits the structure of a symmetric monoidal 
$\infty$-category (via the smash product). This $\infty$-category can be presented (for example) by the 
simplicial model structure on symmetric spectra \cite{HSS}. Other presentations by orthogonal spectra and 
$S$-module spectra are Quillen equivalent to this one \cite[Theorem 0.1]{MMSS}, so these spectra can 
be presented as topological spectra if needed. 
\begin{proposition}[Monoidal structure]
The category of sheaves of spectra $\mathscr{S}{\rm h}_{\infty}(\mathscr{M}{\rm an};\mathscr{S}{\rm p})$ 
admits the structure of a symmetric monoidal $\infty$-category. Moreover, the topological realization functor 
$$
\Pi:\mathscr{S}{\rm h}_{\infty}(\mathscr{M}{\rm an};\mathscr{S}{\rm p})\longrightarrow
\mathscr{S}{\rm p}
$$
preserves products.
\end{proposition}
\theproof
That $\mathscr{S}{\rm h}_{\infty}(\mathscr{M}{\rm an};\mathscr{S}{\rm p})$ carries symmetric monoidal 
structure inherited from $\mathscr{S}{\rm p}$ is \cite[Proposition 1.15]{Lur}. To see that $\Pi$ preserves 
this product, observe that by the proof of \cite[Proposition 1.15]{Lur} (via \cite[Lemma 1.13]{Lur}), the 
smash product of sheaves is computed in presheaves. Hence it suffices to show that $\Pi$ preserves the product 
at the level of presheaves. On presheaves, $\Pi$ is presented by the colimit operation. Since the smash 
product preserves colimits in each variable, we are done.
\endofproof

\subsection{Group completion}

Let $C\mathscr{M}{\rm on}(\sset)$ be the $\infty$-category of commutative $\infty$-monoids and let 
$C\mathscr{G}{\rm r}(\sset)$ be the $\infty$-category of abelian $\infty$-groups. The canonical inclusion 
$$
\xymatrix{
i:C\mathscr{G}{\rm r}(\infty\mathscr{G}{\rm pd}) \; \ar@{^{(}->}[r] &
C\mathscr{M}{\rm on}(\infty\mathscr{G}{\rm pd})
}
$$
admits a left adjoint $\mathscr{K}$. We call this functor the \emph{group completion} functor. The unit of the adjunction $\eta:\mathbb{1}\to i\mathscr{K}$ takes in an $\infty$-monoid and outputs an infinite loop space (i.e. a connected spectrum). The induced functor on the homotopy categories
$$
{\rm Gr}:\mathscr{M}{\rm on}\to C\mathscr{G}{\rm r}\;,
$$
the usual Grothendieck group completion functor. More explicitly ${\rm Gr}$ takes a monoid $M$ and and forms the group whose elements are ordered pairs $(m_+,m_-)$ amd quotients by the equivalence relation 
$$
(a_+,a_-)\simeq (b_+,b_-) \Longleftrightarrow \text{there is\ } n\in M
\ \text{such that\ } a_+ +b_- +n=b_+ +a_-+n \;.
$$


The following is  of main interest in this paper. 
\begin{remark}[Vector bundles with connections]
Let $M$ be a smooth manifold and let ${\rm Vect}_{\nabla}(M)$ be the category of vector bundles with connections on $M$, i.e. 
\item {\bf (i)} the objects are pairs $(V,\nabla)\to M$ with $V$ a vector bundle and $\nabla$ a connection
\item {\bf (ii)}  and the morphisms are bundle maps $f:V\to W$ such that $\nabla(f(s))=f(\nabla(s))$. 

\noindent Let ${\rm Iso}({\rm Vect}_{\nabla}(M))$ be the subcategory on isomorphisms and form
$\pi_0{\rm Iso}({\rm Vect}_{\nabla}(M))$, the isomorphism classes of such objects. 
This set carries a monoidal structure induced from the Whitney sum of vector bundles. Then the elements of the Grothendieck group completion are given by pairs of vector bundles with connection $(V,\nabla_{V})$ and $(W,\nabla_{V})$ on $M$ modulo the equivalence relation 
\begin{align} 
\big((V,\nabla_{V}),(W,\nabla_{W})\big) \sim \big((V^{\prime},\nabla_{V^{\prime}}),&(W^{\prime},\nabla_{W^{\prime}})\big)   \Leftrightarrow \ \exists  (G,\nabla_{G}) \ \text{such that} \nonumber
\\
&  (V,\nabla_{V})\oplus(W^{\prime},\nabla_{W^{\prime}})\oplus (G,\nabla_{G})\cong (W,\nabla_{W})\oplus (V^{\prime},\nabla_{V^{\prime}})\oplus (G,\nabla_{G})\;. \nonumber
\end{align}
\begin{itemize}
\item Addition is performed component-wise via the Whitney sum.
\item The inverse of a pair $((V,\nabla_{V}),(W,\nabla_{W})$ is $((W,\nabla_{W}),(V,\nabla_{V})$. 
\end{itemize}
For this reason, we will usually denote the class of the pair $((V,\nabla_{V}),(W,\nabla_{W})$ as $[(V,\nabla_{V})]-[(W,\nabla_{W})]$. We denote this group by $\op{{\bf K}O}_{\nabla}(M)$. 
\end{remark} 

\begin{example}[Class of the trivial bundle and the reduced theory]
Consider the trivial bundle ${\bf n}:=\RR^n\times M\to M$, equipped with the trivial connection $d$. Then in the Grothendieck group we have $({\bf n},d)=n({\bf 1},d)$. Thus, all trivial vector bundles with connection are generated by $({\bf 1},d)$. After modding out by this copy of $\ZZ$ in the group $\op{{\bf K}O}_{\nabla}(M)$, all trivial bundles equipped with the trivial connection are zero. This corresponds to taking the \emph{reduced} cohomology.
\end{example}


\begin{example}[Virtual reduced classes]
Suppose $V \to M$ is a real Spin vector bundle. Then the class of 
$(V - \op{rank} V)$ in $\op{KO}(M)$ has AH filtration $\geq 4$. This is because 
the classifying map $V \to \ZZ \times \op{BO}$ of $(V - \op{rank} V)$ lifts to $\op{BSpin}$,
whose 3-skeleton is trivial (see \cite[Lemma B.10]{Fr}). Similarly, we can equip $V$ with connection
and form the difference 
$$
[(V,\nabla_{V})]-[({\bf  \op{rank}} V,d)]\;.
$$
\end{example}

\subsection{The hexagon diagram}

In this section, we show how to obtain the differential cohomology hexagon diagram. We will not give 
the proofs of each statement, as these can be found elsewhere (see \cite{Urs}\cite{BNV}\cite{GS3}). 
We will however, provide the proof for the existence of the hexagon 
diagram (which can also be found in \cite{BNV}). 

\medskip 
Note that the category of smooth manifolds admits finite products and every smooth manifold admits a map $x:\ast\to M$, which chooses a point in $M$. By \cite[Proposition 3.4.9]{Urs}, it follows that the $\infty$-category of smooth stacks admits a quadruple adjunction $\Pi\dashv \delta \dashv \Gamma \dashv \delta^{\dagger}$,
\(\label{cohadsset}
\xymatrix{
\mathscr{S}{\rm h}_{\infty}(\mathscr{M}{\rm an})\ar@<-.4em>[rr]|-{\Gamma} \ar@<1.3em>[rr]|-{\Pi} && \infty\mathscr{G}{\rm pd} \ar@<1.3em>[ll]|-{ \delta^{\dagger}} \ar@<-.4em>[ll]|-{\delta} 
}
.
\)
The functor $\delta:\infty\mathscr{G}{\rm pd}\into \mathscr{S}{\rm h}_{\infty}(\mathscr{M}{\rm an})$ takes 
an $\infty$-groupoid and forms the associated constant sheaf. The functor $\Gamma$ evaluates a sheaf 
of spectra on the point manifold. The functor $\Pi$ is called the \emph{topological realization}. 
\footnote{See footnote \ref{footgeom}.}
At the level 
of presheaves, $\Pi$ is left adjoint to the constant diagram functor. Hence, it presents the $\infty$-colimit 
operation over the corresponding presheaf. This functor has the following useful property.
\begin{proposition}
The functor $\Pi$ preserves finite products. 
\end{proposition}
\theproof
The category of smooth manifolds admits finite products and hence $\mathscr{M}{\rm an}^{\rm op}$ 
admits finite coproducts. Thus, $\mathscr{M}{\rm an}^{\rm op}$ is a sifted category. Since $\Pi$ can be 
presented as the colimit functor on presheaves, and colimits over sifted diagrams commute with finite 
products, the claim follows. 
\endofproof

\begin{remark}
The Quillen equivalence $\vert \cdot \vert \dashv {\rm sing}$, with $\vert \cdot \vert :\sset\to \mathscr{T}{\rm op}$ the topological realization of simplicial sets and ${\rm sing}:\mathscr{T}{\rm op}\to \sset$ the right adjoint which takes the singular nerve of a space, allows us to recover a topological space from a smooth stack $X$. Namely, the composite functor 
$$
\vert \Pi \vert:\mathscr{S}{\rm h}_{\infty}(\mathscr{M}{\rm an})\longrightarrow
 \mathscr{T}{\rm op}
$$
takes as input a smooth stack and returns a topological space which is built out of the geometric data encoded by the smooth stack. Throughout the paper, we sometimes refer to the space obtained from a smooth stack $X$ in this way as the \emph{topological realization} of $X$. Other times, we refer to the $\infty$-groupoid $\Pi(X)$ as the topological realization. Whatever the context may be, one can pass between the two models by either geometric realization or the singular nerve construction. 
\end{remark}

According to \cite[Proposition 4.1.9]{Urs}, the adjoints of diagram \eqref{cohadsset} lift to a quadruple 
adjunction in the stable case 
\(\label{cohadsset}
\xymatrix{
\mathscr{S}{\rm h}_{\infty}(\mathscr{M}{\rm an},\mathscr{S}{\rm p})\ar@<-.4em>[rr]|-{\Gamma} \ar@<1.3em>[rr]|-{\Pi} &&  \ar@<1.3em>[ll]|-{ \delta^{\dagger}} \ar@<-.4em>[ll]|-{\delta} \mathscr{S}{\rm p}
}
.
\)
The topological realization has a particularly nice form in the stable case. According to \cite[Proposition 7.6]{BNV}, 
we have an equivalence
\vspace{-2mm} 
\(
\Pi(\E)(M)\simeq {\rm lim}\left\{\vcenter{
\xymatrix{
\mathcal{E}(\Delta^0)\ar@<.05cm>[r]\ar@<-.05cm>[r] &
\E(\Delta^1)\ar@<.1cm>[r]\ar[r]\ar@<-.1cm>[r] & \hdots 
}}
\right\},
\)
where right $\infty$-limit is computed over the simplicial diagram given by evaluating on smooth simplices
$\Delta^n$, viewed as manifolds with corners. Actually, there a much stronger statement was proved -- namely,
that the composite functor $\delta\Pi$ is represented by the `smooth singular sheaf' operation. We will not 
need this here, however. 

\begin{proposition}[Diamond/hexagon for sheaves of spectra]
Let $\E$ be a sheaf of spectra. Set $\mathcal{I}:=\eta_{\E}:\E\to \delta\Pi\E$, where $\eta$ is the unit of the 
adjunction $\eta:\mathbb{1}\to \delta \Pi$, $j:=\epsilon_{\E}:\delta\Gamma\to \mathbb{1}$, where 
$\epsilon:\delta\Gamma\to \mathbb{1}$ is the counit, $\mathcal{R}$ is the canonical map to the cofiber 
and $a$ is the canonical map out of the fiber. We have a homotopy commutative diagram
\vspace{-3mm} 
\(
\xymatrix@R=2em@C=1pt @!C{
&{\rm fib}(\mathcal{I}) \ar[rd]^-{a}\ar[rr] & & {\rm cofib}(j) \ar[rd] &
\\
\Sigma^{-1}{\delta \Pi}{\rm cofib}(j)\ar[ru]\ar[rd] & &
{\E}\ar[rd]^-{\mathcal{I}} \ar[ru]^-{\mathcal{R}}& &  {\delta \Pi}{\rm cofib}(j)\;.
\\
& \delta\Gamma{\E}\ar[ru]^{j}\ar[rr] & & {\delta \Pi}{ \E}\ar[ru] &
}
\)
Moreover, the two squares are homotopy Cartesian, the top and bottom sequences are part of long
fiber sequences and the diagonals are exact. 
\end{proposition}

\newpage

\theproof
In any stable infinity category, a square is Cartesian if and only if the induced map on the fibers (taken at 
the zero point $0:\ast\to \E$) is an equivalence. Moreover, a sequence is a fiber sequences if and only if
it is a cofiber sequence in a stable $\infty$-category. As a left adjoint, the functor $\delta\Pi$ sends
fiber/cofiber sequences to fiber/cofiber sequences. Since $\delta\Pi\delta\Gamma=\delta\Gamma$, 
we see that applying $\delta\Pi$ to the diagonal 
$$
\delta\Gamma\E\overset{j}{\longrightarrow} \E\overset{\mathcal{R}}{\longrightarrow} {\rm cofib}(j)
$$
gives the bottom fiber/cofiber sequence and the right square must be Cartesian.  Arguing 
dually for the left square proves the claim.
\endofproof

\subsection{Ordinary differential cohomology}

In this section we review some of the basic properties of ordinary differential cohomology. 
See \cite{HS}\cite{CS}
\cite{SSu1}\cite{Bun}\cite{Urs}\cite{BNV}. See  also \cite{GS3} and more extensive 
list of references therein.
\begin{definition}\label{rfnhzthr}
A {\rm differential refinement} $\widehat{H}^*(-;\ZZ)$ of integral cohomology $H^*(-;\ZZ)$ consists of the 
following data:
\begin{enumerate} 
\item A functor $\widehat{H}^*(-;\ZZ):Sh_{\infty}(\mathscr{M}{\rm an}_+)^{\rm op} \to \mathscr{A}\mathrm{b}_{gr}$;
\item Three natural transformations:
\begin{enumerate}
\item {\rm Integration:} $\mathcal{I}:\widehat{H}^*(-;\ZZ)\to H^*(-;\ZZ)^*$;
\item {\rm Curvature:} $\mathcal{R}:\widehat{H}^*(-;\ZZ)\to \Omega^*_{\rm cl}$;
\item {\rm Secondary Chern character:} $a:\Omega^{*-1}(-)/{\rm im}(d)\to \widehat{H}^*(-;\ZZ)$;
\end{enumerate}
\end{enumerate}
such that the following axioms hold:
\begin{list}{$\circ$}{}
\item {\rm (Chern-Weil).} We have a commutative diagram 
$$
\xymatrix@R=1.5em{
\widehat{H}^*(-;\ZZ)\ar[rr]^-{R}\ar[d]_{I} && \Omega^*_{\rm cl}\ar[d]^-{q}
\\
H^*(-;\ZZ)\ar[rr]^-{i} && H_*(-;\RR)
}
$$
where $i:\ZZ\into \RR$ is the canonical map. 
\item {\rm (Secondary Chern-Weil).} We have a commutative diagram
$$
\xymatrix@R=1.5em{
\Omega^{*-1}/{\rm im}(d) \ar[rr]^-{d} \ar[dr]_-{a} &&  \Omega^*_{\rm cl}
\\
& \widehat{H}^*(-;\ZZ)^*\ar[ur]_-{\mathcal{R}} & 
}
$$
and an exact sequence
$$
\hdots \longrightarrow H^{*-1}(-;\ZZ) \longrightarrow \Omega^{*-1}/{\rm im}(d)\overset{a}{\longrightarrow} \widehat{H}^*(-;\ZZ) \overset{\mathcal{I}}{\longrightarrow} H^*(-;\ZZ)\longrightarrow \hdots\;.
$$
\end{list}
\end{definition}
There are various models for differential cohomology, such as smooth Deligne cohomology
(see \cite{Urs}\cite{Cech}
\cite{Bun}\cite{FSS1}\cite{FSS2} and references therein), 
Cheeger-Simons differential characters \cite{CS}, and the general construction 
of Hopkins-Singer \cite{HS} (specialized to $\mathscr{H}\ZZ$). We will 
review both perspectives of Deligne cohomology 
and of Hopkins-Singer. 

\begin{proposition}\label{stckdelpr}
Consider the smooth Deligne complex 
$$
\mathcal{D}(n):=\big(\vcenter{\xymatrix{ \hdots \ar[r] & \underline{U(1)}\ar[r]^-{d\log} & \Omega^1\ar[r] 
& \hdots \ar[r] & \Omega^{n-1}}}\big)\;,
$$
where $\Omega^{n-1}$ is in degree zero. The hypercohomology groups $H^0(-;\mathcal{D}(n))$ satisfy the axioms in Definition \ref{rfnhzthr}, restricted to the subcategory of smooth manifolds. More generally, if
$$
\op{DK}:=U\Gamma:\chp\longrightarrow \sab\longrightarrow \sset\;,
$$
is the Dold-Kan functor, then the functor $\pi_0\map(-, \op{DK}(\mathcal{D}(n)))$ defines a presheaf of 
abelian groups on $Sh_{\infty}(\mathscr{M}{\rm an}_+)$ satisfying the properties in Definition \ref{rfnhzthr}.
\end{proposition}
\theproof
The first statement is classical (see for example \cite{Bun} for a review). For the second claim, observe 
that the maps $\mathcal{I}$, $\mathcal{R}$ and $a$ are induced by the morphism of sheaves of chain 
complexes.  By the basic properties of the Dold-Kan correspondence, we have a natural isomorphism 
of presheaves 
\(\label{natisodel}
H^0(-;\mathcal{D}(n))\cong \pi_0\map(-, \op{DK}(\mathcal{D}(n)))\;.
\)
By functoriality of ${\rm DK}$, the maps $\mathcal{I}$,$\mathcal{R}$ and $a$ induce corresponding 
morphisms of smooth stacks. The mapping space on the right is defined for any smooth stack and taking 
$\pi_0$ gives an abelian group, with group structure inherited via the isomorphism \eqref{natisodel}. 
Thus, the second claim follows immediately from the first.
\endofproof
\begin{definition}
We define the moduli stack of $n$-gerbes with connection as the smooth stack presented by the 
simplicial sheaf
$$
\BB^{n-1}U(1)_{\nabla}:=\op{DK}(\mathcal{D}(n))\;.
$$
\end{definition}
From the basic properties of the Dold-Kan correspondence, we have a natural isomorphism 
$$
\pi_0\map(M;\BB^{n-1}U(1)_{\nabla})\cong H^0(M;\mathcal{D}(n))\;,
$$
for each smooth manifold $M$. This gives a presentation of ordinary differential cohomology 
by isomorphism classes of $n$-gerbes with connection (see 
\cite{Urs}\cite{Cech}\cite{FSS1}\cite{FSS2}). Smooth Deligne 
cohomology admits a ring structure via the Deligne-Beilinson cup product 
$$
\cup_{\rm DB}:H^0(M;\mathcal{D}(n))\otimes 
H^0(M;\mathcal{D}(m))\longrightarrow  H^0(M;\mathcal{D}(n+m))\;,
$$
and this product comes from a morphism of sheaves of complexes 
$\cup_{DB}:\mathcal{D}(n)\otimes \mathcal{D}(m)\to \mathcal{D}(n+m)$. 
Using the monoidal properties of ${\rm DK}$, one can show that this product 
induces a morphism of smooth stacks \cite{}
$$
\cup_{\rm DB}:\BB^{n-1}U(1)_{\nabla}\times 
\BB^{m-1}U(1)_{\nabla}\longrightarrow \BB^{n+m-1}U(1)_{\nabla}\;,
$$
which induces a product operation on $U(1)$-gerbes. This product gives ordinary 
differential cohomology a ring structure. 

\medskip
We can also present differential cohomology via the Hopkins-Singer construction. Let 
$H:\mathscr{C}{\rm h}\to \mathscr{S}{\rm p}$ be the Eilenberg-MacLane functor, which 
sends an unbounded chain complex to its corresponding spectrum. Consider the sheaves of 
spectra $\mathscr{H}\ZZ$, $\mathscr{H}(\tau_{\geq 0}\Omega^*)$ and 
$\mathscr{H}(\Omega^*)\simeq \mathscr{H}\RR$. The $\infty$-pullback
$$
\xymatrix{
\widehat{\Sigma^n\mathscr{H}\ZZ}\ar[rr]\ar[d] && \mathscr{H}(\tau_{\geq 0}\Omega^*[n])\ar[d]
\\
\Sigma^n\mathscr{H}\ZZ\ar[r]^-{i} & \Sigma^n\mathscr{H}\RR\ar@{^{(}->}[r]_{\simeq} & \mathscr{H}(\Omega^*[n])
}
$$
gives a sheaf of spectra representing $\widehat{H}^n(-;\ZZ)$. Indeed, unwinding the data described in
\cite[Section 4.4]{BNV}, we immediately see that the cohomology 
$\pi_0\map(M_+,\widehat{\Sigma^n\mathscr{H}\ZZ})$ satisfies the properties of Definition 
\ref{stckdelpr}. In this case, the corresponding differential cohomology hexagon diagram looks as follows
\(\label{ordfcodiam}
\xymatrix@C=1pt @!C{
&\Omega^{*-1}/{\rm im}(d) \ar[rd]^-{a}\ar[rr]^-{d} & & \Omega^*_{\rm cl} \ar[rd] &
\\
H^{*-1}(-;\RR) \ar[ru]\ar[rd] & &
\widehat{H}^*(-;\ZZ)\ar[rd]^-{\mathcal{I}} \ar[ru]^-{\mathcal{R}}& &  H^*(-;\RR)\;.
\\
& H^*(-;U(1))\ar[ru]^{j}\ar[rr] & & H^*(-;\ZZ)\ar[ru] &
}
\)
Differential cohomology is \emph{not} a homotopy invariant theory. In particular, 
$\widehat{H}^*([0,1]\times M;\ZZ)\not \simeq \widehat{H}^*(M;\ZZ)$.
However, the discrepency is measured quite nicely by the homotopy formula
\(\label{hofmdfco}
i^*_0\hat{x}-i_1^*\hat{x}=\int_{[0,1]\times M/M}\mathcal{R}(\hat{x})\;,
\)
where $i_0,i_1:M\into M\times [0,1]$ are the maps which include at the endpoints. 
We refer the reader to \cite{BNV} for the proof in the general case (see also \cite{Bun}). 

\begin{example}[Chern-Simons form]
Given a characteristic form $\omega$ and a pair of connections $\nabla$,$\nabla^{\prime}$ on a 
(real or complex) vector bundle $V\to M$, the associated Chern-Simons form is defined by the transgression
\(\label{csfmtran}
{\rm cs}_{\omega}(\nabla,\nabla^{\prime}):=
\int_{[0,1]\times M/M}\omega(t{\rm pr}_0^*\nabla +(1-t){\rm pr}^*_1\nabla^{\prime})\;.
\)
This defines a class in differential cohomology by application of the natural map $a$ in 
diagram \eqref{ordfcodiam}. 
\end{example}

Consider the natural map $r:H^n(M;\ZZ)\to H^n(M;\RR)$. Given an  integral cohomology 
class $x\in H^n(M;\ZZ)$ and closed differential form  $\omega\in \Omega^n_{\rm cl}(M)$ 
whose underlying de Rham class is equal to $r(x)$,  It is natural to ask whether there  exists 
a unique differential refinement. In fact, this is \emph{not} the case. 

\begin{proposition}[Set of differential refinements]\label{prop-setref}
Let $x\in H^n(M;\ZZ)$ be a class and $\omega\in \Omega^n_{\rm cl}(M)$ be such that $[\omega]=r(x)$. 
Then the set of differential refinements is a torsor for the group 
$a(H^{n-1}(M;\RR))\cong H^{n-1}(M;\RR)/H^{n-1}(M;\ZZ)$. 
\end{proposition}
\theproof
Let $\hat{x}$ be any differential refinement. Then  by the diamond diagram \eqref{ordfcodiam}, for any element $y\in H^{n-1}(M;\RR)$, the class $\hat{x}+a(y)$ satisfies $\mathcal{R}(\hat{x}+a(y))=\mathcal{R}(\hat{x})$. Similarly, by the diagonal exact sequence $\Omega^{*-1}(M)/{\rm im}(d)\overset{a}{\longrightarrow} \widehat{H}^*(M;\ZZ)\overset{\mathcal{I}}{\longrightarrow} H^*(M;\ZZ)$, we have $\mathcal{I}(\hat{x}+a(y))=\mathcal{I}(\hat{x})=x$. This gives an action on the set of lifts. This is transitive, since similar considerations show that any two lifts $\hat{x}$ and $\hat{y}$ the difference $\hat{x}-\hat{y}$ is in the kernel of both $a$ and $\mathcal{R}$. Thus, by the commutativity of the left square in \eqref{ordfcodiam}, there is $z\in H^{n-1}(M;\RR)$ such that $\hat{x}-\hat{y}=a(z)$. The action is clearly free. 
\endofproof



\end{document}